\documentclass{article}
\usepackage{latexsym}
\usepackage{epsfig}

\begin{document}

\title{\huge The Algebraic Proof of the Universality Theorem}
\author{Ai-Ko Liu \footnote{email address: 
akliu@math.berkeley.edu} \footnote{Current Address: 
Mathematics Department of U.C. Berkeley}\footnote{
 HomePage:math.berkeley.edu/$\sim$akliu}}
\maketitle

\newtheorem{main}{Main Theorem}
\newtheorem{theo}{Theorem}
\newtheorem{lemm}{Lemma}
\newtheorem{prop}{Proposition}
\newtheorem{rem}{Remark}
\newtheorem{cor}{Corollary}
\newtheorem{mem}{Examples}
\newtheorem{defin}{Definition}
\newtheorem{axiom}{Axiom}
\newtheorem{conj}{Conjecture}

\section{\bf Preliminary}

 In [LL1] the authors had developed a topological version of
 family Seiberg-Witten invariants. In an attempt to generalize Taubes'
 ``SW=Gr'' to families, it was observed that multiple coverings of
 exceptional curves may show up in the enumeration.
 In [Liu4] a family curve counting scheme
 has been proposed for algebraic families.  
In [L3] the topological/algebraic version of
 family blowup formula was derived. In [L5] the topological/algebraic version
 of family switching formula was derived. It was used in [Liu1] to
 derive the G$\ddot{o}$ttsche-Yau-Zaslow conjecture regarding counting of
 nodal curves on algebraic surfaces, including $K3$.

 The derivation in [Liu1] makes use of the various techniques from 
 differential topology and complex geometry as well as 
 the technique from the theory of 
pseudo-holomorphic curves. In the derivation the
 author had assumed that the reader is reasonably familiar with Taubes proof
of ``SW=Gr''.  On the other hand, a version of algebraic
 family Seiberg-Witten invariant was defined in [Liu3] and its
family blowup formula and switching formula was derived in [Liu3] and
 [Liu5], respectively, based on intersection theory developed in [F]. 

 Thus it is desirable to give a purely algebraic proof of the
 universality theorem, the backbone in deriving the
 G$\ddot{o}$ttsche-Yau-Zaslow formula, based on intersection theory [F] and
 algebraic geometry [Ha].

 The following universality theorem is the main theorem proved in this
 paper.  

\begin{main}\label{main; universal}
Let $\delta\in {\bf N}$ denote the number of nodal singularities.
Let $L$ be a $5\delta-1$ very-ample line bundle on an algebraic surface $M$,
 then the number of $\delta-$node nodal singular curves in a generic 
 $\delta$ dimensional linear sub-system of $|L|$ can be expressed as 
 a universal polynomial (independent to $M$) 
of $c_1(L)^2$, $c_1(L)\cdot c_1(M)$, $c_1(M)^2$, 
$c_2(M)$ of degree $\delta$.
\end{main}

  For simplicity, we take ${\bf C}$ to be our ground field. But the
 same argument works for algebraic closed fields of characteristic
 zero as well. 

  The term ``number of $\delta-$node nodal curves in a linear system'' used here 
is the weighted number of singular algebraic curves with isolated $\delta-$nodal
 singularities defined in [Got]. This concept is
 closed related but not always equal to Gromov-Witten invariant of $c_1(L)$.
  For their potential difference, please consult section \ref{section; compare}.

 In the following we outline the different 
functionalities of the various sections of
 the paper. Because the paper is a bit lengthy--involving quite a few
 new notations not widely used in the community of enumeration
 geometry/algebraic geometry, in subsection
 \ref{subsection; notation} we offer some simple
 advice, a notation table and quite a lot of footnotes, 
which may help the reader to get used to our notations and read 
the paper more fluently.

 In section \ref{section; strata}, we review the basic facts about
 the universal spaces and then introduce the 
relative-universal spaces over a base $B$.
 We also recall the concept of admissible graphs, admissible strata and 
 the stratification of the (relative) universal spaces. 
 To minimize the dependence on [Liu1], we prove a few useful results which 
will be used in the latter sections. As a result, an admissible stratum can be 
 characterized as the locus of co-existence of type $I$ exceptional classes
 attached to it.

\medskip

 In section \ref{section; bundle} and the subsidiary subsections, we develop 
 a technique to construct and identify the quotient bundle ${\bf V}_{quot}$
 of the given obstruction vector bundle ${\bf W}_{canon}|_{Y(\Gamma)\times T(M)}$, 
given some datum of quotient sheaves. 

 In section \ref{section; discrepancy}, we review the
 residual intersection formula of top Chern classes [F] and 
develop an algebraic tool to 
compare the top Chern classes of vector bundles $\sigma:E\mapsto F$
 isomorphic off a closed subset. 
Under the bundle homomorphism $\sigma$, a section $s_0$ of $E$
 induces a section $\sigma(s_0)$ of $F$ and the difference of their
 top Chern classes can be studied by the difference of
 their localized top Chern classes along
 the zero loci $Z(s_0)$, $Z(\sigma(s_0))$, respectively.
We find that after some blowing ups along loci in
 $\overline{Z(\sigma(s_0))-Z(s_0)}$,
 the top Chern class of the residual bundle of $F$ has the same ``numerical
 property'' as the top Chern class of the pull-back of 
$E$. We achieve this goal by a graph 
 construction in the projective space bundle ${\bf P}(E\oplus 1)$.
 This proposition is crucial in controlling the seemingly ``unmanageable''
 blowing ups and relate the modified bundle to better understood objects.

 The canonical cross section $s_{canon}$ of the 
canonical algebraic obstruction bundle \footnote{See definition 5.3 in 
[Liu3] for its definition.}
 ${\bf H}\otimes \pi_X^{\ast}{\bf W}_{canon}$
 defines a zero locus in $X={\bf P}({\bf V}_{canon})$ which contains the  
 sub-locus (closure) $\overline{{\cal M}_{C-{\bf M}(E)E}\times_{M_n}Y_{\gamma_n}}$ 
that we want to study. The proof of the main theorem enables us to 
 attach an invariant to this locus.
 In section \ref{section; main}
 and the sub-sequential sub-sections, we initiate the proof of the main theorem
 by introducing an inductive procedure of blowing up the smooth total space of
 the canonical family algebraic Kuranishi space $X={\bf P}({\bf V}_{canon})$.
 The goal is to apply residual intersection theory inductively and remove
 all the excess contributions.
 In subsection \ref{subsection; modinv}, we define the modified algebraic 
 family Seiberg-Witten invariants attached to the various smooth sub-loci.
 We also address in subsection \ref{subsection; comb} some combinatorial
 questions regarding the sub-loci of ${\cal M}_{C-{\bf M}(E)E}$. In the same
 sub-section, we also
 explain the geometric meanings of several partial orderings $\gg, \succ$,
 $\sqsupset$, etc. which had already appeared in [Liu1], [Liu4], [Liu5] before.

 The section \ref{section; proof} is the core of the current paper. In this
section, we finish up the proof of the main theorem. 
We address in subsection \ref{subsection; independence} some combinatorial
 questions regarding the independence of the localized top Chern classes to the 
 permutations/collapsing of the blowing up orderings. We prove inductively
 that the various excess contributions of the algebraic family invariant
 ${\cal AFSW}_{M_{\delta+1}\times \{t_L\}\mapsto M_{\delta}\times \{t_L\}}(1,
 C-2\sum_{1\leq i\leq \delta}E_i)$ 
can be identified with the various modified family algebraic Seiberg-Witten 
invariants defined in subsection \ref{subsection; modinv}. In subsection 
\ref{subsection; transv}, we enhance G$\ddot{o}$ttsche's argument slightly to get the
 necessary finiteness result on nodal curves that we enumerate.

 Finally in the appendix, section \ref{section; compare}, 
we offer a light-weighted comparison between our family invariants and
 the standard Gromov-Witten invariants of algebraic surfaces.
 Even though family Seiberg-Witten invariants on the universal spaces 
are related to Gromov-Witten invariants of algebraic surfaces, there is 
 some subtle difference between them which may cause confusion to
 the readers who are new to family Seiberg-Witten invariants.

\medskip

 In this paper, we do not attempt to address the issue of identifying these
 universal polynomials in the universality theorem 
algebraically. The discussion about the 
relationship of the above universality theorem with Riemann-Roch formula
 along with some open problems and conjectures related to this theorem will
 be addressed in a separated note [Liu7] elsewhere.

  Finally, it is recommended to use [Liu3], [Liu4], [Liu5] as companions
  reading this paper.

\medskip

\subsection{A Table of Our Notations}\label{subsection; notation}

\bigskip

  In this paper we give an algebraic proof of the universality theorem, and 
most of our notations
 have already appeared in [Liu1], [Liu3], [Liu4], [Liu5], etc. On the
 other hand a few new notations must be introduced in the purely algebraic
argument of the paper. For 
the standard notations in algebraic
 geometry and intersection theory, the reader may consult [Ha], [F].
 In the following, we list the
 frequently used ``global''
\footnote{For its definition, see the next paragraph.}
 notations in the paper. Following the convention of the earlier papers [Liu3], 
[Liu5],
 a locally free sheaf of sections will be denoted by the calligraphic character, 
say ${\cal G}$, 
if and only if the corresponding algebraic vector bundle has been denoted by the
 bold character ${\bf G}$.

 Another convention in this paper is that
the variables/notations defined within the proof of a lemma or a 
proposition is viewed as a ``local'' variable and its ``scope'' is within this 
particular proof.
 We may recycle the same variable/notation in the proofs of the other 
lemma/proposition in a different context. 
On the other hand,
 the variables/notations declared in the definitions are viewed as the 
``global'' variables. Subject to some specialization of values, 
their meanings are fixed 
throughout the paper.
 Finally, the variables/notations declared in the text of the paper
 are ``semi-local'' in 
the sense that
 their scope is the whole section containing the particular text.
 If we refer to this 
variable/notation from 
 a different section, we will indicate to the reader the
 location (page) where it has 
been defined.

\bigskip

  The following is the list of notations widely used throughout the paper.
  Every symbol is leaded by a $\bullet$.

\bigskip

$\bullet$ $adm(n)$--the set of $n$-vertexes admissible graphs satisfying the
 five axioms starting at page \pageref{graphaxiom}.

$\bullet$ $adm_2(n)$--the subset of $adm(n)$ consisting of fan-like admissible 
graphs. See definition \ref{defin; special}
 on page \pageref{defin; special} and fig.4 on page \pageref{fig.4}.

$\bullet$ $C$--the class in $H^{1, 1}(M, {\bf Z})$. In this paper $C$ is assumed to
 satisfy ${\cal R}^1\pi_{\ast}\bigl({\cal E}_C\bigr)=
{\cal R}^2\pi_{\ast}\bigl({\cal E}_C\bigr)=0$.

$\bullet$ $C-{\bf M}(E)E$--In this expression the term 
$-{\bf M}(E)E=-\sum_{1\leq i\leq n}m_iE_i$ 
is interpreted as a cohomology class
 instead of an effective divisor. $C-{\bf M}(E)E$ is a class 
of Hodge type $(1, 1)$. 

$\bullet$ ${\cal C}_{\Gamma}$--the simplicial exceptional cone constant over 
${\cal S}_{\Gamma}\supset Y_{\Gamma}$. 

$\bullet$ $\Delta(n)$--the subset of $adm(n)$ containing admissible graphs 
which satisfy the additional extremal conditions on page 
\pageref{specialcondition}. See definition on page \pageref{defin; special}

$\bullet$ $\Delta_{k_i}$--the (anti)-effective divisor in 
$\tilde{\Xi}_{k_i}\mapsto Y(\Gamma_{e_{k_i}})$ which appears in the
 exact sequence relating 
${\cal R}^1\tilde{\pi}_{\ast}\bigl({\cal O}_{\tilde{\Xi}_{k_i}}\otimes
 {\cal E}_{-{\bf M}(E)E}\bigr)$ and
 ${\cal R}^1\tilde{\pi}_{\ast}\bigl({\cal O}_{\tilde{\Xi}_{k_i}}(e_{k_i})\bigr)$.
 Consult page \pageref{prop; exact} for more details.

$\bullet$ $\delta$--the number of nodes in the main theorem.

$\bullet$ $e_i$--the extremal generators of ${\cal C}_{\Gamma}$, called type 
$I$ exceptional class.

$\bullet$ $E_i$--the effective exceptional divisor $E_{i; n+1}$ in $M_{n+1}$ 
and is viewed as the fiberwise divisor of
 $M_{n+1}\mapsto M_n$.

$\bullet$ $E_{a; b}$--the effective exceptional divisor in $M_n$ corresponding 
to the blowing up along the $(a, b)-$th partial 
diagonal of $M^n$.

$\bullet$ ${\cal E}_C$--the invertible sheaf over $M\times T(M)$ with 
$c_1({\cal E}_C|_{M\times \{t\}})=C$ for $\forall t\in T(M)$.

$\bullet$ ${\cal E}_{C-{\bf M}(E)E}$--the invertible sheaf over 
$M_{n+1}\times T(M)$ corresponding to the cohomology class $C-{\bf M}(E)E$ of Hodge 
type $(1, 1)$.

$\bullet$ ${\cal EC}_b(\underline{C}, Q)$--the type $I$ exceptional 
cone associated with $\underline{C}$ over $b\in M_n$. See definition 5 of [Liu4] for its 
 definition.

$\bullet$ $f_n$--the projection map $M_{n+1}\mapsto M_n$.

$\bullet$ $f_{n-1; k}$--the composition of the projection map 
$f_k\circ f_{k+1}\circ \cdots \circ f_{n-1}:M_n\mapsto M_k$.

$\bullet$ $\Gamma$--A typical element of an admissible graph $\Gamma\in adm(n)$.
 See fig. 1 on page \pageref{fig.1} for an example.

$\bullet$ $\Gamma_{e_i}$--the fan-like admissible graph attached to $e_i$ 
in which the index 
$i$ is the only direct ascendent index and the indexes $j_i$ 
appearing in $e_i=E_i-\sum_{j_i}E_{j_i}$ are
 the direct descendent indexes of $i$. See fig.2 on page \pageref{fig.2}
 for some examples.

$\bullet$ $\gamma_n$ or $\gamma$--the unique admissible graph of $n$ 
vertexes with no one edge. 

$\bullet$ 
${\bf H}$ and ${\cal H}$--the hyperplane bundle and its
 invertible sheaf of sections on ${\bf P}({\cal V}_{canon})=
{\bf P}({\bf V}_{canon}^{\circ})$
 induced by the linear structure ${\bf V}_{canon}$.

$\bullet$ $I_{\Gamma}$--the subset of $\Delta(n)$ which collects all the elements 
smaller than $\Gamma$ under the linear
 ordering $\models$. Consult page \pageref{defin; already} for its definition.

$\bullet$ $\bar{I}_{\Gamma}$--the reduced subset of $I_{\Gamma}$ throwing away 
elements $\in I_{\Gamma}$ which are
 $\ll$ than some other element in $I_{\Gamma}$. 
Consult page \pageref{defin; reduce} for its definition.

$\bullet$ $\bar{I}_{\Gamma}^{\gg}$--The subset of $\bar{I}_{\Gamma}$ collecting 
all the elements in $I_{\Gamma}$ which
 are $\ll \Gamma$. Consult page \pageref{defin; ggset} for its definition.

$\bullet$ $j_i$--the typical direct descendent index of $i$. The subscript $i$
 in $j_i$ indicates $j_i$ is a direct descendent of $i$.

$\bullet$ $k_i$--the subscripts in $\{1, 2, \cdots , n\}$ which corresponds to 
the indexes of type $I$ exceptional classes with $e_{k_i}\cdot (C-{\bf M}(E)E)<0$.

$\bullet$ $\Xi_i$--the ${\bf P}^1$ fibration over $Y(\Gamma_{e_i})$ 
representing the universal curves of the type
 $I$ exceptional class $e_i$. The notation is used starting in section 
\ref{section; bundle}.

$\bullet$ $\tilde{\Xi}_i$--the relative minimal model of $\Xi_i$ which 
has the ${\bf P}^1$ fiber bundle structure over $Y(\Gamma_{e_i})$.
Please consult page \pageref{lemm; sum}, lemma \ref{lemm; sum} in subsection
 \ref{subsection; realcase} for the construction.

$\bullet$ $L$--an line bundle over $M$ with first Chern class 
$C\in H^{1, 1}(M, {\bf Z})$. The bundle
 is assumed to be ``sufficiently very ample'' in this paper.

$\bullet$ $M$--an algebraic surface with irregularity $q=q(M)$ and
 geometric genus $p_g$.

$\bullet$ $M_n$--the $n-$th universal space associated with $M$. See
 section \ref{section; strata} for its construction.

$\bullet$ ${\cal M}_{C-{\bf M}(E)E}$ and 
${\cal M}_{C-{\bf M}(E)E-\sum_{e_i\cdot (C-{\bf M}(E)E)<0}e_i}$--the 
family moduli space associated with $C-{\bf M}(E)E$ and $C-{\bf M}(E)E-
\sum_{e_i\cdot (C-{\bf M}(E)E)<0}e_i$, respectively.
 They can be viewed as the sub-schemes $Z(s_{canon})$ and 
$Z(s_{canon}^{\circ})$ 
of ${\bf P}({\bf V}_{canon})={\bf P}({\bf V}_{canon}^{\circ})$.

$\bullet$ ${\bf M}(E)E$--the multiple covering of $E_i$, 
$\sum_{1\leq i\leq n}m_iE_i$ with non-increasing singular 
multiplicities $0<m_1\leq m_2\leq
\cdots\leq m_n$.

$\bullet$ $\pi$--the projection map $f_n:M_{n+1}\mapsto M_n$ or its 
restriction to the various
 ${\bf P}^1$ fibrations $\Xi_i\mapsto Y(\Gamma_{e_i})$ or the union 
 $\sum_{1\leq i}\Xi_{k_i}\mapsto Y(\Gamma)$. 

$\bullet$ $\tilde{\pi}$--the projection morphism 
$\tilde{\Xi}_i\mapsto Y(\Gamma_{e_i})$ of the
relative minimal model of $\Xi_i$.

$\bullet$ $\pi_I$--the projection map from $M_n$ to $M^{|I|}$ determined by the 
index subset $I\subset \{1, 2, \cdots, n\}$.

$\bullet$ $\check{\pi}_I$--the lifting of $\pi_I$ to $M_n\mapsto M_{|I|}$.

$\bullet$ ${\bf Q}_{k_i}$ and ${\cal Q}_{k_i}$--the line bundle and the corresponding invertible
sheaf appearing in the short exact sequence relating 
${\cal R}^1\tilde{\pi}_{\ast}\bigl({\cal O}_{\tilde{\Xi}_{k_i}}\otimes
 {\cal E}_{-{\bf M}(E)E}\bigr)$ and
 ${\cal R}^1\tilde{\pi}_{\ast}\bigl({\cal O}_{\tilde{\Xi}_{k_i}}(e_{k_i})\bigr)$. Consult
 page \pageref{prop; exact} for more details.

$\bullet$ $s_{canon}$--the canonical section of 
$\pi_{{\bf P}({\bf V}_{canon})}^{\ast}{\bf W}_{canon}\otimes {\bf H}$ induced 
 by the bundle morphism ${\bf V}_{canon}\mapsto {\bf W}_{canon}$.

$\bullet$ $s_{canon}^{\circ}$--the canonical section of 
$\pi_{{\bf P}({\bf V}_{canon})}^{\ast}{\bf W}_{canon}^{\circ}
\otimes {\bf H}$ determined by ${\bf V}_{canon}^{\circ}\mapsto 
{\bf W}_{canon}^{\circ}$.

$\bullet$ ${\cal S}_{\Gamma}$--the subset of $M_n$ over which 
the type $I$ exceptional cone ${\cal EC}_b(C-{\bf M}(E)E; Q)$
 remains constant for $b\in {\cal S}_{\Gamma}$.

$\bullet$ $T(M)$--the connected component of the Picard variety 
of $M$ parametrizing the line bundles with first Chern class 
$C$.

$\bullet$ ${\bf V}_{canon}$--the vector bundle over $M_n\times T(M)$ 
associated with the zero-th derived image sheaf
 ${\cal R}^0\pi_{\ast}\bigl({\cal E}_C\bigr)$.

$\bullet$ ${\bf V}_{canon}^{\circ}$--the vector bundle over $M_n\times T(M)$ 
associated with the zero-th derived image sheaf
 ${\cal R}^0\pi_{\ast}\bigl({\cal E}_C\bigr)$.

$\bullet$ ${\bf V}_{quot}$--the quotient bundle of ${\bf W}_{canon}$ whose 
associated
 locally free sheaf ${\cal V}_{quot}$ is constructed from the 
torsion free summand of a coherent sheaf 
${\cal R}^1\pi_{\ast}\bigl({\cal O}_{\sum_{1\leq i\leq p}\Xi_{k_i}}\otimes
 {\cal E}_{C-{\bf M}(E)E}\bigr)$. 
Consult section \ref{subsection; lf} definition \ref{defin; replace}
 for its definition, and
 proposition \ref{prop; id} 
 on page \pageref{prop; id} for its 
 construction.

$\bullet$ $\tilde{\bf V}_{quot}$--a direct sum of vector bundles which is
 equivalent to ${\bf V}_{quot}$ in the $K$ group. Consult 
\pageref{defin; replace},  definition \ref{defin; replace} for its definition.

$\bullet$ ${\bf W}_{canon}$ and ${\cal W}_{canon}$--${\bf W}_{canon}$ is 
the canonical obstruction bundle associated with
 $C-{\bf M}(E)E$. ${\cal W}_{canon}$ is the locally free sheaf associated with 
${\bf W}_{canon}$. Consult definition 5.3 of [Liu3] for its definition.

$\bullet$ ${\bf W}_{canon}^{\circ}$ and 
${\cal W}_{canon}^{\circ}$--${\bf W}_{canon}^{\circ}$ 
is the canonical obstruction bundle associated with
 $C-{\bf M}(E)E-\sum_{e_i\cdot (C-{\bf M}(E)E)<0}e_i$. 
${\cal W}_{canon}^{\circ}$ is the locally free sheaf associated 
with ${\bf W}_{canon}^{\circ}$. Consult section 5, right in front of lemma 6,
 of [Liu5] for 
 its definition.

$\bullet$ $({\cal X}/B)_n$--the $n-$th relative version of the universal 
space of ${\cal X}\mapsto B$. ${\cal X}/B$ is ${\cal X}$ with $''/B''$ to
 indicate that it has a fiber bundle structure over $B$.

$\bullet$ $Y(\Gamma)$, ${\bf Y}(\Gamma)$--$Y(\Gamma)$ is the closure of the
 admissible strata $Y_{\Gamma}\subset M_n$; 
 ${\bf Y}(\Gamma)$ is the relative version $\subset ({\cal X}/B)_n$.
 Consult section \ref{section; strata} for some of its basic properties.

$\bullet$ $Y_{\Gamma}$, ${\bf Y}_{\Gamma}$--$Y_{\Gamma}$ is the locally
 closed admissible stratum; ${\bf Y}_{\Gamma}$
 is the relative version $\subset ({\cal X}/B)_n$.
Consult section \ref{section; strata} for some of its properties.

\bigskip

$\bullet$ $\gg$--a partial ordering among pairs of the
 form $(\Gamma, \sum_{e_i\cdot (C-{\bf M}(E)E)<0}
e_i)$ with $\Gamma\in \Delta(n)$ which encodes the inclusion
 relationship of 
${\cal M}_{C-{\bf M}(E)E-\sum_{e_i\cdot (C-{\bf M}(E)E)<0}e_i}\times_{M_n}
Y(\Gamma)$. Please consult page \pageref{defin; partial}
 for more details.

$\bullet$ $\sqsupset$--the partial ordering among pairs of the
 form $(\Gamma, \sum_{e_i\cdot (C-{\bf M}(E)E)<0}
e_i)$ with $\Gamma\in \Delta(n)$ which encodes the discrepancy of using
${\bf W}_{canon}^{\circ}$ and $s_{canon}^{\circ}$ to replace
 ${\bf W}_{canon}$ and $s_{canon}$. Please consult page \pageref{defin; sq}
 for more details.

$\bullet$ $\succ$--a partial ordering among $\Gamma\in \Delta(n)$ which 
encodes that
 the type $I$ exceptional cone 
${\cal C}_{\Gamma}$ gets larger under degenerations. Consult
 page \pageref{defin; succ} for more details.

$\bullet$ $\models$--The linear ordering introduced on $\Delta(n)$ and 
therefore on $I_{\Gamma}$ and $\bar{I}_{\Gamma}$. Consult page \pageref{models}
 for more details.

$\bullet$ $\vdash$--the altered linear ordering on $\bar{I}_{\Gamma}$.
 Consult definition \ref{defin; newordering} for more details.

\medskip

\section{A Brief Review about the Admissible Graphs and 
 Admissible Strata} \label{section; strata}

\bigskip

 In [Liu1] we had introduced the concepts of universal spaces (see also
 [V]), admissible
 graphs and  the admissible stratification of the universal spaces $M_n$, 
$n\in {\bf N}$.
  For the convenience of the reader, we extract the basic facts about them
 in this section.

We review the construction of $M_n$ and review the admissible graphs and the
 admissible stratification. Then we generalize it to a relative setting and
 discuss their basic properties and the relationship with type $I$ exceptional classes.
 
 Recall (consult section 3 on page 400 of [Liu1]) that the 
universal space $M_n$ is constructed by an inductive
 procedure. Take $M_0=pt$ and $M_1=M$. Suppose that $M_0, M_1, M_2, \cdots, 
M_{k-1}$ has been constructed and there are natural projection maps 
 $M_{k-1}\mapsto M_{k-2}\mapsto M_{k-3}\cdots \mapsto pt$, then define $M_k$
 to be the blowing up of the relative diagonal 
$\Delta_{M_{k-1}/M_{k-2}}:M_{k-1}\mapsto 
M_{k-1}\times_{M_{k-2}}M_{k-1}$. Then the natural projection map 
$f_{k-1}: M_k\mapsto M_{k-1}$ is a surjection. 
 By mathematical induction, the universal spaces are defined for all $n\in {\bf N}$ and
 there are smooth surjective morphisms $f_k:M_{k+1}\mapsto M_k$ for all $k$.
 As usual the composite map $f_k\circ f_{k+1}\circ\cdots f_n:M_{n+1}\mapsto M_k$
 will be denoted by $f_{n, k}$.

 In [Liu1] we had introduced a concept called admissible graphs. The set of 
$n$-vertex
 admissible graphs, denoted by $adm(n)$, is the set of finite graphs with $n$ 
vertexes
 and a finite number of arrowed one-edges which satisfy five axioms (from
 page 412-413 of
 [Liu1]).

\medskip

\noindent Axiom 1: There is a $1-1$ correspondence between the vertexes 
 of $\Gamma$ and the
 positive 
integers smaller or equal to l. An association of this type is called a
 marking of the graph. More generally, one can mark the graph by any
 finite subset of ${\bf N}$. If ${\bf I}$ is the index set. The graph
 is called ${\bf I}$ admissible. \label{graphaxiom}

\medskip

\noindent Axiom 2:
The one-edges are oriented by arrows from the vertex marked by a smaller integer
 (called a direct ascendent) to the vertex marked by a larger integer (called a direct
 descendent).

\medskip

\noindent Axiom 3:
  The only loops allowed in the graph are triangles
 formed by the three vertexes. Suppose $a<b<c$ are the three different
 vertexes, then $b$, $c$ must be the direct descendents of $a$ while $a$, $b$ must
 be the direct ascendents of $c$. The vertexes $a$, $b$, $c$ form a triangle. 

\medskip

\noindent Axiom 4: Any vertex can have at most two direct ascendents.
 When a vertex has exactly two direct ascendents, it and its two direct 
ascendents form a triangular loop.

\medskip

\noindent Axiom 5:
 Suppose that two adjacent triangles share a common edge, then out of the
 four vertexes in the two triangles
 the ending vertex \footnote{at which the arrow points to}
 of this common edge has exactly one direct descendent among the other
 three vertexes.  
 
\medskip

  The admissible graphs code the combinatorial patterns of the blowing ups.
  We usually denote a typical element in $adm(n)$ by $\Gamma$. 
  The consistency of the above axioms were ensured geometrically in [Liu1]
 by constructing the corresponding 
admissible strata \footnote{Consult [Liu1] or below for the
 definition or construction of $Y_{\Gamma}$ (${\bf Y}_{\Gamma}$).} 
$Y_{\Gamma}\not=\emptyset$ explicitly. 

\begin{defin}\label{defin; codim}
Define the codimension of an admissible graph $\Gamma\in adm(n)$, $codim_{\bf C}\Gamma$,
 to be the number of one-edges in $\Gamma$.
\end{defin}

 The graph in fig.1
 illustrates an example of admissible graphs. 
 
\begin{figure}
\centerline{\epsfig{file=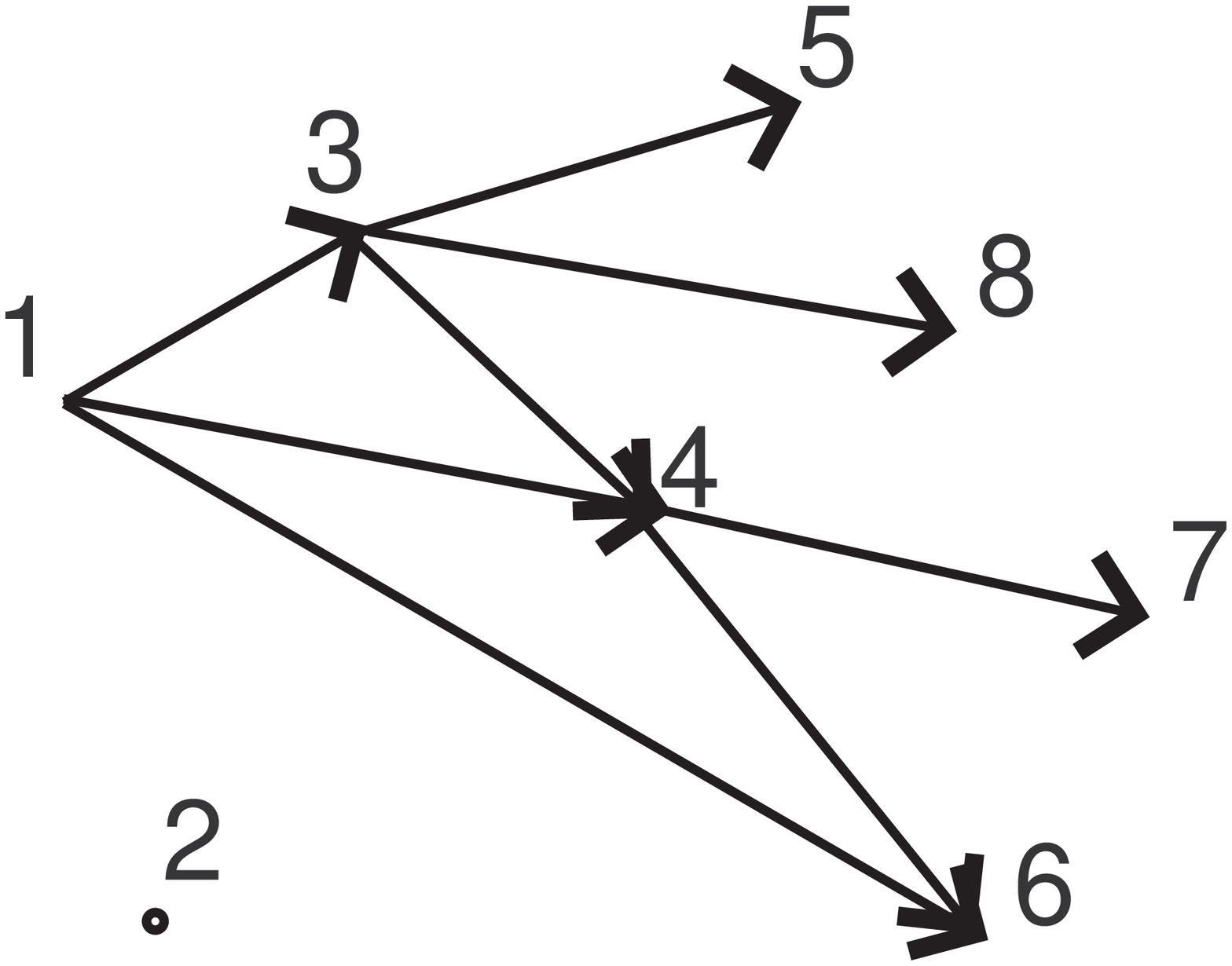,height=4cm}}\label{fig.1}
\centerline{fig.1}
\centerline{An admissible graph with $8$ vertexes.}
\end{figure}

In the set $adm(n)$ there is a special
element $\gamma_n$ (or for simplicity 
skipping the subscript $n$ and denote it by $\gamma$ if it does not cause confusion) 
which consists of $n$ free vertexes
 without one-edges. By construction $codim_{\bf C}\gamma_n=0$ and $\gamma_n$ is the only
 admissible graph in $adm(n)$ which has this special property.

The reason to introduce such a set $adm(n)$ is because that the universal space
 $M_n$ can be stratified by the various 
admissible strata $Y_{\Gamma}$ which have smooth closure $Y(\Gamma)$
 in $M_n$,

\begin{prop}\label{prop; stratification}
 The $n-$th universal space $M_n$ admits an admissible stratification 
$M_n=\coprod_{\Gamma\in adm(n)}Y_{\Gamma}$ into locally closed smooth
 subsets $Y_{\Gamma}$ such that 

(i). The closure of $Y_{\Gamma}$ in $M_n$, $\overline{Y_{\Gamma}}=Y(\Gamma)$ is smooth
 of dimension $dim_{\bf C}M_n-codim_{\bf C}\Gamma$ in $M_n$.

(ii). $Y(\Gamma)$ can be expressed as a union of admissible strata $Y_{\Gamma'}$, 
 $Y(\Gamma)=\coprod_{\Gamma'<\Gamma}Y_{\Gamma'}$.
\end{prop}

  $\Gamma'<\Gamma$ indicates that $Y_{\Gamma'}$ appears in the compactification of 
$Y_{\Gamma}$ into $Y(\Gamma)$ and is said to be a degeneration of the admissible
 graph $\Gamma$. 

 This proposition was proved in proposition 4.2.-4.3. in [Liu1]. In definition 4.8. of
 [Liu1] we had given a combinatorial characterization of the degenerations of
 the admissible graphs.

\begin{rem}\label{rem; interpret}
The fibers of the fiber bundle map $f_n:M_{n+1}\mapsto M_n$ are all $n-$consecutive 
blowing ups from $M$. Therefore the space $M_n$ can be interpreted as the 
`` universal space'' 
parametrizing all the ordered $n-$consecutive pointwise blowing ups from $M$.
\end{rem}

 Let $M_n=\coprod_{\Gamma\in adm(n)}Y_{\Gamma}$ be the admissible stratification of
 $M_n$ into locally closed strata. 
Then the stratum $Y_{\gamma_n}$ is the only top dimensional stratum which 
parametrizes all the ordered distinct $n$ points in $M$. Every
 distinct $n$ points in $M$
 corresponds to $n$ distinct pointwise blowing ups of $M$. The various $Y_{\Gamma}$, 
$\Gamma\not=\gamma_n$, parametrize those $n-$consecutive blowing ups whose blowing up centers
 may lie on the exceptional loci of the previous blowing ups.

Let ${\cal X}\mapsto B$ be a fiber bundle over a base $B$ which
  is smooth of relative dimension two. We use the notation ${\cal X}/B$ to indicate that
 the space ${\cal X}$ has a structure of fiber bundle over $B$.
We define $({\cal X}/B)_0=B$, 
$({\cal X}/B)_1={\cal X}/B$. Suppose that $({\cal X}/B)_1, ({\cal X}/B)_2, \cdots, 
({\cal X}/B)_{k-1}$ have
 been defined and there are natural surjective 
projection maps $({\cal X}/B)_{k-1}\mapsto ({\cal X}/B)_{k-2}
\mapsto \cdots\mapsto B$, define $({\cal X}/B)_k$ to be the blowing up of the relative
 diagonal $\Delta_{({\cal X}/B)_{k-1}/({\cal X}/B)_{k-2}}:({\cal X}/B)_{k-1}\hookrightarrow 
({\cal X}/B)_{k-1}\times_{({\cal X}/B)_{k-2}}({\cal X}/B)_{k-1}$, etc.
 Apparently ${\bf f}_{k-1}:({\cal X}/B)_k\mapsto ({\cal X}/B)_{k-1}$ is surjective.
 By mathematical induction, $({\cal X}/B)_n$ are constructed for all $n\in {\bf N}$ such that
 ${\bf f}_n:({\cal X}/B)_{n+1}\mapsto ({\cal X}/B)_n$ are smooth and surjective.
 
 For a given $n$, ${\bf f}_{n-1, 0}:({\cal X}/B)_n\mapsto B$ is the fiber bundle of the
 $n-$th universal spaces of the fibers.

 We have the following lemma relating different relative universal space constructions, 

\begin{lemm}\label{lemm; fuse}
 For all $n\geq k$, we have the following identity
$$({\cal X}/B)_n/({\cal X}/B)_k=\bigl(({\cal X}/B)_{k+1}/({\cal X}/B)_k\bigr)_{n-k}.$$
\end{lemm}

\noindent Proof: This can be seen by noticing that for $n=k+1$, 
$({\cal X}/B)_{k+1}/({\cal X}/B)_k$
is a fiber bundle projection map. So we may take ${\cal X}'=({\cal X}/B)_{k+1}$ and 
$B'=({\cal X}/B)_k$ and $({\cal X}/B)_{k+1}/({\cal X}/B)_k=({\cal X}'/B')_1$. 
Then by definition 
$({\cal X}/B)_{k+2}\mapsto ({\cal X}/B)_k$ is the blowing up of the 
 relative diagonal of 
$({\cal X}/B)_{k+1}\times_{({\cal X}/B)_k}({\cal X}/B)_{k+1}=
({\cal X}'/B')_1\times_{B'}
 ({\cal X}'/B')_1$. 
So we identify $({\cal X}/B)_{k+2}/({\cal X}/B)_k$ with $({\cal X}'/B')_2$. By a simple
 induction argument and by comparing with the above relative construction of the
 universal spaces, we
 find $({\cal X}/B)_n/({\cal X}/B)_k=({\cal X}'/B')_{n-k}$ for all $n\geq k$. 
Therefore we have the following identity,

$$({\cal X}/B)_n/({\cal X}/B)_k=\bigl(({\cal X}/B)_{k+1}/({\cal X}/B)_k\bigr)_{n-k}.$$
 $\Box$

  If we ignore the base space on the left hand side, we may rewrite the identity as
 $({\cal X}/B)_n=\bigl(({\cal X}/B)_{k+1}/({\cal X}/B)_k\bigr)_{n-k}$.
 By taking $B=pt$ and ${\cal X}=M$, we recover the important special case 
 $M_n=(M_{k+1}/M_k)_{n-k}$.

\medskip

\begin{lemm}\label{lemm; times}
 For all $n\in {\bf N}$, there exists a canonical dominated birational map
$({\cal X}/B)_n\mapsto \times_B^n({\cal X}/B)_1$ from the $n-th$
 relative universal space $({\cal X}/B)_n$ to the $n-th$ fiber products of ${\cal X}/B$. 
\end{lemm}

\noindent Proof: The assertion is apparently true for $n=1$. Suppose that
 the birational map 
$({\cal X}/B)_k\mapsto \times_B^k ({\cal X}/B)_1$ has been constructed, then
 
$$({\cal X}/B)_{k+1}\mapsto ({\cal X}/B)_k\times_{({\cal X}/B)_{k-1}}({\cal X}/B)_k
\mapsto \times_B^k({\cal X}/B)_1\times_{\times_B^{k-1}({\cal X}/B)_1}\times_B^k({\cal X}/B)_1
\cong \times_B^{k+1}({\cal X}/B)_1$$

 is a composition of dominated birational maps. By mathematical induction, the
 lemma is proved. $\Box$

\medskip

   The following proposition and its corollary are about the liftings
 of the projection maps
 to the corresponding (relative) universal spaces.

\begin{prop}\label{prop; lifting}
Let $n$ be a positive integer and let 
$I\subset \{1, 2, \cdots, n\}$ be an index subset. Let $\pi_i:({\cal X}/B)_n\mapsto ({\cal X}/B)$
 be the composite projection map $({\cal X}/B)_n\mapsto \times_B^n({\cal X}/B)\mapsto {\cal X}/B$ 
to the $i-$th direct factor and let $\pi_I=\times_{i\in I}\pi_i:
({\cal X}/B)_n\mapsto \times_B^{|I|}({\cal X}/B)$ be the projection to the fiber products of
 $|I|$ copies of ${\cal X}/B$ indexed by the subset $I$. Then there exists a natural lifting of the map
 $\pi_I$ to $({\cal X}/B)_n\mapsto ({\cal X}/B)_{|I|}$ which makes the following diagram commutative,

\[
 \begin{array}{ccc}
  ({\cal X}/B)_n & \stackrel{\check{\pi}_I}{\longrightarrow} & ({\cal X}/B)_{|I|} \\
  & \searrow &  \Big\downarrow \\
  &   &  \times_B^{|I|}({\cal X}/B)
 \end{array}
\]
 The lifted map $\check{\pi}_I$ is smooth of relative dimension $2n-2|I|$.
\end{prop}

 Notice that the lifted map $\check{\pi_I}$ is usually different 
from the composite projection map
 $({\cal X}/B)_n\mapsto ({\cal X}/B)_{n-1}\mapsto \cdots \mapsto ({\cal X}/B)_{|I|}$, which
 corresponds to the lifting of a very specific $I$.

\begin{rem}\label{rem; forget}
Following remark \ref{rem; interpret}, the lifting map $\check{\pi}_I$ can be
 interpreted as the ``forgetful'' map of forgetting all the blowing ups among
the sequence of $n$-consecutive blowing ups
 marked by indexes in $\{1, 2, \cdots, n\}-I$.
\end{rem}

\medskip

\noindent Proof of the proposition: When $n=1$, the statement holds trivially.

 Let $2\leq n\in {\bf N}$ be a positive integer such that the statement of 
the proposition is known
 to be true for $n-1$, for all the ${\cal X}\mapsto B$ pairs. We would 
like to show that it
 holds for $n$ as well.
 Suppose that $I=\{1, 2, \cdots, n\}$, the statement holds since the lifted map 
 $({\cal X}/B)_n\mapsto ({\cal X}/B)_{|I|}$ is the identity map. So let us 
assume that 
 $|I|<n$. By lemma \ref{lemm; fuse} 
we may rewrite $({\cal X}/B)_n$ as $(({\cal X}/B)_2/({\cal X}/B)_1)_{n-1}$.
 
 $\diamondsuit$ We separate into two cases. (i). $1\not\in I$. (ii). $1\in I$.

 In the first case, the projection to the second factor 
$\pi_2:({\cal X}/B)_2\mapsto ({\cal X}/B)_1$ induces a morphism 
 $({\cal X}/B)_2/({\cal X}/B)_1\mapsto ({\cal X}/B)_1/pt$. This morphism induces an
 morphism on
 the relative $n-1-$th universal spaces
 $(({\cal X}/B)_2/({\cal X}/B)_1)_{n-1}\mapsto (({\cal X}/B)_1/pt)_{n-1}\cong 
({\cal X}/B)_{n-1}$.
  By the assumption $1\not\in I$ we know that $I\subset \{2, 3, \cdots, n\}$.
 Define a new index set $I_{-1}\subset \{1, 2, \cdots, n-1\}$ by subtracting
 $1$ from all the elements of $I$. Then by the induction hypothesis,
 $\pi_{I_{-1}}:({\cal X}/B)_{n-1}\mapsto \times_{B}^{|I_{-1}|}({\cal X}/B)$ can be 
lifted to the
 smooth surjective morphism
 $\check{\pi}_{I_{-1}}:({\cal X}/B)_{n-1}\mapsto ({\cal X}/B)_{|I_{-1}|}$. By
 composing with 

$$\hskip -.6in
({\cal X}/B)_n\cong(({\cal X}/B)_2/({\cal X}/B)_1)_{n-1}\mapsto (({\cal X}/B)_1/pt)_{n-1}
\cong ({\cal X}/B)_{n-1},$$
  
  we get the desired lifting from 
$({\cal X}/B)_n\cong (({\cal X}/B)_2/({\cal X}/B))_{n-1}\mapsto 
\times_B^{|I_{-1}|}({\cal X}/B)=\times_B^{|I|}({\cal X}/B)$
 to $({\cal X}/B)_n\mapsto ({\cal X}/B)_{|I_{-1}|}=({\cal X}/B)_{|I|}$.
 The lifted surjective map is apparently smooth because all the composite factors of maps are 
 smooth and surjective.

 In the second case when $1\in I$,
 we construct an new index set $I'_{-1}$ by subtracting all elements in 
$I-\{1\}$ by $1$, then $I'_{-1}\subset \{1, 2, \cdots, n-1\}$ and we have $|I'_{-1}|=|I|-1$.
 We define a new fiber bundle ${\cal X}'\mapsto B'$  smooth of relative dimension
 two by setting ${\cal X}'=({\cal X}/B)_2, B'=({\cal X}/B)_1$. Then we may rewrite 
$({\cal X}/B)_n$ as $(({\cal X}/B)_2/({\cal X}/B)_1)_{n-1}=
({\cal X}'/B')_{n-1}$. By the induction
 hypothesis, the map $\pi'_{I'_{-1}}:({\cal X}'/B')_{n-1}\mapsto 
\times_{B'}^{|I'_{-1}|}({\cal X}'/B')$ can
 be lifted to $\check{\pi}'_{I'_{-1}}:({\cal X}'/B')_{n-1}\mapsto 
({\cal X}'/B')_{|I'_{-1}|}$. 

 On the other hand $|I'_{-1}|=|I|-1$, we realize by using lemma \ref{lemm; fuse} again that 
$$({\cal X}'/B')_{|I'_{-1}|}=
 ({\cal X}'/B')_{|I|-1}=(({\cal X}/B)_2/({\cal X}/B)_1)_{|I|-1}\cong 
({\cal X}/B)_{|I|-1+1}=({\cal X}/B)_{|I|}.$$

 Then the lifted maps from $({\cal X}/B)_n=({\cal X}'/B')_{n-1}$ to 
$({\cal X}'/B')_{|I|-1}=({\cal X}/B)_{|I|}$ is the
 desired lifting map $\check{\pi}_I$. By inductive hypothesis, it is smooth and surjective.

 Based on mathematical induction, the existence of the lifting is proved. 
 Finally the assertion about the relative dimension $2n-2|I|$ is by a direct comparison of
 the dimensions of the source and the target. 
$\Box$

\begin{cor}\label{cor; lifting}
Let $\pi_i:M_n\mapsto M, 1\leq i\leq n$ 
denote the projection to the $i-$th copy of $M$.
Let $\pi_I=\times_{i\in I}\pi_i:M_n\mapsto M^{|I|}$ is the projection to
 the $|I|$-Cartesian product of $M$ indexed by $I$. 

Then there exists a unique surjective and smooth lifting $\check{\pi}_I:
 M_n\mapsto M_{|I|}$ which makes the
 following diagram commutative, 

\[
\begin{array}{ccc}
  M_n & \stackrel{\check{\pi}_I}{\longrightarrow} & M_{|I|} \\
   &  \searrow &  \Big\downarrow \\
   &  &  M^{|I|}
\end{array}
\]
 
The map $\check{\pi}_I$ is of relative dimension $2n-2|I|$.
\end{cor}

\noindent Proof: By taking $M={\cal X}\mapsto B=pt$ 
in the proposition \ref{prop; lifting}, 
 the corollary is a direct consequence of proposition \ref{prop; lifting}. $\Box$ 

\medskip

 Given an $n\in {\bf N}$, one may prove inductively (see lemma 3.1 and proposition 3.1 
on pages 401-402 of [Liu1]) 
that the birational map $M_{n+1}\mapsto M\times M_n$ can be 
factorized into $n$ codimension-two blowing ups along the cross sections of the
 intermediate fiber bundles $f_{n-1, i}^{\ast}M_{i+1}\mapsto M_n$ induced 
 by the relative diagonals $M_{i+1}\hookrightarrow M_{i+1}\times_{M_i}M_{i+1}$. As 
 a consequence, $M_n$ can be blown up from $M^n$ by ${n(n-1)\over 2}$ 
consecutive codimension-two blowing ups along the partial diagonals. Our
 convention in this paper is that $E_{a; b}, 1\leq a<b\leq n$ denote the (pull-back of)
the exceptional divisor blown up from the strict transforms of the 
$(a, b)-$th partial diagonals. Under
 this convention, the $n$ distinct exceptional divisors, $E_{1; n+1}$, $E_{2; n+1}$,
 $\cdots, E_{n; n+1}$ of the fiber bundle $f_n: M_{n+1}\mapsto M_n$ 
are denoted by $E_1, E_2, \cdots , E_n$ respectively. They will play a special
 role in this paper.
 
\medskip

  Let $\Gamma$ be an admissible graph $\in adm(n)$, then one may attach 
 $n$ distinct type $I$ exceptional classes $e_i$, $1\leq i\leq n$, 
to $\Gamma$. Given an index $i$, with $1\leq i\leq n$,
 let the indexes $j_i$ run through all the direct descendents of $i$ in
 $\Gamma$. Then $e_i=E_i-\sum_{j_i}E_{j_i}$ is the $i-$th type $I$
 exceptional class attached to $\Gamma$. We set $e_i\cdot e_j$ to be
 the fiberwise intersection number of the classes $e_i$ and $e_j$.

 The following proposition will be used frequently in this paper,

\begin{prop}\label{prop; reverse}
 For $1\leq i\leq n$ let $J_i$ be an index subset of $\{1, 2, \cdots, n\}$
 satisfying $inf(J_i)>i$.
 Let $e_1, e_2, \cdots, e_n$ be $n$ divisor classes of the form
 $E_i-\sum_{j\in J_i}E_j$. Suppose that $e_1, e_2, \cdots, e_n$ satisfy
 the condition $e_a\cdot e_b\geq 0$ for all $a\not=b$. Then there exists
 an admissible graph $\Gamma\in adm(n)$ such that $e_1, e_2, e_3, \cdots, e_n$ 
 are the type $I$ exceptional classes associated with $\Gamma$.
\end{prop}

 In other words,
 the locus $Y(\Gamma)\subset M_n$ is the locus of co-existence over which
 $e_1, e_2, \cdots, e_n$ become effective. Please see fig.3 on page \pageref{fig.3}
 for an example that $\Gamma$ is recovered from the fan-like graphs associated with
 these $e_i$.

\noindent Proof: The proposition is proved by an induction 
argument on $n$. The base case $n=1$ is trivial. Suppose that
 the proposition has been proved for $n-1$, we would like to prove the
 existence of such a $\Gamma\in adm(n)$ for $n$. Given $n$ vertexes, construct
 the graph $\Gamma$ by the following rule: Given an $i\leq n$, connect
 an oriented edge from $i$ to any $j>i$ if $j\in J_i$, i.e. if the term 
$-E_j$ appears in the class $e_i$.  We show that $\Gamma$ is an admissible graph, i.e.
 it satisfies the five axioms characterizing admissible graphs.

 Firstly we shift all the indexes by $-1$ temporally 
and denote the new graph marked by the shifted indexes $\{0, 1, 2, \cdots, n-1\}$ as
 $\tilde{\Gamma}$. It is easy to see that $\Gamma$ is admissible
 with respect to $\{1, 2, \cdots, n\}$ iff $\tilde{\Gamma}$ is admissible 
\footnote{relaxing
 the constraint on the index set.}
with respect to the shifted index set $\{0, 1, \cdots, n-1\}$. 
Define $\phi:{\bf Z}\mapsto {\bf Z}$
 by the formula $\phi(i)=i-1$. Define $\tilde{J}_{i}=\phi(J_{i+1})$ for $i\in \{0, 
1, 2, \cdots, n-1\}$. Likewise define $\tilde{E}_i=E_{i+1}$.
  Define accordingly $\tilde{e}_i=\tilde{E}_i-\sum_{j\in \tilde{J}_i}
\tilde{E}_{j}$. It is clear that $\tilde{e}_i=e_{i+1}$ for $0\leq i\leq n-1$
 and their mutual intersection pairings are still non-negative.

 Thus by our induction hypothesis, the classes
\footnote{Notice that the class $\tilde{e}_0$ is excluded here.}
 $\tilde{e}_1, \tilde{e}_2, \cdots, 
\tilde{e}_{n-1}$ satisfy $\tilde{e}_i\cdot \tilde{e}_j=e_{i+1}\cdot e_{j+1}\geq 0$
 for $i\not= j$.
 So there exists an admissible graph $\Gamma'\in adm(n-1)$ such that
 $\tilde{e}_1, \cdots, \tilde{e}_{n-1}$ are the type $I$ exceptional classes
 associated to $\Gamma'$. And by the inductive construction $\Gamma'$ is
 constructed by the datum of $\tilde{e}_1, \cdots, \tilde{e}_{n-1}$ the same
 way we construct $\Gamma$ (and $\tilde{\Gamma}$). So it is clear that
 $\Gamma'$ is a sub-graph of $\tilde{\Gamma}$ by removing the $0-$th vertex
 and all the arrowed edges starting from it. That is because $\tilde{e}_0$ is not used
 in constructing $\Gamma'$. Our inductive assumption
 implies that this sub-graph $\Gamma'$ of $\tilde{\Gamma}$ 
is admissible (with respect to $\{1, 2, \cdots, n-1\}$). 

 Our final task is to show the admissibility of the whole $\tilde{\Gamma}$.
 Among the five axioms which characterize admissible graphs,
 the axiom 1. and axiom 2. are satisfied trivially by construction. 

 To show that axioms
 3., 4., 5., are satisfied, we prove by contradiction.
 
$\diamondsuit$ About axiom 3.: if there is a polygonal loop in $\tilde{\Gamma}$ which 
is not a triangle, we may choose a polygonal loop of vertexes involving
 the least number of vertexes/edges.  Firstly, this 'exotic loop' cannot locate completely 
 within $\Gamma'$ because $\Gamma'$ has been known to be admissible, so
 axiom 3 for $\Gamma'$ rules out this possibility. So the vertex
marked by $0$ must be within this loop. Consider the vertex $v$ marked by the
 largest integer $q$ along the loop. It must be located in $\Gamma'$ because
 $q\not=0$. On the
 other hand, the loop passes through $v$ means that 
there are two edges ending at $v$ (because
 the arrows of the edges always point to vertexes marked by the larger integers).

 If both direct ascendents $v_1, v_2$ of $v$
are not marked by $0$, they are vertexes in $\Gamma'$ as well
 and by the admissibility of $\Gamma'$, axiom 4. for $\Gamma'$ 
implies that $v$, $v_1, v_2$  
 form a triangle. This implies that we can shorten the loop by replacing
 the oriented edges $\overrightarrow{v_1v}$ and $\overrightarrow{v_2v}$ by the
 single edge $\overrightarrow{v_1v_2}$ (or $\overrightarrow{v_2v_1}$, depending on
 which vertex is marked by a larger integer). 
 This violates the assumption that the loop involves the least number of
 vertexes/edges!

If one of the two direct ascendents $v_1, v_2$ 
of $v$, say $v_1$, is marked by $0$, assume that $v_2$ is marked by $p$, with $1\leq 
p\leq n-1$. Then $q\in \tilde{J}_p$ and $q\in \tilde{J}_0$ simultaneously. 
 
 On the other hand, we have assumed that $\tilde{e}_0\cdot \tilde{e}_p\geq 0$. As 
 $(-E_q)^2$ contributes $-1$ to the intersection number, there must
 be a positive counter-term which makes the intersection number non-negative.
 This only occurs when $p\in \tilde{J}_0$ and $(E_p)\cdot (-E_q)=1$ contributes
 positively to the sum. But this implies that the edges linking $v_1, v_2, v$ already
 form a triangle and it is not a polygonal loop. 

 In any case, the axiom 3. holds for $\tilde{\Gamma}$.
 
 $\diamondsuit$ About axiom 4., suppose that a vertex in $\tilde{\Gamma}$ has more than
 two direct ascendents in $\tilde{\Gamma}$. Then we randomly pick three of them, say 
 $v_1, v_2, v_3$ (marked by $p_1<p_2<p_3$, respectively), and derive a contradiction. 

Because $q\in \tilde{J}_{p_1}\cap \tilde{J}_{p_2}\cap \tilde{J}_{p_3}$, 
 by the same argument in checking axiom 3., $\tilde{e}_{p_1}\cdot 
\tilde{e}_{p_2}\geq 0$, $\tilde{e}_{p_1}\cdot \tilde{e}_{p_3}\geq 0$,
 and $\tilde{e}_{p_2}\cdot \tilde{e}_{p_3}\geq 0$ force 
 $v_2, v_3$ to be the direct descendents of $v_1$ while $v_3$ is forced to be a direct
 descendent of $v_2$. Having made such an observation, we recalculate
 $\tilde{e}_{p_1}\cdot\tilde{e}_{p_2}$ again. Because $p_2\in \tilde{J}_{p_1}$, but
 $\{p_3, q\}\subset \tilde{J}_{p_1}\cap \tilde{J}_{p_2}$. This implies that
 $$0\leq 
\tilde{e}_{p_1}\cdot\tilde{e}_{p_2}\leq (-E_{p_2})\cdot E_{p_2}+(-E_{p_3})^2
+(-E_q)^2=-1<0.$$  

 This is absurd! Thus $v$ can have at most two direct descendents in $\tilde{\Gamma}$.

 When $v$ has exactly two direct ascendents $v_1, v_2$ marked by $p_1<p_2$,
 $\tilde{e}_{p_1}\cdot \tilde{p_2}\geq 0$ and $q\in \tilde{J}_{p_1}\cap
\tilde{J}_{p_2}$ imply that $v_2$ is also a direct descendent of $v_1$ and
 $v_1, v_2, v$ form a triangle in $\tilde{\Gamma}$. So the axiom 4. is
 satisfied.

  $\diamondsuit$ Finally about axiom 5.: Suppose that two adjacent triangles are sharing
 a common one-edge. Let $v_1$ and $v$ be the starting and the ending
 vertexes of the one-edge. Suppose that $v_2$, $v_3$ are the other
 two vertexes in these two triangles, 
 let $p_1, p_2, p_3$ and $q$ mark the vertexes $v_1$, $v_2$, $v_3$ and $v$,
 respectively. Because $v$ are in both of the triangles, there must be two different
 one-edges linking $v$ and $v_2$, $v$ and $v_3$, respectively.
 There are three exclusive possibilities. Either

 (i). both $v_2, v_3$ are the direct ascendents of $v$.

 (ii). One of them, say $v_2$, is the direct ascendent of $v$ and the other vertex
 $v_3$ is the direct descendent of $v$.

(iii). Both of $v_2$ and $v_3$ are direct descendents of $v$.

 If the possibility (i) holds, then $v$ has at least three direct ascendents
 in $\tilde{\Gamma}$ and this violates axiom 4 for $\tilde{\Gamma}$, which we have proved already.

 If the possibility (iii). holds, then we argue $v_2, v_3$ are the direct
 descendent of $v_1$ as well. 

By our assumption on the existence of
 adjacent triangles, there must be one-edges between $v_1, v_2$ and $v_1, v_3$.
But $v$ is already known to be a direct descendent of $v_1$. If $v_2, v_3$
 are known to be the direct descendents of $v$, then $v_1$ cannot be
 a direct descendent of either $v_2$ or $v_3$. So the arrows of the oriented one-edges
 must go from $v_1$ to $v_2$ and $v_3$, respectively. Then $v_2, v_3$
 must be the direct descendents of $v_1$.

 Now the vertex $v_1$ has at least three direct descendents $v_2, v_3$ and $v$
 while the vertex $v$ has at least two direct descendents $v_2, v_3$.

 This implies that $\{p_2, p_3\}\subset \tilde{J}_{p_1}\cap \tilde{J}_q$ and
 a direct calculation on $\tilde{e}_{p_1}\cdot \tilde{e}_{q}$ shows that

$$0\leq \tilde{e}_{p_1}\cdot \tilde{e}_{q}\leq (-\tilde{E}_q)\cdot (\tilde{E}_q)
+(-\tilde{E}_{p_2})^2+(-\tilde{E}_{p_3})^2=1-1-1=-1<0.$$

 A contradiction to our assumption!  So the only possibility is $(ii)$ and exactly one of 
$v_2$ or $v_3$ can be the direct descendent of $v$.

 As $\tilde{\Gamma}$ satisfies all five axioms, it is admissible
 with respect to the marking $\{0, 1, 2, \cdots, n-1\}$. Thus the original $\Gamma$
 is admissible with respect to $\{1, 2, \cdots, n\}$. So $\Gamma\in adm(n)$.
$\Box$

\medskip

 Conversely for all $\Gamma\in adm(n)$ the smooth and closed set $Y(\Gamma)\subset 
 M_n$ can be identified to be the transversal intersection 
$\cap_{1\leq i\leq n}Y(\Gamma_{e_i})$, where $\Gamma_{e_i}$ is the
 fan-like admissible graph $\in adm(n)$ such that (i). the vertex marked by
 $i$ is the only direct ascendent among the $n$ vertexes. (ii). The direct
descendents of the vertex marked by $i$ are the direct descendent indexes of $i$
 in $\Gamma$.  Thus $Y(\Gamma)$ can be viewed as the locus over which all the 
 type $I$ exceptional classes $e_i, 1\leq i\leq n$ are simultaneously 
 effective along the fibers of $M_{n+1}\mapsto M_n$. The result has been
 proved in proposition 4.7. of [Liu1], using slightly different
terminologies in terms of pseudo-holomorphic curves.  

 The graphs in fig.2 
are the fan-like sub-graphs from fig.1 on page \pageref{fig.1}. 

\begin{figure}
\centerline{\epsfig{file=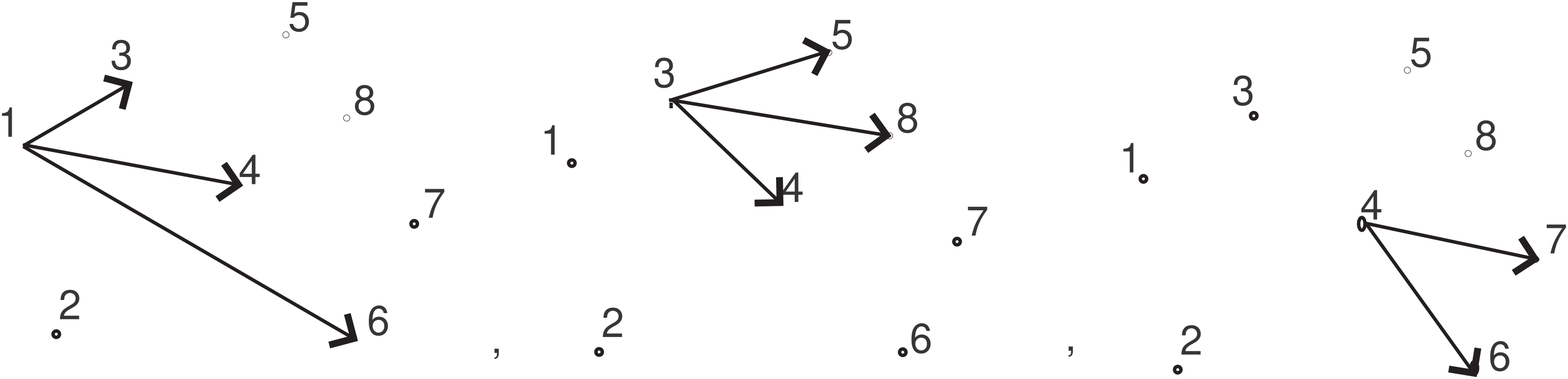,height=4cm}}
\centerline{fig.2}
\centerline{The fan-like subgraphs $\Gamma_{e_1}$, $\Gamma_{e_3}$ and
 $\Gamma_{e_4}$ of the admissible graph $\Gamma$ in fig.1.}\label{fig.2}
\end{figure}

 For the ease of the reader with algebraic geometric background, we
 give an alternative relative construction of $Y(\Gamma)$
 which makes the above property manifest.

 Let $({\cal X}/B)_n$ be the $n-$th relative universal space over $B$.
 The fiber of $({\cal X}/B)_n$ over $b\in B$ is nothing but the $n-$th 
universal space of the fiber of ${\cal X}/B$ above $b$. 

By lemma \ref{lemm; fuse} we have the following canonical 
isomorphism $\bigl(({\cal X}/B)_{i+1}/({\cal X}/B)_i\bigr)_{n-i}=
 ({\cal X}/B)_n$ for each $i\leq n$. Thus we may 
set $B'_i=({\cal X}/B)_i$ for all $i$ and there is a surjection
 $({\cal X}/B)_n\mapsto \times_{B'_i}^{n-i}
\bigl(({\cal X}/B)_{i+1}/B'_i\bigr)$ to the
 $n-i-$fold fiber product of $({\cal X}/B)_{i+1}$ over $B'_i$. 

Parallel to the absolute version, for all $\Gamma\in adm(n)$ 
we may define the relative
 admissible strata ${\bf Y}_{\Gamma}$ (or the closure ${\bf Y}(\Gamma)$) 
\footnote{Here we use bold ${\bf Y}$ to denote the
 relative versions of $Y_{\Gamma}$ or $Y(\Gamma)$}
 to be the union of fiberwise $Y_{\Gamma}$ (or $Y(\Gamma)$)
 over\footnote{See remark \ref{rem; alternative} 
for an outline of an alternative inductive 
definition after proposition \ref{prop; transversal}} $B$.

By the previous 
inductive construction on the relative universal spaces, $({\cal X}/B)_{i+1}/B'_i$
 is the blowing up of the relative diagonal 
$\Delta_{B'_i/B'_{i-1}}:B'_i\mapsto B'_i\times_{B'_{i-1}}B'_i$. Let $D_i\mapsto B'_i$ with 
$D_i\subset ({\cal X}/B)_{i+1}$
 denote the blown up exceptional divisor in $({\cal X}/B)_{i+1}$, which has a structure of 
${\bf P}^1$ bundle over $B'_i$. Let $J_i$ denote the set of direct descendent indexes of $i$ in
 $\Gamma$.
Let $s=|J_i|$ be the cardinality of $J_i$, the number of direct descendents
 of $i$ in $\Gamma_{e_i}$. The number $s$ is also equal to $codim_{\bf C}\Gamma_{e_i}$. 

  The inclusion of the fiber bundle $D_i/B'_i\hookrightarrow ({\cal X}/B)_{i+1}/B'_i$ induces
 the canonical map on the $|J_i|-$th relative universal spaces,
 $$(D_i/B'_i)_{|J_i|}\hookrightarrow (({\cal X}/B)_{i+1}/B'_i)_{|J_i|}\cong 
({\cal X}/B)_{|J_i|+i}.$$

 However, the fiber bundle 
$D_i\mapsto B'_i$ is smooth of relative dimension one. So by a direct check we find 
 (using the fact the codimension one blowing ups are trivial) $(D_i/B'_i)_s\cong 
\times_{B'_i}^sD_i$, the $s$-fold fiber product of $D_i$ over $B'_i$. On the other hand,
 proposition \ref{prop; lifting} implies that for ${\cal X}'/B'_i=({\cal X}/B)_{i+1}/B'_i$,
$\pi_{J_i}:({\cal X}'/B'_i)_{n-i}\mapsto
 \times_{B'_i}^{|J_i|}({\cal X}'/B'_i)$ can be lifted to a smooth and surjective map 
$\check{\pi}_{J_i}:({\cal X}'/B'_i)_{n-i}\mapsto ({\cal X}'/B'_i)_{|J_i|}$. Moreover,
 the isomorphism $({\cal X}/B)_n\stackrel{\psi_{i, n}}{\longrightarrow}
(({\cal X}/B)_{i+1}/B'_i)_{n-i}$ allows us to
 view the $\psi_{i, n}^{-1}$ pre-image of the 
relative \footnote{From the subscript of $\gamma_{n-i}$, 
 one should be able to distinguish ${\bf Y}(\gamma_n)\subset 
({\cal X}/B)_n$ over $B$ and
 ${\bf Y}_{\gamma_{n-i}}\subset \bigl(({\cal X}/B)_{i+1}/B'_i\bigr)_{n-i}$ over 
$B'_i=({\cal X}/B)_i$.} ${\bf Y}_{\gamma_{n-i}}$ over $B'_i$ as a subset of $({\cal X}/B)_n$.

We have the following characterization of 
${\bf Y}(\Gamma_{e_i})$ and ${\bf Y}_{\Gamma_{e_i}}$,

\begin{lemm}\label{lemm; fanlike}
The closed subspace ${\bf Y}(\Gamma_{e_i})\subset ({\cal X}/B)_n$, smooth
of codimension $codim_{\bf C}\Gamma_{e_i}$ in $({\cal X}/B)_n$, is
 the pre-image of $\times_{B'_i}^sD_i\subset (({\cal X}/B)_{i+1}/B'_i)_{|J_i|}$
 under $\check{\pi}_{J_i}^{-1}$. Likewise the locally closed subset 
 ${\bf Y}_{\Gamma_{e_i}}\subset {\bf Y}(\Gamma_{e_i})$
 can be identified with ${\bf Y}(\Gamma_{e_i})\cap \psi_{i, n}^{-1}({\bf Y}_{\gamma_{n-i}})$.
\end{lemm}

\noindent Proof: For all $b\in B$, by remark \ref{rem; interpret} we know that 
$({\cal X}_b)_n$ parametrizes all the
 ordered $n-$consecutive pointwise blowing ups of ${\cal X}_b$. Given 
the fan-like admissible graph $\Gamma_{e_i}$, an ordered $n-$consecutive blowing
 ups from ${\cal X}_b$ lies in the fiberwise $Y(\Gamma_{e_i})$ of
 $({\cal X}/B)_n\mapsto B$ above $b\in B$ iff
 all the $k_1, k_2, \cdots, k_s-$th blown up points, $k_l\in J_i$,
$1\leq l\leq s=codim_{\bf C}\Gamma_{e_i}=|J_i|$, lie
 above the exceptional ${\bf P}^1$ of the $i-$th blown up point.
 Over the relative $i-$th universal space $({\cal X}/B)_i$ which parametrizes the
 first $i-$th blowing ups in the $B$ family, the union of the
$i-$th exceptional ${\bf P}^1$ forms a fiber bundle, which is nothing but 
 $D_i\mapsto ({\cal X}/B)_i$ introduced above. On the one hand 
 all the $k_1, k_2, \cdots, k_s-$th blown up points are allowed to
 move on the fibers of $D_i$ freely. This implies that
 ${\bf Y}(\Gamma_{e_i})\stackrel{\check{\pi}_{J_i}}
{\longrightarrow} \times_{({\cal X}/B)_i}^sD_i$ must be 
surjective. On the other hand, for all $j>i$ which are not the direct
descendent indexes of $i$, the $j-$th blowing up centers within the
 $n-$consecutive pointwise blowing ups are not restricted at all. Therefore
 ${\bf Y}(\Gamma_{e_i})$ can be identified with 
 $\check{\pi}_{J_i}^{-1}(\times_{({\cal X}/B)_i}^sD_i)$.
 
 Inside this smooth space ${\bf Y}(\Gamma_{e_i})$, the sub-locus 
 ${\bf Y}_{\Gamma_{e_i}}$
 corresponds to the set of all the 
$n-$consecutive blowing ups from ${\cal X}_b$ above all $b\in B$ such 
that it is in the fiberwise $Y(\Gamma_{e_i})$ and none of blowing ups marked by 
$\{i+1, i+2, \cdots n\}$ lie above the exceptional loci of one another. 
Thus the space ${\bf Y}_{\Gamma_{e_i}}$
 must map into the relative ${\bf Y}_{\gamma_{n-i}}$ over $B'_i=({\cal X}/B)_i$ 
 under $\psi_{i, n}:({\cal X}/B)_n\mapsto 
\bigl(({\cal X}/B)_{i+1}/({\cal X}/B)_i\bigr)_{n-i}$, as the space ${\bf Y}_{\gamma_{n-i}}$
 parametrizes all the disjoint 
last $n-i-$th consecutive blowing ups of the family $({\cal X}/B)_{i+1}/B'_i$
 of algebraic surfaces. Thus ${\bf Y}_{\Gamma_{e_i}}={\bf Y}(\Gamma_{e_i})\cap 
\psi_{i, n}^{-1}({\bf Y}(\gamma_{n-i}))$ is locally closed.
$\Box$

 By a direct calculation, $codim_{\bf C}\Gamma_{e_i}$, being the
 number of direct descendents of $\Gamma_{e_i}$, is also equal to
 the negation of 
$d_{GT}(e_i)={e_i^2-c_1({\bf K}_{M_{n+1}/M_n})\cdot e_i\over 2}$. This
 implies that the existence locus $\subset ({\cal X}/B)_n$ 
of the fiberwise class $e_i\in {\cal A}_{\cdot}(({\cal X}/B)_{n+1}/({\cal X}/B)_n)$
 (over which $e_i$ becomes effective along the fibers)
is smooth of the expected dimension $dim_{\bf C}({\cal X}/B)_n)+d_{GT}(e_i)$.
 Moreover ${\bf Y}_{\Gamma_{e_i}}$ is the locus over which the type $I$ class $e_i$ is
 effective and irreducible/smooth in the
 fibers of $({\cal X}/B)_{n+1}\mapsto ({\cal X}/B)_n$.

\label{fig.3}
\begin{figure}
\centerline{\epsfig{file=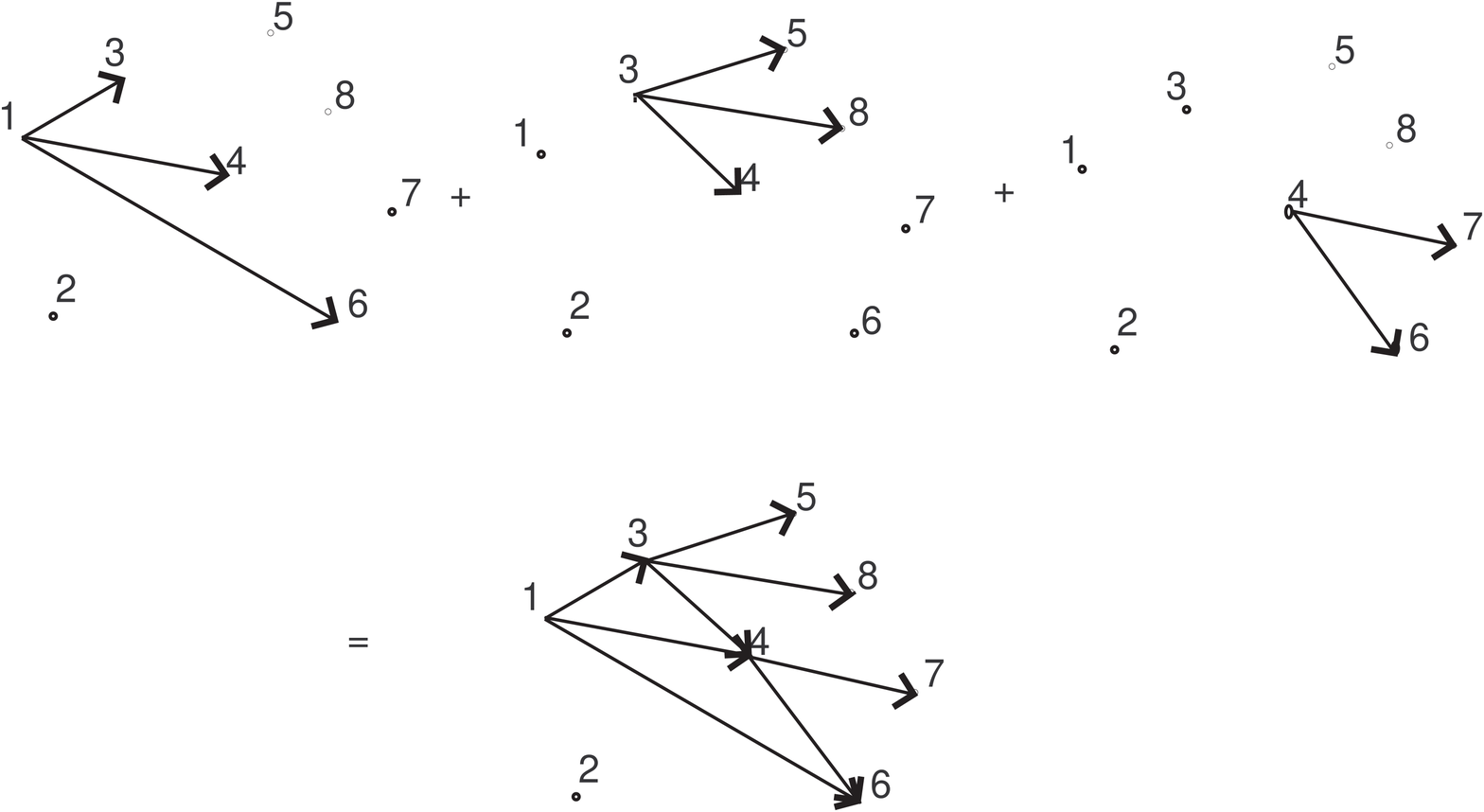,height=6cm}}
\centerline{fig.3}
\centerline{When we fuse the admissible graphs $\Gamma_{e_1}$, 
 $\Gamma_{e_3}$, $\Gamma_{e_4}$ together, we recover the original $\Gamma$}
\end{figure}

 The following proposition characterizes ${\bf Y}(\Gamma)$ in terms of 
 the fan-like graphs $\Gamma_{e_i}$ and can be viewed as 
 the converse of proposition \ref{prop; reverse}, 

\begin{prop}\label{prop; transversal}
 Let $\Gamma\in adm(n)$ be an $n-$vertex admissible graph.
 The smooth and closed subspace ${\bf Y}(\Gamma)\subset ({\cal X}/B)_n$ of 
$codim_{\bf C}\Gamma$ can 
 be identified with the regular intersection $\cap_{1\leq i\leq n}
{\bf Y}(\Gamma_{e_i})$. 
 Likewise the locally closed smooth subspace ${\bf Y}_{\Gamma}\subset 
{\bf Y}(\Gamma)$
 is equal to the intersection $\cap_{1\leq i\leq n}
{\bf Y}_{\Gamma_{e_i}}$.
\end{prop}

\noindent Proof:  

$\diamondsuit$ Auxiliary Statement: Let $K\subset \{1, 2, \cdots, n\}$ 
be an index subset. We claim that the restricted lifted
mapping
$\check{\pi}_{K}:{\bf Y}(\Gamma)\mapsto ({\cal X}/B)_{|K|}$ (or
 its restriction to the subspace ${\bf Y}(\Gamma)$) maps smoothly onto
 ${\bf Y}(\Gamma_K)\subset ({\cal X}/B)_{|K|}$ (or ${\bf Y}(\Gamma_K)$, 
for an admissible $\Gamma_K\in adm(|K|)$ characterized
 as the following: Firstly, the map $\check{\pi}_K$ induces an ordering preserving
 bijection $\phi_K:K\mapsto [1, 2, \cdots, |K|]$ between index sets. Then the
 type $I$ exceptional cycle classes $e_i, 1\leq i\leq n$ along the fibers of 
$({\cal X}/B)_{n+1}\times_{({\cal X}/B)_n}
{\bf Y}(\Gamma)\mapsto {\bf Y}(\Gamma)$ are pushed-forward
 to fiberwise cycle classes
 of $({\cal X}/B)_{|K|+1}\mapsto ({\cal X}/B)_{|K|}$ by the
 following rule: The $e_i\mapsto 0, i\not\in K$; but 
$e_i\mapsto e^K_{\phi_K(i)}$ for $i\in K$. Those 
 $e^K_{\phi_K(i)}$ are constructed from $e_i, i\in K$, by the following 
 substitutions: ${\cal A}_{\cdot}(({\cal X}/B)_{n+1})
\ni E_a\mapsto E_{\phi_K(a)}\in {\cal A}_{\cdot}(({\cal X}/B)_{|K|+1})$, 
$a\in K$, and $E_a\mapsto 0$, $a\not\in
K$. By a simple calculation based on the substitution rules we find that
 $e^K_i\cdot e^K_j\geq e_{\phi_K^{-1}(i)}\cdot e_{\phi_K^{-1}(j)}\geq 0$ for
 all $i\not=j$, $i, j\in [1, \cdots, |K|]$.
Then the desired $\Gamma_K\in adm(|K|)$
 is constructed from the collection of classes $e^K_j, 1\leq j\leq |K|$,
 by applying proposition \ref{prop; reverse}.

 We prove the proposition
 along with the auxiliary statement $\diamondsuit$ 
based on induction arguments on $n$.

 For $n=1$, there is only one admissible graph $\Gamma=\gamma_1\in adm(1)$
 and the proofs of both the statements are trivial. 
Suppose that the statements have been proved for all the pairs ${\cal X}\mapsto B$
 for the natural numbers $\leq n$, we prove them for $n+1$.

 Let $\Gamma\in adm(n+1)$ be an $n+1$-vertex admissible graph.
 As in the proof of proposition \ref{prop; reverse}, 
we subtract all the indexes by $-1$
 and denote the resulting graph by $\tilde{\Gamma}$. Then as before the 
subgraph $\Gamma'$, marked by
 $\{1, 2, \cdots, n\}$, is an admissible graph $\in adm(n)$ in the usual sense.

 By lemma \ref{lemm; fuse} we have 
$({\cal X}/B)_{n+1}=\bigl(({\cal X}/B)_2/({\cal X}/B)_1\bigr)_n$. 
Over the base space
 ${\cal X}/B$ the locus
 ${\bf Y}(\Gamma')$ defines a closed subset of the relative 
$n-$th universal space 
 of $({\cal X}/B)_2\mapsto ({\cal X}/B)_1$, $({\cal X}/B)_{n+1}/({\cal X}/B)_1$.
 By
 induction hypothesis, we may assume ${\bf Y}(\Gamma')=
\cap_{1\leq i\leq n}{\bf Y}(\Gamma'_{e'_i})$
 to be a regular intersection.
 The class $e'_i/\Gamma'_{e'_i}$ are the type $I$
 exceptional classes/fan-like admissible graphs
 associated with the vertexes of $\Gamma'$. On the other hand, $\Gamma'$ is 
the subgraph of
 $\tilde{\Gamma}$ removing the $0-$th vertex and all the arrowed one-edges 
starting from '0'. 
So $e'_i=\tilde{e}_i$ and
 $\Gamma'_{e'_i}=\tilde{\Gamma}_{\tilde{e}_i}$ for $1\leq i\leq n$. 
 As before we define $J_0$ to be the set of all
 direct descendent indexes of $'0'$ in $\tilde{\Gamma}$.
By lemma \ref{lemm; fanlike},
 ${\bf Y}(\tilde{\Gamma}_{\tilde{e}_0})\subset 
\bigl(({\cal X}/B)_2/({\cal X}/B)_1\bigr)_n=({\cal X}/B)_{n+1}$
 is smooth of codimension
 $codim_{\bf C}\tilde{\Gamma}_{\tilde{e}_0}$, the pre-image of
 $\times_{({\cal X}/B)_1}^{codim_{\bf C}\tilde{\Gamma}_{\tilde{e}_0}}D_0$ under
 $\check{\pi}_{J_0}^{-1}$.

 Consider the admissible graph $\Gamma'_{J_0}=\Gamma_{J_0}\in adm(|J_0|)$,
 constructed by the recipe at the beginning of our proof. 
 We claim that it must be a finite union of linear chains 
\footnote{Refer to fig.4 on
 page \pageref{fig.4} for an example.}. This is equivalent to say that
 all the type $I$ exceptional 
classes associated to $\Gamma'_{J_0}$ are either $-1$ or $-2$
classes. If there is a $-k$ class (with $k>2$) 
among the type $I$ exceptional classes
 of $\Gamma'_{J_0}$, then there is an index $a$ with more than one direct
descendent in $\Gamma'_{J_0}$. By the construction of $\Gamma'_{J_0}=\Gamma_{J_0}$
 from
 $\Gamma$, it implies that $\phi_{J_0}^{-1}(a)\in J_0$ and $0$
 share more than one direct descendent in $J_0$. However, this would
\footnote{By a similar calculation as was performed in the proof of
 proposition \ref{prop; reverse}.} imply
 $\tilde{e}_0\cdot \tilde{e}_{\phi_{J_0}^{-1}(a)}\leq -1<0$, violating the
 non-negativity of the intersection numbers between distinct type $I$ 
exceptional classes.

 By the inductive assumption on the auxiliary statement $\diamondsuit$,
 we know that 
${\bf Y}(\Gamma')$ is a smooth fibration over 
${\bf Y}(\Gamma'_{J_0})$ under $\check{\pi}_{J_0}$.
 
   Because the chain-like nature of $\Gamma'_{J_0}$, 
${\bf Y}(\Gamma'_{J_0})\times_{
\bigl(({\cal X}/B)_2/({\cal X}/B)_1\bigr)_{|J_0|}}
\bigl(\times_{({\cal X}/B)_1}^{|J_0|}D_0\bigr)$ is the 
 $codim_{\bf C}\Gamma'_{J_0}$ partial diagonal of the
 $({\bf P}^1)^{|J_0|}$ bundle $\times_{({\cal X}/B)_1}^{|J_0|}D_0$ (i.e. demanding
that the ${\bf P}^1$ coordinates marked by indexes 
within the same connected component of the
chain $\Gamma'_{J_0}$ to be equal). And
 therefore it is a regular intersection (of codimension 
 $codim_{\bf C}\Gamma_0$) in $\times_{({\cal X}/B)_1}^{|J_0|}D_0$ and
 is of codimension $|J_0|=codim_{\bf C}\Gamma'_{J_0}$ 
in ${\bf Y}(\Gamma'_{J_0})$.

 Then by the fact that ${\bf Y}(\tilde{\Gamma}_{\tilde{e}_0})$ is the
 inverse image of $\times_{({\cal X}/B)_1}^{|J_0|}D_0$ under $\check{\pi}_{J_0}$, 
the fiber product (also the pull-back of the smooth fibration
 ${\bf Y}(\Gamma')\mapsto {\bf Y}(\Gamma'_{J_0})$

$$\hskip -1in
{\bf Y}(\Gamma')\times_{ 
\bigl(({\cal X}/B)_2/({\cal X}/B)_1\bigr)_{|J_0|}}
\bigl(\times_{({\cal X}/B)_1}^{|J_0|}D_0\bigr)
={\bf Y}(\Gamma')\cap {\bf Y}(\tilde{\Gamma}_{\tilde{e}_0})=\cap_{0\leq i\leq n}
 {\bf Y}(\tilde{\Gamma}_{\tilde{e}_i})$$

is irreducible, smooth of codimension $codim_{\bf C}{\bf Y}(\Gamma')+
|J_0|=codim_{\bf C}{\bf Y}(\Gamma')+codim_{\bf C}\tilde{\Gamma}_{\tilde{e}_0}=
codim_{\bf C}{\bf Y}(\tilde{\Gamma})=codim_{\bf C}\tilde{\Gamma}$ in
 $({\cal X}/B)_{n+1}=\bigl(({\cal X}/B)_2/({\cal X}/B)_1\bigr)_n$. 

\medskip

Correspondingly by a similar inductive argument the locally closed 
${\bf Y}_{\tilde{\Gamma}}$ is equal to 
 ${\bf Y}_{\Gamma'}\cap {\bf Y}_{\tilde{\Gamma}_{e_0}}=\cap_{0\leq i\leq n}
{\bf Y}_{\tilde{\Gamma}_{\tilde{e}_i}}$, a Zariski dense subset of
 ${\bf Y}(\tilde{\Gamma})$ and therefore a locally closed subset of 
$({\cal X}/B)_{n+1}$.

 By adding $1$ back to all the indexes, we find that 
${\bf Y}(\Gamma)={\bf Y}(\tilde{\Gamma})$ is an irreducible 
 regular intersection $\cap_{1\leq i\leq n+1}{\bf Y}(\Gamma_{e_i})$ of
 codimension $codim_{\bf C}\Gamma$ in $({\cal X}/B)_{n+1}$ and 
 ${\bf Y}_{\Gamma}=\cap_{1\leq i\leq n+1}{\bf Y}_{\Gamma_{e_i}}$ is an open
 subset of ${\bf Y}(\Gamma)$. So we
 have proved the proposition for $n+1$.

 Now let us prove the auxiliary statement $\diamondsuit$ 
on the smoothness of the restricted morphism
 $\check{\pi}_{K}$.
 Let $K$ be an index subset of $\{1, 2, \cdots, n+1\}$. We show that
 $\check{\pi}_K:{\bf Y}(\Gamma)\mapsto ({\cal X}/B)_{|K|}$ maps smoothly
 onto ${\bf Y}(\Gamma_K)$, for the $\Gamma_K\in adm(|K|)$ constructed earlier.

 Firstly, if $K=\{1, \cdots, n+1\}$ itself, then 
$\Gamma_K=\Gamma$ itself and the map is an isomorphism. 
From now on we may assume $|K|\leq n$. Denote $I=\{1, 2, \cdots, n\}$. 
By induction hypothesis ${\bf Y}(\Gamma_I)$ maps
onto ${\bf Y}(\Gamma_{K\cap I})$
 smoothly under $\check{\pi}_{K\cap I}$ and we have
the following commutative diagram,

\[
\begin{array}{ccccccc}
{\bf Y}(\Gamma) & \subset & ({\cal X}/B)_{n+1} & 
\stackrel{\check{\pi}_K}{\longrightarrow}  & ({\cal X}/B)_{|K|} &
 \supset & {\bf Y}(\Gamma_K)\\
  &  & \Big\downarrow\vcenter{%
\rlap{$\scriptstyle{\check{\pi}_I}\,$}} &   & \Big\downarrow
 & & \Big\downarrow\\
{\bf Y}(\Gamma_I) & \subset & ({\cal X}/B)_n & \stackrel{\check{\pi}_{I\cap K}}
{\longrightarrow}
  &  ({\cal X}/B)_{|K\cap I|} & \supset & {\bf Y}(\Gamma_{K\cap I})
\end{array}
\]

 It is easy to see  
that $\Gamma_I\in adm(n)$ can be viewed as
 the admissible sub-graph \footnote{It was denoted by 
$\Gamma(-1)$ in [Liu1].} formed by restricting to 
the first $n$ vertexes (and the one-edges between them) of $\Gamma$.
 Then by remark \ref{rem; interpret} 
the forgetful map (i.e. forgetting the last index $n+1$)
 $\check{\pi}_I:
{\bf Y}(\Gamma)\mapsto {\bf Y}(\Gamma_I)$ is smooth of relative
dimension zero (i.e. isomorphic), dimension one (a ${\bf P}^1$ bundle), or
 dimension two\footnote{If necessary, please consult the inductive 
construction of $Y(\Gamma)$ on page 418-419 of [Liu1] for more details.
 in that construction, a dependence of the fiber bundle structure upon
the number of direct ascendents of $n+1$ was discussed in more details.},
 depending on whether $n+1$ has two, one or no direct ascendent(s)
 in $I$. Thus the map is smooth.
  
 The smoothness of 
 $\check{\pi}_K:{\bf Y}(\Gamma)\mapsto {\bf Y}(\Gamma_K)$
 follows from the smoothness of both
 ${\bf Y}(\Gamma_I)\mapsto {\bf Y}(\Gamma_{K\cap I})$ 
(by the induction hypothesis)
 and ${\bf Y}(\Gamma)\mapsto {\bf Y}(\Gamma_I)$,
 and the commutativity of the above 
diagram. The surjectivity of the map 
$\check{\pi}_K:{\bf Y}(\Gamma)\mapsto 
{\bf Y}(\Gamma_K)$ follows from the fact \footnote{which can be 
checked directly.} that all
the fiberwise smooth and irreducible type $I$ exceptional curves over 
$Y_{\Gamma}$ are mapped to smooth and irreducible type $I$ curves (or points)
 dual to $e^K_i, 1\leq i\leq |K|$ under $({\cal X}/B)_{n+1}/({\cal X}/B)_n \mapsto 
({\cal X}/B)_{|K|+1}/({\cal X}/B)_{|K|}$. By the induction hypothesis of the
 proposition, ${\bf Y}_{\Gamma_K}$ has been the locus over which the
 type $I$ curves representing $e^K_i, 1\leq i\leq |K|$ co-exist as smooth
 curves. As ${\bf Y}_{\Gamma}\mapsto {\bf Y}_{\Gamma_K}$ is onto, the
closure ${\bf Y}(\Gamma)$ has to be mapped onto the closure
$\overline{{\bf Y}(\Gamma_K)}={\bf Y}(\Gamma_K)$. So the inductive proof
 of the auxiliary statement $\diamondsuit$ has been complete.
$\Box$

\medskip

\begin{rem}\label{rem; alternative}
If we desire to minimize the dependence to the reference [Liu1], one may take 
an alternative route.
One may turn lemma \ref{lemm; fanlike} and
 proposition \ref{prop; transversal} into constructions/definitions and
 use them to define the relative admissible strata
${\bf Y}(\Gamma_{e_i})$, ${\bf Y}_{\Gamma_{e_i}}$, 
${\bf Y}(\Gamma)$ and
 ${\bf Y}_{\Gamma}$, etc. 
 Then one may deduce all the basic properties from them. At the end
 we may take $B=pt$ to recover the usual $Y(\Gamma)$ and $Y_{\Gamma}$ as a special case. 
\end{rem}

\medskip

\begin{rem}\label{rem; classexist}
We have remarked that the
 locus $Y(\Gamma_{e_i})\subset M_n$ is the smooth locus over which the type $I$ exceptional
 class becomes effective. So $Y(\Gamma)=\cap_{1\leq i\leq n}Y(\Gamma_{e_i})$ is the
 locus over which all the $e_1, e_2, \cdots, e_n$ become effective.
Likewise, $Y_{\Gamma}\subset Y(\Gamma)\subset M_n$
 is the locally closed locus over which the classes 
$e_1, e_2, \cdots, e_n$ co-exist as smooth and irreducible type $I$ exceptional curves.
\end{rem}

\medskip

\begin{rem}\label{rem; ignore}
 In the intersection $\cap_{1\leq i\leq n}Y(\Gamma_{e_i})$, we may
 ignore all the $i$ such that $e_i^2=-1$. Because each of such $e_i=E_i$ is a
 $-1$ class, the corresponding $i$ has no direct descendent in $\Gamma_{e_i}$.
 We have $\Gamma_{e_i}=\gamma_n$ and $Y(\Gamma_{e_i})=Y(\gamma_n)=M_n$. So the
 intersection with these $Y(\Gamma_{e_i})$ can be skipped.
\end{rem}

 From now on we will make use of this simple observation implicitly.
 In the latter sections, we will use proposition \ref{prop; transversal} 
frequently 
and view $Y(\Gamma)$ as the locus of co-existence of all the type $I$ exceptional
 classes $e_i$, over which they all become effective.
 
\medskip 

\section{\bf The Construction of the Quotient Bundle Based on the
 ${\bf P}^1$ Fibrations of Universal Curves}\label{section; bundle}

\bigskip

In this section, let $\Gamma\in adm(n)$ be an $n$-vertex admissible graph and
 let $Y(\Gamma)$ be the closure of the admissible stratum associated to $\Gamma$, 
as was described in section \ref{section; strata}.
Let $M_{n+1}\times_{M_n}Y(\Gamma)\mapsto Y(\Gamma)$ be the fiber bundle
 of algebraic surfaces over $Y(\Gamma)$ induced by $M_{n+1}\mapsto M_n$ through the
 pull-back map of $Y(\Gamma)\subset M_n$. 
Let $e_1, e_2, \cdots, e_n$ denote the type $I$ exceptional classes associated
 to $Y(\Gamma)$. As usual, we let $e_{k_1}, e_{k_2}, e_{k_3}, \cdot, e_{k_i}, 
\cdots, e_{k_p}$, $k_1<k_2<\cdots k_p$, 
$1\leq i\leq p$, denote the type $I$ exceptional classes which
 pair negatively with $C-{\bf M}(E)E$. Because each $e_i$ is effective and is represented by
 a unique curve over each point of $Y(\Gamma_{e_i})$, 
the notation $\Xi_{i}\mapsto Y(\Gamma_{e_i})$ has been used in [Liu1], [Liu3], [Liu5] to
 denote the ${\bf P}^1$ fibration (embedded in the fiber bundle 
$M_{n+1}\times_{M_n}Y(\Gamma_{e_i})\mapsto 
Y(\Gamma_{e_i})$) representing the universal curves of $e_i$.

 In section 5, proposition 9 of [Liu5], we had
 analyzed the canonical algebraic family Kuranishi models
 of two classes $C-{\bf M}(E)E-\sum_{e_i\cdot (C-{\bf M}(E)E)<0}e_i$ and $C-{\bf M}(E)E$,
 $(\Phi_{{\cal V}_{canon}^{\circ}{\cal W}_{canon}^{\circ}}, {\cal V}_{canon}^{\circ}, 
{\cal W}_{canon}^{\circ})$ and 
$(\Phi_{{\cal V}_{canon}{\cal W}_{canon}}, {\cal V}_{canon}, {\cal W}_{canon})$ under
 the assumption\footnote{We will assume that these conditions hold for the class 
 $C$ throughout the paper. \label{Cample}} that
 ${\cal R}^1\pi_{\ast}\bigl({\cal E}_C\bigr)={\cal R}^2\pi_{\ast}\bigl({\cal E}_C\bigr)=0$.

   We know ${\cal V}_{canon}^{\circ}={\cal V}_{canon}$ but ${\cal W}_{canon}^{\circ}$ and
 ${\cal W}_{canon}$ differ from each other. In fact we have the four term exact sequence
\footnote{by proposition 9 of [Liu5].}, 

$$\hskip -.7in
0\mapsto 
{\cal R}^0\pi_{\ast}\bigl({\cal O}_{\sum_{1\leq i\leq p}\Xi_{k_i}}\otimes 
{\cal E}_{C-{\bf M}(E)E}\bigr)\mapsto {\cal W}_{canon}^{\circ}\mapsto {\cal W}_{canon}
\mapsto {\cal R}^1\pi_{\ast}\bigl({\cal O}_{\sum_{1\leq i\leq p}\Xi_{k_i}}\otimes 
{\cal E}_{C-{\bf M}(E)E}\bigr)\mapsto 0.$$

 In particular their difference in the $K$ group of 
${\bf P}({\bf V}_{canon})\times_{M_n}Y(\Gamma)$ can be
 represented by 
 ${\cal R}^1\pi_{\ast}\bigl({\cal O}_{\sum_{i\leq p} \Xi_{k_i}}\otimes {\cal E}_{C-{\bf M}(E)E}
\bigr)-
{\cal R}^0\pi_{\ast}\bigl({\cal O}_{\sum_{i\leq k} \Xi_{k_i}}\otimes
 {\cal E}_{C-{\bf M}(E)E}\bigr)$.

 Because the canonical sections $s_{canon}^{\circ}$ and $s_{canon}$
 and
 the ${\bf H}-$twisted bundle map $\pi_X^{\ast}{\bf W}_{canon}^{\circ}\otimes {\bf H}
|_{X\times_{M_n}Y(\Gamma)}
\mapsto \pi_X^{\ast}{\bf W}_{canon}\otimes {\bf H}|_{X\times_{M_n}Y(\Gamma)}$ over 
$X\times_{M_n}Y(\Gamma)={\bf P}({\bf V}_{canon})\times_{M_n}Y(\Gamma)$ play important
roles in the paper, it is vital for us to study the map ${\bf W}_{canon}^{\circ}\mapsto 
{\bf W}_{canon}$ in more details.

The ${\bf P}^1$ fibration $\Xi_i\mapsto Y(\Gamma_{e_i})$ may have singular
 fibers which are trees of ${\bf P}^1$ curves.
Despite that the invertible sheaf 
${\cal O}_{\Xi_{k_i}}(-{\bf M}(E)E)$ is of negative relative degree along
 $\Xi_{k_i}\mapsto Y(\Gamma_{e_{k_i}})$, the kernel vector spaces 
${\cal R}^0\pi_{\ast}\bigl({\cal O}_{\sum_{i\leq p} \Xi_{k_i}}(-{\bf M}(E)E)\otimes
 {\cal E}_C\otimes k(y)\bigr)$ may not always be the zero vector spaces
 and the canonical bundle map
 ${\bf W}_{canon}^{\circ}\mapsto {\bf W}_{canon}$ may fail to be injective over
 some sub-locus of $Y(\Gamma)$.

 It is the goal of this section to construct an \label{quot}
 algebraic quotient bundle ${\bf V}_{quot}\mapsto
 {\bf P}({\bf V}_{canon})\times_{M_n}Y(\Gamma)$
 of ${\bf W}_{canon}|_{Y(\Gamma)\times T(M)}$ of rank\footnote{By applying
 curve Riemann-Roch to the fibers of $\sum_{1\leq i\leq p}\Xi_{k_i}$.}
$-p+\sum_{1\leq i\leq p}e_{k_i}\cdot (-{\bf M}(E)E-\sum_{1\leq j<i\leq p}e_{k_j})$ 
and identify its total Chern
 class explicitly.

 To construct ${\bf V}_{quot}$, we consider the torsion free part of the 
coherent sheaf 

 ${\cal R}^1\pi_{\ast}\bigl({\cal O}_{\sum_{i\leq p} 
\Xi_{k_i}}(-{\bf M}(E)E)\otimes
 {\cal E}_C\bigr)$ and show that, 

\medskip

\noindent {\bf Claim}:
The torsion free part of the first right derived image sheaf 
${\cal R}^1\pi_{\ast}\bigl({\cal O}_{\sum_{i\leq p} \Xi_{k_i}}(-{\bf M}(E)E)\otimes
 {\cal E}_C\bigr)$ is locally free.

The proof of the claim will appear in section \ref{subsection; lf} 
 proposition \ref{prop; id}.

 Once we know that the torsion free sheaf is locally free, we denote it by ${\cal V}_{quot}$ and
 the corresponding vector bundle is our desired ${\bf V}_{quot}$.

 The key idea for the explicit determination of its Chern classes 
is to consider the relative minimal model $\tilde{\Xi}_i$
 of $\Xi_i$ (see proposition 5.1 on page 442 of [Liu1]), which has a structure of
 ${\bf P}^1$ bundle over $Y(\Gamma_{e_i})$. The ${\bf P}^1$ fibration
 $\Xi_i\mapsto Y(\Gamma_{e_i})$ can be viewed as some consecutive 
blowing up from $\tilde{\Xi}_i$ along
 some codimension two sub-loci \footnote{Even though it is possible 
to determine the sub-loci, the explicit form of this loci is not crucial to us.}
 determined by the graph $\Gamma$.

 From the brief discussion in subsection \ref{subsubsection; torsion} below on the
 torsion free quotient, we know that there is a canonical surjection, 
 $${\cal R}^1\pi_{\ast}\bigl({\cal O}_{\sum_{i\leq p} \Xi_{k_i}}(-{\bf M}(E)E)\otimes
 {\cal E}_C\bigr)\mapsto 
({\cal R}^1\pi_{\ast}\bigl({\cal O}_{\sum_{i\leq p} \Xi_{k_i}}(-{\bf M}(E)E)\otimes
 {\cal E}_C\bigr))_{torfree}={\cal V}_{quot}.$$ 

 By composing with the surjective
 sheaf morphism ${\cal W}_{canon}|_{Y(\Gamma)\times T(M)}\mapsto 
{\cal R}^1\pi_{\ast}\bigl({\cal O}_{\sum_{i\leq p} \Xi_{k_i}}(-{\bf M}(E)E)\otimes
 {\cal E}_C\bigr)$, we get the surjection ${\cal W}_{canon}|_{Y(\Gamma)\times T(M)}\mapsto {\cal V}_{quot}$.
  Because both sheaves are locally free\footnote{by the claim above.},
 we have the vector bundle short exact sequence over 
$Y(\Gamma)\times T(M)$,  $0\mapsto \underline{\bf W}_{canon}
 \mapsto {\bf W}_{canon}|_{Y(\Gamma)\times T(M)}
 \mapsto {\bf V}_{quot}
 \mapsto 0$, \label{underlinewcanon}
 where $\underline{\bf W}_{canon}$ is defined to be 
the kernel bundle of ${\bf W}_{canon}|_{Y(\Gamma)\times
 T(M)} \mapsto {\bf V}_{quot}$.

\begin{defin}\label{defin; under}
Define $\underline{\bf W}_{canon}$ to be the kernel bundle of 
 ${\bf W}_{canon}|_{Y(\Gamma)\times T(M)}\mapsto {\bf V}_{quot}\mapsto 0$.
\end{defin}

\begin{rem}\label{rem; comment}
 This short exact sequence will plays an essential role in our
 proof of the main theorem
 in section \ref{section; proof}.
\end{rem}

  Because the
 composition ${\bf W}_{canon}^{\circ}|_{Y(\Gamma)\times T(M)}
\mapsto {\bf W}_{canon}|_{Y(\Gamma)\times T(M)}
\mapsto {\bf V}_{quot}$
 is the zero map, ${\bf W}_{canon}^{\circ}|_{Y(\Gamma)\times T(M)}\mapsto 
{\bf W}_{canon}|_{Y(\Gamma)\times T(M)}$ 
factors through the
 kernel $\underline{\bf W}_{canon}$ and,

\begin{lemm}\label{lemm; injective}
The induced bundle map ${\bf W}_{canon}^{\circ}|_{Y(\Gamma)\times T(M)} \mapsto 
\underline{\bf W}_{canon}$ is
 injective over a Zariski open subset $U=Y_{\Gamma}\times T(M)$
 of $Y(\Gamma)\times T(M)$.
\end{lemm}

\noindent Proof of lemma \ref{lemm; injective}: The kernel spaces 
${\cal R}^0\pi_{\ast}\bigl({\cal O}_{\sum_{1\leq i\leq p}\Xi_{k_i}}\otimes 
{\cal E}_{C-{\bf M}(E)E}\otimes k(y)\bigr)$, $y\in Y(\Gamma)\times T(M)$, of 
${\bf W}_{canon}^{\circ}|_{Y(\Gamma)\times T(M)}\mapsto 
{\bf W}_{canon}|_{Y(\Gamma)\times T(M)}$ 
(see proposition 9 of [Liu5]) is supported ``away'' from 
the open dense subset 
$Y_{\Gamma}\times T(M)\subset Y(\Gamma)\times T(M)$ over which the
 fibrations $\Xi_{k_i}\mapsto Y(\Gamma)$ are smooth and 
irreducible for all $1\leq i\leq p$.
 Since the fibers of the restricted fibrations 
$\Xi_{k_i}\times_{Y(\Gamma_{e_{k_i}})}Y_{\Gamma}\mapsto Y_{\Gamma}$, 
$1\leq i\leq p$, remain smooth and irreducible throughout $Y_{\Gamma}\times 
T(M)$, the vanishing of
 ${\cal R}^0\pi_{\ast}\bigl({\cal O}_{\sum_{1\leq i\leq p}\Xi_{k_i}}\otimes 
{\cal E}_{C-{\bf M}(E)E}\bigr)|_{Y_{\Gamma}\times T(M)}$ is due to
 the negative relative degrees of the invertible ${\cal E}_{C-{\bf M}(E)E}$ on
 all the $\Xi_{k_i}$, i.e. $e_{k_i}\cdot (C-{\bf M}(E)E)<0$. So 
${\cal R}^0\pi_{\ast}\bigl({\cal O}_{\sum_{1\leq i\leq p}\Xi_{k_i}}\otimes 
{\cal E}_{C-{\bf M}(E)E}\bigr)$ is a torsion sheaf \footnote{As it is
a torsion sub-sheaf of a locally free sheaf, it vanishes. But the
 sheaf injection ${\cal W}_{canon}^{\circ}\mapsto {\cal W}_{canon}$ does not
induce a bundle injection, because generally speaking 
$\otimes k(y)$ is not left-exact.}
 over $Y(\Gamma)\times T(M)$.
 
On the other hand, the sheaf map ${\cal W}_{canon}^{\circ}|_{Y(\Gamma)\times T(M)}\mapsto 
{\cal W}_{canon}|_{Y(\Gamma)\times T(M)}$ factors through the intermediate 
$\underline{\cal W}_{canon}$. So the
sheaf morphism ${\cal W}_{canon}^{\circ}\mapsto \underline{\cal W}_{canon}$
 is injective over $Y_{\Gamma}\times T(M)$.
$\Box$

\medskip

\subsubsection{A Short Remark about the Torsion Free Sub-sheaves}\label{subsubsection; torsion}

  Let ${\cal F}$ be a coherent sheaf over a smooth and connected scheme $Y$. Let 
$U=Spec(R)$ be an affine open subspace of $Y$. Then $R$ is an integral domain and
 let $K$ denote the field of fractions. Over $U$ the coherent sheaf ${\cal F}$
is the sheaf associated to a finite $R$-module $N$.
Recall that (e.g [Fr]) the
 generic rank of ${\cal F}$ is defined to be the rank of $N\otimes_R K$.

 Define ${\cal F}^{\ast}={\cal HOM}_{{\cal O}_Y}({\cal F}, {\cal O}_Y)$ to be
 the dual sheaf. Then there is a natural morphism ${\cal F}\mapsto 
{\cal F}^{\ast\ast}$. Define ${\cal F}_{torfree}$ to be the image of
 ${\cal F}$ in ${\cal F}^{\ast\ast}$, and it is a torsion free sheaf.

 On the other hand by corollary 21 on page 44 of [Fr] 
we have the following short exact sequence,
 $$0\mapsto {\cal F}_{tor}\mapsto {\cal F}\mapsto {\cal F}_{torfree}\mapsto 0,$$

 where the cokernel ${\cal F}_{tor}$ is the torsion sub-sheaf of ${\cal F}$ 
over $Y$. From now on, we call ${\cal F}_{torfree}$ the torsion free quotient
(part) of ${\cal F}$.

In general the inclusion ${\cal F}_{torfree}\subset {\cal F}^{\ast\ast}$ is
 not always an
 equality. A torsion free 
sheaf ${\cal F}\cong {\cal F}^{\ast\ast}$ under the injection is
 called a reflexive sheaf. It is well known that locally free sheaves are
 reflexive. 

\begin{lemm}\label{lemm; unique}
Let $Y$ be a reduced, smooth and connected scheme.
Let ${\cal F}$ be coherent and let ${\cal E}$ be locally free such that
 generic rank of ${\cal F}$=rank of ${\cal E}$.
Let ${\cal F}\mapsto {\cal E}\mapsto 0$ be a sheaf surjection, then
 ${\cal F}_{torfree}\cong {\cal E}$.
\end{lemm}

\noindent Proof of lemma \ref{lemm; unique}: Consider the following
 commutative diagram, 

\[
\begin{array}{ccccc}
{\cal F}& \longrightarrow &  {\cal E}  & \longrightarrow & 0\\
 \Big\downarrow &  &  \Big\downarrow\vcenter{%
\rlap{$\scriptstyle{\mathrm{\cong}}\,$}} & & \\
{\cal F}^{\ast\ast} & \longrightarrow & {\cal E}^{\ast\ast} & & 
\end{array}
\]

The surjectivity of ${\cal F}\mapsto {\cal E}\cong{\cal E}^{\ast\ast}$
 implies that the image of ${\cal F}$ in ${\cal F}^{\ast\ast}$,
 ${\cal F}_{torfree}$, maps 
 surjectively onto ${\cal E}^{\ast\ast}$. Because ${\cal F}_{torfree}$ 
and ${\cal E}^{\ast\ast}$ are of the same generic rank, the kernel
 of ${\cal F}_{torfree}\mapsto {\cal E}^{\ast\ast}\mapsto 0$
 has to be a torsion sheaf. But a torsion sheaf can never maps injectively
 into a torsion free sheaf ${\cal F}_{torfree}$. So it must vanish and
 ${\cal F}_{torfree}\cong {\cal E}^{\ast\ast}$ and therefore
 ${\cal F}_{torfree}$ is isomorphic to ${\cal E}$ itself. $\Box$

 In other words, ${\cal E}$ is isomorphic to the torsion free quotient of 
 ${\cal F}$.

\begin{lemm}\label{lemm; split}
 Let $Y$ be a smooth, connected and reduced scheme.
 Let ${\cal F}$ be coherent and let ${\cal E}$ be locally free and
 generic rank of ${\cal F}$=rank of ${\cal E}$.
 Let $0\mapsto {\cal E}\mapsto {\cal F}$ be a sheaf injection such that
 $\bullet \otimes k(y)$ is left exact for all the closed points $y\in Y$.
 Then ${\cal E}\cong {\cal F}_{torfree}$ under the composition
 ${\cal E}\mapsto {\cal F}\mapsto {\cal F}_{torfree}$.
\end{lemm}

 In such a case, ${\cal E}\cong {\cal F}_{torfree}$
 induces a morphism ${\cal F}_{torfree}\mapsto {\cal F}$
 through ${\cal E}$ and  we may write ${\cal F}={\cal F}_{torfree}
\oplus {\cal F}_{tor}$ and
 call ${\cal F}_{torfree}$ and ${\cal F}_{tor}$ the torsion free and the 
torsion summands of ${\cal F}$.

\noindent Proof:  Consider the following commutative diagram,

\[
\begin{array}{ccc}
{\cal E} & \longrightarrow & {\cal F}  \\
\Big\downarrow\vcenter{%
\rlap{$\scriptstyle{\mathrm{\cong}}\,$}} &  &  \Big\downarrow \\
{\cal E}^{\ast\ast} & \longrightarrow & {\cal F}^{\ast\ast} 
\end{array}
\]

To show that ${\cal E}\cong {\cal F}_{torfree}$, it suffices to show
that ${\cal E}\cong {\cal E}^{\ast\ast}\cong {\cal F}^{\ast\ast}$.

  Because ${\cal E}$ and ${\cal F}$ are of the same generic rank, the
 cokernel of the short exact sequence 
 $0\mapsto {\cal E}\mapsto {\cal F}\mapsto {\cal R}\mapsto 0$ is
 torsion. Applying the
 contravariant left exact functor ${\cal HOM}_{{\cal O}_Y}(\bullet, {\cal O}_Y)$
 to this short exact sequence and observing 
${\cal HOM}_{{\cal O}_Y}({\cal R}, {\cal O}_Y)=0$ because of its torsion nature,
we have

 $$0\mapsto {\cal F}^{\ast} \mapsto {\cal E}^{\ast}\mapsto {\cal R}'\mapsto 0.$$
 The cokernel ${\cal R}'$ appears as the 
${\cal HOM}_{{\cal O}_Y}(\bullet, {\cal O}_Y)$ functor is not always right exact.
 By proposition 24 on page 45 of [Fr], ${\cal F}^{\ast}$ is reflexive and
 torsion free. On the other hand the $k(y)$ vector space morphism 
${\cal F}^{\ast}\otimes k(y)\mapsto
{\cal E}^{\ast}\otimes k(y)$ is equivalent to
 $HOM_{k(y)}({\cal F}\otimes k(y), k(y))\mapsto 
HOM_{k(y)}({\cal E}\otimes k(y), k(y))$ and is always surjective since
 ${\cal E}\otimes k(y)\mapsto {\cal F}\otimes k(y)$ is always injective for
 all the closed points $y\in Y$. As $\bullet \otimes k(y)$ is right
 exact, the $rank_{k(y)}{\cal R}'\otimes k(y)=0$ for all the closed points $y$
 of $Y$.
  Then by exercise II.5.8 on page 125 of [Ha] and the fact that $Y$ is reduced, 
${\cal R}'$ is locally free of rank
 $0$, so ${\cal R}'$ vanishes. 

 Therefore ${\cal F}^{\ast}\cong {\cal E}^{\ast}$. By dualizing this equality
again we get the desired ${\cal E}^{\ast\ast}\cong {\cal F}^{\ast\ast}$.
 The lemma is proved.
$\Box$

\medskip

  The reader may consult page 42-46 of [Fr] for some basic knowledge about 
 torsion free sheaves.

  In our paper we consider the torsion free quotient of coherent sheaves
 which are the derived image sheaves of invertible sheaves along a union of 
${\bf P}^1$ fibrations or along finite morphisms. 
In the proof of proposition \ref{prop; id}, we show that these torsion
 free sheaves are in fact locally free. When this situation occurs, we denote 
 ${\cal F}_{torfree}$ by an alternative notation ${\cal F}_{free}$ (or
 $({\cal F})_{free}$) to indicate that it is not only torsion free, but is
 actually locally free.

\medskip

\subsection{The Construction of ${\bf P}^1$ fiber bundles 
$\tilde{\Xi}_i\mapsto Y(\Gamma_{e_i})$}
\label{subsection; realcase}

\medskip

 Let $\Gamma\in adm(n)$ be an $n$-vertex admissible graph and let 
$Y(\Gamma)\subset M_n$ be the smooth closure of the corresponding admissible
 stratum $Y_{\Gamma}$. As usual let $m_1, m_2, \cdots, m_n$ be the
multiplicities satisfying $0<m_a\leq m_b$ whenever $1\leq a\leq b\leq n$. 
We assume that
such a multiple covered ${\bf M}(E)E=\sum_{1\leq i\leq n}m_i E_i$ has been fixed.

  Let $\pi:\Xi_i\mapsto Y(\Gamma_{e_i})$ be
 \footnote{The construction of $\Xi_i$ will be outlined
 in the proof of lemma \ref{lemm; sum}.}
 the ${\bf P}^1$ fibration 
  in $M_{n+1}\times_{M_n}Y(\Gamma_{e_i})$ representing the universal 
exceptional curves dual to $e_i$. 
 In the following we would like to construct smooth
 ${\bf P}^1$ fiber bundles 
$\tilde{\pi}:\tilde{\Xi}_i\mapsto Y(\Gamma_{e_i})$ birational to $\Xi_i$
for $1\leq i\leq n$. To simplify our notations,
 we would like to  drop the restriction symbol\footnote{The reader should
be able to recover the restriction notation from the base space we are using.} and 
denote their restrictions to the sub-locus 
$Y(\Gamma)=\cap_{1\leq i\leq p}Y(\Gamma_{e_i})\subset
 Y(\Gamma_{e_i})$, $\Xi_i|_{Y(\Gamma)}$ or
 $\tilde{\Xi}_i|_{Y(\Gamma)}$ by the same notations.

 Given a subscript $1\leq i\leq n$, we define $I_i$ in the following,

\begin{defin}\label{defin; indexset}
Define the index set $I_i$ to be the set of all the subscripts of $E$ appearing in
 $e_i=E_i-\sum_{j_i}E_{j_i}$. I.e. the union of $\{i\}$ and all the direct
descendent indexes of $i$ in $\Gamma$.
\end{defin}

 Given an index subset $I\subset \{1, 2, \cdots, n\}$, by 
corollary \ref{cor; lifting} in section \ref{section; strata} there exists
 the canonical lifting 
$\check{\pi}_I: M_n\mapsto M_{|I|}$ of $\pi_I:M_n\mapsto M^{|I|}$ and it induces 
the canonical map $Y(\Gamma)\mapsto M_{|I|}$ by composing $Y(\Gamma)\hookrightarrow M_n$
 and $M_n\mapsto M_{|I|}$.

 By taking $I=\{1, 2, \cdots, i-1\}\cup I_i$ in the above setting, 
we may construct the total space of the ${\bf P}^1$ fibration 
$\tilde{\Xi}_i$ as a divisor in the fiber product 
$M_{i+|I_i|}\times_{M_{i-1+|I_i|}}Y(\Gamma_{e_i})$.

\begin{lemm}\label{lemm; sum}
There exists an ${\bf P}^1$ fibration over $Y(\Gamma_{e_i})$, $\tilde{\pi}:
\tilde{\Xi}_i\mapsto Y(\Gamma_{e_i})$ such that 

(i). $\tilde{\Xi}_i$ is pulled back from a ${\bf P}^1$ sub-fibration of
 $M_{i+|I_i|}\mapsto M_{i-1+|I_i|}$ by $Y(\Gamma_{e_i})\stackrel{\check{\pi}_I}{\mapsto} 
M_{|I|}=M_{i-1+|I_i|}$.

\medskip

(ii). $\tilde{\Xi}_i\mapsto Y(\Gamma_{e_i})$ has a structure of ${\bf P}^1$ fiber 
bundle over $Y(\Gamma_{e_i})$.

\medskip

(iii).
The space $\Xi_i$ maps birationally onto $\tilde{\Xi}_i$ and the birational map $\Xi_i\mapsto \tilde{\Xi}_i$
is a consecutive blowing ups along codimension two smooth centers.
\end{lemm}

 Thus $\tilde{\Xi}_i$ is the relative minimal model of $\Xi_i$.

\medskip
 
\noindent Proof of lemma \ref{lemm; sum}: Because $I=\{1, 2, \cdots, i-1\}\cup I_i$, 
 $|I|=i-1+|I_i|$.
Consider an ordering preserving bijection
 $\phi:I\mapsto \{1, 2, \cdots , |I|\}$. Then $\phi(j)=j$ for $j\leq i$.
Consider an $|I|-$vertex fan-like admissible graph
 $\Gamma_i\in adm(|I|)$ with one-edges from $i-$th vertex to all the vertexes marked in
 $\phi(I_i-\{i\})$. Then the $i-$th vertex is the direct
 ascendent of all the other vertexes in $\Gamma_i$ marked by $\phi(I_i-\{i\})$ 
and it is the only direct ascendent vertex in $\Gamma_i$. So $Y(\Gamma_{e_i})$
 is mapped onto $Y(\Gamma_i)$ under $Y(\Gamma_{e_i})\mapsto M_{|I|}$. This can be seen
 by the construction of ${\bf Y}(\Gamma_{e_i})$ in lemma \ref{lemm; fanlike} of section
 \ref{section; strata} as the pre-image \footnote{Please refer to the proof of
 lemma \ref{lemm; fanlike} for more details.}
 of ${\bf Y}(\Gamma_i)\cong
 \times_{M_i/M_{i-1}}^{|I_i|-1} D_i$ under\footnote{The space $D_i\mapsto M_i$ is the 
exceptional divisor by blowing up along $\Delta_{M_i}:M_i\subset M_i\times_{M_{i-1}}M_i$ and has 
 a ${\bf P}^1$ bundle structure over $M_i$.} $\check{\pi}_{I}$. 
Under the surjection $Y(\Gamma_{e_i})\mapsto Y(\Gamma_i)$, 
 the open subset $Y_{\Gamma_{e_i}}$ is mapped surjectively onto $Y_{\Gamma_i}$.

 As a subspace of the space $M_{|I|}$, $Y(\Gamma_i)$ is characterized as the existence locus
 of the type $I$ exceptional class $E_i-\sum_{|I_i|\geq j\geq 2}E_{j+i-1}$. So there exists
 a ${\bf P}^1$ fibration of universal curves 
over $Y(\Gamma_i)\subset M_{|I|}$, ${\cal C}_i\mapsto Y(\Gamma_i)$,
 whose fiber over $b\in 
 Y(\Gamma_i)\subset M_{|I|}$ is the type $I$ exceptional curve representing
 $E_i-\sum_{|I_i|\geq j\geq 2}E_{j+i-1}$ in the algebraic surface $M_{|I|+1}|_b$. So
 ${\cal C}_i\mapsto Y(\Gamma_i)$ can be viewed as a sub-fibration of $M_{|I|+1}\times_{M_|I|}
Y(\Gamma_i)$.  By pulling back ${\cal C}_i\mapsto Y(\Gamma_i)$ by $Y(\Gamma_{e_i})\mapsto Y(\Gamma_i)$,
 we define $\tilde{\Xi}_i$ to be ${\cal C}_i\times_{Y(\Gamma_i)}Y(\Gamma_{e_i})$. Then 
 the condition (i). holds by our construction.

\medskip

  To prove (ii)., it suffices to show 
\footnote{An alternative way to achieve this is to check that the type $I$
 class $E_i-\sum_{2\leq j\leq |I_i|}E_{j+i-1}$ can not be broken into
 two distinct type $I$ classes.}
that ${\cal C}_i\mapsto Y(\Gamma_i)$ is a ${\bf P}^1$ 
fiber bundle. By a special case of lemma \ref{lemm; fuse}, we have $M_{|I|+1}\cong (M_{i+1}/M_i)_{|I|-i+1}=
(M_{i+1}/M_i)_{|I_i|}$.
  On the other hand, the exceptional divisor $D_i\subset M_{i+1}/M_i$ blown up from the
 relative diagonal $M_i\mapsto M_i\times_{M_{i-1}}M_i$ has a ${\bf P}^1$ bundle structure over $M_i$.
 So we have $(D_i/M_i)_{|I_i|}\subset (M_{i+1}/M_i)_{|I_i|}$. On the other hand, we have
 the commutative diagram,

\[
\begin{array}{ccccc}
(D_i/M_i)_{|I_i|}& \longrightarrow & (M_{i+1}/M_i)_{|I_i|} & \cong & M_{|I_i|+i}\\
\Big\downarrow & & \Big\downarrow & & \\
(D_i/M_i)_{|I_i|-1} &\longrightarrow & (M_{i+1}/M_i)_{|I_i|-1} & \cong & M_{|I_i|+i-1} 
\end{array}
\]

 Because $D_i\mapsto M_i$ is smooth of relative dimension one\footnote{Here we are using the
 fact that the codimension one blowing ups are trivial.}, the projection of fiber products
 $\times_{M_i}^{|I_i|}D_i \cong (D_i/M_i)_{|I_i|}\mapsto  (D_i/M_i)_{|I_i|-1}\cong 
\times_{M_i}^{|I_i|-1}D_i\cong Y(\Gamma_i)$
 has a ${\bf P}^1$ fiber bundle structure. It suffices to identify 
 $(D_i/M_i)_{|I_i|}$ with ${\cal C}_i$.

 By induction it is easy to see that the fiber bundle 
$M_{|I|+1}\mapsto M_{|I|}$ can be constructed from the trivial
 bundle $M_1\times_{M_0}M_{|I|}=
M\times M_{|I|}\mapsto M_{|I|}$ by $|I|-$ consecutive blowing ups along the cross 
sections \footnote{They are pull-backs of the relative diagonals $\Delta_{M_{k+1}/M_k}:M_{k+1}\mapsto 
 M_{k+1}\times_{M_k}M_{k+1}$ by $f_{|I|-1, k+1}:M_{|I|}\mapsto M_{k+1}$.} of the intermediate fiber bundles
 $M_{k+1}\times_{M_k}M_{|I|}\mapsto M_{|I|}$, for $0\leq k\leq |I|-1$.
  Consider $f_{|I|-1, i}:M_{|I|}\mapsto M_i$ and the pulled-back fiber bundle
$f_{|I|-1, i}^{\ast}D_i\subset M_{i+1}\times_{M_i}M_{|I|}$ is isomorphic to  
the exceptional divisor $E_i$ of the
 $i-$th intermediate fiber bundle. Thus
 the projection $M_{|I|+1}/M_{|I|}\mapsto M_{i+1}\times_{M_i}M_{|I|}/M_{|I|}$ to the $i-$th
 intermediate space can be constructed by
 $|I|-i=|I_i|-1-$consecutive blowing ups along 
cross sections of the intermediate fiber bundles. When we restrict to the locus 
 $Y(\Gamma_i)\subset M_{|I|}$, the $j-$th cross section,
 for all $1\leq j\leq |I_i|-1$, maps into the sub-bundle
 $f_{|I|-1, i}^{\ast}D_i\mapsto M_{|I|}$
 and becomes a cross section of $f_{|I|-1, i}^{\ast}D_i|_{Y(\Gamma_i)}\mapsto
 Y(\Gamma_i)$.

 On the other hand, by Chapter II corollary 7.15. of [Ha], 
 the strict transform of the restriction of the exceptional divisor $E_i\times_{M_{|I|}}Y(\Gamma_i)\subset 
 M_{i+1}\times_{M_i}Y(\Gamma_i)$ inside $M_{|I|+1}\times_{M_{|I|}}Y(\Gamma_i)$ 
is nothing but the $|I_i|-1-$consecutive blowing ups
 of $f_{|I|-1, i}^{\ast}D_i|_{Y(\Gamma_i)}$ along the $|I_i|-1$ distinct cross sections.
 Because the ${\bf P}^1$ fiber bundle is smooth of relative dimension one, all the
 blowing ups along cross sections are trivial. Thus
 its strict transform ${\cal C}_i$, 
representing $E_i-\sum_{|I_i|\geq j\geq 2}E_{j+i-1}$ in 
$M_{|I|+1}\times_{M_{|I|}}Y(\Gamma_i)$, is still isomorphic to 
 $f_{|I|-1, i}^{\ast}D_i|_{Y(\Gamma_i)}$. 
 The condition (ii). is proved.

 In the following, we derive the conclusion (iii) based on a similar argument as above.
 Consider the projection map $f_{n-1, i}:M_n\mapsto M_i$ and the induced ${\bf P}^1$ bundle
$f_{n-1, i}^{\ast}D_i\subset M_{i+1}\times_{M_i}M_n$ is the exceptional divisor $E_i$ of the
 $i-$th intermediate fiber bundle in-between $M_{n+1}$ and the trivial product $M\times M_n$.
 Similar to the above argument 
the map $M_{n+1}\mapsto M_{i+1}\times_{M_i}M_n$ can be constructed by $n-i-$consecutive blowing ups
 along cross sections of the $n-i-$intermediate fiber bundles. Similar to the above discussion to
 ${\cal C}_i$,
 $\Xi_i$ is the strict transform of $E_i\times_{M_n}Y(\Gamma_{e_i})\subset M_{i+1}\times_{M_i}M_n$ under
 these consecutive blowing ups. Again by Chapter II corollary 7.15. of [Ha], $\Xi_i$ can be
 identified with the $n-i-$
consecutive blowing ups of $F_0=f_{n-1, i}^{\ast}D_i\times_{M_n}Y(\Gamma_{e_i})$ along the 
 intersections (of the intermediate blown up spaces from $F_0$) 
with the various cross sections in the intermediate fiber bundles 
\footnote{Pulled-back from $M_{k+i+1}\mapsto M_{k+i}$ by $M_n\mapsto M_{k+i}$.}
$M_{k+i+1}\times_{M_{k+i}}M_n\mapsto M_n$. Denote $C_0\subset F_0$ to be the first blowup
 center and inductively define $F_k=BlowUp_{C_{k-1}}F_{k-1}$. 
Denote $C_k\subset F_k$ to be
 the $k-$th blowup center, for $k$ ranging in $0\leq k\leq n-i-1$. 
At the end we have $F_{n-i}=\Xi_i$ and it suffices to show that
 all the blowup centers $C_k$ ($0\leq k\leq n-i-1$) 
are smooth of codimension two/one in $F_k$.

 Because $C_k$ is the intersection of 
$F_k$ with a cross section of the ambient fiber bundle
 $M_{k+i+1}\times_{M_{k+i}}M_n\mapsto M_n$, 
the projection $C_k\subset F_k \mapsto Y(\Gamma_{e_i})$
 induces an isomorphism onto the image of the intersection locus.
Suppose that $C_k$ maps onto $Y(\Gamma_{e_i})$, then
$C_k$ must be a cross section of $F_k\mapsto Y(\Gamma_{e_i})$ and therefore 
\footnote{By proposition \ref{prop; stratification} the space $Y(\Gamma_{e_i})$ is smooth.}
 is smooth. This can only occur when $k+i+1$ is a direct descendent index 
of $i$ throughout $Y(\Gamma_{e_i})$,
 which happens only when $k+i+1\in I_i$. If $k+i+1\not\in I_i$, then $k+i+1$ is not a direct
descendent index of $i$ in $\Gamma_{e_i}$. Consider a particular 
degeneration $\Gamma_{e_i; k+i+1}$ of $\Gamma_{e_i}$ 
by adding a single one-edge from $i$ to $k+i+1$. Then by
 proposition \ref{prop; stratification}, $Y(\Gamma_{e_i; k+i+1})\subset Y(\Gamma_{e_i})$ is
 a smooth divisor in $Y(\Gamma_{e_i})$. On the other hand, $C_k$ is the 
 intersection of $F_k$, i.e. the strict transform of $F_0$,
 with the cross section of $M_{k+i+1}\times_{M_{k+i}}M_n\mapsto M_n$ induced
 by the relative diagonal. So at the location where $F_k$ intersects the
 cross section, the $k+i+1-$th blowing up in $F_k$ determined by the
 intersection locus is located in the strict transform of the 
exceptional locus $E_i$ of the
 $i-$th blowing up. By the interpretation of remark \ref{rem; interpret} 
this occurs exactly when $k+i+1$ becomes a (direct or indirect)
 descendent of $i$ and so
$C_k$ maps onto $Y(\Gamma_{e_i; k+i+1})$. Therefore $C_k$ is 
isomorphic to $Y(\Gamma_{e_i; k+i+1})\subset Y(\Gamma_{e_i})$ and 
by proposition \ref{prop; stratification} it is smooth. When this occurs $C_k$ is of codimension
 two in $F_k$. Then by induction $\Xi_i=F_{n-k}$ is an $n-k-$consecutive blowing up of $F_0$ along 
codimension two smooth centers.

  Because $f_{|I|-1, i}\circ f_{n-1, |I|}=f_{n-1, i}$, we have 
$\tilde{\Xi}_i=(f_{n-1, |I|}|_{Y(\Gamma_{e_i})})^{\ast}{\cal C}_i=
(f_{n-1, |I|}|_{Y(\Gamma_{e_i})})^{\ast}(f_{|I|-1, i})^{\ast}D_i=
(f_{n-1, i}|_{Y(\Gamma_{e_i})})^{\ast}D_i=F_0$. So 
$\Xi_i=F_{n-k}$ projects onto $F_0=\tilde{\Xi}_i$ and $\tilde{\Xi}_i$ is the
 relative minimal model of $\Xi_i$. This finishes the proof of $(iii)$.   $\Box$

\medskip

\begin{rem}\label{rem; specialfiber}
Because some of these $C_k$ are not cross sections and are not dominating $Y(\Gamma_{e_i})$, the blowing ups along those $C_k$ cause the special fibers of $\Xi_i\mapsto Y(\Gamma_{e_i})$ to become
 a finite tree of normal-crossing ${\bf P}^1$.
\end{rem}

\medskip

 For a fixed $i$ one may re-write
 the cohomology class $C-{\bf M}(E)E=C-\sum_{1\leq a\leq n}m_aE_a$
 as $C-\sum_{a\in I_i}m_aE_a-\sum_{a\not\in I_i}m_aE_a$ and there is a canonical
 (up to rescaling of ${\bf C}^{\ast}$) sheaf morphism ${\cal E}_{C-{\bf M}(E)E}=
{\cal E}_{C-\sum_{a\in I_i}m_aE_a-\sum_{a\not\in I_i}m_aE_a}\mapsto 
 {\cal E}_{C-\sum_{a\in I_i}m_aE_a-\sum_{a<i}m_aE_a}$ by tensoring with the defining sections of
 $\sum_{i<a\not\in I_i}m_aE_a$ on ${\cal M}_{n+1}$.

 The main reason that we introduce the fibrations $\tilde{\Xi}_{k_i}$, $1\leq i\leq p$ for
 $e_{k_i}$ is because of the following,

\begin{lemm}\label{lemm; vanish}
  The sheaf  ${\cal R}^0\tilde{\pi}_{\ast}\bigl({\cal O}_{\tilde{\Xi}_{k_i}}
\otimes {\cal E}_{C-\sum_{a\in I_i}m_aE_a-\sum_{a<i}m_aE_a}\bigr)$ is the zero sheaf.
  The first derived image sheaf ${\cal R}^1\tilde{\pi}_{\ast}\bigl({\cal O}_{\tilde{\Xi}_{k_i}}
\otimes {\cal E}_{C-\sum_{a\in I_i}m_aE_a-\sum_{a<i}m_aE_a}\bigr)$ is locally free.
\end{lemm}

\medskip

\noindent Proof:  Based on the condition $e_{k_i}\cdot (C-{\bf M}(E)E)<0$,
 the relative degree of the invertible sheaf 
${\cal F}={\cal E}_{C-\sum_{a\in I_i}m_aE_a-\sum_{a<i}m_aE_a}$ along $\tilde{\Xi}_{k_i}\mapsto
 Y(\Gamma_{e_{k_i}})$ is negative. Since by lemma \ref{lemm; sum}, 
$\tilde{\Xi}_{k_i}\mapsto Y(\Gamma_{e_{k_i}})$ 
is a ${\bf P}^1$ fiber bundle. The negativity of the relative degree
 implies the vanishing of ${\cal R}^0\tilde{\pi}_{\ast}\bigl({\cal O}_{\tilde{\Xi}_{k_i}}
\otimes {\cal E}_{C-\sum_{a\in I_i}m_aE_a-\sum_{a<i}m_aE_a}\bigr)$. On the other hand,
 the vanishing of the zero-th derived image sheaf implies $h^0(y, {\cal F})=0$ for 
$y\in Y(\Gamma_{e_{k_i}})\times T(M)$ and
 by curve Riemann-Roch formula it implies the constancy of $h^1(y, {\cal F})$ throughout
 $Y(\Gamma_{e_{k_i}})\times T(M)$. Then by chapter II, corollary 12.9 of [Ha] the sheaf 
${\cal R}^1\tilde{\pi}_{\ast}\bigl({\cal E}_{C-\sum_{a\in I_i}m_aE_a-\sum_{a<i}m_aE_a}\bigr)$
 is locally free. $\Box$

\medskip

  The locally freeness of the derived image sheaves along the various 
$\tilde{\Xi}_{k_i}\mapsto Y(\Gamma_{e_{k_i}})$ enables us
 to find an explicit representative \footnote{For the definition
 of ${\cal V}_{quot}$, please consult page \pageref{quot}.}
 of $[{\cal V}_{quot}]\in K_0(Y(\Gamma)\times T(M))$ in
 the following subsection.

\medskip

\subsection{The Locally Freeness of 
${\cal V}_{quot}$ and its Explicit Representative in the $K$ Group}
\label{subsection; lf}

\bigskip

 In this subsection, we would like to prove the locally freeness of the torsion
 free quotient
 ${\cal V}_{quot}$ of ${\cal R}^1\pi_{\ast}\bigl({\cal O}_{\sum_{1\leq i\leq p}\Xi_{k_i}}\otimes 
 {\cal E}_{C-{\bf M}(E)E}\bigr)$ over $Y(\Gamma)\times T(M)$ and 
 we also
  give an explicit identification of $[{\cal V}_{quot}]\in K_0(Y(\Gamma)\times T(M))$.

\begin{lemm}\label{lemm; freepart}
The torsion free quotient of the coherent sheaf 
 ${\cal R}^1\pi_{\ast}\bigl({\cal O}_{\Xi_{k_p}}\otimes
 {\cal E}_{C-\sum_{a\leq n}m_aE_a}\bigr)$ is locally free and is isomorphic to
\footnote{The symbol $E_{a; b}$, $a<b$, denote the exceptional divisor in 
 $M_n$ by blowing up the the $(a, b)-$th partial diagonal.} 
 ${\cal R}^1\tilde{\pi}_{\ast}\bigl({\cal O}_{\tilde{\Xi}_{k_i}}\otimes
 {\cal E}_{C-\sum_{a\in I_{k_i}}m_aE_a}\bigr)\otimes \tilde{\pi}^{\ast}
{\cal O}(-\sum_{a<k_i}m_aE_{a; k_i})$.
\end{lemm}

\noindent Proof of the lemma: The above sheaves are of the same generic 
rank, by a direct curve Riemann-Roch calculation along smooth fibers above
 $Y_{\Gamma}\times T(M)$

The locally freeness of 
 ${\cal R}^1\tilde{\pi}_{\ast}\bigl({\cal O}_{\tilde{\Xi}_{k_i}}\otimes
{\cal E}_{C-\sum_{a\in I_{k_i}}m_aE_a-\sum_{a<k_i}m_aE_a}\bigr)$ has been proved in
 lemma \ref{lemm; vanish}. The isomorphism 

$$\hskip -1.3in
{\cal R}^1\tilde{\pi}_{\ast}\bigl({\cal O}_{\tilde{\Xi}_{k_i}}\otimes
{\cal E}_{C-\sum_{a\in I_{k_i}}m_aE_a-\sum_{a<k_i}m_aE_a}\bigr)\cong
{\cal R}^1\tilde{\pi}_{\ast}\bigl({\cal O}_{\tilde{\Xi}_{k_i}}\otimes
{\cal E}_{C-\sum_{a\in I_{k_i}}m_aE_a}\bigr)\otimes \tilde{\pi}^{\ast}
{\cal O}(-\sum_{a\in k_i}m_aE_{a; k_i})$$

 follows from the projection formula (exercise II.8.3. on page 253 of [Ha]) and the
 fact that ${\cal O}(-E_a)|_{\tilde{\Xi}_{k_i}}=\tilde{\pi}^{\ast}
{\cal O}(E_{a; k_i})|_{Y(\Gamma)}$, for $a<k_i$.

 To prove the lemma, by lemma \ref{lemm; unique} it suffices to prove that
 ${\cal R}^1\pi_{\ast}\bigl({\cal O}_{\Xi_{k_p}}\otimes
 {\cal E}_{C-\sum_{1\leq a\leq n}m_aE_a}\bigr)\mapsto 
{\cal R}^1\tilde{\pi}_{\ast}\bigl({\cal O}_{\tilde{\Xi}_{k_i}}\otimes
 {\cal E}_{C-\sum_{a\in I_{k_i}}m_aE_a-\sum_{a<k_i}m_aE_a}\bigr)$ is surjective.

 Firstly by lemma \ref{lemm; sum} (iii). $\Xi_{k_i}\mapsto \tilde{\Xi}_{k_i}$ is a composite
blowing down map, the push-forward of 
${\cal O}_{\Xi_{k_i}}(-\sum_{n\geq a\geq k_i; a\not\in I_{k_i}}m_aE_a)$
 to $\tilde{\Xi}_{k_i}$ defines an ideal sheaf \footnote{The subscript $t$ of the 
notation $Z_t$ stands for ``torsion'' because $Z_t$ is closely related to the
 torsion part of a sheaf.} of the sub-scheme ${\cal I}_{Z_t}\subset
 {\cal O}_{\tilde{\Xi}_{k_i}}$. 
 To show the surjectivity of the original sheaf map, it suffices to show that $Z_t\mapsto
 Y(\Gamma_{e_{k_i}})$ is at most of relative dimension zero. I.e. the fibers of $Z_t\mapsto 
 Y(\Gamma_{e_{k_i}})$ are either empty or are zero dimensional.
In fact $\tilde{\Xi}_{k_i}$
 is a ${\bf P}^1$ fiber bundle and all the fibers are smooth and irreducible. On the other 
 hand $Z_t$ supports over the image of $\sum_{n\geq a\geq k_i; a\not\in I_{k_i}}m_aE_a$
 in $\tilde{\Xi}_{k_i}$.
  So the only chance for $Z_t|_y$, $y\in Y(\Gamma)$ to be one dimensional
 is when $Z_t|_y$ supports over the whole $\tilde{\Xi}_{k_i}|_y$. This implies that
 the defining section of ${\cal O}(\sum_{n\geq a\geq k_i; a\not\in I_{k_i}}m_aE_a)|_y$ is
 divisible by the defining section of 
$(E_{k_i}-\sum_{j_{k_i}}E_{j_{k_i}})|_y$. This implies the existence of
 an admissible graph $\Gamma', \Gamma'<\Gamma$ and $y\in Y_{\Gamma'}$,
 such that $k_i$ is a descendent of
 some $a\geq k_i, a\not\in I_{k_i}$. This is absurd because the axiom 2. of admissible 
graphs forbids 
 $k_i$ with $k_i<a$ to be the descendent of $a$.
 
Therefore
 $Z_t\mapsto Y(\Gamma)$ is generically empty and is 
at most relative dimension zero.  Then by tensoring the defining
 exact sequence
 $0\mapsto {\cal I}_{Z_t}\mapsto {\cal O}_{\tilde{\Xi}_{k_i}}\mapsto 
{\cal O}_{Z_t}\mapsto 0$ with ${\cal E}_{C-\sum_{a\in I_{k_i}}m_aE_a-\sum_{a<k_i}m_aE_a}$, 
 we get the desired surjectivity from a portion of its derived long exact sequence,

$$
{\cal R}^1\tilde{\pi}_{\ast}\bigl({\cal I}_{Z_t}\otimes
{\cal E}_{C-\sum_{a\in I_{k_i}}m_aE_a-\sum_{a<k_i}m_aE_a}\bigr)\mapsto
{\cal R}^1\tilde{\pi}_{\ast}\bigl({\cal O}_{\tilde{\Xi}_{k_i}}\otimes 
{\cal E}_{C-\sum_{a\in I_{k_i}}m_aE_a-\sum_{a<k_i}m_aE_a}\bigr)$$
$$\mapsto 
 {\cal R}^1\tilde{\pi}_{\ast}\bigl({\cal O}_{Z_t}\otimes
{\cal E}_{C-\sum_{a\in I_{k_i}}m_aE_a-\sum_{a<k_i}m_aE_a}\bigr)$$

 and the vanishing of ${\cal R}^1\tilde{\pi}_{\ast}\bigl({\cal O}_{Z_t}\otimes
{\cal E}_{C-\sum_{a\in I_{k_i}}m_aE_a-\sum_{a<k_i}m_aE_a}\bigr)$
  by the fact that $Z_t\mapsto Y(\Gamma)$ is 
``at most relative dimension zero''.

 Once the surjectivity
 has been achieved, this surjection induces
 an isomorphism between the torsion free quotients of 
${\cal R}^1\pi_{\ast}\bigl({\cal O}_{\Xi_{k_i}}\otimes
 {\cal E}_{C-\sum_{a\leq n}m_aE_a}\bigr)$, $\cong {\cal R}^1\tilde{\pi}_{\ast}\bigl({\cal I}_{Z_t}\otimes
{\cal E}_{C-\sum_{a\in I_{k_i}}m_aE_a-\sum_{a<k_i}m_aE_a}\bigr)$,
 with ${\cal R}^1\tilde{\pi}_{\ast}\bigl({\cal O}_{\tilde{\Xi}_{k_i}}\otimes
 {\cal E}_{C-\sum_{a\in I_{k_i}}m_aE_a-\sum_{a<k_i}m_aE_a}\bigr)$, by the
 argument in lemma \ref{lemm; unique}. $\Box$

\medskip

 The following lemma will be used frequently in the following discussion.

\begin{lemm}\label{lemm; tor}
 Let ${\cal G}_i, 0\leq i\leq 4$ be five coherent sheaves over a smooth, connected 
and reduced scheme $Y$ and let
 $ {\cal G}_0\mapsto {\cal G}_1\mapsto 
{\cal G}_2\mapsto {\cal G}_3\mapsto
 {\cal G}_4$ be a sheaf exact sequence. Suppose that both of ${\cal G}_0$ and ${\cal G}_4$
 are torsion sheaves and $({\cal G}_2)_{torfree}$ is locally free.
  Suppose additionally that the induced morphism 
$({\cal G}_1)_{torfree}\otimes k(y)\mapsto ({\cal G}_2)_{free}\otimes k(y)$
 is injective for all the closed points $y\in Y$, then the torsion free quotients
 of ${\cal G}_1, {\cal G}_2, {\cal G}_3$ are all locally free and they
 form a short exact sequence of locally free sheaves,

$$0\mapsto ({\cal G}_1)_{free}\mapsto ({\cal G}_2)_{free}\mapsto 
({\cal G}_3)_{free}\mapsto 0.$$

\end{lemm}

\noindent Proof:  It is well known that
 any morphism from a torsion sheaf to a torsion free sheaf is
 trivial and any morphism from a torsion free sheaf to a torsion sheaf
 cannot be injective. By taking the double-duals of the original sequence,
 we get 

\[
\begin{array}{ccccc}
{\cal G}_1 & \mapsto & {\cal G}_2 & \mapsto  & {\cal G}_3\\
\Big\downarrow & & \Big\downarrow & & \Big\downarrow \\
{\cal G}_1^{\ast\ast} & \mapsto & {\cal G}_2^{\ast\ast} & \mapsto  & 
{\cal G}_3^{\ast\ast}
\end{array}
\]

 The second row is acyclic and it induces an acyclic sequence on
 the torsion free parts $({\cal G}_i)_{torfree}$, $1\leq i\leq 3$,

$$\hskip -.5in ({\cal G}_1)_{torfree} \mapsto ({\cal G}_2)_{torfree}\mapsto 
({\cal G}_3)_{torfree} \mapsto 0.$$

 The above sequence is right-exact because by the commutative diagram 
$({\cal G}_2)_{tor}$ is in the kernel of
 the composite surjection 
${\cal G}_2\mapsto {\cal G}_3\mapsto ({\cal G}_3)_{torfree}$.

 To show that it is also left exact and 
exact in the middle, notice that the acyclicity of the above
 sequence implies the surjection 
$({\cal G}_2)_{torfree}/Im(({\cal G}_1)_{torfree})\mapsto 
({\cal G}_3)_{torfree}\mapsto 0$. 

 Consider the sequence 
$$({\cal G}_1)_{torfree}\mapsto ({\cal G}_2)_{torfree}
\mapsto ({\cal G}_2)_{torfree}/Im(({\cal G}_1)_{torfree})\mapsto 0.$$

 By the assumption of our lemma,
 $({\cal G}_1)_{torfree}\otimes k(y)\mapsto ({\cal G}_2)_{torfree}\otimes k(y)$ 
is injective for all closed points $y$.
Because 
 both $({\cal G}_1)_{torfree}$ and 
$({\cal G}_2)_{torfree}/Im(({\cal G}_1)_{torfree})$ are coherent, 
 by exercise II.5.8(a) of [Ha], both $rank_{k(y)}({\cal G}_1)_{torfree}
\otimes k(y)$ and 
 $rank_{k(y)}\bigl(({\cal G}_2)_{torfree}/Im(({\cal G}_1)_{torfree})
\bigr)\otimes 
k(y)=rank_{k(y)}\bigl(
({\cal G}_2)_{torfree}\otimes k(y)/Im(({\cal G}_1)_{torfree}\otimes
 k(y))\bigr)$
 are upper semi-continuous.
 But by assumption $({\cal G}_2)_{torfree}$ is known to be locally free, so by
 exercise II.5.8(b) of [Ha] and the connectivity of $Y$, 
 $rank_{k(y)}({\cal G}_2)_{torfree}\otimes k(y)$ is constant throughout the 
connected scheme $Y$. 
This forces
 $rank_{k(y)}({\cal G}_1)_{torfree}\otimes k(y)$ and 
 $rank_{k(y)}({\cal G}_2)_{torfree}/Im(({\cal G}_1)_{torfree})\otimes k(y)=
rank_{k(y)}({\cal G}_2)_{free}\otimes k(y)-rank_{k(y)}({\cal G}_1)_{torfree}\otimes k(y)$ to be lower semi-continuous and to be constant.
 Therefore by exercises 3.17, 5.7-5.8 of chapter II of [Ha], 
$({\cal G}_1)_{torfree}$ and the quotient 
$({\cal G}_2)_{torfree}/Im(({\cal G}_1)_{torfree})$ are also locally free.
 In particular, $({\cal G}_1)_{torfree}\mapsto ({\cal G}_2)_{torfree}$ induces
 a bundle injection 
\footnote{Because of the injectivity of the $\otimes k(y)$ version of morphisms.} 
and it has to be a sheaf injection. 

On the other hand, the surjection 

$$({\cal G}_2)_{torfree}/Im(({\cal G}_1)_{torfree})\mapsto
 ({\cal G}_3)_{torfree}\mapsto 0$$

 implies that $({\cal G}_3)_{torfree}$ is the quotient of a locally
 free sheaf of the same generic rank.
 This implies that the kernel sheaf of this
 surjection must be a torsion sheaf. As there is
 no sheaf injection from a torsion sheaf into a locally free sheaf, 
 the surjection is in fact a sheaf isomorphism and thus 
$({\cal G}_3)_{torfree}$ is also locally free. 

  So we may replace the sheaves in the original short exact sequence 
$$0\mapsto ({\cal G}_1)_{torfree}\mapsto ({\cal G}_2)_{torfree}\mapsto 
({\cal G}_2)_{torfree}/Im(({\cal G}_1)_{torfree})\mapsto 0$$
 by $({\cal G}_1)_{free}$, $({\cal G}_2)_{free}$, and $({\cal G}_3)_{free}$, 
respectively, and the proof of this lemma is complete. $\Box$

\bigskip

 In the following proposition we prove the locally freeness of the
 sheaf ${\cal V}_{quot}=({\cal R}^1
\pi_{\ast}\bigl({\cal O}_{\sum_{1\leq i\leq p}\Xi_{k_i}}\otimes
 {\cal E}_{C-{\bf M}(E)E}\bigr))_{torfree}$ and
 identify the equivalence class of ${\cal V}_{quot}$ in the $K$ group explicitly.

\begin{prop}\label{prop; id}
The torsion free quotient of the coherent sheaf 
${\cal R}^1\pi_{\ast}\bigl({\cal O}_{\sum_{1\leq i\leq p}\Xi_{k_i}}\otimes
 {\cal E}_{C-{\bf M}(E)E}\bigr)$ is locally free and
 it is equivalent to the direct sum of locally free sheaf 
 $\oplus_{1\leq l\leq p}{\cal R}^1
\tilde{\pi}_{\ast}\bigl({\cal O}_{\tilde{\Xi}_{k_l}}
\otimes {\cal E}_{C-\sum_{b\in I_{k_l}}m_bE_b-\sum_{p\geq d>l}
e_{k_d}}\bigr)\otimes  \tilde{\pi}^{\ast}{\cal O}(-\sum_{1\leq a<k_l}
m_aE_{a; k_l})$ in $K_0(Y(\Gamma)\times T(M))$.
\end{prop}

\noindent Proof of the proposition: For $p=1$, the sum of the divisors
$\sum_{i\leq p}\Xi_{k_i}$ collapses to a single $\Xi_{k_1}$. 
The conclusion of locally freeness and the identity in $K_0(Y(\Gamma)\times T(M))$
 are direct consequences of lemma \ref{lemm; freepart}. We prove the general
 case based on induction upon $p$. 

 By induction hypothesis, we know that the ``locally free'' quotient of  
${\cal R}^1\pi_{\ast}\bigl({\cal O}_{\sum_{1\leq i\leq p-1}\Xi_{k_i}}
\otimes {\cal E}_{C-{\bf M}(E)E-e_{k_p}}\bigr)$ is equivalent to 
$\oplus_{1\leq l\leq p-1}{\cal R}^1
\tilde{\pi}_{\ast}\bigl({\cal O}_{\tilde{\Xi}_{k_l}}
\otimes {\cal E}_{C-\sum_{b\in I_{k_l}}m_bE_b-\sum_{p\geq d>l}e_{k_d}}\bigr)\otimes 
 \tilde{\pi}^{\ast}{\cal O}(-\sum_{1\leq a<k_l}m_aE_{a; k_l})$. To prove 
 the $p-$th version of our proposition it suffices to prove the 
 locally freeness of the torsion free sheaf and then prove the existence of
 a short exact sequence of
 locally free sheaves,

$$\hskip -.8in
0\mapsto ({\cal R}^1\pi_{\ast}\bigl({\cal O}_{\sum_{i\leq p-1}\Xi_{k_i}}(-\Xi_{k_p})
\otimes {\cal E}_{C-\sum_{1\leq a\leq n}m_aE_a}\bigr))_{free}\mapsto 
({\cal R}^1\pi_{\ast}\bigl({\cal O}_{\sum_{i\leq p}\Xi_{k_i}}
\otimes {\cal E}_{C-\sum_{1\leq a\leq n}m_aE_a}\bigr))_{free}$$
$$\mapsto 
({\cal R}^1\pi_{\ast}\bigl({\cal O}_{\Xi_{k_p}}\otimes 
{\cal E}_{C-\sum_{1\leq a\leq n}m_aE_a}\bigr))_{free} \mapsto 0.$$

By pushing forward 
the short exact sequence

$$\hskip -1.2in 0\mapsto {\cal O}_{\sum_{1\leq i\leq p-1}\Xi_{k_i}}(-\Xi_{k_p})
\otimes {\cal E}_{C-\sum_{1\leq a\leq n}m_aE_a}\mapsto {\cal O}_{\sum_{1\leq i\leq p}\Xi_{k_i}}
\otimes {\cal E}_{C-\sum_{1\leq a\leq n}m_aE_a}\mapsto {\cal O}_{\Xi_{k_p}}\otimes 
{\cal E}_{C-\sum_{1\leq a\leq n}m_aE_a}\mapsto 0,$$

 we get a long exact sequence

$$ {\cal R}^0\pi_{\ast}\bigl({\cal O}_{\Xi_{k_p}}\otimes 
{\cal E}_{C-\sum_{1\leq a\leq n}m_aE_a}\bigr)\mapsto 
{\cal R}^1\pi_{\ast}\bigl({\cal O}_{\sum_{1\leq i\leq p-1}\Xi_{k_i}}
\otimes {\cal E}_{C-\sum_{1\leq a\leq n}m_aE_a-e_{k_p}}\bigr)$$
$$\mapsto
{\cal R}^1\pi_{\ast}\bigl({\cal O}_{\sum_{i\leq p}\Xi_{k_i}}
\otimes {\cal E}_{C-\sum_{1\leq a\leq n}m_aE_a}\bigr)\mapsto 
{\cal R}^1\pi_{\ast}\bigl({\cal O}_{\Xi_{k_p}}\otimes 
{\cal E}_{C-\sum_{1\leq a\leq n}m_aE_a}\bigr)\mapsto 0.$$

 The above sequence is right exact because $\sum_{1\leq i\leq p-1}\Xi_{k_i}\mapsto Y(\Gamma)$ is of 
relative dimension one over the base $Y(\Gamma)$ and so the corresponding
 second derived image sheaf along $\sum_{1\leq i\leq p-1}\Xi_{k_i}\mapsto Y(\Gamma)$ vanishes.

 Because of the degree constraint on $\sum_{1\leq a\leq n}m_aE_a$, 
${\cal R}^0\pi_{\ast}\bigl({\cal O}_{\Xi_{k_p}}\otimes 
{\cal E}_{C-\sum_{1\leq a\leq n}m_aE_a}\bigr)$ vanishes on the Zariski-open 
subset $Y_{\Gamma}\times T(M)\subset Y(\Gamma)\times T(M)$
 and is a torsion sheaf over $Y(\Gamma)\times T(M)$. By lemma \ref{lemm; freepart} and by the
 induction hypothesis the
 torsion free sheaves $({\cal R}^1\pi_{\ast}\bigl({\cal O}_{\Xi_{k_p}}\otimes 
{\cal E}_{C-\sum_{1\leq a\leq n}m_aE_a}\bigr))_{torfree}$
 and $({\cal R}^1\pi_{\ast}\bigl({\cal O}_{\sum_{1\leq i\leq p-1}\Xi_{k_i}}
\otimes {\cal E}_{C-\sum_{1\leq a\leq n}m_aE_a-e_{k_p}}\bigr))_{torfree}$ are known to be
 locally free.
 
 Then the desired short exact sequence of locally free sheaves is constructed from
 the acyclic sequence formed by
 the torsion free quotients
 of the above long exact sequence, after we have shown the
 torsion free quotient of the middle factor ${\cal V}_{quot}$ is locally free.

 The rest of the proposition is devoted to derive the locally freeness 
 of $({\cal R}^1\pi_{\ast}\bigl({\cal O}_{\sum_{1\leq i\leq p}\Xi_{k_i}}\otimes
 {\cal E}_{C-\sum_{1\leq a\leq n}m_aE_a}\bigr))_{torfree}$ and the exactness
 of the above acyclic sequence.

\bigskip 

\noindent {\bf Step I}: The locally freeness of the torsion free quotient (part).

 \medskip

 The invertible sheaf ${\cal O}_{M_{n+1}}(-\sum_{1\leq a\leq n}m_aE_a)$
 pulls back to an invertible sheaf on $\sum_{1\leq i\leq p}\Xi_{k_i}$, 
denoted by ${\cal O}_{\sum_{1\leq i\leq p}\Xi_{k_i}}(-\sum_{1\leq 
a\leq n}m_aE_a)$. The invertible sheaf fails to be a sub-sheaf of 
 ${\cal O}_{\sum_{1\leq i\leq p}\Xi_{k_i}}$, therefore it is not an ideal
 sheaf on $\sum_{1\leq i\leq p}\Xi_{k_i}$. We point out the main cause.
Let $P$ be the union of $\{k_1, k_2, k_3, \cdots, k_p\}$ union with
 their direct and indirect ascendents in $\Gamma$. 
The canonically defined sheaf morphism 
${\cal O}_{M_{n+1}}(-\sum_{a\in P}m_aE_a)
\mapsto {\cal O}_{M_{n+1}}$ 
 vanishes on the whole sub-scheme $\sum_{1\leq i\leq p}\Xi_{k_i}\subset M_{n+1}$
 because the sections defining the
 divisors $E_a$, $a\in P$
  vanish \footnote{This is why we want $P$ to contain $k_i$s or their ascendents.
 when $i$ is a ascendent of $j$ in $\Gamma$, 
the defining section of $E_i$ is divisible by the defining section of $E_j$
 above $Y(\Gamma)$.}
on $\Xi_{k_i}$ for all $1\leq i\leq p$. Thus,
 the defining section of 
${\cal O}_{M_{n+1}}(-\sum_{1\leq a\leq n}m_aE_a)
\mapsto {\cal O}_{M_{n+1}}$ vanishes on $\sum_{1\leq i\leq p}\Xi_{k_i}$ as well.

Consider the fiber product (intersection) of 
$\sum_{1\leq a\leq n}m_aE_a\subset M_{n+1}$ with
 $\sum_{1\leq i\leq p}\Xi_{k_i}\subset M_{n+1}$. Even though the fiber product
 is not a sub-scheme of $\sum_{1\leq i\leq p}\Xi_{k_i}$, it still
 contains a maximal sub-scheme $Z$ as a divisor in
$\sum_{1\leq i\leq p}\Xi_{k_i}$.

Since $Z$ is a divisor,
 the ideal sheaf defining $Z$, ${\cal I}_Z$ is locally free $\cong 
{\cal O}_{\sum_{1\leq i\leq p}\Xi_{k_i}}(-Z)$. 

 Then we may write ${\cal O}_{\sum_{1\leq i\leq p}\Xi_{k_i}}(-\sum_{1\leq a\leq n}
m_aE_a)$ as ${\cal I}_Z\otimes {\cal J}$, with ${\cal J}$ being locally free. 

Tensoring the defining short exact sequence 
$0\mapsto {\cal I}_Z\mapsto {\cal O}_{\sum_{i\leq p}\Xi_{k_i}}
\mapsto {\cal O}_Z\mapsto 0$ by 
${\cal J}_C={\cal J}\otimes {\cal E}_C$
and taking the derived long exact sequence
 along $\pi:\sum_{i\leq p}\Xi_{k_i}\mapsto Y(\Gamma)$, we find 

$$ 
{\cal R}^0\pi_{\ast}\bigl({\cal O}_{\sum_{i\leq p}\Xi_{k_i}}\otimes 
{\cal J}_C\bigr)
\mapsto {\cal R}^0\pi_{\ast}\bigl({\cal O}_Z\otimes 
{\cal J}_C\bigr)$$
$$\stackrel{\delta}{\mapsto} {\cal R}^1\pi_{\ast}\bigl(
{\cal O}_{\sum_{i\leq p}\Xi_{k_i}}
\otimes {\cal E}_{C-\sum_{1\leq a\leq n}m_aE_a}\bigr)\mapsto 
{\cal R}^1\pi_{\ast}\bigl({\cal O}_{\sum_{i\leq p}\Xi_{k_i}}\otimes 
{\cal J}_C\bigr).$$

We know that ${\cal R}^1\pi_{\ast}\bigl({\cal O}_{\sum_{i\leq p}\Xi_{k_i}}\otimes 
{\cal J}_C\bigr)$ is a torsion sheaf, since $deg_{\Xi_{k_i}/Y(\Gamma)}{\cal J}_C
=m_i>0$ by a calculation shown on page \pageref{degree}.

Now we get the following short exact sequence on 
the torsion free quotients of the
 above sequence (based on lemma \ref{lemm; tor}),

$$0\mapsto ({\cal R}^0\pi_{\ast}\bigl({\cal O}_{\sum_{i\leq p}\Xi_{k_i}}\otimes 
{\cal J}_C\bigr))_{free}
\mapsto ({\cal R}^0\pi_{\ast}\bigl({\cal O}_Z\otimes 
{\cal J}_C\bigr))_{free}$$
$$\mapsto 
({\cal R}^1\pi_{\ast}\bigl({\cal O}_{\sum_{i\leq p}\Xi_{k_i}}
\otimes {\cal E}_{C-\sum_{1\leq a\leq n}m_aE_a}\bigr))_{free}\mapsto 0,$$

 once one shows (i). 
$({\cal R}^0\pi_{\ast}\bigl({\cal O}_Z\otimes 
{\cal J}_C\bigr))_{torfree}$ is locally free.

(ii). the injectivity  $({\cal R}^0\pi_{\ast}\bigl({\cal O}_{\sum_{i\leq p}\Xi_{k_i}}\otimes 
{\cal J}_C\bigr))_{torfree}\otimes k(y)\mapsto 
({\cal R}^0\pi_{\ast}\bigl({\cal O}_Z\otimes 
{\cal J}_C\bigr))_{torfree}\otimes k(y)$ for all the closed
 points $y\in Y(\Gamma)\times T(M)$.

\bigskip

\noindent The proof of condition (i):
Take \footnote{The subscripts $f$ and $t$ of $Z_f$ and $Z_t$ correspond to
 the keywords ``free'' and ``torsion'' as they are closely related to the 
locally free part and
 the torsion part of ${\cal R}^0\pi_{\ast}\bigl({\cal O}_Z\otimes 
{\cal J}_C\bigr)$.}
 $Z_f$ to be the union of the irreducible components of the
 divisor $Z\subset \sum_{1\leq i\leq p}\Xi_{k_i}$ which dominates $Y(\Gamma)$. 
Take $Z_t$ be the union of the irreducible components of divisors in 
$Z$ which do not dominate $Y(\Gamma)$. 
Because $Z$ is the divisor induced by $\sum_{1\leq a\leq n}m_aE_a$, 
for a fixed $i$ with $1\leq i\leq p$, 
the restriction of 
$\sum_{j_{k_i}}m_{j_{k_i}}E_{j_{k_i}}$ to $\Xi_{k_i}$ defines a sub-scheme in
 $\Xi_{k_i}$,
 a union (with multiplicities) of cross-sections of $\Xi_{k_i}\mapsto Y(\Gamma)$.
So the 
 map $Z_f\mapsto Y(\Gamma)$ is a finite morphism. On the other hand, $Z_t\subset 
\sum_{1\leq i\leq p}\Xi_{k_i}$ maps onto a union of divisors in $Y(\Gamma)$.
 Because $Z=Z_f\cup Z_f$ is a union of divisors, 
we have the following short exact sequence of divisors
 in $\sum_{i\leq p}\Xi_{k_i}$, 
 
 $$0\mapsto {\cal O}_{Z_f}(-Z_t)\mapsto {\cal O}_Z\mapsto {\cal O}_{Z_t}\mapsto 0.$$

Since $Z_f\mapsto Y(\Gamma)$ is a finite morphism, 
${\cal R}^1\pi_{\ast}\bigl({\cal O}_{Z_f}(-Z_t)\otimes 
{\cal J}_C\bigr)$ vanishes and 
${\cal R}^0\pi_{\ast}\bigl({\cal O}_{Z_f}(-Z_t)\otimes {\cal J}_C
\bigr)$ is locally free\footnote{Its rank is nothing but the relative 
length of $Z_f\mapsto Y(\Gamma)$.}. 
So the ${\cal J}_C$ twisted 
derived long exact sequence of the above short exact sequence truncates
 to a sheaf short exact sequence 

$$\hskip -.4in 
0\mapsto {\cal R}^0\pi_{\ast}\bigl({\cal O}_{Z_f}(-Z_t)\otimes 
{\cal J}_C\bigr)\mapsto
{\cal R}^0\pi_{\ast}\bigl({\cal O}_Z\otimes {\cal J}_C\bigr)
\mapsto {\cal R}^0\pi_{\ast}\bigl({\cal O}_{Z_t}\otimes 
{\cal J}_C\bigr)\mapsto 0.$$

 Because $Z_t$ is mapped into a proper sub-scheme of $Y(\Gamma)$ under $\pi$, 
its intersection with the generic fibers of 
$\pi:\sum_{i\leq p}\Xi_{k_p}\mapsto Y(\Gamma)$
 must be empty. Thus ${\cal R}^0\pi_{\ast}\bigl({\cal O}_{Z_t}\otimes 
{\cal J}_C\bigr)$ is a torsion sheaf.

 One the other hand, it is easy to check \footnote{By comparing with
the corresponding short exact sequence of the fiber above $y$.
 Consult the next argument below
 and the next footnote for a similar argument.} that
 ${\cal R}^0\pi_{\ast}\bigl({\cal O}_{Z_f}(-Z_t)\otimes 
{\cal J}_C\bigr)\otimes k(y)\mapsto 
{\cal R}^0\pi_{\ast}\bigl({\cal O}_Z\otimes 
{\cal J}_C\bigr)\otimes k(y)$ is 
injective for all the closed points $y$ on $Y(\Gamma)\times T(M)$. 
 lemma \ref{lemm; split} the torsion free quotient of 
${\cal R}^0\pi_{\ast}\bigl({\cal O}_Z\otimes 
{\cal J}_C \bigr)$ 
is isomorphic to ${\cal R}^0\pi_{\ast}\bigl({\cal O}_{Z_f}(-Z_t)\otimes 
{\cal J}_C \bigr)$ and is known to be locally free.

 So we know that 
$({\cal R}^0\pi_{\ast}\bigl({\cal O}_Z\otimes 
{\cal J}_C\bigr))_{torfree}$ is locally free and the 
condition (i). has been proved.

\medskip

\noindent The proof of condition (ii): By the derivation of condition (i), the 
original condition (ii) is equivalent to  
the injectivity  $({\cal R}^0\pi_{\ast}\bigl({\cal O}_{\sum_{i\leq p}\Xi_{k_i}}
\otimes {\cal J}_C
\bigr))_{torfree}\otimes k(y)\mapsto 
{\cal R}^0\pi_{\ast}\bigl({\cal O}_{Z_f}(-Z_t)\otimes 
{\cal J}_C \bigr)\otimes k(y)$ for all the closed
 points $y\in Y(\Gamma)\times T(M)$. On the other hand, the $k(y)$-twisted 
zero-th derived image sheaves map into the zero-th fiberwise sheaf cohomologies 
injectively\footnote{This follows from taking the global sections of a twisted 
version of
 the exact sequence $0\mapsto {\cal I}_{\sum_{i\leq p}\Xi_{k_i}|_y}\mapsto 
{\cal O}_{\sum_{i\leq p}\Xi_{k_i}}\mapsto {\cal O}_{\sum_{i\leq p}\Xi_{k_i}|_y}
\mapsto 0$
 for all $y\in Y(\Gamma)$.}, so
 it suffices to check the following injection
$H^0(\sum_{i\leq p}\Xi_{k_i}\times_{Y(\Gamma)\times T(M)}\{y\}, 
{\cal J}_C|_y)\mapsto 
H^0(Z_f\times_{Y(\Gamma)\times T(M)}\{y\}, {\cal O}_{Z_f}(-Z_t)\otimes 
{\cal J}_C|_y)$ for all the closed points $y$ 
in $Y(\Gamma)\times T(M)$.
 By composing with the natural morphism induced by 
${\cal O}_{Z_f}(-Z_t)\mapsto {\cal O}_{Z_f}$, it suffices to check the 
injectivity of 
$H^0(\sum_{i\leq p}\Xi_{k_i}\times_{Y(\Gamma)\times T(M)}\{y\}, 
{\cal J}_C|_y)\mapsto 
H^0(Z_f\times_{Y(\Gamma)\times T(M)}\{y\}, {\cal O}_{Z_f}\otimes 
{\cal J}_C|_y)$. But this map is nothing but the 
restriction map
 of fiberwise global sections above $y$
to the finite sub-scheme $Z_f\times_{Y(\Gamma)\times T(M)}\{y\}$. 
If this map is not injective, there
 must be a non-zero global section of 
${\cal J}_C|_y$ on
 the fiber $\sum_{i\leq p}\Xi_{k_i}\times_{Y(\Gamma)\times T(M)}\{y\}$ which 
vanishes along
 the finite sub-scheme $Z_f\times_{Y(\Gamma)\times T(M)}\{y\}$. In particular,
 we may
 restrict this fiberwise global 
section to one particular $\Xi_{k_l}\times_{Y(\Gamma)\times T(M)}\{y\}$ (for some
 $1\leq l\leq p$)
 over which the section does not vanish identically. We derive the contradiction by
 computing the degree of the invertible sheaf in two different ways.
 
 Consider the index subset $P\subset \{1, 2, \cdots, n\}$ collecting $k_l$
 and all its direct or indirect ascendents in $\Gamma$.

On the one hand, the degree of the
 invertible sheaf \footnote{The degree can be calculated along any smooth
 fiber of $\Xi_{k_l}$ over $Y_{\Gamma}$. moreover it is easy to see that
the inverse image of the invertible sheaf ${\cal J}$ to the sub-locus 
$\Xi_{k_l}\times_{M_n} Y_{\Gamma}$
 is isomorphic to ${\cal O}_{\Xi_{k_l}\times_{M_n}
Y_{\Gamma}}(-\sum_{a\in P}m_aE_a)$ over $Y_{\Gamma}$. it is because
 the defining section of $E_a$, $a\in P$ vanishes on $\Xi_{k_i}|_y$ over $y$.}, 
$deg_{\Xi_{k_l}\times_{Y(\Gamma)\times T(M)}\{y\}}
{\cal J}_C|_y$,  \label{degree}
 is  $deg_{\Xi_{k_l}\times_{Y(\Gamma)\times T(M)}\{y\}}
{\cal O}_{\sum_{1\leq i\leq p}\Xi_{k_i}}(-\sum_{a\in P}m_aE_a)$
$=e_{k_l}\cdot (-\sum_{a\in P}m_aE_a)=
(E_{k_l}-\sum_{j_{k_l}}E_{j_{k_l}})\cdot
 (-\sum_{a\in P}m_aE_a)= -m_{k_l}E_{k_l}^2=m_{k_l}$ (since $k_l\in P$ and
 the descendent $j_{k_l}$ of $k_l$ can never be an ascendent of $k_l$,  
 any term of the form $(-E_{j_{k_l}})\cdot (-m_aE_a)$
contributes trivially
 to the sum). On the 
other hand, the same degree must be no less than the relative length, 
$length((Z_f\cap \Xi_{k_l})
\times_{Y(\Gamma)\times T(M)}\{y\})$, of the finite scheme $(Z_f\cap \Xi_{k_l})
\times_{Y(\Gamma)\times T(M)}\{y\}$ along which the section vanishes.
 However the length of this finite scheme is nothing but the sum of 
 multiplicities
 for all $m_{j_{k_l}}E_{j_{k_l}}$ along\footnote{The index 
$j_{k_l}$ is a typical
 direct descendent index of $k_l$ in the admissible graph $\Gamma$.}
 $\Xi_{k_l}$ and is equal to $\sum_{j_{k_l}}m_{j_{k_l}}$. 
 
 Combining these observations together we get an inequality 
$m_{k_l}\geq \sum_{j_{k_l}}m_{j_{k_l}}$, which
 implies $e_{k_l}\cdot (C-{\bf M}(E)E)\geq 0$. This violates our 
choices of $e_{k_l}$ of 
 making $e_{k_l}\cdot (C-{\bf M}(E)E)<0$!
 So the original $k(y)$-vector space morphisms must be injective for all $y$.
 The proof of condition (ii). is 
 finished.

\bigskip

\noindent {\bf Step II}: The proof of exactness of the acyclic sequence.

\medskip

  We plan to derive it by adopting
 the commutative diagram chasing argument.

  Recall that the inclusion 
$\sum_{i\leq p-1}\Xi_{k_i}\subset \sum_{i\leq p}\Xi_{k_i}$ of ${\bf P}^1$ 
fibrations
 (removing the last $\Xi_{k_p}$) induces an inverse image of the ideal sheaf 
${\cal I}_Z$ 
 in $\sum_{i\leq p-1}\Xi_{k_i}$, denoted by ${\cal I}_{Z'}$. Then the 
restriction of $Z$,
 $Z'$, can be viewed
 as a sub-scheme of $Z$ by removing the intersection $Z''=\Xi_{k_p}\cap Z$. 
Similarly,
 both $Z_f$ and $Z_t$ are restricted to $Z'_f$ and $Z'_t$, respectively and 
likewise we also have
 $Z''=Z''_f+Z''_t=\Xi_{k_p}\cap Z_f+\Xi_{k_p}\cap Z_t$.  To summarize, we have the
 following commutative diagram of three rows of short exact sequences,

\[
\hskip -.2in
\begin{array}{ccccc}
 {\cal O}_{Z'_f}(-Z'_t-Z'') & \mapsto & {\cal O}_{Z'}(-Z'') & 
\mapsto & {\cal O}_{Z'_t}(-Z'')\\
 \Big\downarrow &  & \Big\downarrow &  & \Big\downarrow\\
{\cal O}_{Z_f}(-Z_t)  & \mapsto & {\cal O}_Z & \mapsto & {\cal O}_{Z_t} \\
\Big\downarrow &  & \Big\downarrow &  & \Big\downarrow\\
 {\cal O}_{Z''_f}(-Z''_t) & \mapsto & {\cal O}_{Z''} & \mapsto & {\cal O}_{Z''_t}
\end{array}
\] 

All of the three rows and the first two 
columns in this above diagram are short exact sequences constructed from tensoring
 invertible sheaves with a divisorial exact sequence of the following
 form, $0\mapsto {\cal O}_B(-A)\mapsto {\cal O}_{A+B}\mapsto {\cal O}_B\mapsto 0$.

 Moreover we have the following commutative diagram \footnote{We have
skipped the inverse image notation for ${\cal J}_C$ to the various
 sub-schemes, in order to make the formula less complicated.}
 (after substituting by the 
short hand
 notations $\underline{C}=C-\sum_{1\leq a\leq n}m_aE_a=C-{\bf M}(E)E$). 

\[
\hskip -1.5in
\begin{array}{ccccc}
 ({\cal R}^0\pi_{\ast}\bigl({\cal O}_{\sum_{i\leq p-1}\Xi_{k_i}}(-\Xi_{k_p})\otimes
 {\cal J}_C\bigr))_{free} & \mapsto &
 {\cal R}^0\pi_{\ast}\bigl({\cal O}_{Z'_f}(-Z'_t-Z'')\otimes
 {\cal J}_C\bigr) & \stackrel{\delta}{\mapsto} & 
({\cal R}^1\pi_{\ast}\bigl({\cal O}_{\sum_{i\leq p-1}\Xi_{k_i}}
\otimes {\cal E}_{\underline{C}-e_{k_p}}\bigr))_{free}\\
\Big\downarrow & & \Big\downarrow & &  \Big\downarrow \\
({\cal R}^0\pi_{\ast}\bigl({\cal O}_{\sum_{i\leq p}\Xi_{k_i}}\otimes 
{\cal J}_C\bigr))_{free} & \mapsto & 
{\cal R}^0\pi_{\ast}\bigl({\cal O}_{Z_f}(-Z_t)\otimes 
{\cal J}_C \bigr)  & \stackrel{\delta}{\mapsto}& 
 ({\cal R}^1\pi_{\ast}\bigl({\cal O}_{\sum_{i\leq p}\Xi_{k_i}}\otimes 
{\cal E}_{\underline{C}}\bigr))_{free}\\
\Big\downarrow & & \Big\downarrow & &  \Big\downarrow\\
({\cal R}^0\pi_{\ast}\bigl({\cal O}_{\Xi_{k_p}}\otimes
{\cal J}_C\bigr))_{free}  & \mapsto  & 
{\cal R}^0\pi_{\ast}\bigl({\cal O}_{Z''_f}(-Z''_t)\otimes 
{\cal J}_C\bigr) & \stackrel{\delta}{\mapsto} & 
({\cal R}^1\pi_{\ast}\bigl({\cal O}_{\Xi_{k_p}}\otimes 
{\cal E}_{\underline{C}}\bigr))_{free}
\end{array}
\]

 $\diamondsuit$ 
We claim that all three rows and the first two columns in this commutative 
diagram 
are short exact sequences:
 By the earlier discussion based on lemma \ref{lemm; tor}, the sheaves 
in the second rows are all locally free (this justifies the
 usage of the subscript $(\cdot)_{free}$ above). In the same argument
the second row has been shown to be short exact. We may adopt a
parallel argument
 on $\sum_{i\leq p-1}\Xi_{k_i}$ applied to $Z'=Z'_f\cup Z'_t$ or on
 $\Xi_{k_p}$ applied to $Z''=Z''_f\cup Z''_t$, thus 
the sheaves in the first and the 
third rows
 are all locally free and both of the first and the third rows are short exact 
as well.
  The first column is the locally (torsion) free summand of 
 a derived long exact sequence, it is exact based on lemma \ref{lemm; tor} and 
we argue
 as the following:  The locally freeness of the factor in the 
 middle is already known. The injectivity of the $\bullet\otimes k(y)$ version
 of the first column in the above diagram

$$\hskip -.8in
({\cal R}^0\pi_{\ast}\bigl({\cal O}_{\sum_{i\leq p-1}\Xi_{k_i}}(-\Xi_{k_p})\otimes
 {\cal J}_C\bigr))_{free}\otimes k(y)\mapsto 
({\cal R}^0\pi_{\ast}\bigl({\cal O}_{\sum_{i\leq p}\Xi_{k_i}}\otimes 
{\cal J}_C\bigr))_{free}\otimes k(y)$$

 is a direct consequence of the injectivity of the fiberwise zero-th cohomologies,

$$\hskip -1.3in
H^0(\sum_{i\leq p-1}\Xi_{k_i}\times_{Y(\Gamma)\times T(M)}\{y\}, 
{\cal O}(-\Xi_{k_p}\times_{Y(\Gamma)\times T(M)}\{y\})\otimes 
{\cal J}_C|_y)
\mapsto 
H^0(\sum_{i\leq p}\Xi_{k_i}\times_{Y(\Gamma)\times T(M)}\{y\}, 
{\cal J}_C|_y).$$

 The desired injectivity of the $H^0$ morphism has been 
 the direct consequence of the following short exact sequence on the fiber 
above $y$,
$$\hskip -.9in
0\mapsto {\cal O}_{\sum_{i\leq p-1}\Xi_{k_i}
\times_{Y(\Gamma)\times T(M)}\{y\}}(-\Xi_{k_p}\times_{Y(\Gamma)\times 
T(M)}\{y\})\mapsto 
{\cal O}_{\sum_{i\leq p}\Xi_{k_i}
\times_{Y(\Gamma)\times T(M)}\{y\}}\mapsto {\cal O}_{\Xi_{k_p}
\times_{Y(\Gamma)\times T(M)}\{y\}}\mapsto 0.$$

 So the exactness of the first column is ensured.

Finally the second column is short exact as it is the derived short exact 
sequence (remembering that $Z'_f\mapsto Y(\Gamma)$ is
 a finite morphism) of a twisted short exact sequence in the the first column 
of the 
previous commutative diagram
 on $Z$, $Z_f$, $Z_t$, $Z''$, $Z''_f$, $Z''_t$, $Z'$, $Z'_f$, $Z'_t$.

  The third column has been known to be acyclic.
  Then its exactness follows from the standard diagram-chasing 
technique.
$\Box$

\medskip

 At the end of the subsection, we define a short-hand notation, 

\begin{defin}\label{defin; replace} 
Define the locally free sheaf $\tilde{\cal V}_{quot}$ to be 

$$\tilde{\cal V}_{quot}=\oplus_{1\leq l\leq p}{\cal R}^1
\tilde{\pi}_{\ast}\bigl({\cal O}_{\tilde{\Xi}_{k_l}}
\otimes {\cal E}_{C-\sum_{b\in I_{k_l}}m_bE_b-\sum_{p\geq d>l}e_{k_d}}\bigr)
\otimes 
 \tilde{\pi}^{\ast}{\cal O}(-\sum_{1\leq a<k_l}m_aE_{a; k_l}).$$
\end{defin}

 Following our convention, the corresponding vector bundle 
will be denoted by $\tilde{\bf V}_{quot}$.

\medskip

\section{\bf The Localized Chern Classes and the Discrepancy of the 
Top Chern Classes} \label{section; discrepancy}

\bigskip

 Let us consider the following general set up. Let $X$ be
\footnote{The notation $X$ used 
in this section is not the same $X$ used in the 
the proof of our main theorem.} a purely $m$ 
dimensional scheme and let ${\bf E}\mapsto X$, $F\mapsto X$ be vector
 bundles over $X$
of the same rank, say $e$, and let $\sigma: {\bf E}\mapsto {\bf F}$ be a bundle 
homomorphism on $X$ exact off a closed subset $Z$. Then in principle 
the difference
of Chern classes of ${\bf F}$ and ${\bf E}$ should be expressible as cycle classes 
'localized' in $Z$. In particular, when $Z=\emptyset$, the map
 $\sigma$ induces an isomorphism between ${\bf F}$ and ${\bf E}$ and their Chern 
classes
 coincide.

 For the convenience of the reader, we review the construction of
 [F], page 348 (c). as a special case of the graph construction of MacPherson.
 We have changed a few notations from the original notations of [F].

\medskip

\begin{prop}\label{prop; graph}
Let $cl({\bf E}), cl({\bf F})$ denote a polynomial of Chern classes of ${\bf E}$
 and ${\bf F}$, 
respectively. Let $G=Grass_e({\bf E}\oplus {\bf F})$
 be the $e$-plane Grassmanian bundle 
over $X$ and let $\zeta\mapsto G$ be the universal rank $e$ bundle over $G$. 
Then there exists a cycle $\sum_{i\geq 1}n_iV_i\subset G$ supporting
 over $Z$, $\eta_i:V_i\mapsto Z$ the projection map, such that

$$cl({\bf E})\cap [X]-cl({\bf F})\cap [X]=\sum_{V_i\mapsto Z}n_i\eta_{i\ast}
\{cl(\zeta)\cap
[V_i]\}.$$
\end{prop}

\medskip

 Define $\phi: X\times {\bf A}\mapsto G\times {\bf P}^1$ by
 sending $(x, t)$ to the graph of $t\sigma(x)\times (1, t)$. Then define $W$ to be
 the closure of the image in $G\times {\bf P}^1$. Set 
$W_{\infty}=i_{\infty}^{\ast}[W]=\sum_{i\geq 0}n_iV_i$ to be the
 fiber of $W\mapsto {\bf P}^1$ above $\{\infty\}\subset {\bf P}^1$.

 Because $\sigma$ is bijective on $X-Z$, there exists a special 
component, say $V_0$, 
with $n_0=1$, and is birational to $X$ through the projection $G\mapsto X$.
 In fact it is isomorphic to $X$. Excluding this component from $W_{\infty}$,
it turns out that the remaining $\sum_{i>0}n_iV_i$, $\eta_i:V_i\mapsto Z$
 has the desired property. For the full details, consult page 340-348, section
 18.1 the graph construction, of [F].

\medskip

 The graph construction allows us to express the difference of 
characteristic classes of ${\bf E}$ and ${\bf F}$ by cycle classes in $Z$, 
constructed 
 from $\sigma$ through the limiting process. In particular, it 
 implies the following identity

$$\{c_e({\bf E})-c_e({\bf F})\}\cap [X]=
\sum_{V_i\mapsto Z}n_i\eta_{i\ast}\{c_e(\eta)\cap [V_i]\}.$$

 But it might be technical to identify these $[V_i]$ explicitly.

 Suppose that one is given additionally 
a global section $s_0\in\Gamma(X, {\bf E})$, then
 $s=\sigma(s_0)\in \Gamma(X, {\bf F})$ is a global section of ${\bf F}$.

 The localized top Chern class construction on page 244, section 14.1 of [F]
 defines localized classes ${\bf Z}(s_0)\in {\cal A}_{m-e}(Z(s_0))$, 
${\bf Z}(s)\in {\cal A}_{m-e}(Z(s))$ and their 
push-forward into $X$ are equal to $c_e({\bf E}), c_e({\bf F})$, respectively.

Thus the datum of the sections $s_0, s=\sigma(s_0)$ may be used to express the
 difference $c_e({\bf E})-c_e({\bf F})$ as geometric cycles relating to $s_0$ and
 $\sigma$.

 The first hint to such a possibility is the following proposition,

\medskip

\begin{prop}\label{prop; section}
 Let $\sigma:{\bf E}\mapsto {\bf F}$, $X$, $s_0:X\mapsto {\bf E}$ be as above and 
let
 $s_{\bf E}:X\mapsto {\bf E}$ denote the zero section of ${\bf E}$ and let 
$\pi_{\bf E}:{\bf E}\mapsto X$ be the
 bundle projection, then
 the kernel $Ker(\sigma)$ determines an algebraic sub-cone ${\bf C}_{\rho}$ 
of the total space
 of ${\bf E}$ and there exists an scheme theoretical equality 
 $Z(s)=Z(s_0)\cup  \pi_{\bf E}(({\bf C}_{\rho}-s_{\bf E}(X))\cap s_0(X))$ between
 the zero loci. 
\end{prop}

\medskip

 In general we may write ${\bf C}_{\rho}=\cup_{i\geq 0}{\bf C}_{\rho_i}$, where 
 ${\bf C}_{\rho_0}$ is the zero section cone $s_{\bf E}(X)$ and 
$\cup_{i>0}{\bf C}_{\rho_i}$
 is the union of the remaining irreducible components supporting inside $Z$.
 Because the proposition is parallel to the discussion in proposition 12
 in [Liu5], we only give a sketch of the proof:

\medskip

\noindent Sketch of the proof: Let ${\cal E}$, ${\cal F}$ be the 
 locally free sheaves over $X$ associated to ${\bf E}$, ${\bf F}$, respectively. 

 The sheaf morphism ${\cal E}\mapsto {\cal F}$ induces a dual
 morphism ${\cal F}^{\ast}\mapsto {\cal E}^{\ast}$ with cokernel sheaf
 ${\cal R}$.

Consider the
 ${\cal O}_X$ algebra ${\bf S}^{\cdot}$ generated by ${\bf S}^1={\cal R}$,
 then ${\bf C}_{\rho}=Spec({\bf S}^{\cdot})$ defines a sub-cone in the vector 
bundle cone of ${\bf E}$. By tensoring with $k(x)$ (which is right exact) 
for all $x\in X$ and taking 
the left exact contravariant functor $HOM_{k(x)}(\cdot, k(x))$, one may
see easily that this cone is the kernel sub-cone of $\sigma: {\bf E}\mapsto 
{\bf F}$.

  One may observe $Z(s)=Z(s_0)\cup 
\pi_{\bf E}(({\bf C}_{\rho}-s_{\bf E}(X))\cap s_0(X))$ on the set theoretical 
level rather
 easily. The equality as schemes follows from a parallel
 discussion as in proposition 12/corollary 3 of [Liu5]. We omit the details 
 here. $\Box$

 One notices that besides the unique component ${\bf C}_{\rho_0}$ 
which is equal to the zero section $s_{\bf E}(X)$,
 the union of the remaining cones, 
${\bf C}_{\rho}-s_{\bf E}(X)=\cup_{i>0}{\bf C}_{\rho_i}$, supports exactly on $Z$.

\medskip

 On the one hand, we have the following residual intersection theory formula
 on the (localized) top Chern class (see page 245, example 14.1.4. of [F])

\medskip

\begin{prop}\label{prop; rif}
Let ${\bf F}\mapsto X$ be a rank $e$ vector bundle over a purely $m$-dimensional 
scheme
 $X$ and let $s:X\mapsto {\bf F}$ be a global section.
  Let $D$ be an effective Cartier divisor contained in $Z(s)$, then
 there exists a section $s'$ of ${\bf F}\otimes {\cal O}(-D)$ such that

\noindent (i). ${\bf F}\otimes {\cal O}(-D)\mapsto {\bf F}$ maps $s'$ to $s$.

\noindent (ii). ${\bf Z}(s)={\bf Z}(s')+\sum_{1\leq i\leq e}(-1)^{i-1}
c_{e-i}({\bf F})
\cap D^{i-1}\cap [D]$.
\end{prop}

\medskip

  It makes sense to combine proposition \ref{prop; graph}, \ref{prop; section} 
and \ref{prop; rif} and unify these observations together.

 Firstly, let $\sigma: {\bf E}\mapsto {\bf F}$ be isomorphic off $Z\subset X$ 
as before.
 Consider the sub-scheme $Z(s)\cap Z\subset X$. One may blow it up into
 an exceptional Cartier divisor, denoted as $D$ in the blown up scheme $\tilde{X}$.
 From the general construction
 of blowing up coherent sheaves of ideals, ( see page 163-169 of [Ha] and
 B.6 page 435-437 of [F]), $D$ is isomorphic to the
 projectified normal cone ${\bf P}({\bf C}_{Z(s)\cap Z}X)$. 
 Then one may apply proposition \ref{prop; rif} to $\tilde{X}$, 
$g^{\ast}{\bf F}\mapsto \tilde{X}$, where $g:\tilde{X}\mapsto X$ denotes the 
blowing
 down map with the exceptional divisor $D$.

 The following simple lemma identifies the term 
$\sum_{1\leq i\leq e}(-1)^{i-1}c_{e-i}({\bf F})\cap D^{i-1}\cap [D]$ for us,

\medskip

\begin{lemm}\label{lemm; local} 
The cycle class
$g_{\ast}\{\sum_{1\leq j\leq e}(-1)^{j-1}c_{e-j}({\bf F})\cap D^{j-1}\cap [D]\}\in 
{\cal A}_{m-e}(Z(s)\cap Z)$ is equal to the localized
 contribution of top Chern class of 
 $Z(s)\cap Z$,
${\bf Z}_{Z(s)\cap Z}(s)
=\{c({\bf F}|_{Z(s)\cap Z})\cap s(Z(s)\cap Z, X)\cap [X]\}_{m-e}$. (see definition
 1 of in section 5 of [Liu5])
\end{lemm}

\medskip

\noindent Proof of the lemma: Recall from page 71 of [F] 
that the Segre class of a cone\footnote{Do not confuse this cone $C$ with
 the cohomology class $C\in H^{1, 1}(M, {\bf Z})$.}
 $C\mapsto Y$ over $Y$ is defined to be 
$$s(C)=q_{\ast}(\sum_{j\geq 0}c_1({\cal O}(1))^j\cap [{\bf P}(C\oplus 1)])
\in {\cal A}_{\ast}(Y),$$

where $q:{\bf P}(C\oplus 1)\mapsto Y$ is the projection map.

 In our context we take $C$ to be the normal cone ${\bf C}_{Z(s)\cap Z}X$ and 
 $Y=Z(s)\cap Z$, and so $D={\bf P}(C_{Z(s)\cap Z}X)$.
On the one hand, ${\bf P}({\bf C}_{Z(s)\cap Z}X)$ is a divisor in
 ${\bf P}({\bf C}_{Z(s)\cap Z}X\oplus 1)$ defined by $c_1({\cal O}(1))$.
 On the other hand,
$D={\bf P}({\bf C}_{Z(s)\cap Z}X)$ is an exceptional divisor
 in $\tilde{X}$, thus $D=-c_1({\cal O}(1))$ and we have

$$g_{\ast}\{\sum_{1\leq j\leq e}(-1)^{j-1}c_{e-j}(g^{\ast}{\bf F})\cap D^{j-1}
\cap [D]\}
=g_{\ast}\{\sum_{1\leq j\leq e}c_{e-j}(g^{\ast}{\bf F})\cap c_1({\cal O}(1))^{j-1} 
\cap [{\bf P}({\bf C}_{Z(s)\cap Z}X)]\}$$
$$=\{c({\bf F}|_{Z(s)\cap Z})\cap g_{\ast}(\sum_{1\leq j\leq e}
c_1({\cal O}(1))^j\cap 
 [{\bf P}({\bf C}_{Z(s)\cap Z}\oplus 1)])\}_{m-e}$$

$$=\{c({\bf F}|_{Z(s)\cap Z})g_{\ast}(\sum_{1\leq j}c_1({\cal O}(1))^j\cap 
[{\bf P}({\bf C}_{Z(s)\cap Z}\oplus 1)])\}_{m-e}$$
$$=\{c({\bf F}|_{Z(s)\cap Z})\cap s(Z(s)\cap Z, X)\}_{m-e}.$$   
$\Box$

\medskip

 The expression $\{c({\bf F}|_{Z(s)\cap Z})\cap s(Z(s)\cap Z, X)\}_{m-e}$ is
 nothing but the localized contribution of the top Chern class ${\bf Z}_Z(s)$
 discussed in section 6 of [Liu5].

   Thus, the identity in proposition \ref{prop; rif} (ii). can be re-written
 as ${\bf Z}(s)={\bf Z}(s')+{\bf Z}_{Z(s)\cap Z}(s)$. When $\sigma:{\bf E}
\mapsto {\bf F}$ is
 exact off $Z$, we compare it to the top Chern class identity 
$$c_e({\bf F})\cap [X]=c_e({\bf E})\cap [X]-\sum_{V_i\mapsto Z}n_i
\eta_{i\ast}\{cl(\zeta)
\cap [V_i]\}.$$

\medskip

 Introducing the following equivalent relationship which will be
 essential to the invariant enumeration in section \ref{section; proof}.

\begin{defin}\label{defin; numerical}
 Let $\eta_1, \eta_2$ with $\eta_1\cap, \eta_2\cap:
{\cal A}_k(X)\mapsto {\cal A}_0(X)$ be two grade-$k$ characteristic 
classes on an $m$ dimensional complete scheme $X$. The classes $\eta_1$ and
 $\eta_2$ are 
said to be numerically equivalent, denoted as
 $\eta_1\stackrel{n}{=}\eta_2$ if for all $\alpha\in {\cal A}_{m-k}(X)$,
 $\int_X\eta_1\cap \alpha=\int_X\eta_2\cap \alpha$. In other words, 
 $\int_X\eta_1\cap \cdot $ and $\int_X\eta_2\cap \cdot $ define identical 
integral operations
 from ${\cal A}_{m-k}(X)$ to ${\cal A}_0(pt)$.
\end{defin}

  It makes sense to ask the following question,

\medskip

\noindent {\bf Question}: Are $i_{\ast}{\bf Z}(s')\stackrel{n}{=}
c_e(g^{\ast}{\bf E})\cap [\tilde{X}]$ and 
 ${\bf Z}_{Z(s)\cap Z}(s)\stackrel{n}{=}
-\sum_{V_i\mapsto Z}n_i\eta_{i\ast}\{c_e(\zeta)\cap [V_i]\}$?

\medskip

The following proposition answers the question affirmatively.
 In order to apply the current discussion to the explicit enumeration
 problem in section \ref{section; proof},
 we generalize the setting slightly.

\medskip

 Let $X$ be a purely $m$ dimensional reduced complete scheme as above. 
$Z\subset X$ is
 a closed sub-scheme of $X$. 
  
\begin{prop}\label{prop; equivalent}
 Let ${\bf E}, {\bf F}$ be two rank $e$ vector bundles on $X$.
 Suppose that $\sigma:{\bf E}\mapsto {\bf F}$ is a bundle morphism on $X$
 isomorphic off
 $Z$, and that $s_0:X\mapsto {\bf E}$ is a global section of ${\bf E}$ inducing
 the global section of ${\bf F}$, 
$s=\sigma(s_0):X\mapsto {\bf F}$.  According to proposition \ref{prop; section},
 there exists a union of irreducible cones $\cup_{i\in I=\{1, 2\cdots, n\}}
 {\bf C}_{\rho_i}$ 
supported over $Z$ such that $Z(s)=Z(s_0)\cup \cup_{i>0}
\pi_{\bf E}({\bf C}_{\rho_i}\cap s_0(X))$.

 Let $I=\coprod_{1\leq p\leq r}I_p$ be a partition of the index set $I$
 into disjoint 
subsets $I_p\subset I$. Consider the $r$-consecutive scheme theoretical
 blowing ups
 of $X$ along the strict transformations of 
$\cup_{i\in I_p}\pi_{\bf E}({\bf C}_{\rho_i}\cap s_0(X))$,
 and denote the resulting blown up scheme by $\tilde{X}$.
Let $f:\tilde{X}\mapsto X$ to be the $r$-compositions of
 blowing down projection maps.

 Let $D=f^{-1}(Z(s)\cap Z)
\subset \tilde{X}$ be the exceptional Cartier divisor in $\tilde{X}$.
 Let $s':\tilde{X}:f^{\ast}{\bf F}\otimes {\cal O}(-D)$ be the residual section
 which maps to $s$ through $f^{\ast}{\bf F}\otimes {\cal O}(-D)\mapsto 
f^{\ast}{\bf F}$.
 Let $i:Z(s')\mapsto \tilde{X}$ be the inclusion map, then
 $$c_e(f^{\ast}{\bf E})\cap [\tilde{X}]\stackrel{n}{=}
i_{\ast}{\bf Z}(s')\in {\cal A}_{m-e}(\tilde{X}),$$

 i.e. they define the same cap product operation
 from ${\cal A}_e(\tilde{X})$ to ${\cal A}_0(pt)\cong {\bf Z}$.
\end{prop}

\medskip

 Even though $s'$ is not directly related to $s_0$ and ${\bf E}$, the cycle 
 ${\bf Z}(s')$ still defines a version of localized top Chern class 
"localized" in $Z(f^{\ast}s)$ away from $D$. The proposition implies that
 its image under the push-forward morphism $i_{\ast}$ is numerically
 equivalent to $c_e(f^{\ast}{\bf E})$.

\medskip

\noindent Proof of proposition \ref{prop; equivalent}: The main idea of the 
proof is to construct an ambient space containing $\tilde{X}$, some
 auxiliary vector bundles and sections which are used to relate both sides of the
 equality.
Define 
$Y={\bf P}({\bf E}\oplus {\bf C})$ and let $\pi_Y:Y\mapsto X$ denote the projection
 map. Through the map ${\bf v}\mapsto ({\bf v}, 1)$ the
 total space of the vector bundle ${\bf E}$ can be viewed as an open
 subspace of $Y$, which is 
the complement of the closed hypersurface ${\bf P}({\bf E})\subset Y$ at
 infinity. Thus,
 ${\bf P}({\bf E})$ can be viewed as the compactification at infinity of 
 ${\bf E}\subset {\bf P}({\bf E}\oplus {\bf C})$.

 Consider the hyperplane line bundle on $Y$, denoted as ${\bf H}$. 
 Then the projection map ${\bf E}\oplus {\bf C}\mapsto {\bf C}$ over $X$ induces
 a section of ${\bf H}$ vanishing exactly at ${\bf P}({\bf E})$.
 Then $[{\bf P}({\bf E})]\in {\cal A}_{\cdot}(Y)$ is equal to $c_1({\bf H})\cap 
[Y]$. 
 On the other hand, the zero section $s_{\bf E}(X)$ embedded in ${\bf E}\subset Y$
 can be 
viewed as the zero locus of a canonical 
section of $\pi_Y^{\ast}{\bf E}\otimes {\bf H}$ determined 
by the bundle map ${\bf E}\oplus {\bf C}\mapsto {\bf E}$. Thus 
$[s_{\bf E}(X)]=c_e(\pi_Y^{\ast}{\bf E}\otimes {\bf H})\cap [Y]$.

 The composition of ${\bf E}\oplus {\bf C}\mapsto {\bf E}$ and $\sigma:{\bf E}
\mapsto {\bf F}$
 induces a tautological section $\underline{s}$ of $\pi_Y^{\ast}{\bf F}\otimes 
{\bf H}$ on $Y$.
 The following lemma characterizes its zero locus $Z(\underline{s})$.

\medskip

\begin{lemm}\label{lemm; ch}
 Let ${\bf C}_{\rho}$ denote the algebraic sub-cone of ${\bf E}$ corresponding to
 $Ker(\sigma)$. Then ${\bf C}_{\rho}$ can be identified canonically 
with a locally closed sub-scheme of
 $Y$ and $Z(\underline{s})\subset Y$ is the closure of ${\bf C}_{\rho}$, 
${\bf P}({\bf C}_{\rho}\oplus 1)$, in $Y$.
\end{lemm}

\medskip

\noindent Proof of the lemma: A point in $Y={\bf P}({\bf E}\oplus {\bf C})$ is
 inside $Z(\underline{s})$ if and only if the corresponding ray in
 ${\bf E}\oplus {\bf C}$ maps to zero under ${\bf E}\oplus {\bf C}\mapsto 
{\bf E}\mapsto {\bf F}$.
 In other words, when the ray is in the direction
 inside the cone ${\bf C}_{\rho}\oplus 1$ corresponding to the kernel of
 ${\bf E}\oplus {\bf C}\mapsto {\bf F}$. So $Z(\underline{s})=
{\bf P}({\bf C}_{\rho}\oplus 1)$.
 $\Box$

On the other hand, the space $Y$ can be viewed as $1$-plane Grassmanian bundle of
 ${\bf E}\oplus {\bf C}$ over $X$. Viewed as a universal object, one may use it 
 as our playground to prove proposition \ref{prop; equivalent}.

Firstly write ${\bf C}_{\rho}$ as $\cup_{0\leq i\leq n}{\bf C}_{\rho_i}$
 with ${\bf C}_{\rho_0}=s_{\bf E}(X)$. Then 
$\cup_{i\in I}{\bf C}_{\rho_i}={\bf C}_{\rho}-s_{\bf E}(X)$ is a union of
 irreducible sub-cones and $G={\bf P}((\cup_{i\in I}{\bf C}_{\rho_i})\oplus 1)$
 defines a closed sub-scheme of $Y$.

 Notice that $G=\cup_{i\in I}{\bf P}({\bf C}_{\rho_i}\oplus 1)$ and
 we may set $G_l=\cup_{i\in I_l}{\bf P}({\bf C}_{\rho_i}\oplus 1)$. Then
 we may write 
$G=\cup_{1\leq l\leq r}G_l$.  It is obvious that $Z(\underline{s})=s_{\bf E}(X)\cup
 G$.

\medskip  

Secondly one blows up $Y$ consecutively along (the strict transformations 
under previous blowing ups of)
$G_p, 1\leq p\leq r$, following the exactly the same blowing
 up orders to construct $\tilde{X}$ from $X$. Denote the resulting
 scheme \footnote{We have skipped the dependence of $\tilde{Y}$ on the
 choices of the blowing ups.}
$\tilde{Y}$ and denote the union of the 
resulting exceptional Cartier divisors $D_Y$.
 Denote the composite blowing down map $\tilde{Y}\mapsto Y$ by $f_Y$.
 Because
 $G\subset Z(\underline{s})$ and the sub-scheme $G$ has been
 blown up consecutively to get $\tilde{Y}$, the pull-back section 
$(f_Y)^{\ast}\underline{s}$ is divisible by the defining section of  
$D_Y$. Let $\underline{s}'$ denote the residual section in
 $(f_Y)^{\ast}(\pi_Y^{\ast}{\bf F}\otimes {\bf H})\otimes {\cal O}(-D_Y)$.
 Then by proposition \ref{prop; rif} 
$f_Y^{\ast}\underline{s}$ is the image of $\underline{s}'$ under
 $(f_Y)^{\ast}(\pi_Y^{\ast}{\bf F}\otimes {\bf H})\otimes {\cal O}(-D_Y)\mapsto
 (f_Y)^{\ast}(\pi_Y^{\ast}{\bf F}\otimes {\bf H})$.

\medskip

 Consider the closure of $s_{\bf E}(X)-G$ in $\tilde{Y}$, which is nothing but the
 strict (proper)
 transformation of $s_{\bf E}(X)\subset Y$ under the composite blowing ups. We
 denote
 the resulting scheme by $R$. Because $\underline{s}'$ is the residual section
 of $(f_Y)^{\ast}\underline{s}$ vanishing on $R\cup D_Y$,
it is clear that the zero locus of $\underline{s}'$ in $\tilde{Y}$, 
$Z(\underline{s}')$, is equal to $R$, and is of codimension $e$ in
 $\tilde{Y}$. The section 
$\underline{s}'$ may not be regular since $R$ may not be always smooth. 
Nevertheless by example 14.3.1. on 
 page 251 of [F], when $[Z(\underline{s}')]=\sum_i {\bf m}_i [\Omega_i]$, 
 we know that $[{\bf Z}(\underline{s}')]=\sum_i {\bf e}_i [\Omega_i]$
 with ${\bf e}_i\leq {\bf m}_i$.
 But its zero locus 
$R=Z(\underline{s}')$ is birational to $s_{\bf E}(X)\cong X$, the initial base 
space. 
 Because $X$ is reduced, so is $R$, then $[R]=m[R_{red}]$ with $m=1$. 
 Thus we may still conclude that 
 ${\bf Z}(\underline{s}')=[R]$ without the regularity assumption
\footnote{Nevertheless, $R$ is still of the right codimension and is regular
 on a dense open subset.} on $R$.
 \label{step; second}
 Let $i_R$ denote the inclusion $i_R:R\mapsto \tilde{Y}$. 
Then $i_{R\ast}[R]=i_{R\ast}[Z(\underline{s}')]=i_{R\ast}{\bf Z}(\underline{s}')
=c_e(f_Y^{\ast}(\pi_Y^{\ast}{\bf F}\otimes 
{\bf H})\otimes {\cal O}(-D_Y))\cap [\tilde{Y}]$.

\medskip

Thirdly the bundle map ${\bf E}\oplus {\bf C}\mapsto {\bf E}$ on $X$ induces 
a tautological regular 
section $s_{tauto}$ of $\pi_Y^{\ast}{\bf E}\otimes {\bf H}$ on $Y$ vanishing
 at $s_{\bf E}(X)\subset Y$. The pull-back section $(f_Y)^{\ast}s_{tauto}$ of
 $(f_Y)^{\ast}(\pi_Y^{\ast}{\bf E}\otimes {\bf H})$ defines a zero locus
 $Z((f_Y)^{\ast}s_{tauto})=(f_Y)^{-1}(Z(s_{tauto}))=(f_Y)^{-1}(s_{\bf E}(X))$. 
Because
 none of the sub-cones $G_i$ we blow up is contained in $s_{\bf E}(X)=
{\bf C}_{\rho_0}$, 
 the sub-scheme $f_Y^{-1}(s_{\bf E}(X))$ can be identified with the closure
 of $s_{\bf E}(X)-G$ in $\tilde{Y}$, which is nothing but $R$. By the same 
reasoning
 as  above, we have

\label{step; third}
$$i_{R\ast}[R]=
i_{R\ast}[Z((f_Y)^{\ast}s_{tauto})]=i_{R\ast}
{\bf Z}((f_Y)^{\ast}s_{tauto})=
c_e((f_Y)^{\ast}(\pi_Y^{\ast}{\bf E}\otimes {\bf H}))\cap [\tilde{Y}].$$

\medskip

 Fourthly the section $s_0(X)\subset {\bf E}\subset {\bf P}({\bf E}\oplus 1)$ 
can be viewed as a sub-scheme in $Y$, 
denoted by the same symbol. Because $\pi_{\bf E}|_{s_0(X)}:s_0(X)\mapsto X$ induces
 an isomorphism and 
$s_0(X)\cap {\bf P}({\bf E})=\emptyset$ (${\bf P}({\bf E})$ is at infinity), 
the hyperplane line bundle ${\bf H}\mapsto Y$ is trivialized
 over $s_0(X)$ by its cross section which vanishes exactly on ${\bf P}({\bf E})$.
 So
${\bf H}|_{s_0(X)}\cong {\bf C}$, and $\pi_Y^{\ast}{\bf E}\otimes 
{\bf H}|_{s_0(X)}\cong
 {\bf E}$. Then under the Gysin homomorphism the $s_0$ 
pull-back of the formula $c_e(\pi_Y^{\ast}{\bf E}\otimes {\bf H})\cap [Y]=
[s_{\bf E}(X)]\in 
 {\cal A}_{\cdot}(Y)$
 by $s_0^{\ast}:{\cal A}_{\cdot}(Y)\mapsto {\cal A}_{\cdot-e}(X)$ 
 becomes $c_e({\bf E})\cap [X]=s_0^{\ast}[X]\in {\cal A}_{\cdot}(X)$.

 Next we construct an embedding of the blown up scheme $\tilde{X}$ into the
 blown up projective bundle $\tilde{Y}$.

 One notices that the intersection 
$s_0(X)\cap {\bf P}({\bf C}_{\rho_i}\oplus 1)$
 lies inside the cone ${\bf C}_{\rho_i}$ and is equal to 
$s_0(X)\cap {\bf C}_{\rho_i}$.
 By corollary 7.15. on page 165 of [Ha] and the subsequent definition,
 the strict transform of a closed sub-scheme of a scheme theoretical blowing up 
along
 a given blowing up center can be identified to be the blowing up 
 of this sub-scheme along its intersection with the given blowing up center.
  Thus one finds that the closure of $s_0(X)-\cup_{i\in I}{\bf C}_{\rho_i}$
 inside $\tilde{Y}$ is isomorphic to $\tilde{X}$, the consecutive blowing ups of
 $X$ along strict transforms of the forms 
$\pi_{\bf E}((\cup_{i\in I_l}{\bf C}_{\rho_i})\cap s_0(X))$ with $l=1, 2, 
\cdots, r$.
 By abusing the notation slightly, we
 fix such an identification and still denote 
the resulting sub-scheme of $\tilde{Y}$ by the same symbol $\tilde{X}$.

 We have established the following crucial facts after identifying the strict 
transformation of $s_0(X)$ in $\tilde{Y}$ with $\tilde{X}$,

\noindent (i). ${\bf H}|_{\tilde{X}}\cong {\bf C}$.

\medskip

\noindent (ii). $D_Y|_{\tilde{X}}=D$.

\medskip

\noindent (iii). $(f_Y)^{\ast}(\pi_Y^{\ast}{\bf F}\otimes {\bf H})\otimes 
{\cal O}(-D_Y)|_{\tilde{X}}=f^{\ast}{\bf F}\otimes {\cal O}(-D)$.

and 

\medskip

\noindent (iv). The sections $\underline{s}'$, $s'$ of the vector 
bundles in (iii) are compatible. Namely, $\underline{s}'|_{\tilde{X}}=s'$.

\medskip

\medskip

Set $i_{\tilde{X}}:\tilde{X}\mapsto \tilde{Y}$ and set $i:Z(s')
\mapsto \tilde{X}$ to be the inclusion maps.

Then we may conclude that for all $\alpha\in {\cal A}_e(\tilde{X})$ the
 following identification argument:
 By using (ii)., (iii)., the projection formula of Chern classes (see page
 3.2.(c), page 50 of [F]) and 
the relationship between the global and the
 localized top Chern classes of $(f_Y)^{\ast}(\pi_Y^{\ast}{\bf F}\otimes
 {\bf H}\otimes {\cal O}(-D_Y))$,

$$i_{\tilde{X}\ast}\{c_e(f^{\ast}{\bf F}\otimes {\cal O}(-D))\cap \alpha\}=
i_{\tilde{X}\ast}(c_e(i_{\tilde{X}}^{\ast}(f_Y)^{\ast}(\pi_Y^{\ast}{\bf F}\otimes
 {\bf H}\otimes {\cal O}(-D_Y)))\cap \alpha)$$
$$=c_e((f_Y)^{\ast}(\pi_Y^{\ast}{\bf F}\otimes
 {\bf H}\otimes {\cal O}(-D_Y)))\cap [\tilde{Y}]\cap i_{\tilde{X}\ast}\alpha
=i_{R\ast}{\bf Z}(\underline{s}')\cap i_{\tilde{X}\ast}\alpha.$$

And by the defining formula of the localized top Chern class of $\underline{s}'$ 
 and the concluding equality of the second and the third 
statements on page \pageref{step; second},

$$=i_{R\ast}(\underline{s}')^{!}[\tilde{Y}]\cap i_{\tilde{X}\ast}\alpha
=i_{R\ast}[R]\cap i_{\tilde{X}\ast}\alpha
=(c_e(f_Y^{\ast}(\pi_Y^{\ast}{\bf E}\otimes {\bf H}))
\cap [\tilde{Y}])\cap i_{\tilde{X}\ast}\alpha.$$

 Then by projection formula of Chern classes again, 

$$=c_e((f_Y)^{\ast}(\pi_Y^{\ast}{\bf E}\otimes {\bf H}))
\cap i_{\tilde{X}\ast}\alpha=i_{\tilde{X}\ast}\{c_e(i_{\tilde{X}}^{\ast}
(f_Y)^{\ast}(\pi_Y^{\ast}{\bf E}\otimes {\bf H}))\cap \alpha\}.$$

 Then by (i). above we know ${\bf H}|_{\tilde{X}}={\bf C}$, so finally the
 above expression 

$$\hskip -1.in=i_{\tilde{X}\ast}\{c_e(f^{\ast}({\bf E})\cap i_{\tilde{X}\ast}
\alpha\}.$$

Because $\tilde{X}\mapsto pt$ factorizes as 
$\tilde{X}\stackrel{i_{\tilde{X}}}{\longrightarrow}\tilde{Y}\mapsto pt$, 
this implies that 
 $$\int_{\tilde{X}}c_e(f^{\ast}{\bf F}\otimes {\cal O}(-D))\cap \alpha=
\int_{\tilde{X}}c_e(f^{\ast}{\bf E})\cap i_{\tilde{X}\ast}\alpha\in 
{\cal A}_0(pt),$$
 for all $\alpha\in {\cal A}_e(\tilde{X})$.

Therefore $c_e(f^{\ast}{\bf F}\otimes {\cal O}(-D))\cap [\tilde{X}]
\stackrel{n}{=}c_e(f^{\ast}{\bf E})\cap
[\tilde{X}]$. The proposition is proved. $\Box$

\medskip

\begin{rem}\label{rem; num}
 In the proof of this proposition, if we have the knowledge of the
 injectivity $0\mapsto {\cal A}_{m-e}(\tilde{X})\mapsto {\cal A}_{m-e}(\tilde{Y})$,
 then our argument implies a stronger result that 
$c_e(f^{\ast}{\bf F}\otimes {\cal O}(-D))\cap [\tilde{X}]=
c_e(f^{\ast}{\bf E})\cap [\tilde{X}]\in {\cal A}_{m-e}(\tilde{X})$.
 In our paper's main application to the algebraic family Seiberg-Witten
 invariants, the top Chern classes are paired with other cycle classes and 
then pushed-forward to
 ${\cal A}_0(pt)$ to form
 algebraic family Seiberg-Witten invariants. Our result ensures that one may
 replace $c_e(f^{\ast}{\bf F}\otimes {\cal O}(-D))$ by $c_e(f^{\ast}{\bf E})$ 
whenever
 $c_e(f^{\ast}{\bf F}\otimes {\cal O}(-D))$ appears in an 
integration of the top intersection pairing.
 For this purpose, one may view them as ``equal'' without causing potential
 confusion.
\end{rem}

\medskip

 There are many different blowing up sequences which can bring 
 $\cup_{i>0}{\bf P}({\bf C}_{\rho_i}\oplus 1)$ into an exceptional Cartier divisor.
 If one chooses to blow up the whole $\cup_{i>0}{\bf P}({\bf C}_{\rho_i}\oplus 1)$
 in ${\bf P}({\bf E}\oplus {\bf C})$ all at once, the proposition implies that
 $c_e({\bf F})\cap [X]-c_e({\bf E})\cap [X]$ is equal to the push-forward of 
the local contribution of top Chern class from $Z=\cup_{i>0}supp({\bf C}_{\rho_i})$
 into $X$. If one chooses to group various irreducible 
${\bf P}({\bf C}_{\rho_i}\oplus 1)$
 into different sub-schemes and blow up consecutively, one may
 apply proposition \ref{prop; rif} and lemma \ref{lemm; local} inductively
 and get a sum of cycles supported in $Z=supp(\cup_{i>0}{\bf C}_{\rho_i})$.
 It is natural to wonder if the result is 
invariant to the various choices of the orders of 
the blowing ups.

 A corollary of proposition \ref{prop; equivalent} is the following,

\medskip

\begin{cor}\label{cor; equal}
 With the same notations as in proposition \ref{prop; equivalent}, 
the expression $\sum_{1\leq i\leq e}(-1)^{i-1}c_{e-i}({\bf F}) \cap D^{i-1}
\cap [D]$
 is numerically equivalent
 to $\{c({\bf F}|_{Z(s)\cap Z})\cap s(Z(s)\cap Z, X)\cap [X]\}_{m-e}$ and
 is therefore independent to the choices of blowing up processes.
\end{cor}

\medskip

\noindent Proof: 
Let 
$I=\coprod_{1\leq l\leq r}I_l$ be a partition of the index set $I$ and
 let 

$$\tilde{X}=\tilde{X}_r\mapsto \tilde{X}_{r-1}\mapsto \tilde{X}_{r-2}\mapsto
 \cdots \tilde{X}_1\mapsto \tilde{X}_0=X$$

be the sequence of blowing up processes where the $l-$th blowing up
 $\tilde{X}_l\mapsto \tilde{X}_{l-1}$ is centered at the strict transform
 of $\pi_{\bf E}({\bf s}_0\cap (\cup_{i\in I_l}{\bf C}_{\rho_i}))$ under
 $\tilde{X}_{l-1}\mapsto \tilde{X}_0$. 

First we notice that the section $s_0$ does not intersect with the infinity of
 the projective bundle ${\bf P}({\bf E})\subset {\bf P}({\bf E}\oplus 1)$ 
and therefore 
$s_0\cap (\cup_{i\in I_l}{\bf C}_{\rho_i})
=s_0\cap (\cup_{i\in I_l}{\bf P}({\bf C}_{\rho_i}\oplus
 1))$ for all $1\leq l\leq r$.

Let us fix a few notations. Let $\tilde{h}_l:\tilde{X}\mapsto \tilde{X}_l$ be the 
blowing down map from the final ($r-$th blowing up) to the $l-$th 
intermediate blowing up of $X$.
Let $D_l$, $1\leq l\leq r$ denote
 the exceptional divisor of $\tilde{X}_l\mapsto \tilde{X}_{l-1}$ 
indexed by the subscript
 $l$. Let $\tilde{D}_l$ denote its pre-image $\tilde{h}_l^{-1}(D_l)
\subset \tilde{X}$. 

 Then by an induction argument, proposition \ref{prop; rif} implies the following
 identities on the Chern classes,

$$c_e(\tilde{h}_0^{\ast}{\bf F}\otimes 
\otimes_{j\leq l-1}\tilde{h}_j^{\ast}{\cal O}(-D_j))\cap [\tilde{X}_r]-
c_e(\tilde{h}_0^{\ast}{\bf F}\otimes \otimes_{j\leq l}\tilde{h}_j^{\ast}
{\cal O}(-D_j))\cap [\tilde{X}_r]$$
$$=i_{\tilde{D}_l\ast}\sum_{1\leq i\leq e}(-1)^{i-1}
c_{e-i}(\tilde{h}_0^{\ast}{\bf F}\otimes 
\otimes_{j\leq l-1}{\cal O}(-\tilde{D}_j)) \cap 
\tilde{D}_l^{i-1}\cap [\tilde{D}_l],$$

 for $1\leq l\leq r$.
 If we sum up all these $r$ equations, a simple cancellation of the
 intermediate terms leads to
 the final equation

$$\hskip -.9in 
c_e(\tilde{h}_0^{\ast}{\bf F})\cap [\tilde{X}]-
c_e(\tilde{h}_0^{\ast}{\bf F}\otimes 
\otimes_{1\leq l\leq r}{\cal O}(-\tilde{D}_l))\cap
[\tilde{X}]=\sum_{1\leq l\leq r}
i_{\tilde{D}_l\ast}\sum_{1\leq i\leq e}(-1)^{i-1}c_{e-i}(\tilde{h}_0^{\ast}
{\bf F}\otimes 
\otimes_{j\leq l-1}{\cal O}(-\tilde{D}_j)) \cap \tilde{D}_l^{i-1}\cap 
[\tilde{D}_l].$$

  By realizing ${\cal O}(D)=\otimes_{1\leq l\leq r}{\cal O}(\tilde{D}_l)$,
 the left hand side of the identity is $c_e(\tilde{h}_0^{\ast}{\bf F})\cap 
[\tilde{X}]-
c_e(\tilde{h}_0^{\ast}{\bf F}\otimes {\cal O}(-D))\cap [\tilde{X}]$. Thus 
 the right hand side can be identified with
$\sum_{1\leq i\leq e}(-1)^{i-1}c_{e-i}(f^{\ast}{\bf F}) \cap D^{i-1}\cap [D]$ 
by the well known property
 on the top Chern class.
 On the other hand, proposition \ref{prop; equivalent} has implied that
 the push-forward of the left hand side under
 $\tilde{X}\mapsto X$ is numerically equivalent
 to $c_e({\bf F})\cap [X]-c_e({\bf E})\cap [X]$ and 
therefore is
 independent to all the grouping and the ordering 
choices involved in the blowing ups of $\pi_{\bf E}(s_0\cap 
(\cup_{i>0}{\bf C}_{\rho_i}))\subset X$.
 In particular,
 we may take $I=I_1, r=1$ to be the single partition of $I$ and the 
$f:\tilde{X}\mapsto X$ is
 constructed from $X$ by a single blowing up centered at
 $\pi_{\bf E}({\bf s}_0\cap (\cup_{i>0}
{\bf C}_{\rho_i}))$. In this case, the identification of
 $f_{\ast}i_{D\ast}\sum_{1\leq i\leq e}(-1)^{i-1}c_{e-i}(f^{\ast}{\bf F}) 
\cap D^{i-1}\cap [D]$ with the local contribution of top Chern class
 $\{c({\bf F}|_{Z(s)\cap Z})\cap s(Z(s)\cap Z, X)\cap [X]\}_{m-e}$ 
is the direct consequence of lemma \ref{lemm; local}. $\Box$

\medskip

\begin{rem}\label{rem; segreiden}
 An alternative way to prove corollary \ref{cor; equal} and show that
 it is independent to the ordering of the blowing ups is to notice that
 $f(D)$ is always equal to $Z(s)\cap Z$ no matter which blowing up
 sequence we choose. Then by realizing 
$\sum_{i>0}(-1)^{i-1}D^{i-1}=s(D, \tilde{X})$ and by using proposition 
4.2.(a) of\footnote{Cited in proposition \ref{prop; birational}.} [F], 

$$f_{\ast}s(D, \tilde{X})=deg(\tilde{X}/X)s(f(D), f(\tilde{X}))=1\cdot
 s(Z(s)\cap Z, X),$$ 
 one may identify $f_{\ast}\bigl(c(f^{\ast}{\bf F})\cap 
(\sum_{i>0}(-1)^{i-1}D^{i-1}[D])
\bigr)$ with $(c({\bf F})\cap s(Z(s)\cap Z, X)$.
\end{rem}

 The observation in remark \ref{rem; segreiden} will be used in proving
 proposition \ref{prop; identical} in section \ref{subsection; independence}.

\medskip

\subsection{\bf Some Observation about Residual Intersection Formula of Top Chern 
Classes}\label{subsection; subbundle}

\bigskip

 In this subsection, we consider the following geometric
 setting. Let $X$ be a purely
$m$ dimensional scheme and let ${\bf E}$ be a rank $e$ vector bundle over $X$.
Suppose that ${\bf E}_0$ is a rank $e_0$ 
sub-bundle of ${\bf E}$ with a section $s_0:X\mapsto {\bf E}_0$,
 and suppose that we have the following bundle short exact sequence,

$$0\mapsto {\bf E}_0\mapsto {\bf E}\mapsto {\bf E}/{\bf E}_0\mapsto 0.$$

The section $s_0$ and the 
bundle injection ${\bf E}_0\mapsto {\bf E}$ induces a section $s:X\mapsto {\bf E}$
 and we know $Z(s_0)=Z(s)$.

 We raise the following question:

\noindent {\bf Question}: How are the residual intersection 
formulae of the top Chern classes of
 ${\bf E}_0$ and ${\bf E}$ related to each other?

 More precisely, let $Z_1, Z_2, \cdots, Z_k$ be a finite number of
 closed proper sub-schemes of $X$. One may blow up the strict transforms 
(under the previous blowing ups) of 
$Z(s_0)\cap Z_i$, $1\leq i\leq k$ consecutively and get a residual 
intersection formulae
 of top Chern classes of ${\bf E}_0$. On the other hand, we may blow up the
 strict transforms of $Z(s)\cap Z_i$, $1\leq i\leq k$ consecutively and
get the residual top Chern classes of ${\bf E}$. Because $Z(s_0)=Z(s)$, we expect
 these two residual intersection 
formulae to be closed related. This is the content of the
 following proposition,

\medskip

\begin{prop}\label{prop; cap}
 Let $\tilde{X}$ denote the scheme repeatedly blown up from $X$ centered at the
 strict transforms of 
 $Z(s_0)\cap Z_i$, $1\leq i\leq k$ and let $f:\tilde{X}\mapsto X$ denote the
 composite blow down projection map.  Let $D=f^{-1}(\cup_i Z_i\cap Z(s_0))$ be the
 exceptional Cartier divisor above $\cup_iZ_i\cap Z(s_0)=\cup_iZ_i\cap Z(s)$.
 Let $(f^{\ast}s_0)'$ and
 $(f^{\ast}s)'$ denote the residual sections in 
$f^{\ast}{\bf E}_0\otimes {\cal O}(-D)$ and $f^{\ast}{\bf E}\otimes {\cal O}(-D)$ of
 $f^{\ast}s_0\in \Gamma(\tilde{X}, f^{\ast}{\bf E}_0)$ and $f^{\ast}s\in
 \Gamma(\tilde{X}, f^{\ast}{\bf E})$, 
 respectively. 

 By proposition \ref{prop; rif} there is a residual intersection formula
 of the localized top Chern class of $f^{\ast}{\bf E}_0$,

$${\bf Z}(f^{\ast}s_0)={\bf Z}((f^{\ast}s_0)')+
\sum_{1\leq i\leq e_0}(-1)^{i-1}c_{e_0-i}({\bf E}_0|_D)
\cap D^{i-1}\cap [D].$$

Suppose we cap the above formula
 with the top Chern class 
$c_{e-e_0}(f^{\ast}({\bf E}/{\bf E}_0)|_D)$, one gets the corresponding 
residual intersection
 formula of the localized top Chern class of $f^{\ast}s$,

$${\bf Z}(f^{\ast}s)={\bf Z}((f^{\ast}s)')+
\sum_{1\leq i\leq e}(-1)^{i-1}c_{e-i}({\bf E}|_D)
\cap D^{i-1}\cap [D].$$
\end{prop}

\noindent Proof of proposition \ref{prop; cap}: Recall that (see [F] page 244
 and proposition 6.1.(a) page 94), ${\bf Z}(f^{\ast}s_0)$ is equal to
$\{c(f^{\ast}{\bf E}_0)\cap s(Z(f^{\ast}s_0), \tilde{X})\}_{m-e_0}$ 
(we have skipped the bundle restriction notation of $f^{\ast}{\bf E}_0$ to 
$Z(f^{\ast}s_0)$ to simplify the notation). 
By proposition 13 of [Liu5],
$\{c(f^{\ast}{\bf E}_0)\cap s(Z(f^{\ast}s_0), \tilde{X})\}_{m-e_0+r}=0$ for 
all negative $r\in -{\bf N}$, i.e. the localized contribution of top Chern 
 class is the lowest degree term of the
 cycle class formed by the total Chern/Segre classes. 
Because $Z(f^{\ast}s_0)=Z(f^{\ast}s)$, we may
 compare it with\footnote{We have changed $r$ to $r'=r-(e-e_0)$ in the
 second line of the equalities.} 

$$\hskip -1.3in
\{c(f^{\ast}{\bf E})\cap s(Z(f^{\ast}s), \tilde{X})\}_{m-e}
=\{c(f^{\ast}{\bf E}_0)\cap c(f^{\ast}({\bf E}/{\bf E}_0)) 
\cap s(Z(f^{\ast}s), \tilde{X})\}_{m-e}
=\{c(f^{\ast}{\bf E}_0)\cap c(f^{\ast}({\bf E}/{\bf E}_0)) 
\cap s(Z(f^{\ast}s_0), \tilde{X})\}_{m-e}$$
$$\hskip -1.3in 
=\sum_{r\geq 0}\{c(f^{\ast}{\bf E}_0)
\cap s(Z(f^{\ast}s_0), \tilde{X})\}_{m-e_0-(e-e_0)+r}\cap 
c_r(f^{\ast}({\bf E}/{\bf E}_0)) 
=\sum_{r'\geq -(e-e_0)}\{c(f^{\ast}{\bf E}_0)
\cap s(Z(f^{\ast}s_0), \tilde{X})\}_{m-e_0+r'}\cap 
c_{r'+(e-e_0)}(f^{\ast}{\bf E}/{\bf E}_0).$$

 Because ${\bf E}/{\bf E}_0$ is of rank $e-e_0$, 
$c_{r'+(e-e_0)}(f^{\ast}({\bf E}/{\bf E}_0))=0$ for all
 $r'\in {\bf N}$. So by the vanishing of terms with grades $m-e_0+r'$ for
 $r'<0$ in $\{c(f^{\ast}{\bf E}_0)\cap s(Z(f^{\ast}s_0), \tilde{X})
\}$ the above sum is reduced to a single term and the result
is nothing but the cap product of
 $\{c(f^{\ast}{\bf E}_0)
\cap s(Z(f^{\ast}s_0), \tilde{X})\}_{m-e_0}$ with 
$c_{e-e_0}(f^{\ast}{\bf E}/{\bf E}_0)$.

 The same discussion can be applied to ${\bf Z}((f^{\ast}s_0)')$ and
 ${\bf Z}((f^{\ast}s)')$ as well. Then the correspondence of the two formulae
 follows as the correspondence under capping with  
 $c_{e-e_0}({\bf E}/{\bf E}_0)$ has been shown to hold for 
two out of the three terms. $\Box$

\medskip

\begin{rem}\label{rem; finish}
 When there is only one $Z=Z_1$ and $\tilde{X}$ is the blown up of $X$ along 
$Z_1\cap Z(s)$,
 the direct proof of the equality of $c_{e-e_0}(f^{\ast}({\bf E}/{\bf E}_0)|_D)
\cap 
\sum_{1\leq i\leq e_0}(-1)^{i-1}c_{e_0-i}({\bf E}_0|_D)
\cap D^{i-1}\cap [D]$ and $\sum_{1\leq i\leq e}(-1)^{i-1}c_{e-i}({\bf E}|_D)
\cap D^{i-1}\cap [D]$ follows from the identification with the
 localized top Chern class 
contribution of $Z\cap Z(s)$ as was done in lemma \ref{lemm; local}.
 When more than one $Z_i$ is present, one may adopt the similar argument in
 the proof of corollary \ref{cor; equal} and lemma \ref{lemm; local}
 to identify them directly.
 As there is no new idea involved, we leave the details to
 the readers.
\end{rem}

\medskip

\section{\bf Residual Intersection Formula and Inductive Blowing Ups
 of $X={\bf P}({\bf V}_{canon})$}\label{section; main}

\bigskip

\medskip

 We follow the same notations as in [Liu1], [Liu3] and [Liu5]. Let
 ${\bf M}(E)E=\sum_{i\leq n}m_iE_i$ be the sum of the exceptional divisors
 $E_i$ with multiplicities $m_i$, $1\leq i\leq n$. To simplify our discussion,
 we require that $m_i\leq m_j$ for all $i<j$.

 As in [Liu3], we take 
$(\Phi_{{\bf V}_{canon}{\bf W}_{canon}}, 
{\bf V}_{canon}, {\bf W}_{canon})$ over $M_n\times T(M)$
to be the canonical algebraic Kuranishi model of the class $C-{\bf M}(E)E$
 with respect to $f_n:M_{n+1}\mapsto M_n$. We take the initial total space 
 of our discussion to be $X={\bf P}_{M_n\times T(M)}({\bf V}_{canon})$. This
 space parametrizes all the curves in the non-linear system of $C$ along 
the family $M_n\times T(M)$.

 Recall from section \ref{section; strata} that the universal space 
$M_n$ allows a stratification by the admissible strata
 $Y_{\Gamma}$, $\Gamma\in adm(n)$. 
We consider the finite set of admissible strata 
satisfying the {\bf Special Condition} first stated in [Liu5] section 6.1:

For all the type $I$ exceptional classes $e_1, e_2, \cdots, e_n$ of $Y_{\Gamma}$, 

$\diamondsuit$ either \label{specialcondition}

\noindent (i). $(C-{\bf M}(E)E)\cdot e_i<0$, i.e. ${\bf M}(E)E\cdot e_i>0$. 

$\diamondsuit$ or 

\noindent (ii). the condition 
$e_i^2=-1$ holds, i.e. $e_i$ is a type $I$ $-1$ class. 

\medskip

 This maximality special condition means there is no ``redundant'' type $I$ classes
 which pair non-negatively with $C-{\bf M}(E)E$.

 Recall from section \ref{section; strata} 
that the notation $adm(n)$ denotes the finite set of all $n$-vertex
admissible graphs $\Gamma$. We introduce a few subsets of $adm(n)$ here.

\begin{defin}\label{defin; special}
Let $\Delta(n)\subset adm(n)$ denote the subset of $adm(n)$
 consisting of all $n$-vertex admissible graphs satisfying the 
 {\bf special condition} above.
 Let $adm_2(n)\subset adm(n)$ denote the 
 subset of $n$-vertex admissible graphs satisfying the
 condition that each vertex has at most
 one direct descendent.
\end{defin}

 The graphs \footnote{See fig.4 for an example.} 
in $adm_2(n)$ may have more than one connected component.
 Each component looks like a chain of vertexes connected by a chain of arrows.
 We will refer to them as chain-like admissible graphs in the following
 discussion.

\begin{figure}
\centerline{\epsfig{file=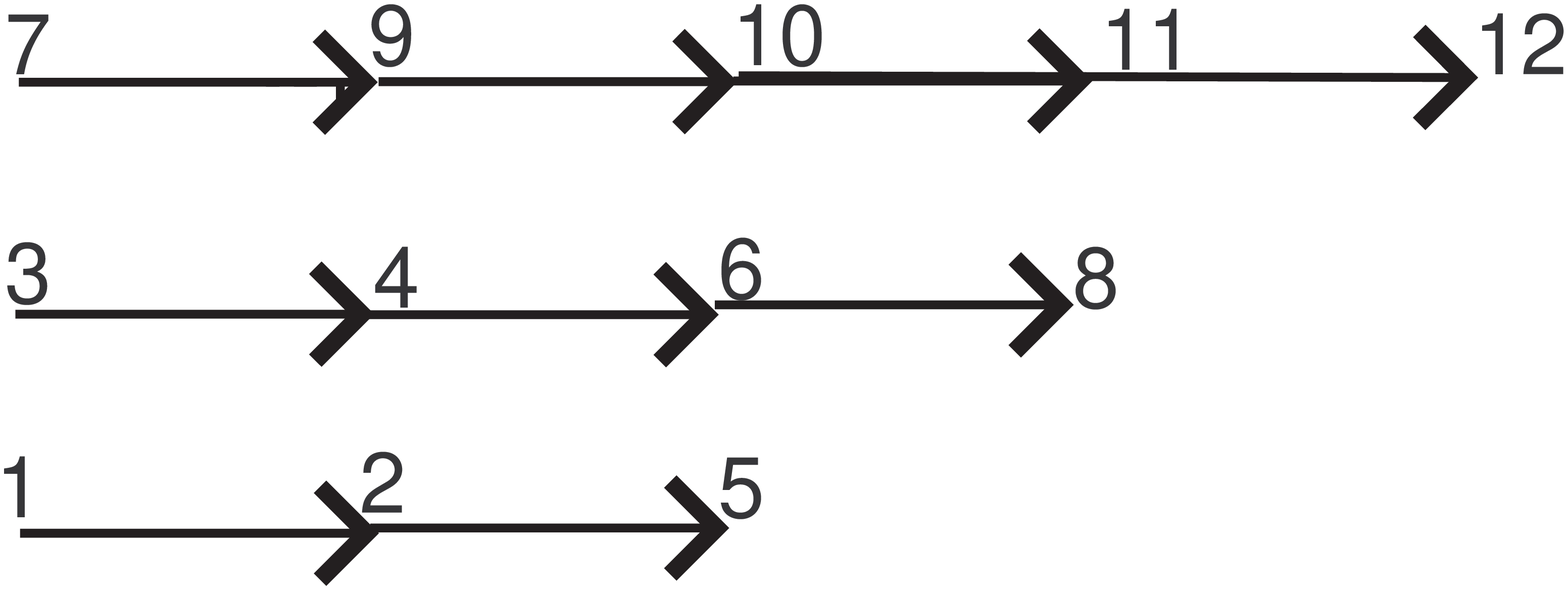,height=4cm}}
\centerline{fig.4}
\centerline{a chain-like admissible graph with three connected components 
 and 12 vertexes}\label{fig.4}
\end{figure}

 Firstly we point out that the union of the closure
 of all such strata $Y_{\Gamma}$, $\Gamma\in \Delta(n)$, form 
a large closed subset of $M_n$. 

\begin{prop}\label{prop; large}
The closed subset $\cup_{\Gamma\in \Delta(n)-\{\gamma_n\}}Y({\Gamma})$ is 
of at least complex codimension
 two in $M_n$. Its complement can be expressed as 
$\coprod_{\Gamma\in P} Y_{\Gamma}$ for some $P\subset adm_2(n)$.
\end{prop}

\medskip

\noindent Proof: Firstly, we know that $\Delta(n)\not=\emptyset$ because
 the admissible graph $\gamma_n$ with no one-edges is in $\Delta(n)$.
 
 We know that each $Y(\Gamma)$ is closed, so 
$\cup_{\Gamma\in \Delta(n), \Gamma\not=\gamma_n}Y({\Gamma})$ is a finite
 union of closed sets and is closed.

We observe that if for some $i\leq n$, the $i$-th vertex of $\Gamma$
 contains more than one direct descendent vertex, then the
 intersection pairing of $e_i=E_i-\sum_{j_i}E_{j_i}$ and $C-{\bf M}(E)E$
is negative because the assumption $m_i\leq m_j$, for $i\leq j$,

$$e_i\cdot (C-{\bf M}(E)E)=(E_i-\sum_{j_i}E_{j_i})\cdot (C-\sum_km_kE_k)
=m_i-\sum_{j_i}m_{j_i}<0.$$

 If the index $i$ has exact one direct descendent in $\Gamma\in adm(n)$, 
denoted as $j$, then the
 intersection pairing of this $-2$ class $e_i=E_i-E_j$ and 
 $C-{\bf M}(E)E$ is $m_i-m_j\leq 0$.

\bigskip

 Take an arbitrary
 $\Gamma\not\in adm_2(n)$, then there must be some index $i\leq n$ such that
 the $i-$th vertex contains more than one direct descendent in $\Gamma$. 
  Given the graph $\Gamma$, one may construct the type $I$ exceptional
 classes associated to it following the recipe in section \ref{section; strata}:
 Given any index $i, 1\leq i\leq n$, let
 $j_i$ be the indexes of all the direct descendents of $i$ in the graph $\Gamma$.
 Then take $e_i=E_i-\sum_{j_i}E_{j_i}$.

Among all such $e_i, 1\leq i\leq n$, consider
all the type $I$ classes $e_{k_i}, 1\leq i\leq p$ (for some $p$ depending on
 both $\Gamma$ and ${\bf M}(E)E$), attached to
 $\Gamma$ which have negative pairings with 
the given $C-{\bf M}(E)E$. Each $e_{k_i}$ is represented by an effective fiberwise 
 divisor in the fiber algebraic surfaces
 of $M_{n+1}\mapsto M_n$ over a smooth locus $Y(\Gamma_{e_{k_i}})$ with a complex 
codimension\footnote{$=codim_{\bf C}\Gamma$.}
 equal to the number of direct descendents of $k_i$ in $\Gamma$.
 This locus is usually called the 
existence locus of $e_{k_i}$.

By proposition \ref{prop; reverse} in section 
\ref{section; strata} one can construct
 an admissible graph $\Gamma_0$ from $e_{k_i}$, $1\leq i\leq p$, satisfying
 $e_{k_i}\cdot e_{k_j}\geq 0$ for $i\not=j$. Then by the corollary of
 proposition \ref{prop; transversal}
 and remark \ref{rem; ignore} we have
 $Y(\Gamma_0)=\cap_{1\leq i\leq p}Y(\Gamma_{e_{k_i}})$ and $Y(\Gamma_0)$ is
 the locus in $M_n$ over which $e_{k_1}, e_{k_2}, \cdots, e_{k_p}$ co-exist
 as effective curve/divisor classes.

 Because these $e_{k_i}$ are effective over $Y_{\Gamma}$, so we know
 that $\Gamma$ is a degeneration of $\Gamma_0$, i.e. $\Gamma<\Gamma_0$,
 and $Y_{\Gamma}\subset Y(\Gamma_0)$. By the construction of $\Gamma_0$, 
 all these $e_{k_i}, 1\leq i\leq p$ are the only type $I$ exceptional
 classes associated to $\Gamma_0$
 which are not $-1$ classes. 
Therefore it satisfies the {\bf special condition} and 
 $\Gamma_0\subset \Delta(n)$.

 Thus we conclude that every stratum $Y_{\Gamma}$, $\Gamma\not\in adm_2(n)$,
 must be contained in $\cup_{\Gamma\in \Delta(n)-\{\gamma_n\}}Y(\Gamma)$.

Secondly, we separate into two cases 
to determine the dimension of $\cup_{\Gamma\in 
\Delta(n)-\{\gamma_n\}}Y(\Gamma)$. If all the singular multiplicities 
 are equal $m_1=m_2=\cdots=m_n$, then
 $C-{\bf M}(E)E$ has vanishing pairings with all the $-2$ type $I$ exceptional 
classes of
 the general form $E_i-E_j$. The admissible graphs in $adm_2(n)$ are exactly
 those graphs whose associated type $I$ exceptional classes are either $-1$
 or $-2$ classes. In this case, the codimension of 
$\cup_{\Gamma\in \Delta(n)-\{\gamma_n\}}Y(\Gamma)$ is exactly two because all the 
codimension one admissible strata are parametrized by some very simple 
graphs \footnote{Indeed, graphs with a single one-edge.} in $adm_2(n)$. 
The index set $P$ can be taken to be the whole $adm_2(n)$.

 If there exists a pair of singular multiplicities 
$m_i<m_j$ for $1\leq i < j\leq n$, then 
 $C-{\bf M}(E)E\cdot (E_i-E_j)<0$. 
 Consider the admissible graph $\Gamma_{i, j}$
 with a unique one-edge starting at the
 $i-$th vertex and ending at the $j-$th vertex. Apparently it belongs to
 $adm_2(n)$. Because $m_i<m_j$, this graph $\Gamma_{i, j}$ also belongs to 
$\Delta(n)$. In this case, codimension of 
$\cup_{\Gamma\in \Delta(n)-\{\gamma_n\}}Y(\Gamma)$ 
is one. The index subset $P$ is then chosen to be the
 proper subset of $adm_2(n)$, removing all the chain-like admissible graphs 
 of the type $\Gamma_{i, j}$ with single
 edges from $i$ to $j$, for all the pairs $m_i<m_j$, $i<j$. $\Box$

\medskip
 
In the paper [Liu4], we have considered 
the set $Q$ and describe a curve-counting scheme
 based on $Q$. In our current setting of type $I$ exceptional classes,
 define $Q$ to be the finite set of all the classes $E_i-\sum_{i<j}E_j$ which
 pair negatively with $C-{\bf M}(E)E$. The Set $Q$ encodes all the possible
 type $I$ exceptional classes which can appear above the family $M_{n+1}\mapsto 
M_n$.

 Recall the definition of the type $I$ exceptional cone over a point $b\in M_n$,

\begin{defin}\label{defin; econe}
Let $b$ be an arbitrary point in $M_n$, define the type $I$ exceptional
 cone over $b$,  ${\cal EC}_b(C-{\bf M}(E)E; Q)$,
 to be the real cone in $H^{1, 1}(M_{n+1}|_b, {\bf R})$ 
 generated by all the type $I$ exceptional classes effective over $b$ which
 pair with $C-{\bf M}(E)E$ negatively.
\end{defin}

 According to proposition 4 of [Liu4], the cone is always simplicial, and we call the
 primitive generators at the 1-edge the extremal generators of the cone. 

  Because the fiber bundle $M_{n+1}\mapsto M_n$ has no non-trivial monodromy,
  one can discuss about the change of the cone un-ambiguously.
  The variation of ${\cal EC}_b(C-{\bf M}(E)E; Q)$ with respect to $b$
 gives us important information about how to organize the admissible
 strata $Y_{\Gamma}, \Gamma\in \Delta(n)$.

 Given a $Y(\Gamma)$, $\Gamma\in \Delta(n)$, the type $I$
 exceptional cone ${\cal EC}_b(C-{\bf M}(E)E; Q)$ may vary when $b$ specializes
 to the boundary points $\coprod_{\Gamma'<\Gamma}
 Y_{\Gamma'}$.

 Let ${\cal C}_{\Gamma}$ denote the type $I$ exceptional cone
 ${\cal EC}_b(C-{\bf M}(E)E; Q)$ constant for all $b\in Y_{\Gamma}$.
 There is a distinguished 
 locally closed subset 
 $S_{\Gamma}\subset Y(\Gamma)$, $S_{\Gamma}$ (should be
 denoted by $S_{{\cal C}_{\Gamma}}$ if we follow the notation in [Liu4])
 over which the exceptional cone 
${\cal EC}_b(C-{\bf M}(E)E; Q)\equiv {\cal C}_{\Gamma}$ remain
 unchanged.

 Because for all $b\in Y_{\Gamma'}$, their 
${\cal EC}_b(C-{\bf M}(E)E; Q)$ remain constant \footnote{Which may
 be different from ${\cal C}_{\Gamma}$ though.}, $S_{\Gamma}$ itself
 is a union of admissible strata and one may write $S_{\Gamma}$ formally as
 $Y_{\Gamma}\coprod_{Y_{\Gamma'}\cap S_{\Gamma}\not=\emptyset} Y_{\Gamma'}$.
 
\medskip

\begin{lemm}\label{lemm; union3}
 The union $\cup_{\Gamma\in \Delta(n)}Y(\Gamma)=Y_{\gamma_n} 
\cup_{\Gamma\in \Delta(n)-\{\gamma_n\}}Y(\Gamma)\subset M_n$ is equal to
 the disjoint union $\coprod_{\Gamma\in \Delta(n)}S_{\Gamma}$.
\end{lemm}

\medskip

\noindent Proof: Following the proof of proposition \ref{prop; large},
 for all $Y_{\Gamma}$ over which at least one type $I$ exceptional class
 pairs negatively with $C-{\bf M}(E)E$, 
$Y_{\Gamma}\subset S_{\Gamma_0}$ for
 some unique $\Gamma_0$ constructed (in the proof of proposition 
 \ref{prop; large} through the usage of proposition \ref{prop; reverse})
 by the co-existence of different type $I$ curves.

 On the other hand, an effective type $I$ exceptional class pairing
 negatively with $C-{\bf M}(E)E$ above $Y_{\Gamma}$ 
still remains effective and pairs negatively
 with $C-{\bf M}(E)E$ over the boundary $\partial Y_{\Gamma}=Y(\Gamma)-Y_{\Gamma}$.
 But it may break into more than one irreducible 
component. This implies that over any such
 degenerated stratum $Y_{\Gamma'}\subset
 Y(\Gamma)$ there must still exist at least one type $I$ class pairing negatively
 with $C-{\bf M}(E)E$. 

 This shows that
$\cup_{\Gamma\subset \Delta(n)}Y(\Gamma)\subset \cup_{\Gamma
\in \Delta(n)} S_{\Gamma}$. The opposite inclusion follows from 
 the inclusion $S_{\Gamma}\subset Y(\Gamma)$ for $\Gamma\in \Delta(n)$. 
 Finally $S_{\Gamma}\cap S_{\Gamma'}=\emptyset$ if $\Gamma\not=\Gamma'$ in
 $\Delta(n)$. It is because ${\cal C}_{\Gamma}\not={\cal C}_{\Gamma'}$ and
 by definition of $S_{\Gamma}$, they can not overlap. So the union is
 a disjoint union. $\Box$

\medskip

Among all such $S_{\Gamma}$, $\Gamma\in \Delta(n)$, 
one may introduce a partial ordering $\succ$, 
 as has been done in [Liu5] for a slightly more general setting.
 The partial ordering induces a partial ordering on the corresponding 
graphs in $\Delta(n)$,
 denoted by the same symbol.

\begin{defin}\label{defin; succ}
  Let $\Gamma_1, \Gamma_2\in \Delta(n)$. The graph $\Gamma_1$ is said to be
 greater than $\Gamma_2$ under the partial ordering $\succ$, denoted as
 $\Gamma_1\succ \Gamma_2$, if ${\cal C}_{\Gamma_1}\subset {\cal C}_{\Gamma_2}$.
\end{defin}

Please refer to page \pageref{fig.7}, fig.7 for an example. In that example,
 the smaller cone is generated by 
 $E_1-E_2-E_3-E_7$, $E_2$, $E_3$, $E_4$, $E_5$, $E_6$, $E_7$.
The larger cone is generated by 
 $E_1-E_2-E_3-E_4-E_7$, $E_2-E_5-E_6$, $E_3$, $E_4$, $E_5$, $E_6$, $E_7.$

A sufficient condition to check whether
 $\Gamma_1\succ \Gamma_2$ is the following,

\begin{lemm}\label{lemm; suff}
Suppose that $\overline{S_{\Gamma_1}}=Y(\Gamma_1)$ intersects
 with $S_{\Gamma_2}$ non-trivially, then $\Gamma_1\succ \Gamma_2$.
\end{lemm}

\noindent Proof: The follows from the fact that the cones get larger under
 degenerations of points from $b\in S_{\Gamma_1}$ to $b\in S_{\Gamma_2}$. $\Box$

\medskip

  Our goal is to study the local contribution of the algebraic family
Seiberg-Witten invariant over $\cup_{\Gamma\in \Delta(n)}Y(\Gamma)$ and decompose
 the algebraic family Seiberg-Witten 
invariant ${\cal AFSW}_{M_{n+1}\times T(M)\mapsto 
 M_n\times T(M)}(1, C-{\bf M}(E)E)$ or the restricted version
 ${\cal AFSW}_{M_{n+1}\times \{t_L\}\mapsto 
 M_n\times \{t_L\}}(1, C-{\bf M}(E)E)$, for some $t_L\in T(M)$, 
into the various excess local contributions
 from $\cup_{\Gamma\in \Delta(n)-\{\gamma_n\}}Y(\Gamma)$ 
and the residual contribution
 from $M_n-\cup_{\Gamma\in \Delta(n)-\{\gamma_n\}}Y(\Gamma)$ based on the following
 two simple but fundamental observations,

\medskip

\noindent {\bf Observation 1}: The family algebraic Seiberg-Witten invariant
 is defined by the push-forward into ${\cal A}_0(pt)$ of 
the cap product of a certain power of $c_1({\bf H})$ (determined
 by the dimension formula) with the top Chern class of the canonical obstruction
bundle 
$c_{top}({\bf H}\otimes \pi_{{\bf P}({\bf V}_{canon})}^{\ast}{\bf W}_{canon})$.

\medskip

\noindent {\bf Observation 2}: The residual intersection formula of top
 Chern class allows us to decompose the total invariant contribution into the
 local contribution to some closed subset of 
$X={\bf P}_{M_n\times T(M)}({\bf V}_{canon})$ and the residual contribution.
 The residual contribution again involves the top Chern class of a modified
 bundle over a blown up space. 

\medskip

  A direct but probably naive approach is to set 
$Z=\pi_X^{-1}(\cup_{\Gamma\in \Delta(n)-\{\gamma_n\}}
Y(\Gamma))\subset X$, the pre-image in the ${\bf P}({\bf V}_{canon})$ of 
 $\cup_{\Gamma\in \Delta(n)-\{\gamma_n\}}
Y(\Gamma)$, and then apply the residual
 intersection formula (proposition \ref{prop; rif} and lemma \ref{lemm; local})
 to the vector bundle 
$E={\bf H}\otimes \pi_X^{\ast}{\bf W}_{canon}$, the section 
$s=s_{canon}$ and $Z\cap Z(s_{canon})$. The apparent drawback
 of this approach is that the local contribution of the set $Z(s_{canon})\cap Z$ 
to the family invariant
  is very hard to enumerative directly, due to the complicated geometric 
structure of $\cup_{\Gamma\in \Delta(n)-\{\gamma_n\}}
Y(\Gamma)$, $Z$ and therefore $Z(s_{canon})\cap Z$.

  Instead we construct a more refined consecutive blowing ups of sub-schemes in
 $X$ and make use of $Y(\Gamma)$
 as the co-existence locus of all the type $I$ exceptional classes $e_i$, $1\leq
 i\leq n$ effectively. 
In section 6.1, theorem 4 of [Liu5] we have demonstrated that (under some
additional special assumption \footnote{See theorem 4 of [Liu5] for details.}) 
the local contribution of
 the family invariant to $X\times_{M_n}Y(\Gamma)$ can be identified 
 with the mixed family invariant of $C-{\bf M}(E)E-\sum_{1\leq i\leq p}e_{k_i}$ 
 over $Y(\Gamma)$.  This motivates us to consider the following refined approach
 in section \ref{subsection; blowup}.

\medskip

\subsection{\bf The Repeated Blowing Ups of Sub-Loci in $X$}
\label{subsection; blowup}

\bigskip

 The element $\gamma_n$ is apparently the largest element under $\succ$
 in $\Delta(n)$. Over the open top stratum the type I exceptional classes $e_i$
 are the $-1$ classes $E_i, 1\leq i\leq n$ and the exceptional cone 
\footnote{Generated by $E_1$, $E_2$, $E_3$, $\cdots$, $E_n$.}
${\cal C}_{\gamma_n}$ it
 generates is the smallest. 

 List all the $\Gamma\in \Delta(n)-\{\gamma_n\}$
 and they form a finite graph (each $\Gamma\in \Delta(n)-\{\gamma_n\}$
 being a vertex in the graph) under the
 partial ordering $\succ$. For all $\Gamma\in \Delta(n)-\{\gamma_n\}$, we consider
 the fiber product $X\times_{M_n}Y(\Gamma)$.

 By definition the family moduli space ${\cal M}_{C-{\bf M}(E)E}$ 
over $M_n\times T(M)$ 
of curves dual to $C-{\bf M}(E)E$ collects all the
 curves within the fibers of $M_{n+1}\times T(M)\mapsto M_n\times T(M)$ which
 are dual to $C-{\bf M}(E)E$. When we use the canonical algebraic 
Kuranishi model, $Z(s_{canon})={\cal M}_{C-{\bf M}(E)E}$ for
 $s_{canon}\in \Gamma(X, {\bf H}\otimes \pi_X^{\ast}{\bf W}_{canon})$ (for
the definitions of $s_{canon}$, ${\bf W}_{canon}$, please consult 
section 5.1 of [Liu3] and section 5, proposition 9, 10 of [Liu5]).
 So ${\cal M}_{C-{\bf M}(E)E}$ can be viewed as a sub-scheme
 of $X$ and
 the inclusion $Z(s_{canon})\subset X$ induces the natural projection map to 
$M_n$. 
The schemes 
$Z(s_{canon})\cap (X\times_{M_n}Y(\Gamma))=Z(s_{canon})\times_{M_n} Y(\Gamma)$,
 $\Gamma\in \Delta(n)-\{\gamma_n\}$,
 are sub-schemes of $Z(s_{canon})$ and the ultimate goal is to
 enumerate the residual contribution of ($t_L\in T(M)$)

$$\hskip -1.2in 
{\cal AFSW}_{M_{n+1}\times \{t_L\}\mapsto M_n\times \{t_L\}}(1, C-{\bf M}(E)E)
=c_1({\bf H})^{p_g+rank_{\bf C}{\bf V}_{canon}-rank_{\bf C}{\bf W}_{canon}}
\cap c_{top}({\bf H}\otimes \pi_X^{\ast}{\bf W}_{canon}),$$

or 

$$\hskip -1.2in 
{\cal AFSW}_{M_{n+1}\times T(M)\mapsto M_n\times T(M)}(1, C-{\bf M}(E)E)
=c_1({\bf H})^{p_g+rank_{\bf C}{\bf V}_{canon}-rank_{\bf C}{\bf W}_{canon}}
\cap c_{top}({\bf H}\otimes \pi_X^{\ast}{\bf W}_{canon}),$$

outside $\cup_{\Gamma\in \Delta(n)-\{\gamma_n\}}Z(s_{canon})\times_{M_n}Y(\Gamma)$
 and show that under some additional assumption on $t_L$, the residual contribution
 localizes to lie above the open sub-space $X\times_{M_n}Y_{\gamma_n}$.

 We achieve this by blowing up inductively along the various 
loci $Z(s_{canon})\times_{M_n}Y(\Gamma)$ (or more precisely
 their strict transforms under
 the previous blowing ups), 
starting from the minimal $\Gamma$ under $\succ$ 
and running in the reversed orders of
 $\succ$. We will discuss extensively regarding the ambiguities involved in 
the orders of the blowing ups and how ``doesn't''
 it affect the enumeration process.

  Suppose that $Y(\Gamma_2)\subset 
Y(\Gamma_1)$, $\Gamma_1, \Gamma_2\in \Delta(n)$, and 
then $Z(s_{canon})\times_{M_n} Y(\Gamma_2)\subset Z(s_{canon})\times_{M_n} 
 Y(\Gamma_1)$,
 ${\cal C}_{\Gamma_2}\supset {\cal C}_{\Gamma_1}$. If ${\cal C}_{\Gamma_1}$
 is a proper sub-cone of ${\cal C}_{\Gamma_2}$, 
 then $Y(\Gamma_2)=\overline{S_{\Gamma_2}}$
 can never intersect $S_{\Gamma_1}$ non-trivially. Otherwise at any of 
the intersection
 points the type $I$ exceptional cone is ${\cal C}_{\Gamma_1}$, by the
 definition of $S_{\Gamma_1}$. But this intersection
 point is also in $Y_{\Gamma_2}$, or can
 by degenerated from points in $Y_{\Gamma_2}$, i.e. it is in 
$Y(\Gamma_2)-Y_{\Gamma_2}$. Thus
 ${\cal C}_{\Gamma_2}\subset {\cal C}_{\Gamma_1}$ (by degenerations of cones) and
 is impossible by our assumption $\Gamma_1>\Gamma_2$, or equivalently
 $Y(\Gamma_1)\supset Y(\Gamma_2)$. 

  Therefore in such a situation the graph $\Gamma_2$ can never $\succ\Gamma_1$.
In fact lemma \ref{lemm; suff} and $\overline{S_{\Gamma_1}}\cap S_{\Gamma_2}=
Y(\Gamma_1)\cap S_{\Gamma_2}=S_{\Gamma_2}\not=\emptyset$ implies $\Gamma_1\succ
 \Gamma_2$.

 On the other hand, we have the following lemma regarding the possible relationship
 between two admissible graphs in $\Delta(n)$,

\begin{lemm}\label{lemm; common}
Let $\Gamma_1, \Gamma_2\in \Delta(n)$ be two distinct admissible graphs
 satisfying the special maximality conditions (on page \pageref{specialcondition}).
If $Y(\Gamma_1)\cap Y(\Gamma_2)\not=\emptyset$, there
 are three mutually exclusive possibilities.

(a). $\Gamma_1\succ \Gamma_2$.

(b). $\Gamma_2\succ \Gamma_1$.

(c). Neither $\Gamma_1\succ \Gamma_2$ nor $\Gamma_2\succ \Gamma_1$. But
there exists a ``refined'' 
$\Gamma_3\in \Delta(n)$ such that $\Gamma_i\succ \Gamma_3$ for
 both $i=1, 2$. I.e. $\Gamma_3$ is smaller than $\Gamma_1, 
\Gamma_2$ simultaneously.
\end{lemm}

\medskip

\noindent Proof of lemma \ref{lemm; common}:
 This can be shown 
by contradiction easily. Assuming that neither (a). nor (b). holds, then
 lemma \ref{lemm; suff} implies that both 
$S_{\Gamma_1}\cap Y(\Gamma_2)=S_{\Gamma_2}\cap Y(\Gamma_1)=\emptyset$. Along
 with the fact that $S_{\Gamma_1}\cap S_{\Gamma_2}=\emptyset$ for
 ${\cal C}_{\Gamma_1}\not={\cal C}_{\Gamma_2}$, it 
implies that 

$$Y(\Gamma_1)\cap Y(\Gamma_2)\subset (Y(\Gamma_1)-S_{\Gamma_1})\cap 
(Y(\Gamma_2)-S_{\Gamma_2}).$$

 Let $b\in Y(\Gamma_1)\cap Y(\Gamma_2)$. Because
 $b\not\in S_{\Gamma_1}\cup S_{\Gamma_2}$ but $b\in \overline{S_{\Gamma_1}}\cap
 \overline{S_{\Gamma_2}}$, ${\cal EC}_b(C-{\bf M}(E)E; Q)$
 contains both ${\cal C}_{\Gamma_1}, {\cal C}_{\Gamma_2}$ as proper sub-cones.
 Let $e_{k_i}$, $1\leq i\leq p$ be the primitive type $I$ generators of 
 the simplicial cone 
${\cal EC}_b(C-{\bf M}(E)E; Q)$. Because $e_{k_i}$ are represented by 
 irreducible curves $e_{k_i}\cdot e_{k_j}\geq 0$ for $i\not=j$, then 
by proposition \ref{prop; reverse}, one can construct a $\Gamma_3\in adm(n)$ associated
 with these $e_{k_i}$, $1\leq i\leq p$. By proposition \ref{prop; transversal}
 the co-existence of these type $I$ exceptional
 curves dual to $e_{k_i}, 1\leq i\leq p$, characterizes the admissible stratum 
$Y_{\Gamma_3}\subset Y(\Gamma_3)$ over
 which $e_{k_i}$, $1\leq i\leq p$,
 are represented as smooth type $I$ exceptional curves in the
 fibers of $M_{n+1}\times_{M_n}Y_{\Gamma_3}\mapsto Y_{\Gamma_3}$ and all
 $e_j, j\not\in \{k_1, k_2, \cdots, k_p\}$ are $-1$ classes. By definition 
 of ${\cal EC}_b(C-{\bf M}(E)E, Q)$,
 $e_{k_i}\cdot (C-{\bf M}(E)E)<0$, $1\leq i\leq p$ and so $\Gamma_3\in \Delta(n)$.

 By our construction of $\Gamma_3$, we have $b\in S_{\Gamma_3}$ since 
${\cal C}_{\Gamma_3}={\cal EC}_b(C-{\bf M}(E)E; Q)$. It is apparent that
 $b\in S_{\Gamma_3}\cap Y(\Gamma_i)=S_{\Gamma_3}\cap \overline{S_{\Gamma_i}}$, 
for $i=1, 2$. Thus by lemma \ref{lemm; suff},
$\Gamma_1\succ \Gamma_3$, $\Gamma_2\succ \Gamma_3$ simultaneously. $\Box$

\medskip

\begin{rem}\label{rem; ambiguity}
Because the ambiguity of choices of 
the point $b$, the graph $\Gamma_3$ constructed in the
proof may not be unique.
\end{rem}

\medskip

 To apply the residual intersection theory to $Z(s_{canon})\times_{M_n}Y(\Gamma)$,
$\Gamma\in \Delta(n)-\{\gamma_n\}$ 
inductively, each $Z(s_{canon})\times_{M_n}Y(\Gamma)\subset
 X={\bf P}({\bf V}_{canon})$ determines a blowup center and in the
 following we decide the
 order of blowing ups by upgrading the partial ordering $\succ$ into a linear
 ordering called $\models$.

\bigskip
$\diamondsuit$ Definition of $\models$:\label{models}

\medskip

 Initially define the {\bf current index set} to be $\Delta(n)$.
 Because $(\Delta(n), \succ)$ is a partial ordered finite set, there must be 
 some (maybe non-unique) minimal elements in $\Delta(n)$ which are not larger
than other element in $\Delta(n)$ under $\succ$.

(1). List all the minimal elements in the {\bf current index set} 
$\Delta(n)$ under the partial ordering
 $\succ$. Select one of them (this introduces some ambiguity if the minimal
elements are not unique\footnote{This ambiguity does not affect the
 result of our enumeration. see proposition \ref{prop; identical} on page
 \pageref{prop; identical}.}).

(2). Remove the selected element from the {\bf current index set}
 and list all the minimal
elements from the residual set. Define the new {\bf current index set} to be the
 residual set.  Select one of the minimal elements again.

(3). Go back to step (2)., then repeat the above process and iterate. 

\bigskip

(4). After a finite number of times, one will exhaust the whole $\Delta(n)$ and
 determine a sequence of blowing up centers.

In this way we have determined a linear ordering on $\Delta(n)$, denoted by
 $\models$.

 The discussion right in front of lemma \ref{lemm; common}
 indicates the following: Suppose that we blow up the strict transforms of
 $Z(s_{canon})\times_{M_n}Y(\Gamma)$ following the reversed ordering of
 $\models$ starting from the smallest element in $\Delta(n)-\{\gamma_n\}$.
After blowing up the strict transform of $Z(s_{canon})\times_{M_n}
 Y(\Gamma)$, in our set up we will never 
blow up any sub-locus completely lying inside $Z(s_{canon})\times_{M_n}
 Y(\Gamma)$. In fact, the strict transformation of any such sub-locus 
 will be blown up prior to the blowing up
 of the strict transformation of $Z(s_{canon})\times_{M_n}
 Y(\Gamma)$, due to the fact they are smaller under the partial ordering $\succ$, 
 and therefore the linear ordering $\models$.

 Let us blow up $X$ inductively along the strict transforms of the various
 $Z(s_{canon})\times_{M_n}Y(\Gamma)$, $\Gamma\in \Delta(n)-\{\gamma_n\}$. 
Let $D_{\Gamma}\subset X_{\Gamma}$, 
with $\Gamma\in \Delta(n)$, denote the exceptional Cartier divisor
 blown up from the strict transform of $Z(s_{canon})\times_{M_n}Y(\Gamma)$ 
 and denote the intermediate blown up scheme by $X_{\Gamma}$.
 Define $\tilde{X}$ to be the resulting scheme after blowing up
 all the (strict transforms of $Z(s_{canon})\times_{M_n}Y(\Gamma)$, $\Gamma\in 
\Delta(n)-\{\gamma_n\}$) and the projection map
 $\tilde{X}\mapsto X$ can be factorized into the compositions of
 the various intermediate blowing down map. 

In the following discussion, we may pull back ${\cal O}(D_{\Gamma})$ from
 $X_{\Gamma}$ to
 $\tilde{X}$ from the various birational models (intermediate blowing ups) of $X$. 
 To avoid complicated notations involving the line bundle or divisor 
pull-backs, we {\bf skip the
 pull-back notations} consistently. The reader should be able to judge from the 
 context of the formula and restore the pull-back notations accordingly.

\medskip

 At the end of this subsection, we 
introduce an index set $I_{\Gamma}\subset \Delta(n)$ collecting those
 $\Gamma'$ smaller than $\Gamma\in \Delta(n)$ under $\models$.

\begin{defin}\label{defin; already}
 Let $\Gamma\in \Delta(n)$. The 
linear ordering $\models$ among all the 
$\Gamma\in \Delta(n)$ determines the ordering of the blowing ups to
construct $\tilde{X}$ from $X$.
Define $I_\Gamma$ to be the subset of $\Delta(n)$ satisfying
 $I_{\Gamma}=\{\Gamma'|\Gamma\models \Gamma', \Gamma'\in \Delta(n)\}$.  
\end{defin}

The index set $I_{\Gamma}, \Gamma\in\Delta(n)-\{\gamma_n\}$ 
collects all the admissible graphs $\Gamma'$ in 
$\Delta(n)-\{\gamma_n\}$ whose associated zero loci 
$Z(s_{canon})\times_{M_n}Y(\Gamma')$ 
(or more accurately their strict transformations) are blown up prior to 
$Z(s_{canon})\times_{M_n}Y(\Gamma)$.

We notice that $I_{\gamma_n}=\Delta(n)-\{\gamma_n\}$.
 The collection of index sets 
$I_{\Gamma}$for $\Gamma\in \Delta(n)$
 will be used to define the modified algebraic family
 invariant in the next sub-section. 

\medskip

\subsection{The Definition of Modified Algebraic Family Invariants}
\label{subsection; modinv}

\bigskip

 In this subsection we define a version of modified algebraic family 
Seiberg-Witten invariant associated each $\Gamma\in \Delta(n)$. Recall that
 in section 5.3 on page 448 of [Liu1], a version of modified family Seiberg-Witten
 invariant has been defined in the differentiable category. The 
 modified algebraic family Seiberg-Witten invariant we are going to define
  is its algebraic analogue.

 The first step is to define a class $\tau_{\Gamma}\in K_0(Y(\Gamma)\times T(M))$, 
representable by a locally free sheaf (vector bundle) on the connected space 
$Y(\Gamma)\times T(M)$.

As usual we use $e_1, e_2, e_3, \cdots, e_n$ to denote the $n$ type $I$
exceptional classes over $Y_{\Gamma}$. Let $e_{k_i}, 1\leq i\leq p$,
 be the type $I$ exceptional classes over
 $Y_{\Gamma}$ which
 pair negatively with the class $C-{\bf M}(E)E$.

 As usual let $\Gamma_{e_{k_i}}$ denote the fan-like admissible graph such that
 the type $I$ exceptional class $e_{k_i}$ is effective and smooth/irreducible 
over the locally closed 
$Y_{\Gamma_{e_{k_i}}}$ (consult section \ref{section; strata} for more
details).

\begin{prop}\label{prop; exact}
Let $\tilde{\Xi}_{k_i}\mapsto Y(\Gamma_{e_{k_i}})$ be the relatively minimal 
${\bf P}^1$ fiber bundle associated with the type $I$ class 
$e_{k_i}=E_{k_i}-\sum_{j_{k_i}} E_{j_{k_i}}$.  

Suppose that $e_{k_i}^2<e_{k_i}\cdot (C-{\bf M}(E)E)<0$, then there exists
 an invertible sheaf ${\cal Q}_{k_i}$ over $\tilde{\Xi}_{k_i}$, 
pulled-back from $Y(\Gamma_{e_{k_i}})$,
 an effective relative divisor $\Delta_{k_i}\subset
 \tilde{\Xi}_{k_i}\mapsto Y(\Gamma_{e_{k_i}})$ of relative degree 
$-e_{k_i}\cdot ({\bf M}(E)E+e_{k_i})$
and the following short exact sequence of
 locally free sheaves,

$$\hskip -.5in
0\mapsto {\cal R}^0
\tilde{\pi}_{\ast}\bigl({\cal O}_{\Delta_{k_i}}(-m_{k_i}E_{k_i}-
\sum_{j_{k_i}}m_{j_{k_i}}E_{j_{k_i}})\bigr)\otimes
{\cal Q}_{k_i}
\mapsto {\cal R}^1\tilde{\pi}_{\ast}\bigl({\cal O}_{\tilde{\Xi}_{k_i}}(E_{k_i}-\sum
_{j_{k_i}}E_{j_{k_i}})\bigr)\otimes {\cal Q}_{k_i}$$
$$\hskip -.5in 
\mapsto {\cal R}^1\tilde{\pi}_{\ast}
\bigl({\cal O}_{\tilde{\Xi}_{k_i}}(-m_{k_i}E_{k_i}-
\sum_{j_{k_i}}m_{j_{k_i}}E_{j_{k_i}})\bigr)\otimes {\cal O}_{Y(\Gamma_{e_{k_i}})}(
-\sum_{1\leq l<k_i} m_lE_{l; k_i})
\mapsto 0.
$$\footnote{The symbol $E_{a; b}$, $a<b$, denotes the exceptional divisor in 
 $M_n$ by blowing up the strict transform of the $(a, b)-$th partial diagonal.}

 Suppose that 
$e_{k_i}^2\geq e_{k_i}\cdot (C-{\bf M}(E)E)$, then there 
exists an invertible sheaf ${\cal Q}_{k_i}$ pulled-back from $Y(\Gamma_{e_{k_i}})$,
 an effective relative divisor 
 $\Delta_{k_i}\subset \tilde{\Xi}_{k_i}\mapsto Y(\Gamma_{e_{k_i}})$ of
 relative degree $e_{k_i}\cdot ({\bf M}(E)E+e_{k_i})$ and the following
short exact sequence of locally free sheaves on $Y(\Gamma_{e_{k_i}})$,

$$\hskip -.5in
0\mapsto {\cal R}^0\tilde{\pi}_{\ast}\bigl({\cal O}_{\Delta_{k_i}}(E_{k_i}-\sum
_{j_{k_i}}E_{j_{k_i}})\bigr)\otimes {\cal Q}_{k_i}
\mapsto {\cal R}^1\tilde{\pi}_{\ast}
\bigl({\cal O}_{\tilde{\Xi}_{k_i}}(-m_{k_i}E_{k_i}-
\sum_{j_{k_i}}m_{j_{k_i}}E_{j_{k_i}})\bigr)\otimes {\cal O}_{Y(\Gamma_{e_{k_i}})}(
-\sum_{1\leq l<k_i} m_lE_{l; k_i})$$
$$\mapsto {\cal R}^1\tilde{\pi}_{\ast}
\bigl({\cal O}_{\tilde{\Xi}_{k_i}}(E_{k_i}-\sum
_{j_{k_i}}E_{j_{k_i}})\bigr)\otimes {\cal Q}_{k_i} \mapsto 0.$$
\end{prop}

\noindent Proof: The proof is almost identical to the proof of proposition 15
 in [Liu5] and the reader can consult the cited paper for its derivation.
 Instead of using the ${\bf P}^1$ fibrations $\Xi_{k_i}\mapsto
 Y(\Gamma_{e_{k_i}})$, we use the
  ${\bf P}^1$ fiber bundles $\tilde{\Xi}_{k_i}$. Because the
 details of the argument is almost identical, we omit it here. $\Box$

\medskip

 The following lemma characterizes the significance of the locally free sheaf
 ${\cal R}^1\tilde{\pi}_{\ast}\bigl({\cal O}_{\tilde{\Xi}_{k_i}}(E_{k_i}-
\sum_{j_{k_i}}E_{j_{k_i}})\bigr)$.

\begin{lemm}\label{lemm; normal}
Let ${\cal N}_{Y(\Gamma_{e_{k_i}})}$ denote the normal sheaf of 
$Y(\Gamma_{e_{k_i}})\subset M_n$. Then there exists a canonical isomorphism
 ${\cal R}^1\tilde{\pi}_{\ast}\bigl({\cal O}_{\tilde{\Xi}_{k_i}}(E_{k_i}-
\sum_{j_{k_i}}E_{j_{k_i}})\bigr)\mapsto {\cal N}_{Y(\Gamma_{e_{k_i}})}$.
\end{lemm}

\medskip

\noindent Outline of the Proof: The key idea is to study the canonical
algebraic family Kuranishi model of $e_{k_i}$, 

$$0\mapsto {\cal R}^0\pi_{\ast}
\bigl({\cal O}_{M_{n+1}}(E_{k_i}-\sum_{j_{k_i}}E_{j_{k_i}})\bigr)\mapsto 
{\cal R}^0\pi_{\ast}{\cal O}_{M_{n+1}}(E_{k_i})$$
$$\mapsto {\cal R}^0\pi_{\ast}\bigl(
{\cal O}_{\sum_{j_{k_i}}E_{j_{k_i}}}(E_{k_i})\bigr)\mapsto
 {\cal R}^1\pi_{\ast}
\bigl({\cal O}_{M_{n+1}}(E_{k_i}-\sum_{j_{k_i}}E_{j_{k_i}})\bigr)\mapsto 0.$$
Then the lemma is a direct consequence of 
 lemma 9 in section 6.1 of [Liu5], once we realize that the
 birational projection (see lemma \ref{lemm; sum})
 $\Xi_{k_i}\mapsto \tilde{\Xi}_{k_i}$ induces an isomorphism
 ${\cal R}^1\pi_{\ast}\bigl({\cal O}_{\Xi_{k_i}}(E_{k_i}-\sum_{j_{k_i}}
 E_{j_{k_i}})\bigr)\stackrel{\cong}{\mapsto} {\cal R}^1\tilde{\pi}_{\ast}\bigl(
{\cal O}_{\tilde{\Xi_{k_i}}}(E_{k_i}-\sum_{j_{k_i}}
 E_{j_{k_i}})\bigr)$.  $\Box$

\medskip

If at least one of $e_{k_i}$ satisfies $e_{k_i}^2<e_{k_i}\cdot (C-{\bf M}(E)E)<0$,
 then we define $\tau_{\Gamma}\equiv 0\in K_0(Y(\Gamma)\times T(M))$ following 
the rationale of theorem \label{tau0}
 3. of [Liu5] and case II of the proof of proposition \ref{prop; =mo} starting
 at page \pageref{case2}.
From now on we may assume that $0>e_{k_i}^2\geq e_{k_i}\cdot (C-{\bf M}(E)E)$ for
 all $1\leq i\leq p$.

Recall that it was defined in subsection \ref{subsection; realcase} definition
 \ref{defin; indexset}
 that the index subset $I_{k_l}$ collects the
 indexes in $\{1, 2, \cdots, n\}$
 occurring in $e_{k_l}=E_{k_l}-\sum_{j_{k_l}}E_{j_{k_l}}$, i.e.
 $k_l$ and all its direct descendent indexes in the graph $\Gamma$.

We define $\tau_{\Gamma}\in K_0(Y(\Gamma)\times T(M))$ as the following,

\begin{defin}\label{defin; tau}
Let $e_{k_i}, 1\leq i\leq p$ denote the type $I$ exceptional classes among
 $e_1, e_2, \cdots, e_n$ over
$Y(\Gamma)$ which have negative pairings with $C-{\bf M}(E)E$. Suppose that 
 $e_{k_i}^2\geq e_{k_i}\cdot (C-{\bf M}(E)E)$ for all $1\leq i\leq p$, define
 $\tau_{\Gamma}\equiv [\oplus_{1\leq l\leq p}{\cal R}^1\tilde{\pi}_{\ast}\bigl(
{\cal O}_{\tilde{\Xi}_{k_l}}\otimes {\cal E}_{
C-\sum_{a\in I_{k_l}}m_aE_a-\sum_{p\geq a>l}e_{k_a}}
\bigr)\otimes {\cal O}_{Y(\Gamma)}(-\sum_{1\leq r<k_l}m_rE_{r; k_l})
-\oplus_{1\leq l\leq p}{\cal R}^1\tilde{\pi}_{\ast}\bigl(
{\cal O}_{\tilde{\Xi}_{k_l}}(E_{k_l}-\sum_{j_{k_l}}E_{j_{k_l}})\bigr)
\otimes {\cal Q}_{k_l}\otimes {\cal E}_C]\in K_0(Y(\Gamma)\times T(M))$.
Otherwise\footnote{If some $e_{k_l}$ satisfies
 $e_{k_l}^2\leq e_{k_l}\cdot (C-{\bf M}(E)E)<0$.}, set $\tau_{\Gamma}$ to be zero.
\end{defin}

Please compare the definition of $\tau_{\Gamma}$ with definition 
\ref{defin; replace} of $\tilde{\cal V}_{quot}$ on page \pageref{defin; replace}.

\medskip

\begin{lemm}\label{lemm; free}
 The element $\tau_{\Gamma}$ can be represented by a locally free sheaf of 
rank  $\sum_{1\leq l\leq p}e_{k_l}\cdot (e_{k_l}+{\bf M}(E)E+\sum_{1\leq 
j<l}e_{k_j})$,

$$\hskip -1.3in
\oplus_{1\leq l\leq p}({\cal R}^0\tilde{\pi}_{\ast}\bigl(
{\cal O}_{\Delta_{k_l}}(E_{k_l}-\sum_{j_{k_l}}E_{j_{k_l}})\bigr)\otimes 
{\cal Q}_{k_l}\otimes {\cal E}_C
\oplus_{1\leq l<t\leq p}
{\cal R}^0\tilde{\pi}_{\ast}\bigl({\cal O}_{\tilde{\Xi}_{k_l}\cap
 \tilde{\Xi}_{k_t}}(-\sum_{a\in I_{k_l}}m_aE_a) \otimes {\cal E}_C\bigr))
\otimes {\cal O}_{Y(\Gamma)}(-\sum_{1\leq r<k_l}m_rE_{r; k_l}).$$
\end{lemm}

The symbol $\tilde{\Xi}_{k_l}\cap \tilde{\Xi}_{k_t}$ used here denotes the
 cross section of $\tilde{\Xi}_{k_l}|_{Y(\Gamma)}\mapsto 
Y(\Gamma)$ (or $\tilde{\Xi}_{k_t}|_{Y(\Gamma)}\mapsto Y(\Gamma)$ induced
 by $E_{k_t}$ (or $E_{k_l}$) whenever $e_{k_l}\cdot e_{k_t}=1$. 
It is taken to be the empty set when $e_{k_l}\cdot e_{k_t}=0$.

\noindent Proof: 
The lemma is a direct consequence of a collection of
${\cal O}_{Y(\Gamma)}(-\sum_{1\leq r<k_l}m_rE_{r; k_l})$
 twisted version of short exact sequences for different $l$,

$$
0\mapsto \oplus_{p\geq i>l\geq 1}{\cal R}^0\tilde{\pi}_{\ast}\bigl(
{\cal O}_{\tilde{\Xi}_{k_l}\cap \tilde{\Xi}_{k_i}}\otimes {\cal E}_{C-
\sum_{a\in I_{k_l}}m_aE_a}\bigr)
\mapsto {\cal R}^1\tilde{\pi}_{\ast}\bigl(
{\cal O}_{\tilde{\Xi}_{k_l}}\otimes 
{\cal E}_{C-\sum_{a\in I_{k_l}}m_aE_a-\sum_{p\geq t>l}e_{k_t}}\bigr)$$
$$\mapsto {\cal R}^1\tilde{\pi}_{\ast}\bigl({\cal O}_{\tilde{\Xi}_{k_l}}\otimes
 {\cal E}_{C-\sum_{a\in I_{k_l}}m_aE_a}\bigr)
\mapsto 0$$

for all $l$ ranging in 
$1\leq l\leq p$ and the ${\cal E}_C$-twisted versions of 
the short exact sequences in proposition \ref{prop; exact}. The above sheaf short
 exact sequences are the derived exact sequences of sheaf short exact sequences
 on $\tilde{\Xi}_{k_l}$, $1\leq l\leq p$, of the divisors 
 $\cup_{p\geq t>l}\tilde{\Xi}_{k_l}\cap \tilde{\Xi}_{k_t}\subset 
 \tilde{\Xi}_{k_l}$.
 They truncate to short exact sequences because
 ${\cal R}^0\tilde{\pi}_{\ast}\big({\cal O}_{\tilde{\Xi}_{k_l}}\otimes
 {\cal E}_{C-\sum_{a\in I_{k_l}}m_aE_a}\bigr)=0$, due to
 the negativity of the relative degrees $e_{k_l}\cdot (C-{\bf M}(E)E)=
e_{k_l}\cdot (C-\sum_{a\in I_{k_l}}m_aE_a)$ on the ${\bf P}^1$ fiber bundles. 

 The calculation on its rank follows from 
$deg\bigl(\Delta_{k_i}/Y(\Gamma_{e_{k_i}})\bigr)
=e_{k_i}\cdot ({\bf M}(E)E+e_{k_i})=
e_{k_i}\cdot (\sum_{a\in I_{k_i}}m_aE_a+e_{k_i})$ and
 $deg(\bigl(\sum_{p\geq i>l\geq 1}
\tilde{\Xi}_{k_i}\bigr)\cap \tilde{\Xi}_{k_l}/Y(\Gamma))
=\sum_{p\geq i>l\geq 1}e_{k_i}\cdot e_{k_l}$. 

The explicit representative of $\tau_{\Gamma}$ is locally free because each
 summand is a zero-th derived image sheaf along a finite morphism onto $Y(\Gamma)$.
$\Box$

\medskip

 The inductive 
definitions of the modified algebraic family Seiberg-Witten invariants are
 parallel to the induction procedure in enumerating the local contributions 
of the family invariants. The reader who wants to find out the geometric 
motivation for our definition may consult section \ref{section; proof}
for the parallelism. On the other hand, the current inductive scheme
 is also parallel to the 
 definition of modified family Seiberg-Witten invariants in the differentiable
 category. The reader may consult subsection 5.3 of [Liu1] for more details.

Recall (see section 6.5 of [Liu5]) the following definition of
 the partial ordering $\gg$ among the pairs 
$(\Gamma, \sum_{e_i\cdot (C-{\bf M}(E)E)<0}e_i)$.

\begin{defin}\label{defin; partial}
Let $\Gamma>\Gamma'$ be two $n$-vertex admissible graphs and let $e_i$, $e_i'$, 
$1\leq i\leq n$ denote the type $I$ exceptional classes associated with
 $Y_{\Gamma}$, 
$Y_{\Gamma'}$, respectively. 
The pair $(\Gamma, \sum_{e_i\cdot (C-{\bf M}(E)E)<0}e_i)$ is
said to be greater than $(\Gamma', \sum_{e_i'\cdot (C-{\bf M}(E)E)<0}e_i')$ under
$\gg$, denoted as 

$$(\Gamma, \sum_{e_i\cdot (C-{\bf M}(E)E)<0}e_i)\gg 
(\Gamma', \sum_{e_i'\cdot (C-{\bf M}(E)E)<0}e_i'),$$

 if the following conditions hold.

\noindent (A). For all the indexes $i$, $1\leq i\leq n$ such that the 
type $I$ classes over $Y_{\Gamma}$, 
 $e_i$, satisfy $e_i\cdot (C-{\bf M}(E)E)<0$, then $e_i'=e_i$.

 and

\noindent (B). There $\exists$ at least one 
$e_j'$ over $Y_{\Gamma'}$, $e_j'\cdot (C-{\bf M}(E)E)<0$ but 
the corresponding $e_j$ over $Y_{\Gamma}$ with the same subscript $j$ 
 satisfies $e_j\cdot (C-{\bf M}(E)E)\geq 0$.
\end{defin}

 The conditions (A). and (B). imply that some new 
type $I$ class which pairs negatively with $C-{\bf M}(E)E$ shows up
 above the sub-locus $Y_{\Gamma'}\subset Y_{\Gamma}$ 
while the original negative
 $C-{\bf M}(E)E-$paired $e_i$ persists to be irreducible over $Y_{\Gamma'}$ 
and do {\bf not break up}.

 Refer to fig.5 below for an example\footnote{To simplify 
the notations, the $\sum_{e_i\cdot (C-{\bf M}(E)E)<0}e_i$ part has been
 skipped. this causes no problem when a ${\bf M}(E)E$ is fixed throughout the
 discussion.}. 


\label{fig5}

\begin{figure}
\centerline{\epsfig{file=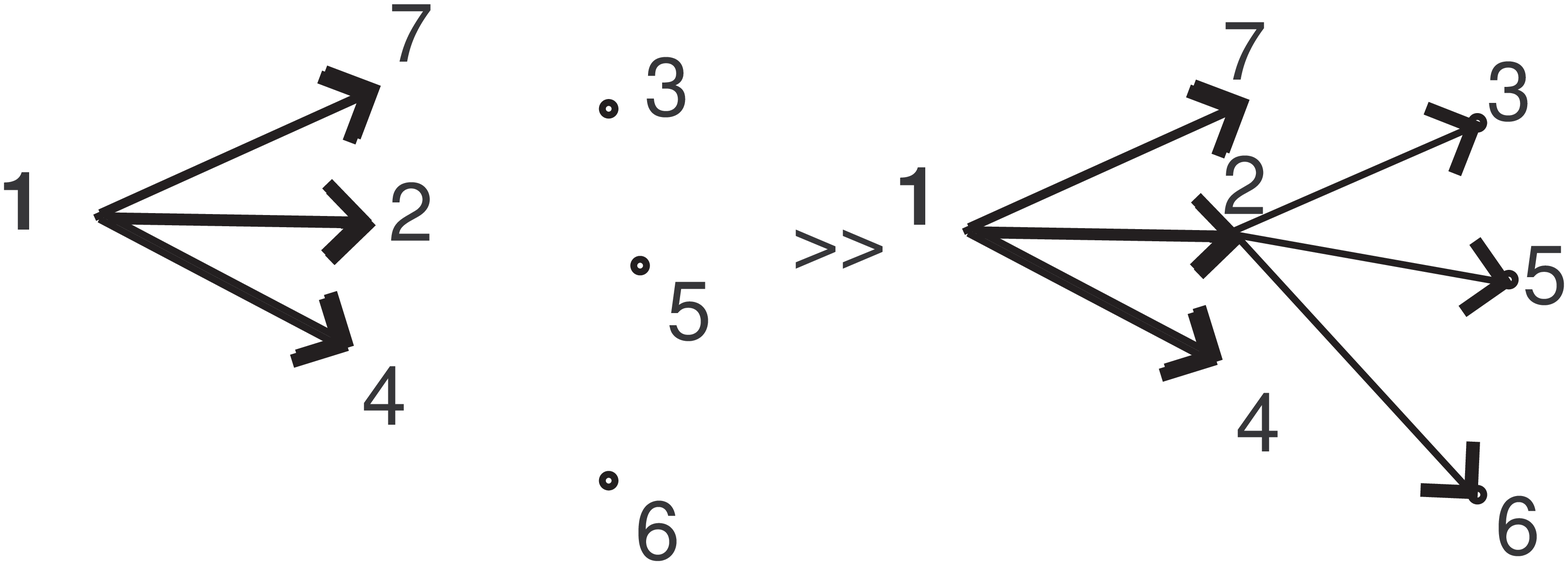,height=4cm}}
\centerline{fig.5}
\centerline{A pair of admissible graphs related by $\gg$, the graph on the
 right hand side is a degeneration} 
\centerline{from the left hand side such that
 $E_2$ breaks into $E_2-E_3-E_5-E_6$ and the union of $E_3$, $E_5$, $E_6$.}
\end{figure}

We abbreviate the above partial order relationship by
 $\Gamma\gg \Gamma'$ if a multiplicity function ${\bf M}(E)E$
has been fixed throughout the discussion.

 We have the following simple observation regarding $\gg$ and $\succ$.

\begin{lemm}\label{lemm; compare}
Let $\Gamma \gg 
\Gamma'$, then $\Gamma\succ \Gamma'$.
\end{lemm}

\noindent Proof:  The cone ${\cal C}_{\Gamma}$ are generated by
 the effective type $I$ classes 
$e_i$ such that $e_i\cdot (C-{\bf M}(E)E)<0$ and some other type $I$ 
$-1$ classes.
  By our assumption
 on $\Gamma\gg \Gamma'$,
 these $e_i$ persist to become $e_i'$ over $Y_{\Gamma'}$ 
with $e_i'\cdot (C-{\bf M}(E)E)<0$. As these $e_i'$ are a subset
 of the generators of 
${\cal C}_{\Gamma'}$, 
this implies that ${\cal C}_{\Gamma'}\supset {\cal C}_{\Gamma}$.
 Therefore $\Gamma\succ \Gamma'$. $\Box$

\medskip

The set $\Delta(n)$ is a finite set. Therefore there must be minimal elements
(may be non-unique) under the partial ordering $\gg$. 

\begin{defin}\label{defin; initial}
Let $\Gamma\in \Delta(n)$ be a minimal element under $\gg$. Define the modified
 algebraic family invariant
 ${\cal AFSW}^{\ast}_{M_{n+1}\times_{M_n}Y(\Gamma)\times T(M)\mapsto
 Y(\Gamma)\times T(M)}(c_{total}(\tau_{\Gamma}), C-{\bf M}(E)E-
\sum_{e_i\cdot (C-{\bf M}(E)E)<0}e_i)$ to be
 ${\cal AFSW} _{M_{n+1}\times_{M_n}Y(\Gamma)\times T(M)\mapsto
 Y(\Gamma)\times T(M)}(c_{total}(\tau_{\Gamma}), C-{\bf M}(E)E-
\sum_{e_i\cdot (C-{\bf M}(E)E)<0}e_i)$, where $c_{total}(\tau_{\Gamma})$ is the
 total Chern class of $\tau_{\Gamma}\in Y(\Gamma)\times T(M)$ defined in 
definition \ref{defin; tau}.
\end{defin}

Let $\Gamma$ be in $\Delta(n)$. Suppose that for all the elements in $\Delta(n)$
 smaller than $(\Gamma, \sum_{e_i\cdot (C-{\bf M}(E)E)<0}e_i)$ under $\gg$, 
the modified
 algebraic family invariants have been defined already. Set the modified 
invariant attached to $Y_{\Gamma}$,
${\cal AFSW}^{\ast}_{M_{n+1}\times_{M_n}Y(\Gamma)\times T(M)\mapsto
 Y(\Gamma)\times T(M)}(c_{total}(\tau_{\Gamma}), C-{\bf M}(E)E-
\sum_{e_i\cdot (C-{\bf M}(E)E)<0}e_i)$ to be,

\begin{defin}\label{defin; inductive}
Let $\tau_{\Gamma}, \tau_{\Gamma'}$ 
be the $K_0$ theory classes by definition \ref{defin; tau} 
associated with $\Gamma, \Gamma'$, respectively.
Define
${\cal AFSW}^{\ast}_{M_{n+1}\times_{M_n}Y(\Gamma)\times T(M)\mapsto
 Y(\Gamma)\times T(M)}(c_{total}(\tau_{\Gamma}), C-{\bf M}(E)E-
\sum_{e_i\cdot (C-{\bf M}(E)E)<0}e_i)$ to be

$${\cal AFSW}_{M_{n+1}\times_{M_n}Y(\Gamma)\times T(M)\mapsto
 Y(\Gamma)\times T(M)}(c_{total}(\tau_{\Gamma}), C-{\bf M}(E)E-
\sum_{e_i\cdot (C-{\bf M}(E)E)<0}e_i)$$
$$-\sum_{\Gamma\gg \Gamma'}
{\cal AFSW}^{\ast}_{M_{n+1}\times_{M_n}Y(\Gamma')\times T(M)\mapsto
 Y(\Gamma')\times T(M)}(c_{total}(\tau_{\Gamma'}), C-{\bf M}(E)E$$
$$-\sum_{e_i'\cdot (C-{\bf M}(E)E)<0}e_i').$$
\end{defin}

\medskip

One may argue easily that the procedure of the inductive definition
 always continues until all the elements in $\Delta(n)$ are exhausted.
Suppose that the process halts before exhausting the elements in $\Delta(n)$.
 Namely, there exists no $\Gamma\in \Delta(n)$
such that the modified algebraic family invariants have been defined for all
 $(\Gamma', \sum_{e_i'\cdot (C-{\bf M}(E)E)<0}e_i')\ll 
(\Gamma, \sum_{e_i\cdot (C-{\bf M}(E)E)<0}e_i)$. If it is the case, then 
for any $\Gamma\in \Delta(n)$ that the modified algebraic family Seiberg-Witten
 invariant is not defined yet, one must be able to find at least one
 $\Gamma'\in \Delta(n)$ such that 
$(\Gamma', \sum_{e_i'\cdot (C-{\bf M}(E)E)<0}e_i')\ll 
(\Gamma, \sum_{e_i\cdot (C-{\bf M}(E)E)<0}e_i)$ and
the modified algebraic 
family invariant is not defined for $(\Gamma', \sum_{e_i'\cdot 
(C-{\bf M}(E)E)<0}e_i')$, either. Then one may trace along the smaller and
 smaller elements under $\gg$ (the procedure involves choices and may not be
 canonical) and it has to stop after a finite number of steps since $\Delta(n)$
 is a finite set. But such a terminal graph $\Gamma$ has to be a minimal
element under $\gg$ and its modified algebraic family invariant has been
defined in definition \ref{defin; initial} already. This generates a
 contradiction and thus the above procedure never
halts unless all the elements in $\Delta(n)$ has been exhausted.

After a finite number of steps and the defining process has to terminate at
 $\gamma_n\in \Delta(n)$. In this case one
 defines ${\cal AFSW}^{\ast}_{M_{n+1}\times T(M)\mapsto M_n\times T(M)}(1, 
C-{\bf M}(E)E)$ by the following recipe,

\begin{defin}\label{defin; generic}
Define ${\cal AFSW}^{\ast}_{M_{n+1}\times T(M)\mapsto M_n\times T(M)}(1, 
C-{\bf M}(E)E)$ to be

$$\hskip -1.2in
{\cal AFSW}_{M_{n+1}\times T(M)\mapsto M_n\times T(M)}(1, C-{\bf M}(E)E)-
\sum_{\Gamma\in \Delta(n)-\{\gamma_n\}}
{\cal AFSW}^{\ast}(c_{total}(\tau_{\Gamma}), C-{\bf M}(E)E-
\sum_{e_i\cdot (C-{\bf M}(E)E)<0}e_i).$$
\end{defin}

 Definition \ref{defin; generic} can be viewed as an extension of
 definition \ref{defin; inductive} once we realize that for $\gamma\in adm(n)$,
 $Y(\gamma_n)=M_n$, $\gamma_n>\Gamma$ for all $\Gamma\in \Delta(n)-\{\gamma_n\}$
 and $(\gamma_n, 0)\gg (\Gamma, \sum_{e_i\cdot (C-{\bf M}(E)E)<0}e_i)$ for all
 $\Gamma\in \Delta(n)-\{\gamma_n\}$.

\medskip

\begin{rem}\label{rem; tL}
Suppose that $t_L\in T(M)$ is a point of the connected component $T(M)$ (determined
 by the first Chern class $C$) of the Picard variety $Pic(M)$.
 There is a corresponding version of ``$t_L$ restricted'' modified algebraic
 family Seiberg-Witten invariants defined by inserting 
$[t_L]\in {\cal A}_0(T(M))$ into each terms in the definition. The resulting
 modified invariant is denoted by 
 
$${\cal AFSW}^{\ast}_{M_{n+1}\times_{M_n}Y(\Gamma)\times T(M)\mapsto 
Y(\Gamma)\times 
T(M)}(c_{total}(\tau_{\Gamma})\cap [t_L], 
C-{\bf M}(E)E-\sum_{e_i\cdot (C-{\bf M}(E)E)<0}e_i),$$

 or equivalently

$${\cal AFSW}^{\ast}_{M_{n+1}\times_{M_n}Y(\Gamma)\times \{t_L\}\mapsto 
Y(\Gamma)\times \{t_L\}}(c_{total}(\tau_{\Gamma}), 
C-{\bf M}(E)E-\sum_{e_i\cdot (C-{\bf M}(E)E)<0}e_i).$$
\end{rem}

 When the point 
$t_L\in T(M)$ 
is determined by an algebraic line bundle $L\mapsto M$ with $c_1(L)=C$, the
 ``$t_L$-restricted'' 
modified algebraic family 
invariants enumerate the curves dual to $C-{\bf M}(E)E$ resolved
 from the linear subsystem of $|L|$.

\medskip

\begin{rem}\label{rem; noCon}
In the earlier paper [Liu5], we had shown that under the
 {\bf Special Condition}, the dominated localized top Chern class contribution
 of $Y(\Gamma)$ is nothing but the mixed algebraic family Seiberg-Witten
 invariant $${\cal AFSW}_{M_{n+1}\times_{M_n}Y(\Gamma)\times T(M)\mapsto
 Y(\Gamma)\times T(M)}(c_{total}(\tau_{\Gamma}), C-{\bf M}(E)E-
\sum_{e_i\cdot (C-{\bf M}(E)E)<0}e_i).$$
 When the {\bf Special Condition} is not met, our definitions of
 the modified invariants indicates that there are correction terms captured by
 the partial ordering $\gg$ besides the dominated term. 
\end{rem}

 A key proposition in proving the main theorem of the paper is the following,

\begin{prop}\label{prop; polynomial}
Let $\Gamma\in \Delta(n)$, then
 the modified algebraic family Seiberg-Witten invariant 
${\cal AFSW}^{\ast}_{M_{n+1}\times_{M_n}Y(\Gamma)\times T(M)\mapsto
 Y(\Gamma)\times T(M)}(c_{total}(\tau_{\Gamma}), C-{\bf M}(E)E-
\sum_{e_i\cdot (C-{\bf M}(E)E)<0}e_i)$ can be expressed as a homogeneous
 universal
 polynomial of $C^2=C\cdot C$, $C\cdot c_1({\bf K}_M)$, $c_1({\bf K}_M)^2$ and
 $c_2(M)$ of degree $n$ multiplied by ${\cal ASW}(C)$. The universal
 polynomial depends on the graph $\Gamma$ and the singular multiplicities
 ${\bf M}(E)E$ but does not depend upon the algebraic surface $M$.
\end{prop}

 When $\Gamma=\gamma_n$, we set $c_{total}(\tau_{\gamma_n})=1$.

\noindent Proof: As all of the modified invariants 
${\cal AFSW}^{\ast}_{M_{n+1}\times T(M)
\mapsto M_n\times T(M)}(1, C-{\bf M}(E)E)$ and 
${\cal AFSW}^{\ast}_{M_{n+1}\times_{M_n}Y(\Gamma)\times T(M)\mapsto
 Y(\Gamma)\times T(M)}(c_{total}(\tau_{\Gamma}), C-{\bf M}(E)E-
\sum_{e_i\cdot (C-{\bf M}(E)E)<0}e_i)$, $\Gamma\in \Delta(n)-\{\gamma_n\}$,
 are defined by
 inductive procedures based on the mixed algebraic invariants, we
 prove that ${\cal AFSW}_{M_{n+1}\times T(M)
\mapsto M_n\times T(M)}(1, C-{\bf M}(E)E)$ and 
 all the ${\cal AFSW}_{M_{n+1}\times_{M_n}Y(\Gamma)\times T(M)\mapsto
 Y(\Gamma)\times T(M)}(c_{total}(\tau_{\Gamma}), C-{\bf M}(E)E-
\sum_{e_i\cdot (C-{\bf M}(E)E)<0}e_i)$, $\Gamma\in \Delta(n)$,
 can be expressed as 
 universal (independent of $M$) homogeneous polynomials of degree $n$ of
 $C\cdot C$, $C\cdot c_1({\bf K}_M)$, $c_1({\bf K}_M)\cdot c_1({\bf K}_M)$,
 and $c_2(M)$, multiplied by the algebraic Seiberg-Witten invariant 
${\cal ASW}(C)$ of $C$.
 
 We present the detailed argument 
 for ${\cal AFSW}_{M_{n+1}\times_{M_n}Y(\Gamma)\times T(M)\mapsto
 Y(\Gamma)\times T(M)}(c_{total}(\tau_{\Gamma}), C-{\bf M}(E)E-
\sum_{e_i\cdot (C-{\bf M}(E)E)<0}e_i)$, $\Gamma\in \Delta(n)-\{\gamma_n\}$ and the
 proof for the case of ${\cal AFSW}_{M_{n+1}\times T(M)
\mapsto M_n\times T(M)}(1, C-{\bf M}(E)E)$ is essentially parallel.

 Recall\footnote{Consult lemma 3.1 on page 401 of [Liu1].}
 that the fiber bundle projection map 
$f_n:M_{n+1}\mapsto M_n$ can be constructed from 
 $M_n\times M\mapsto M_n$ through $n$ consecutive blowing ups along 
 (codimension two) cross sections of the intermediate fiber bundles.
 This implies that its pull-back to $Y(\Gamma)$,
$M_{n+1}\times_{M_n}Y(\Gamma)\mapsto Y(\Gamma)$, can be constructed from 
 the Cartesian projection 
$Y(\Gamma)\times M\mapsto Y(\Gamma)$ through $n$ consecutive blowing ups along
 cross sections of the intermediate blown up spaces.
 Schematically this implies that we may apply the
 family blowup formula of the algebraic family Seiberg-Witten
 invariants [Liu3] to relate
  ${\cal AFSW}_{M_{n+1}\times_{M_n}Y(\Gamma)\times T(M)\mapsto
 Y(\Gamma)\times T(M)}(c_{total}(\tau_{\Gamma}), C-{\bf M}(E)E-
\sum_{e_i\cdot (C-{\bf M}(E)E)<0}e_i)$ and the mixed algebraic family Seiberg-Witten 
invariant ${\cal AFSW}_{Y(\Gamma)\times T(M)\times M\mapsto 
Y(\Gamma)\times T(M)}(c_{total}(\tau_{\Gamma})\cap c_{tatal}(
{\bf U}_{{\bf M}(E)}), C)$. The bundle ${\bf U}_{{\bf M}(E)}$ appearing in
 the identity is the 
relative obstruction bundle 

$${\bf U}_{{\bf M}(E)E}=
\oplus_{1\leq l\leq n} {\cal E}_{C-\sum_{1\leq a\leq l-1}m_aE_a}\otimes
{\bf S}^{m_l-1}({\bf C}\oplus (f_{n-1; l}|_{Y(\Gamma)}^{\ast}
{\bf T}^{\ast}_{M_l/f^{\ast}_{l-1}M_{l-1}}))
$$

, gotten from applying the algebraic family blowup formula $n$ times from
 $M_n\times M\mapsto M_n$ to $M_{n+1}\mapsto M_n$. 
The map 
$f_{n-1; l}|_{Y(\Gamma)}:Y(\Gamma)
\mapsto M_l$ is the composition 
$Y(\Gamma)\subset M_n\stackrel{f_{n-1; l}}{\longrightarrow}M_l$.  

We have the following simple factorization lemma regarding the family invariants,

\begin{lemm}\label{lemm; product}
Let $M\times B\mapsto B$ be a product algebraic fiber bundle over
 a complete and smooth
 base $B$ and let $\underline{C}$
 be a $(1, 1)$ class on the algebraic surface $M$, then
${\cal AFSW}_{B\times M\mapsto B}(\eta, \underline{C})=0$ for
 $\eta\not\in {\cal A}_0(B)$ and is equal to ${\cal ASW}(\underline{C})
\cdot (\int_B\eta)$ for 
$\eta\in {\cal A}_0(B)$.
\end{lemm}

\noindent Proof of lemma \ref{lemm; product}: For simplicity
 we assume that ${\cal E}_{\underline{C}}$ has a vanishing second
 derived image sheaf\footnote{When we apply this lemma to the
 concrete situation below, this additional assumption is satisfied.}
 over $T(M)\times B$. 
 Consider the algebraic family moduli space of $\underline{C}$, 
${\cal M}_{\underline{C}}$, over $B$. Because the fiber bundle of 
 algebraic surfaces is a trivial product, the space ${\cal M}_{\underline{C}}$ 
is also a trivial product over $B$ and the algebraic family Kuranishi models
of $\underline{C}$ are pulled back from $T(M)$ to $T(M)\times B$.
Let $({\bf V}, {\bf W}, \Phi_{{\bf V}{\bf W}})$ be one algebraic family 
Kuranishi model of $\underline{C}$, where ${\bf V}, {\bf W}$ are vector bundles
 over $T(M)\times B$ pulled-back from $T(M)$.

 Then for $\eta\in {\cal A}_k(B)$, $k\leq dim_{\bf C}B$, 
the mixed family invariant 
${\cal AFSW}_{B\times M\mapsto B}(\eta, \underline{C})$ can be expressed as
 the push-forward of

$$\hskip -.7in
\int_{{\bf P}({\bf V})}c_1({\bf H})^{rank_{\bf C}({\bf V}-{\bf W})-1+q(M)+k}
\cap c_{top}({\bf W}\otimes {\bf H})\cap\eta\cap [{\bf P}({\bf V})]\in
 {\cal A}_0(pt)\cong {\bf Z},$$
 
 into ${\cal A}_0(pt)\cong {\bf Z}$.

Because ${\bf P}({\bf V})$ is also a trivial product over $B$, 
$c_1({\bf H})^{rank_{\bf C}({\bf V}-{\bf W})+q(M)-1+k}
\cap c_{top}({\bf W}\otimes {\bf H})\cap [{\bf P}({\bf V})]=0$ for all
 $k>0$. On the other hand, when $\eta\in {\cal A}_0(B)$ the mixed invariant
can be expressed as (for some $b\in B$)

$$\int_{{\bf P}({\bf V}|_{T(M)\times \{b\}})}
c_1({\bf H})^{rank_{\bf C}({\bf V}-{\bf W})-1+q(M)}\cap
 c_{top}({\bf W}|_{T(M)\times \{b\}}\otimes {\bf H})\cap 
{\bf P}({\bf V}|_{T(M)\times \{b\}})\cdot \int_B\eta=
{\cal ASW}(\underline{C})\cdot
 \int_B\eta.$$ 

The case when the second derived image sheaf of ${\cal E}_{\underline{C}}$
 is not vanishing can be discussed similarly and we omit the details here. $\Box$

 By applying lemma \ref{lemm; product} to our context, the above mixed invariant
is equal to ${\cal ASW}(C)\cdot \int_{M_n}\{c_{total}(\tau_{\Gamma})\cap c_{total}(
{\bf U}_{{\bf M}(E)})\cap [Y(\Gamma)]\}_0$.

 Our goal is to show that the intersection number 
$\int_{M_n}\{c_{total}(\tau_{\Gamma})\cap c_{total}({\bf U}_{{\bf M}(E)})\cap [
Y(\Gamma)]\}_0$ is a universal homogeneous polynomial of degree $n$ in terms of
 the Chern numbers $C^2[M]$, $C\cdot c_1({\bf K}_M)[M]$, $c_1({\bf K}_M)^2[M]$ and
 $c_2(M)[M]$.  If for at least one type $I$ class $e_{k_i}$, $1\leq i\leq p$,
  the inequality 
$e_{k_i}^2<e_{k_i}\cdot (C-{\bf M}(E)E)$ holds, 
then $\tau_{\Gamma}\equiv 0$ and the
 mixed algebraic family invariant over $Y(\Gamma)$ has been defined to be zero. 

 So we may assume that $e_{k_i}^2\geq e_{k_i}\cdot (C-{\bf M}(E)E)$ for 
 all $1\leq i\leq p$. Then by lemma \ref{lemm; free} the explicit expression 
of the class $\tau_{\Gamma}$ enables to conclude that,

\begin{lemm}\label{lemm; total}
 Let $\pi_i:M_n\mapsto M$ be the composite projection map 
$M_n\mapsto M^n\mapsto
 M$ to the $i$-th copy of $M$, $1\leq i\leq n$ and 
let\footnote{They were
denoted as $E_{a}(b)\in H^2(M_n, {\bf Z})$ at  page 402, proposition 3.1 of [Liu1], in the 
topological category.} $E_{a; b}$ (for $a<b$) 
 denote the exceptional divisor of $M_n\mapsto M^n$ associated with the 
$(a, b)$-th partial diagonal of $M^n$. 

 The image of the 
total Chern class $c_{total}(\tau_{\Gamma})$ under 
${\cal A}_{\cdot}(Y(\Gamma))\mapsto {\cal A}_{\cdot}(M_n)$ 
can be expressed as a universal degree $rank_{\bf C}\tau_{\Gamma}$
polynomial of the cycle classes 
$\pi_i^{\ast}C\in {\cal A}_{2n-2}(M_n)$, $E_{a; b}\in 
{\cal A}_{2n-2}(M_n)$, for $a<b$. 
\end{lemm}

\medskip

\noindent Proof: By the locally free representative of $\tau_{\Gamma}$ in lemma 
\ref{lemm; free}, one may write

 $$\hskip -1in
c_{total}(\tau_{\Gamma})=c_{total}(\oplus_{1\leq l\leq p}({\cal R}^0
\tilde{\pi}_{\ast}\bigl({\cal O}_{\Delta_{k_l}}(E_{k_l}-
\sum_{j_{k_l}}E_{j_{k_l}})\bigr)\otimes {\cal Q}_{k_l}\otimes {\cal E}_C)
\oplus_{p\geq t>l\geq 1}
{\cal R}^0\tilde{\pi}_{\ast}\bigl({\cal O}_{\tilde{\Xi}_{k_l}\cap
 \tilde{\Xi}_{k_t}}(-\sum_{a\in I_{k_l}}m_aE_a)\otimes {\cal E}_C\bigr))).$$

 It suffices to show that the total Chern class of each of the locally free 
sheaves ${\cal R}^0\tilde{\pi}_{\ast}\bigl({\cal O}_{\Delta_{k_l}}(E_{k_l}-
\sum_{j_{k_l}}E_{j_{k_l}})\bigr)\otimes {\cal Q}_{k_l}\otimes {\cal E}_C$ and
 ${\cal R}^0\tilde{\pi}_{\ast}\bigl({\cal O}_{\tilde{\Xi}_{k_l}\cap
 \tilde{\Xi}_{k_t}}(-\sum_{a\in I_{k_l}}m_aE_a)
\otimes {\cal E}_C\bigr)$ is a polynomial 
in terms of all the flat pull-back $\pi_i^{\ast}C$ and all the $E_{i;j}$, $i<j$.

 Firstly recall how the invertible sheaf 
${\cal E}_C\mapsto T(M)\times M_{n+1}$, $n\in {\bf N}$, has been constructed.

Choose a point $t_{L_0}\in T(M)$ which corresponds to a invertible sheaf
 ${\cal L}_0\mapsto M$ with $c_1({\cal L}_0)=C$. 
Consider the universal invertible sheaf 
 ${\cal L}_{univ}\mapsto T(M)\times M$, then ${\cal L}_{univ}\otimes
 \pi_{T(M)\times M\mapsto M}^{\ast}{\cal L}_0$ defines the invertible 
sheaf ${\cal E}_C$ over $T(M)\times M$.  
To pull it back to $T(M)\times M_{n+1}$,
 we consider the projection $T(M)\times M_{n+1}\mapsto T(M)\times M^{n+1}$
 to the trivial product. By
 composing it with $T(M)\times M^{n+1}\mapsto T(M)^{n+1}\times M^{n+1}\cong 
(T(M)\times M)^{n+1}$ induced by 
$T(M)\ni \{t\}\mapsto \{t\}\times \{t_{L_0}\}\times 
\cdots \{t_{L_0}\}\in T(M)^{n+1}$, the pulled-back invertible sheaf is what
 we denote as ${\cal E}_C$ throughout this paper. It is easy to see that
 the construction is independent to the choices of $t_{L_0}\in T(M)$.

 Fix the ${\bf P}^1$ fiber bundle 
$\tilde{\Xi}_{k_l}\mapsto Y(\Gamma_{e_{k_l}})$
 for some $1\leq l\leq p$, the relative divisors 
 $\Delta_{k_l}\mapsto Y(\Gamma_{e_{k_l}})$ and 
$\tilde{\Xi}_{k_l}\cap \tilde{\Xi}_{k_t}=E_{k_t}|_{\tilde{\Xi}_{k_l}}$ 
are (multiples of) cross sections of
 the given ${\bf P}^1$ fiber bundle $\tilde{\Xi}_{k_l}\mapsto 
Y(\Gamma_{e_{k_l}})$. 
For every direct descendent index $j_{k_l}$ of $k_l$ in the admissible
 graph $\Gamma$, 
the exceptional divisor $E_{j_{k_l}}$ determines a cross-section
 of  $\tilde{\Xi}_{k_l}\mapsto Y(\Gamma_{e_{k_l}})$ and the restriction of 
the invertible 
sheaf ${\cal O}_{\tilde{\Xi}_{k_l}}(E_j)$, $j\not=j_{k_l}$,
 to this cross section determined by $E_{j_{k_l}}$ is isomorphic to the pull-back of
 ${\cal O}_{Y(\Gamma_{e_{k_l}})}(E_{min(j_{k_l}. j); max(j_{k_l}. j)})$ from the base.

 It is easy to see by a simple induction argument that the derived image sheaf 
${\cal R}^0\tilde{\pi}_{\ast}\bigl({\cal O}_{mE_{j_{k_l}}|_{\tilde{\Xi}_{k_l}}}\bigr)
\cong\oplus_{0\leq i\leq m-1}({\cal O}_{Y(\Gamma_{e_{k_l}})}(E_{k_l; j_{k_l}}))^{\otimes -i}$ 
for any direct descendent index $j_{k_l}$ of $k_l$. On the other hand, 
 $c_1({\cal E}_C|_{E_{j_{k_l}}\cap\tilde{\Xi}_{k_l}})=c_1(\pi_{k_l}^{\ast}
{\cal L}_0|_{Y(\Gamma_{e_{k_l}})})=\pi_{k_l}^{\ast}C|_{Y(\Gamma_{e_{k_l}})}$. 
By combining these ingredients, the total Chern class of $\tau_{\Gamma}$ can be 
determined and must be an $M$-independent polynomial in terms of the various 
 $\pi_i^{\ast}C, E_{i; j}$. etc.
$\Box$

\begin{rem}\label{rem; difference}
  If one pulls back the invertible sheaf from ${\cal L}_{univ}\otimes {\cal L}_0
\mapsto T(M)\times M$ by $T(M)\times M_{n+1}\mapsto (T(M)\times M)^{n+1}$ 
 which instead factors through the diagonal embedding $T(M)\times M^{n+1}
\stackrel{\Delta_{T(M)\times id_{M^{n+1}}}}{\longrightarrow} 
T(M)^{n+1}\times M^{n+1}$, then this invertible sheaf differs from our
 ${\cal E}_C$ by an invertible sheaf pulled-back from the base $T(M)\times M_n$.
 If one adopts this alternative invertible sheaf and calculates its first Chern
 class, it will involves not only $c_1({\cal L}_0)$ but also 
$c_1({\cal L}_{univ})$. On the other hand, the algebraic family Kuranishi model
 it determines can be gotten from $({\cal V}_{canon}, {\cal W}_{canon}, 
\Phi_{{\cal V}_{canon}{\cal W}_{canon}})$ by twisting the invertible sheaf 
 pulled-back from $T(M)\times M_n$. It is easy to see that the
 final answer of 
${\cal AFSW}_{T(M)\times M_{n+1}\mapsto T(M)
\times M_n}^{\ast}(1, C-{\bf M}(E)E)$ is independent of the twisting on
 the algebraic family Kuranishi models. Our choice of ${\cal E}_C$ 
 has the benefit of separating the contribution of $c_1({\cal L}_{univ})$ 
 to the family invariant through the factor ${\cal ASW}(C)$ (see remark \ref{rem; poincare}).
\end{rem}

\medskip

 Let $\Gamma\in \Delta(n)$. Then the type $I$ exceptional classes $e_{k_i}$, 
$1\leq i\leq p$,
 are effective over $Y(\Gamma)$ and the type $I$ exceptional classes $e_j$, 
$j\not\in \{k_1, 
k_2, k_3, \cdots, k_p\}$ are $-1$ classes. 
Recall from proposition \ref{prop; transversal} of section \ref{section; strata} 
that the smooth locus $Y(\Gamma)$ is the transversal intersection of
 $Y(\Gamma_{e_{k_i}})$, $1\leq i\leq p$, where $\Gamma_{e_{k_i}}$ is the fan-like
 admissible graph associated with $e_{k_i}$ (see section \ref{section; strata}).
 Then $Y(\Gamma)=\cap_{1\leq i\leq p}Y(\Gamma_{e_{k_i}})$ and in 
 ${\cal A}_{dim_{\bf C}Y(\Gamma)}(M_n)$ we have the equality of cycle classes
$[Y(\Gamma)]=\cap_{i\leq p} [Y(\Gamma_{e_{k_i}})]$.

 Recall $codim_{\bf C}\Gamma$ is the number of one-edges in $\Gamma$.
 Each $[Y(\Gamma_{e_{k_i}})]$ is an algebraic cycle class of dimension
 $dim_{\bf C}M_n+{e_{k_i}\cdot e_{k_i}-c_1({\bf K}_{M_{n+1}/M_n})\cdot e_{k_i}\over 2}=
 2n-codim_{\bf C}\Gamma_{e_{k_i}}$.  To calculate the cycle class explicitly, there are
 essentially two equivalent methods. Either we consider the canonical algebraic family
 Kuranishi model of $e_{k_i}=E_{k_i}-\sum_{j_{k_i}}E_{j_{k_i}}$ as was done in 
section 6.2, lemma 9 of [Liu5], and $Y(\Gamma_{e_{k_i}})$ is the regular zero locus of the
 canonical obstruction bundle, so $[Y(\Gamma_{e_{k_i}})]\in {\cal A}_{2n-
codim_{\bf C}\Gamma_{e_{k_i}}}(M_n)$ represents the top Chern class of the 
obstruction bundle and can be determined explicitly. 
 Or one may apply the algebraic family blow up formula to the class 
 $e_{k_i}=E_{k_i}-\sum_{j_{k_i}}E_{j_{k_i}}$, one may find the top Chern class of the
 obstruction bundle inductively. By either means the answer of $[Y(\Gamma_{e_{k_i}})]$ 
 is expressible as

$$\bigcap_{1\leq s\leq codim_{\bf C}\Gamma_{e_{k_i}}} 
(E_{k_i; j^s_{k_i}}-\sum_{r\leq s-1}
E_{j^r_{k_i}; j^s_{k_i}})\in {\cal A}_{2n-codim_{\bf C}
\Gamma_{e_{k_i}}}(M_n),$$

 where $j_{k_i}^1<j_{k_i}^2<j_{k_i}^3<\cdots<j_{k_i}^{codim_{\bf C}\Gamma_{e_{k_i}}}$ are
 the direct descendent indexes of $k_i$.

  By combining these calculations together, we find that the 
 mixed family invariant ${\cal AFSW}_{M_{n+1}\times_{M_n}Y(\Gamma)\times T(M)
\mapsto Y(\Gamma)\times T(M)}(c_{total}(\tau_{\Gamma}), C-{\bf M}(E)E-\sum_{e_{k_i}\cdot
(C-{\bf M}(E)E)<0}e_{k_i})$ can be expressed as ${\cal ASW}(C)$ times
the $\int_{M_n}$ of a polynomial expression of the various 
$E_{i;j}$, $1\leq i<j\leq n$, $c_{total}({\bf T}M_n)$, 
$c_{total}(f_{n-1, k}^{\ast}{\bf T}M_k)$, and the various $\pi_i^{\ast}C$. 

 To determine the push-forward of the zero cycle class 
under ${\cal A}_0(M_n)\mapsto {\cal A}_0(pt)$, 
observe that the morphism factors through ${\cal A}_0(M_n)\mapsto {\cal A}_0(M^n)
\mapsto {\cal A}_0(pt)\cong {\bf Z}$. We notice the following facts: 

(i). the projection $M_n\mapsto M^n$ 
 can be factorized as ${n(n-1)\over 2}$ consecutive codimension-two blowing down maps,
 
and  

(ii). the well known blowup formula of Chern classes(see page 298,
 section 15.4 of [F]), 

 and the fact that

(iii). the exceptional divisor of a codimension-two 
blowing up along a smooth center has
 the structure of a ${\bf P}^1$ fiber bundle, the projectification of the
 normal bundle of the blowing up center,

and 

(iv). the restriction of the exceptional divisor to itself is equal to
 the first Chern class of the tautological line bundle, and its various 
self-intersections can be expressed by the Chern classes of the normal bundle 
(see page 47-51 of [F] or page 270 of [BT] for a corresponding
 statement in the cohomology ring).

 By combining (i)-(iv)., we reduce the intersection numbers of 
$E_{i; j}$, $c_{total}({\bf T}M_n)$, $c_{total}(f_{n-1, k}^{\ast}{\bf T}M_k)$, 
$\pi_i^{\ast}C$ in ${\cal A}_0(M_n)$ to the intersection numbers of
 $\pi_i^{\ast}c_{total}({\bf T}M)$ and $\pi_i^{\ast}C$, $1\leq i\leq n$
 in ${\cal A}_0(M^n)$.  As the only non-vanishing pairings among these
 classes can be expressed as polynomials in terms of 
 $\pi_i^{\ast}C^2$, $\pi_i^{\ast}C\cap c_1({\bf T}M)$, 
$\pi_i^{\ast}c_1({\bf T}M)^2$ and $\pi_i^{\ast}c_2({\bf T}M)$, $1\leq i\leq n$,
 the integral valued intersection number is a degree $n$ homogeneous polynomial of 
 $C^2\cap [M]$, $C\cap c_1(M)\cap M$, $c_1(M)^2\cap [M]$, $c_2(M)\cap [M]$. 
 $\Box$

\medskip

\begin{rem}\label{rem; poincare}
 When the irregularity $q=0$, ${\cal ASW}(C)=1$ because it is the top intersection
 pairing of $c_1({\bf H})$ on a projective space ${\bf P}({\bf V})$. 
 When $q>0$, the 
${\cal ASW}(C)=\sum_{a+2b=q, a, b\in {\bf N}\cup \{0\}}
(-1)^a\int_{T(M)}{ch_1^a\over a!}\cap {ch_2^b\over b!}$, 
$ch_1=\pi_{T(M)\ast}({c_1({\cal L}_{univ})^2\over 2}\cap
 (2C+c_1({\bf T}M)))$, $ch_2=\pi_{T(M)\ast}{c_1({\cal L}_{univ})^4\over 4!}$, 
 depends on the top intersection pairing on $T(M)$ and 
was calculated in [LL1], [LL2] in the
topological category. Over here $\pi_{T(M)}:M\times T(M)\mapsto T(M)$ denote the 
Cartesian projection to $T(M)$.
\end{rem}

\medskip

\begin{rem}\label{rem; insert}
 In the above discussion, the mixed family invariant enumerates all curves
 within the family $M_{n+1}\times_{M_n}Y(\Gamma)\mapsto Y(\Gamma)$ 
 and dual to $C-{\bf M}(E)E$. We do not require the image curve in $M$ to
 lie within a particular complete linear system associated with
 a holomorphic line bundle over $M$.
 In case we restrict the holomorphic structure to a $t_L\in T(M)$, one may
 insert the zero cycle class $\{t_L\}$ into the family invariant. This has the
 effect of reducing the torus $T(M)$ to a single point and the modified 
 mixed invariant associated to $Y(\Gamma)$ is of the form
 $${\cal AFSW}_{M_{n+1}
\times_{M_n}Y(\Gamma)\times \{t_L\}\mapsto Y(\Gamma)\times \{t_L\}}^{\ast}
(c_{total}(\tau_{\Gamma}), C-{\bf M}(E)E-\sum_{e_i\cdot (C-{\bf M}(E)E)<0}
e_i)$$ and the generic modified family invariant is of the form
 $${\cal AFSW}_{M_{n+1}\times \{t_L\}\mapsto M_n\times \{t_L\}}^{\ast}
(1, C-{\bf M}(E)E).$$ 

An analogue of proposition 
\ref{prop; polynomial} holds while we replace ${\cal ASW}(C)$ by 
${\cal ASW}([t_L], C)=1$.
\end{rem}

\medskip

\subsection{\bf The Combinatorics Involved in the Enumerations}
\label{subsection; comb}

\bigskip

In this subsection, we address the combinatorial issues regarding the
 linear ordering $\models$ and the partial ordering $\gg, \sqsupset$,
 involved in the
 blowing up construction and the inclusion relation on the various 
 restricted family moduli spaces ${\cal M}_{C-{\bf M}(E)E-\sum_{e_i\cdot
(C-{\bf M}(E)E)<0}e_i}\times_{M_n}Y(\Gamma)$, $\Gamma\in \Delta(n)$.
 As will be demonstrated later, it has significant implications on the
 enumeration problem.

 Let us start by noticing that,

\begin{lemm}\label{lemm; trivial}
The localized top Chern class contribution along $D_{\Gamma}$,

$$\sum_{1\leq i\leq e}(-1)^{i-1}
c_{e-i}(E\otimes_{\Gamma'\in I_{\Gamma}}
{\cal O}(-D_{\Gamma'})|_{D_{\Gamma}})
D_{\Gamma}^{i-1}[D_{\Gamma}]$$

can be identified with

$$\sum_{1\leq i\leq e}(-1)^{i-1}
c_{e-i}(E\otimes_{\Gamma'\in I_{\Gamma}; Y(\Gamma)\cap Y(\Gamma')\not=
\emptyset}{\cal O}(-D_{\Gamma'})|_{D_{\Gamma}})
D_{\Gamma}^{i-1}[D_{\Gamma}].$$

I.e. in evaluating the localized contribution of the family invariant along
 $D_{\Gamma}$, one may remove those ${\cal O}(D_{\Gamma'})$ with
 $Y(\Gamma')\cap Y(\Gamma)=\emptyset$.
\end{lemm}

\medskip

\noindent Proof: The defining section of ${\cal O}(D_{\Gamma'})$ vanishes
 exactly on $D_{\Gamma'}$. If $Y(\Gamma')\cap Y(\Gamma)=\emptyset$,
 $D_{\Gamma'}$ is totally disjoint from $D_{\Gamma}$ in the space
 $\tilde{X}$. Thus, 
the line bundle ${\cal O}(D_{\Gamma'})|_{D_{\Gamma}}$ is isomorphic to
the trivial bundle ${\bf C}|_{D_{\Gamma}}$. Therefore they can be 
removed from the expression involving Chern classes of 
${\cal O}(D_{\Gamma'})|_{D_{\Gamma}}$. $\Box$

\medskip

 Because of lemma \ref{lemm; trivial}, in identifying the localized
 top Chern class contribution from $D_{\Gamma}$ we may
 discard all the $Y(\Gamma')$, $\Gamma'\in I_{\Gamma}$, which do not intersect
 $Y(\Gamma)$ at all. Consider all the $Y(\Gamma')$, $\Gamma'\in I_{\Gamma}$,
 and digest briefly the geometric structure of 
${\cal M}_{C-{\bf M}(E)E}\times_{M_n}Y(\Gamma)$ relative to the various
${\cal M}_{C-{\bf M}(E)E}\times_{M_n}Y(\Gamma')$ which touch it non-trivially.

\medskip

\begin{prop}\label{prop; union2}
 Let $e_i$, and $e_i'$, $1\leq i\leq n$ be the type 
$I$ exceptional classes over $Y_{\Gamma}$ and $Y_{\Gamma'}$, respectively.
The restriction of the family moduli space ${\cal M}_{C-{\bf M}(E)E}$ to
 $Y(\Gamma)$,  ${\cal M}_{C-{\bf M}(E)E}\times_{M_n}Y(\Gamma)
=Z(s_{canon})\times_{M_n}Y(\Gamma)$, can be identified as the scheme
 theoretical union of the images
 of the primary component
 ${\cal M}_{C-{\bf M}(E)E-\sum_{e_i; e_i\cdot (C-{\bf M}(E)E)<0}e_i}\times_{M_n}
Y(\Gamma)$ and of the union of secondary components,

$$\cup_{\Gamma'\in I_{\Gamma}}{\cal M}_{C-{\bf M}(E)E-\sum_{e_j; e_j'\cdot 
(C-{\bf M}(E)E)<0}e_j'}\times_{M_n}(Y(\Gamma)\cap Y(\Gamma'))$$

 under the natural inclusions ${\cal M}_{C-{\bf M}(E)E-\sum_{e_i; e_i\cdot 
(C-{\bf M}(E)E)<0}e_i}\subset {\cal M}_{C-{\bf M}(E)E}$ and 
${\cal M}_{C-{\bf M}(E)E-\sum_{e_i'; e_i'\cdot 
(C-{\bf M}(E)E)<0}e_i'}\subset {\cal M}_{C-{\bf M}(E)E}$, respectively.
\end{prop}

\medskip

\noindent Proof of proposition \ref{prop; union2}: Firstly we identify
 them on the set theoretical level.
Let $z\in {\cal M}_{C-{\bf M}(E)E}\times_{M_n}Y(\Gamma)$. 
The point $z$ represents an 
 algebraic curve dual to $C-{\bf M}(E)E$ above a point in $Y(\Gamma)$.
 Suppose that
 $z$ is in the subspace ${\cal M}_{C-{\bf M}(E)E}\times_{M_n}S_{\Gamma}$, then
 $z$ is above a point $b\in S_{\Gamma}$ and 
${\cal EC}_b(C-{\bf M}(E)E; Q)={\cal C}_{\Gamma}$. Then 
$e_i$ with
 $e_i\cdot (C-{\bf M}(E)E)<0$ are exactly the generators of 
${\cal C}_{\Gamma}$. This implies that the effective curve dual to $C-{\bf M}(E)E$
represented by $z$ must contain irreducible components dual to
 each of the $e_i\in {\cal C}_{\Gamma}$. Thus, $z$ is in the image of
 ${\cal M}_{C-{\bf M}(E)E-\sum_{e_i; e_i\cdot 
(C-{\bf M}(E)E)<0}e_i}\subset {\cal M}_{C-{\bf M}(E)E}$.

 If $z\in {\cal M}_{C-{\bf M}(E)E}\times_{M_n}(Y(\Gamma)-S_{\Gamma})$, then
 $z$ is above a point \footnote{The $subset$ is a consequence of
 lemma \ref{lemm; union3} and lemma \ref{lemm; suff}.}
 $b\in Y(\Gamma)-S_{\Gamma}\subset 
\cup_{\Gamma'\in \Delta(n); \Gamma'\prec \Gamma}S_{\Gamma'}$.
 In particular, $b\in Y_{\Gamma'}$ for some $\Gamma'\in \Delta(n)$.  Let
 $e_i'$ with $e_i'\cdot (C-{\bf M}(E)E)<0$
 be the type $I$ exceptional classes over $Y_{\Gamma'}$ which generate
 the simplicial cone ${\cal C}_{\Gamma'}$. 
 Then a similar
argument applies as well and $z$ is in the image of 
${\cal M}_{C-{\bf M}(E)E-\sum_{e_i'; e_i'\cdot 
(C-{\bf M}(E)E)<0}e_i'}\subset {\cal M}_{C-{\bf M}(E)E}$. 

 To identify them on the scheme theoretical level, simply realize that
 the difference of $Z(s_{canon})\times_{M_n}Y(\Gamma)$ from 
$Z(s_{canon}^{\circ})\times_{M_n}Y(\Gamma)$ can be
 analyzed by the various 
analogues of $s_{canon}^{\circ}$ and ${\bf W}_{canon}^{\circ}$,
 involving $\sum_{e_i'\cdot (C-{\bf M}(E)E)<0}e_i'$ over $Y_{\Gamma'}$,
 $Y(\Gamma')\cap Y(\Gamma)\not=\emptyset$.
 By induction, we may get the equality on the sub-schemes of $X$.
$\Box$

\medskip

\begin{rem}\label{rem; gg}
In the special case when $(\Gamma, \sum_{e_i\cdot (C-{\bf M}(E)E)<0}e_i)\gg
 (\Gamma', \sum_{e_j'\cdot (C-{\bf M}(E)E)<0}e_j')$ (consult definition
 \ref{defin; partial} for
 its definition), we have $Y(\Gamma')\subset Y(\Gamma)$ and 
${\cal M}_{C-{\bf M}(E)E-\sum_{e_j'\cdot (C-{\bf M}(E)E)<0}e_j'}
\times_{M_n}Y(\Gamma')$
 is naturally embedded into 
${\cal M}_{C-{\bf M}(E)E-\sum_{e_i\cdot (C-{\bf M}(E)E)<0}e_i}
\times_{M_n}Y(\Gamma)$ as a sub-scheme\footnote{And therefore into
 ${\cal M}_{C-{\bf M}(E)E}\times_{M_n}Y(\Gamma)$ automatically.}.
 Thus we may ignore such 
${\cal M}_{C-{\bf M}(E)E-\sum_{e_j'\cdot (C-{\bf M}(E)E)<0}e_j'}
\times_{M_n}Y(\Gamma')$ in the above union. 
\end{rem}

\medskip

$\diamondsuit$
 Heuristically we may interpret the blowing ups of (the strict transform of)
 ${\cal M}_{C-{\bf M}(E)E}\times_{M_n}Y(\Gamma')$ into $D_{\Gamma'}$ 
and factorizing 
 $s_{canon}$ by $D_{\Gamma'}$ as a mean of removing the 
 contribution of the top Chern class of $\pi_X^{\ast}{\bf W}_{canon}\otimes 
{\bf H}$ from the image scheme of ${\cal M}_{C-{\bf M}(E)E-\sum_{e_i\cdot
 (C-{\bf M}(E)E)<0}e_i}\times_{M_n}Y(\Gamma)$ inside
${\cal M}_{C-{\bf M}(E)E}\times_{M_n}Y(\Gamma)$. This explains why the localized
 contribution is expected to be related to some sort of 
mixed family invariant attached to
 ${\cal M}_{C-{\bf M}(E)E-\sum_{e_i\cdot (C-{\bf M}(E)E)<0}e_i}$ over
 $Y(\Gamma)$, as had been demonstrated in theorem 4 of [Liu5]. 

\medskip

It is vital to reflect and ask the following question,

\noindent $\diamondsuit$ {\bf Question}: How to enumerate/identify the exact 
localized contribution of top Chern class along $D_{\Gamma}$? 
 Do we expect the answer to be expressible as a mixed invariant or do we
 expect to get additional ``correction terms''?

\medskip

 If there are additional correction terms, we have to understand where
 do these terms come from!

In fact when we blow up the various $D_{\Gamma'}$, $\Gamma'\in \Delta(n)$,
 $\Gamma\models \Gamma'$, we have also blown up along the image of
${\cal M}_{C-{\bf M}(E)E-\sum_{e_j'\cdot (C-{\bf M}(E)E)<0}e_j'}
\times_{M_n}Y(\Gamma')$ in ${\cal M}_{C-{\bf M}(E)E}\times_{M_n}Y(\Gamma)$
, with $\Gamma\gg \Gamma'$ as well.

This
 suggests that we have also removed ``accidentally'' some additional 
contributions of the family invariant
 from 
${\cal M}_{C-{\bf M}(E)E-\sum_{e_i\cdot {\bf M}(E)E<0}e_i}\times_{M_n}Y(\Gamma)$
as well\footnote{by the observation in remark \ref{rem; gg}.}.
Therefore, we should expect to 
identify the localized contribution with
 a ``modified'' object instead of a normal mixed family invariant
 of $C-{\bf M}(E)E-\sum_{e_i\cdot {\bf M}(E)E<0}e_i$ over $Y(\Gamma)$.

This explains why in the definitions of the modified invariants, there
 are ``correction terms'' from $\Gamma'\ll \Gamma$ besides the dominated
 leading term.  While these correction terms
 appear naturally in our current setting of the blowing up
 construction and the residual intersection theory, it is also needed to 
avoid the troublesome problem of over-subtracting (see the beginning of 
section 6.5 of [Liu5] for an explanation).

The complete answer to the above question will be given in section \ref{section; 
proof} where we identify the localized algebraic family invariant contribution
 along $D_{\Gamma}$ with ${\cal AFSW}^{\ast}$ inductively. Before we
 present the proof, some additional knowledge about the geometric/combinatorial 
structure of ${\cal M}_{C-{\bf M}(E)E}$ is essential.

\bigskip

Consider the partial ordering $\sqsupset$,

\begin{defin}\label{defin; sq}
 Let $\Gamma, \Gamma'\in \Delta(n)$. 
 The pair $(\Gamma, \sum_{e_i\cdot (C-{\bf M}(E)E)<0}e_i)$ is greater than
 $(\Gamma', \sum_{e_i'\cdot (C-{\bf M}(E)E)<0}e_i')$ under $\sqsupset$ if
 
(i). $Y(\Gamma)\cap Y(\Gamma')\not=\emptyset$.

(ii). The combination of type $I$ classes
 $\sum_{e_i\cdot (C-{\bf M}(E)E)<0}e_i-\sum_{e_i'\cdot (C-{\bf M}(E)E)<0}e_i'$
 is semi-effective over $Y(\Gamma)\cap Y(\Gamma')$. I.e. the combination is
 either identically zero or is represented by effective curves over 
 $Y(\Gamma)\cap Y(\Gamma')$.
\end{defin}

 Geometrically the partial ordering $\sqsupset$ signalizes that 
 some of the type $I$ curves dual to  $e_i, e_i\cdot (C-{\bf M}(E)E)<0$, break into
 more than one component over $Y(\Gamma)\cap Y(\Gamma')$ 
and these $e_j', e_j'\cdot (C-{\bf M}(E)E)<0$ are dual to certain
 components among them.  When $\Gamma\sqsupset \Gamma'$, 
 the semi-effectiveness of $\sum_{e_i\cdot (C-{\bf M}(E)E)<0}e_i-\sum_{e_i'\cdot
 (C-{\bf M}(E)E)<0}e_i'$ over $Y(\Gamma)\cap Y(\Gamma')$ implies that
 ${\cal M}_{C-{\bf M}(E)E-\sum_{e_i\cdot (C-{\bf M}(E)E)<0}e_i}\times_{M_n}
Y(\Gamma)\cap Y(\Gamma')$ is embedded into
${\cal M}_{C-{\bf M}(E)E-\sum_{e_i'\cdot (C-{\bf M}(E)E)<0}e_i'}\times_{M_n}
Y(\Gamma)\cap Y(\Gamma')$. In other words, the image of 
${\cal M}_{C-{\bf M}(E)E-\sum_{e_i\cdot (C-{\bf M}(E)E)<0}e_i}\times_{M_n}
Y(\Gamma)$ in ${\cal M}_{C-{\bf M}(E)E}\times_{M_n}Y(\Gamma)$ does not
 capture ${\cal M}_{C-{\bf M}(E)E}$ above $Y(\Gamma)\cap Y(\Gamma')$ accurately 
and the two loci $Z(s_{canon}^{\circ})={\cal M}_{C-{\bf M}(E)E-
\sum_{e_i\cdot (C-{\bf M}(E)E)<0}e_i}$ and $Z(s_{canon})={\cal M}_{C-{\bf M}(E)E}$
 may differ over
 $Y(\Gamma)\cap Y(\Gamma')$.
 
\medskip

 Refer to fig.6 for an example\footnote{To simplify
the notations, the $\sum_{e_i\cdot (C-{\bf M}(E)E)<0}e_i$ parts have been
 skipped.}.


\begin{figure}
\centerline{\epsfig{file=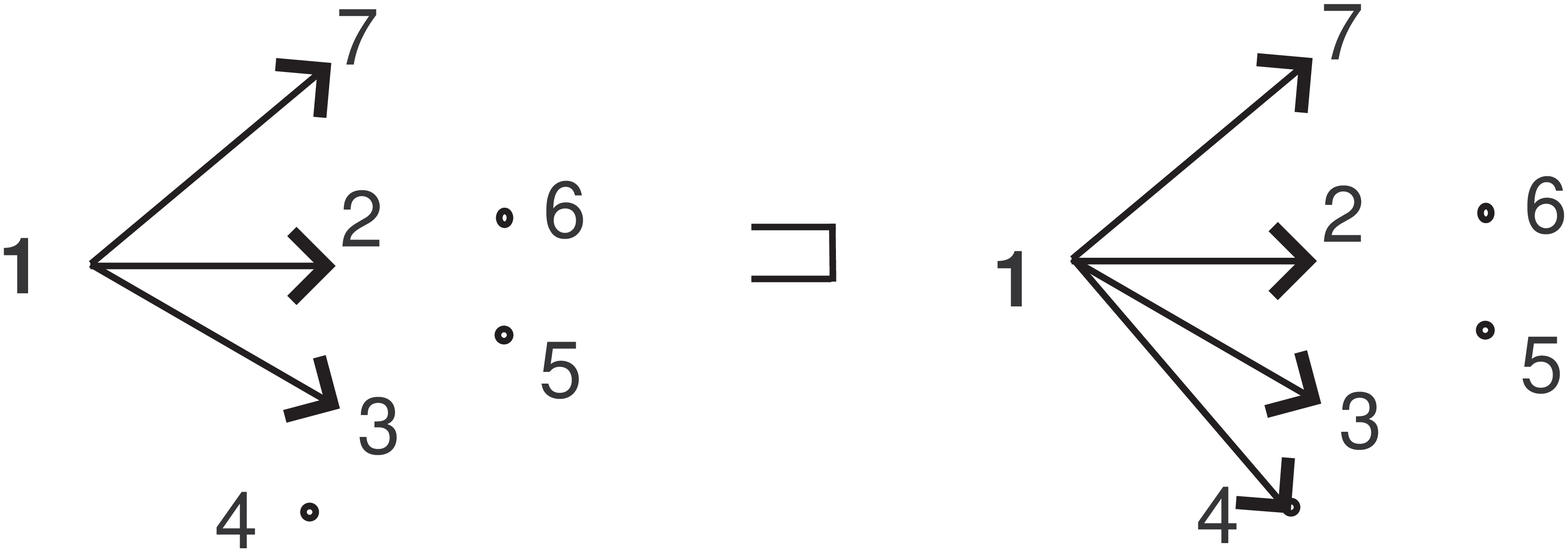,height=4cm}}
\centerline{fig.6}
\centerline{a pair of admissible graphs in $adm(6)$ 
related by the $\sqsupset$ partial ordering.}
\end{figure}

\medskip

\begin{rem}\label{rem; exclude}
These two partial orderings $\gg$ and $\sqsupset$ are exclusive in the
 following sense that if $(\Gamma, \sum_{e_i\cdot (C-{\bf M}(E)E)<0}e_i)\gg
(\Gamma', \sum_{e_i'\cdot (C-{\bf M}(E)E)<0}e_i')$, then the expression
 $\sum_{e_i\cdot (C-{\bf M}(E)E)<0}e_i-\sum_{e_i'\cdot (C-{\bf M}(E)E)<0}e_i'$
 is anti-effective over $Y(\Gamma)\cap Y(\Gamma')=Y(\Gamma')$. Thus,
 the $(\Gamma, \sum_{e_i\cdot (C-{\bf M}(E)E)<0}e_i)$ cannot be greater than
 $(\Gamma', \sum_{e_i'\cdot (C-{\bf M}(E)E)<0}e_i')$ under $\sqsupset$.
\end{rem}

\medskip

 Consider the following setting among three admissible graphs.
Fix a $\Gamma\in \Delta(n)$ and let $\Gamma_1\in I_{\Gamma}$ satisfy
 $Y(\Gamma)\cap Y(\Gamma_1)\not=\emptyset$. Let $\Gamma_2\in adm(n)$ with 
$\Gamma_2<\Gamma$, $\Gamma_2\leq \Gamma_1$,
 be an admissible graph. 
 This implies that $Y_{\Gamma_2}\subset Y(\Gamma)\cap Y(\Gamma_1)$.

  In the following proposition we discuss the few possibilities which can occur, 

\begin{prop}\label{prop; relation}
Suppose that $\Gamma_2\not\in \Delta(n)$, then $Y_{\Gamma_2}\in S_{\Gamma'}$ 
for some $\Gamma'\in \Delta(n)$. If $\Gamma_2\in \Delta(n)$, take 
 $\Gamma'=\Gamma_2$ itself.

  As usual, let $e_i'$, $1\leq i\leq n$,
 denote the type $I$ exceptional classes over $Y_{\Gamma'}$.

 Then either 

(a). $(\Gamma, \sum_{e_i\cdot (C-{\bf M}(E)E)<0}e_i)\gg 
(\Gamma', \sum_{e_i'\cdot (C-{\bf M}(E)E)<0}e_i')$.

or 
(b). $\exists$ an intermediate $\Gamma''\in \Delta(n)$, $\Gamma''\not=\Gamma$,
 such that 
 $(\Gamma, \sum_{e_i\cdot (C-{\bf M}(E)E)<0}e_i)\sqsupset 
(\Gamma'', \sum_{e_i''\cdot (C-{\bf M}(E)E)<0}e_i'')$ and 
 $(\Gamma'', \sum_{e_i''\cdot (C-{\bf M}(E)E)<0}e_i'')\gg
(\Gamma', \sum_{e_i\cdot (C-{\bf M}(E)E)<0}e_i')$.
\end{prop}

 Suppose that 
$\Gamma'$ is as described in the statement of the proposition, then 
$Y(\Gamma)\cap S_{\Gamma'}\not=\emptyset$.
 It implies that $\Gamma\succ \Gamma'$ and then $\Gamma'\in I_{\Gamma}$.

The figure 7, 8 illustrate how proposition \ref{prop; relation}
 holds in a special case.


\begin{figure}\label{fig.7}
\centerline{\epsfig{file=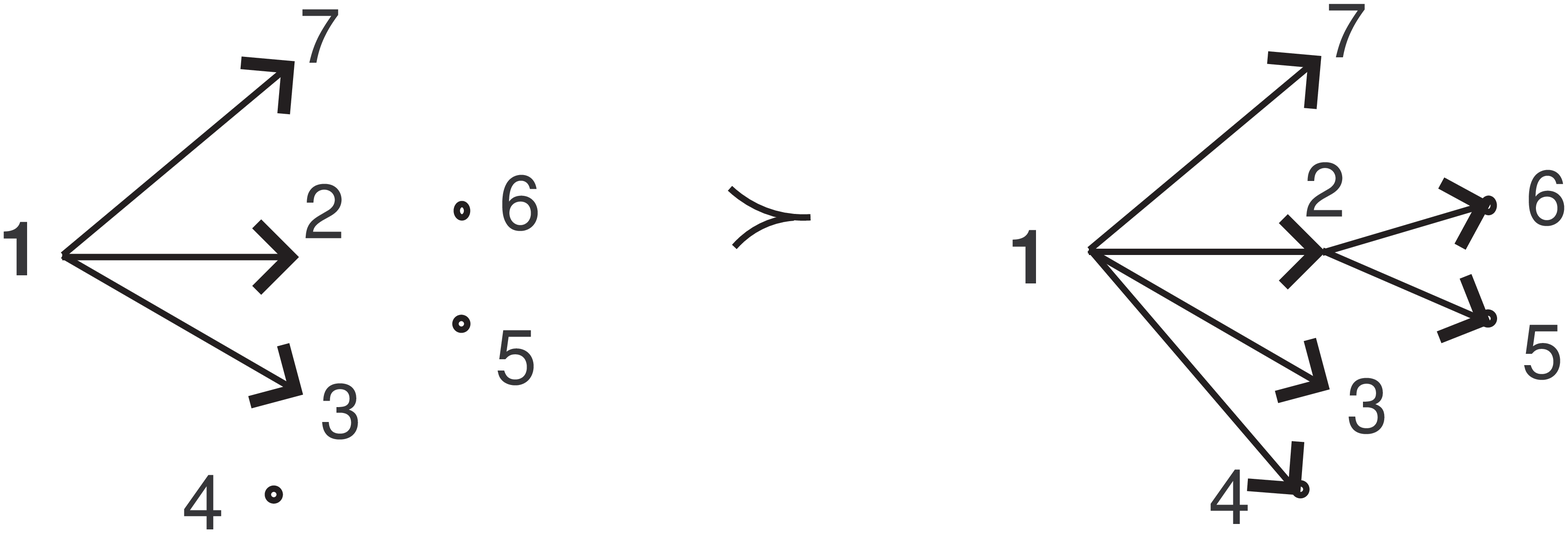,height=4cm}}
\centerline{fig.7}
\centerline{a pair of admissible graphs in $adm(7)$ 
related by the $\succ$ partial ordering.}
\end{figure}


\begin{figure}
\centerline{\epsfig{file=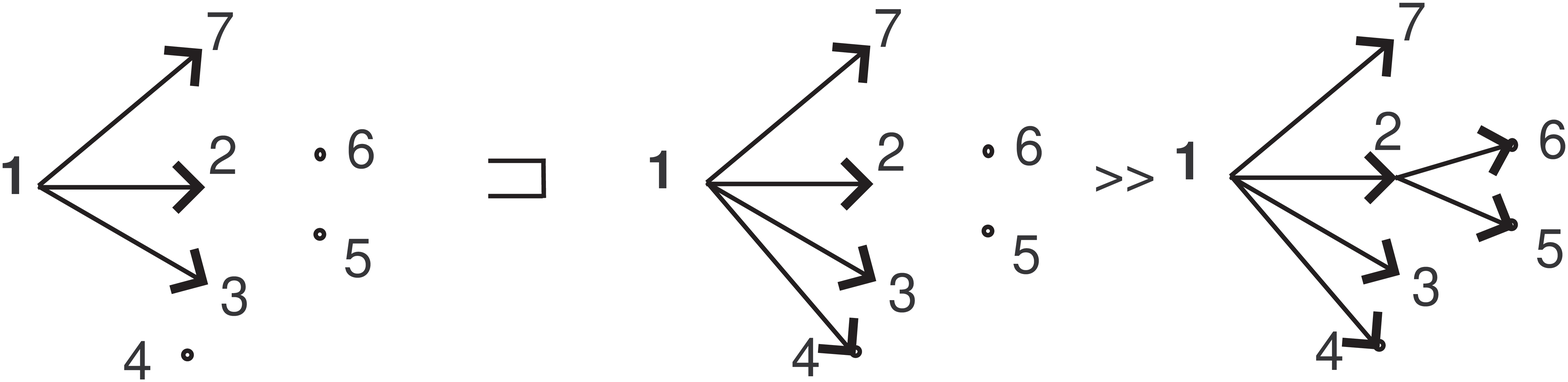,height=4cm}}
\centerline{fig.8}
\centerline{following fig.7 above, the inserted admissible graph $\in adm(7)$ 
 in the middle}
\centerline{is $\sqsubset$ than the admissible graph on the 
 left hand side, but $\gg$ than the admissible graph}
\centerline{on the right hand side. this situation 
corresponds to {\bf proposition} \ref{prop; relation} 
 case (b).}
\end{figure}

\medskip

\noindent Proof of proposition \ref{prop; relation}:
 If $\Gamma_2\in \Delta(n)$, take $\Gamma'=\Gamma_2$. 
If $\Gamma_2\not\in \Delta(n)$, by proposition \ref{prop; large} and 
lemma \ref{lemm; union3}, $Y_{\Gamma_2}\subset
 S_{\Gamma'}$ for some $\Gamma'\in \Delta(n)$. We know that
 $\Gamma'$ cannot be $\Gamma$ itself. Otherwise it implies immediately
 that ${\cal C}_{\Gamma}\supset {\cal C}_{\Gamma_1}$ because now
 we have $Y_{\Gamma_2}\subset
 S_{\Gamma'}=S_{\Gamma}$ 
(where the exceptional cone is equal to ${\cal C}_{\Gamma}$) and
because $Y_{\Gamma_2}\subset Y(\Gamma_1)$. But this implies that
 $\Gamma_1\succ \Gamma$ and then such a 
$\Gamma_1$ can {\bf NOT} be in $I_{\Gamma}$.
 A contradiction to our assumption!  
 
 Therefore one may replace $\Gamma_2$ by some $\Gamma'\in \Delta(n)$.
 In any case we still have
 $Y(\Gamma')\cap Y(\Gamma)\supset Y(\Gamma_2)\cap Y(\Gamma)\not=\emptyset$.

Consider the cone ${\cal C}_{\Gamma'}$. Because 
$Y_{\Gamma_2}\subset S_{\Gamma'}$, the extremal generators
\footnote{I.e. primitive generators of one-edges (extremal rays) at the
 boundary of the cone.} 
of ${\cal C}_{\Gamma'}$
 are exactly the type $I$ exceptional 
classes over $Y_{\Gamma_2}$ which pair negatively
 with $C-{\bf M}(E)E$. On the other hand, $Y_{\Gamma_2}\subset Y(\Gamma)$ 
(therefore $S_{\Gamma'}\cap Y(\Gamma)\not=\emptyset$ and
 this implies that ${\cal C}_{\Gamma}\subset {\cal C}_{\Gamma'}$.

 We separate our discussion into a few cases.

 \noindent (A). Suppose that 
all the extremal generators (among the type $I$ exceptional classes)
 of the simplicial cone ${\cal C}_{\Gamma}$
 remain to be the extremal generators of the simplicial cone ${\cal C}_{\Gamma'}$:
 As we know ${\cal C}_{\Gamma}\not={\cal C}_{\Gamma'}$, there must be
 additional extremal generators of ${\cal C}_{\Gamma'}$ away from the
 boundary of ${\cal C}_{\Gamma}$. Then
 $Y(\Gamma)$ is the locus of co-existence of the type $I$ classes 
$e_i$, $e_i\cdot (C-{\bf M}(E)E)<0$. Then $Y(\Gamma')$ is characterized
 (by proposition \ref{prop; transversal}) as the
 co-existence loci of all $e_i, e_i\cdot (C-{\bf M}(E)E)<0$ and some other type 
$I$ exceptional classes. So we know $Y(\Gamma')\subset Y(\Gamma)$ and we have,

 $$(\Gamma, \sum_{e_i\cdot (C-{\bf M}(E)E)<0}e_i)\gg 
(\Gamma', \sum_{e_i'\cdot (C-{\bf M}(E)E)<0}e_i').$$

\noindent (B). Suppose that at least one
 of the extremal generators (among the type $I$ exceptional
 classes) fails to be an extremal generator of ${\cal C}_{\Gamma'}$. Then
 the curve it represents must break into more than one irreducible component above
 $Y(\Gamma_2)$. Define an index set $P\subset \{1, 2, \cdots, n\}$ to consist of
 all the subscripts $i$ of $e_i$ such that the extremal generator $e_i$
 of the cone ${\cal C}_{\Gamma}$ fails to be an extremal generator of 
${\cal C}_{\Gamma'}$.

 By proposition 4 of [Liu4], each $e_i, i\in P$, can be expressed uniquely
 as an effective 
integral combination of the extremal generators of the simplicial cone
 ${\cal C}_{\Gamma'}$. Let the index set $P''\subset \{1, 2, \cdots, n\}$ 
 consists of all the subscripts $j$ so that $e_j'$ involve in the 
 integral effective combination of at least one $e_i, i\in P$.
 Then we may write $e_i=\sum_{j\in P''}c_{i, j}e_j'$, with $c_j\geq 0$. We know
 that $\sum_{i\in P}c_{i, j}\geq 1$ for all $j\in P''$.

 Then the collection of type $I$ classes $e_i', i\in P''$, generate a 
simplicial sub-cone ${\cal C}$
 of ${\cal C}_{\Gamma'}$.  Because $e_i'\cdot e_j'\geq 0$ for 
$i\not=j$ in $P''$, by proposition \ref{prop; reverse} of section \ref{section;
strata}, $\exists$ an admissible graph $\Gamma''$ such that 
 $e_i''=e_i'$ for $i\in P''$ and $e_i''$ are $-1$ classes for $i\not\in P''$. 
The co-existence locus of $e_i', i\in P''$ characterizes
 the closure of the admissible stratum $Y(\Gamma'')$ 
with ${\cal C}_{\Gamma''}={\cal C}$ and since (by definition)
 $e_j''=e_j', j\in P''$, they are exactly the type $I$ exceptional classes above 
$Y_{\Gamma''}$ which pair negatively with $C-{\bf M}(E)E$.

 Apparently by the construction of $\Gamma''$ we have\footnote{Because
 non $-1$ classes $e_i''$ are selected from $e_j'$.}, 
$$(\Gamma'', \sum_{e_i''\cdot (C-{\bf M}(E)E)<0}e_i'')\gg
(\Gamma', \sum_{e_i'\cdot (C-{\bf M}(E)E)<0}e_i').$$

We may rewrite $\sum_{e_i\cdot (C-{\bf M}(E)E)<0}e_i-
\sum_{e_j''\cdot (C-{\bf M}(E)E)<0}e_j''$ as 

$$\sum_{e_i\cdot (C-{\bf M}(E)E)<0; i\not\in P}e_i+
\sum_{i\in P}e_i-
\sum_{j\in P''}e_j''=\sum_{e_i\cdot (C-{\bf M}(E)E)<0; i\not\in P}e_i+
\sum_{i\in P}\sum_{j\in P''}c_{i, j}e_j''-
\sum_{j\in P''}e_j''$$
$$=\sum_{e_i\cdot (C-{\bf M}(E)E)<0; i\not\in P}e_i+
\{\sum_{j\in P''}\bigl((\sum_{i\in P}c_{i, j})-1\bigr)e_j''\}.$$

 The first sum in the final expression is semi-effective over $Y(\Gamma)$ because
 it is a semi-effective combinations of some $e_i, i\not\in P$ which are effective
 over $Y(\Gamma)$. On the other hand, the second sum in the final expression
 is semi-effective over $Y(\Gamma'')$ because (i). Each $e_j'', j\in P''$ is 
effective over $Y(\Gamma'')$ and $P''$ is defined to be the index set containing
 all the $j$ such that $e_j''=e_j'$ are used in expressing $e_i, i\in P$.

 So for any fixed $j\in P''$, the sum $\sum_{i\in P}c_{i, j}\geq 1$ and
 therefore $(\sum_{i\in P}c_{i, j})-1\geq 0$.  So the sum of these
 two expressions is still semi-effective over the intersection
 $Y(\Gamma)\cap Y(\Gamma'')$.  So the existence of such an
 intermediate $\Gamma''$ has been proved. $\Box$ 

\medskip

 As a consequence of proposition \ref{prop; relation}, the image of the
 moduli space ${\cal M}_{C-\sum_{e_i'\cdot (C-{\bf M}(E)E)<0}e_i'}\times_{M_n}
Y(\Gamma')$ into 
 ${\cal M}_{C-{\bf M}(E)E}\times_{M_n}Y(\Gamma)$ is either contained 
directly in
 the image of ${\cal M}_{C-\sum_{e_i\cdot (C-{\bf M}(E)E)<0}e_i}
\times_{M_n}Y(\Gamma)$
 in ${\cal M}_{C-{\bf M}(E)E}\times_{M_n}Y(\Gamma)$ itself or is contained in
 the image of some intermediate 
${\cal M}_{C-{\bf M}(E)E-\sum_{e_i''\cdot (C-{\bf M}(E)E)<0}e_i''}\times_{M_n}
Y(\Gamma'')$. This gives an additional hierarchical 
 structure among the various
${\cal M}_{C-{\bf M}(E)E-\sum_{e_i\cdot (C-{\bf M}(E)E)<0}e_i}$.

 When we attempt to identify the localized algebraic family 
invariant contribution to $D_{\Gamma}$, 
the above observation motivates us to
 reduce the list of blowing ups preceding $\Gamma$, 
$I_{\Gamma}$, by the following recipe:

\begin{defin}\label{defin; reduce}
Define the subset $\bar{I}_{\Gamma}\subset I_{\Gamma}$ by
removing all those 
 $\Gamma'\in I_{\Gamma}$ such that
  (i). $S_{\Gamma'}\cap Y(\Gamma)=\emptyset$ or
 (ii). $\exists \Gamma''\in I_{\Gamma}$ with
 $(\Gamma'', \sum_{e_i''\cdot (C-{\bf M}(E)E)<0}e_i'')
\gg (\Gamma', \sum_{e_i'\cdot (C-{\bf M}(E)E)<0}e_i')$.
\end{defin}

The reduced index 
set $\bar{I}_{\Gamma}$ collects all the $\Gamma'\in I_{\Gamma}$ which
 satisfy either $(\Gamma, \sum_{e_i\cdot (C-{\bf M}(E)E)<0}e_i)\gg
 (\Gamma', \sum_{e_i'\cdot (C-{\bf M}(E)E)<0}e_i')$ or
 $(\Gamma, \sum_{e_i\cdot (C-{\bf M}(E)E)<0}e_i)\sqsupset 
 (\Gamma', \sum_{e_i'\cdot (C-{\bf M}(E)E)<0}e_i')$.

The set $\bar{I}_{\Gamma}$ inherits a linear ordering from $I_{\Gamma}$, still
 denoted by $\models$.

\medskip

 For all $\Gamma\in \Delta(n)$, we can define the corresponding 
$\bar{I}_{\Gamma}$ by the recipe of definition \ref{defin; reduce}.
 The next lemma characterizes the relationship between the 
index set $\bar{I}_{\overline{\Gamma}}$ and $\bar{I}_{\Gamma}$ when 
$\overline{\Gamma}\ll \Gamma$.

\begin{lemm}\label{lemm; restriction}
Let $\Gamma\in \Delta(n)$ and let $\overline{\Gamma}\in \bar{I}_{\Gamma}$ satisfies
 $(\Gamma, \sum_{e_i\cdot (C-{\bf M}(E)E)<0}e_i)\gg  
(\overline{\Gamma}, \sum_{\overline{e}_i\cdot (C-{\bf M}(E)E)<0}\overline{e}_i)$. 
Then the reduced index set 
$\bar{I}_{\overline{\Gamma}}$ of $\overline{\Gamma}$
 is the set of elements $\Gamma'$ 
in $\bar{I}_{\Gamma}$ satisfying,

\medskip

(i). $Y(\overline{\Gamma})\cap S_{\Gamma'}\not=\emptyset$,

\medskip

(ii).  $\Gamma'$ is smaller than
$\overline{\Gamma}$ under the linear ordering \footnote{Consult 
 page \pageref{models} for its definition.} $\models$ in $\bar{I}_{\Gamma}$.
\end{lemm}

\noindent Proof: 
It is apparent that when $Y(\overline{\Gamma})\subset Y(\Gamma)$,
 the condition $Y(\overline{\Gamma})\cap S_{\Gamma'}\not=\emptyset$ (the 
direct consequence of $\Gamma'\in \bar{I}_{\overline{\Gamma}}$) implies 
 $Y(\Gamma)\cap S_{\Gamma'}\not=\emptyset$. 
On the other hand $\overline{\Gamma}\ll \Gamma$ and $\Gamma''\sqsubset \Gamma$
 imply (see the next footnote below on page \pageref{explain})
 $\overline{\Gamma}\sqsupset \Gamma''$.
Thus to show that
 the condition (ii) is satisfied it suffices to show that
 for all $\Gamma'\in I_{\Gamma}-\bar{I}_{\Gamma}$, either it is due to  

(1). $\Gamma'$ is 
larger than or equal to 
$\overline{\Gamma}$ under the linear ordering $\models$, or 

(2). $Y(\overline{\Gamma})\cap S_{\Gamma'}=\emptyset$,
 or 

(3). $\Gamma'\in I_{\overline{\Gamma}}-\bar{I}_{\overline{\Gamma}}$.

 Suppose that $\Gamma'\in I_{\Gamma}-\bar{I}_{\Gamma}$, 
 but the conditions 
(1). and (2). do not hold. We plan to show that (3). has to hold.
Because of the violating of (1). an (2).,
 $\Gamma'$ is smaller than $\overline{\Gamma}$ under the linear ordering 
$\models$ and
 $Y(\overline{\Gamma})\cap S_{\Gamma'}\not=\emptyset$.
  
 Then by the definition of $\bar{I}_{\Gamma}$, the assumption that
 $\Gamma'$ is not in $\bar{I}_{\Gamma}$ and by proposition \ref{prop; relation} (
notice that $\Gamma'$ satisfies the assumption of this proposition),
 there must exist some intermediate
$\Gamma''\in I_{\Gamma}\subset \Delta(n)$ such that

$$(\Gamma, \sum_{e_i\cdot (C-{\bf M}(E)E)<0}e_i)\sqsupset 
(\Gamma'', \sum_{e_i''\cdot (C-{\bf M}(E)E)<0}e_i'')\gg
 (\Gamma', \sum_{e_i'\cdot (C-{\bf M}(E)E)<0}e_i').$$
 
 The condition $(\Gamma, \sum_{e_i\cdot (C-{\bf M}(E)E)<0}e_i)\sqsupset 
(\Gamma'', \sum_{e_i''\cdot (C-{\bf M}(E)E)<0}e_i'')$ is equivalent to
 the semi-effectiveness of 
 $\sum_{e_i\cdot (C-{\bf M}(E)E)<0}e_i-\sum_{e_i''\cdot (C-{\bf M}(E)E)<0}e_i''$
 over $Y(\Gamma'')\cap Y(\Gamma)\not=\emptyset$. 

 By our assumption on $\overline{\Gamma}$ in the statement of this lemma,
 $(\Gamma, \sum_{e_i\cdot (C-{\bf M}(E)E)<0}e_i)\gg  
(\overline{\Gamma}, \sum_{\overline{e}_i\cdot (C-{\bf M}(E)E)<0}\overline{e}_i)$.
 This implies
that for all the type $I$ exceptional classes above $Y_{\Gamma}$ satisfying 
$e_i\cdot (C-{\bf M}(E)E)<0$, $\overline{e}_i=e_i$ is the 
corresponding
 type $I$ exceptional class over $Y_{\overline{\Gamma}}$ and there are
some additional $\overline{e}_s$ with $\overline{e}_s\cdot (C-{\bf M}(E)E)<0$
 other than those $e_i$, $e_i\cdot (C-{\bf M}(E)E)<0$.

 This implies that the class $\sum_{\overline{e}_j\cdot (C-{\bf M}(E)E)<0}
\overline{e}_j-
\sum_{e_i\cdot (C-{\bf M}(E)E)<0}e_i$ has to be effective over 
$Y(\overline{\Gamma})=\overline{Y(\overline{\Gamma})}$, 
which is a subset of $Y(\Gamma)$.

Thus, 
$$\hskip -1.3in
\{\sum_{\overline{e}_j\cdot (C-{\bf M}(E)E)<0}\overline{e}_j-
\sum_{e_i\cdot (C-{\bf M}(E)E)<0}e_i\}+
\{\sum_{e_i\cdot (C-{\bf M}(E)E)<0}e_i-\sum_{e_i''\cdot (C-{\bf M}(E)E)<0}e_i''\}
=\sum_{\overline{e}_j\cdot (C-{\bf M}(E)E)<0}\overline{e}_j-\sum_{e_i''
\cdot (C-{\bf M}(E)E)<0}e_i''$$

is effective over the intersection of 
 $Y(\Gamma'')\cap Y(\Gamma)$ and $Y(\overline{\Gamma})$, $Y(\Gamma'')
\cap Y(\Gamma)\cap Y(\overline{\Gamma})=Y(\Gamma'')\cap 
Y(\overline{\Gamma})$. In particular, the final expression is semi-effective over
\label{explain} \footnote{The argument essentially shows
 that $\Gamma\gg \overline{\Gamma}$,
and $\Gamma\sqsupset \Gamma''$ imply $\overline{\Gamma}\sqsupset \Gamma''$.}
 $Y(\Gamma'')\cap Y(\overline{\Gamma})$ (since $Y(\overline{\Gamma})\subset
 Y(\Gamma)$.

This implies that 

$$(\overline{\Gamma}, \sum_{\overline{e}_i\cdot (C-{\bf M}(E)E)<0}
\overline{e}_i)\sqsupset 
 (\Gamma'', \sum_{e_i''\cdot (C-{\bf M}(E)E)<0}e_i'')\gg 
(\Gamma', \sum_{e_i\cdot (C-{\bf M}(E)E)<0}e_i)$$
 as well (the $\gg$ inequality within the formula is already known).

 So the element 
$\Gamma'$ must be removed from $I_{\overline{\Gamma}}$ in forming 
 $\bar{I}_{\overline{\Gamma}}$ and 
is not in $\bar{I}_{\overline{\Gamma}}$, either. $\Box$

\medskip

 Fixing a $\Gamma''$ such that 
 $(\Gamma, \sum_{e_i\cdot (C-{\bf M}(E)E)<0}e_i)\sqsupset 
  (\Gamma'', \sum_{e_i''\cdot(C-{\bf M}(E)E)<0}e_i'')$,
 the reduction from $I_{\Gamma}$ to $\bar{I}_{\Gamma}$ enables us to group
 the family moduli spaces above $Y(\Gamma')$, 
 ${\cal M}_{C-{\bf M}(E)E}\times_{M_n}Y(\Gamma')$, of all the $\Gamma'$, 
 satisfying $(\Gamma'', \sum_{e_i''\cdot(C-{\bf M}(E)E)<0}e_i'')\gg (\Gamma',
 \sum_{e_i'\cdot (C-{\bf M}(E)E)<0}e_i')$ together as sub-moduli spaces
 of ${\cal M}_{C-{\bf M}(E)E}\times_{M_n}Y(\Gamma'')$. 
 Instead of blowing up all these ${\cal M}_{C-{\bf M}(E)E}\times_{M_n}Y(\Gamma')$
 individually, we blow up the whole 
${\cal M}_{C-{\bf M}(E)E}\times_{M_n}Y(\Gamma'')$ at once
\footnote{In section \ref{subsection; independence}, 
proposition \ref{prop; identical}
 implies that the re-grouping of the restricted family moduli spaces like this 
does not affect the localized contribution of 
the family invariant along $D_{\Gamma}$, thanks to the birational invariance
 of Segre classes of normal cones.}.

 For all $\Gamma'\in \bar{I}_{\Gamma}$, it satisfies either
 $(\Gamma, \sum_{e_i\cdot (C-{\bf M}(E)E)<0}e_i)\gg 
 (\Gamma', \sum_{e_i'\cdot (C-{\bf M}(E)E)<0}e_i')$, or 
 $(\Gamma, \sum_{e_i\cdot (C-{\bf M}(E)E)<0}e_i)\sqsupset 
 (\Gamma', \sum_{e_i'\cdot (C-{\bf M}(E)E)<0}e_i')$.  For the enumeration purpose,
 we would like to rearrange the blowing up orderings (by $\models$) such that
 those $\Gamma'$, with $(\Gamma, \sum_{e_i\cdot (C-{\bf M}(E)E)<0}e_i)\gg 
 (\Gamma', \sum_{e_i'\cdot (C-{\bf M}(E)E)<0}e_i')$ are blown up later than
 those related to $\Gamma$ by $\sqsupset$.

 To achieve this goal, introduce a new linear ordering among
 $\bar{I}_{\Gamma}$, denoted by $\vdash$,

\begin{defin}\label{defin; newordering}
 Let $\Gamma_1\in \bar{I}_{\Gamma}$. Define the ``accumulation'' 
$A_{\Gamma_1}=\{\Gamma_1\}$ if
 $\Gamma_1\ll \Gamma$. Define $A_{\Gamma_1}=\{\Gamma'|\Gamma'\ll \Gamma_1, 
 \Gamma'\in I_{\Gamma}\}$ if $\Gamma_1\sqsubset \Gamma$.
 For all $\Gamma_1\in \bar{I}_{\Gamma}$, each accumulation 
$A_{\Gamma_1}$ has a unique
 smallest element under $\models$.

 Define a new linear ordering $\vdash$ on $\bar{I}_{\Gamma}$ by the following
 recipe:

(i). Suppose that both $\Gamma_1$ and $\Gamma_2$ are simultaneously 
$\sqsubset \Gamma$ or $\ll \Gamma$, define $\Gamma_1$ is greater
 than $\Gamma_2$ under $\vdash$ if the smallest element within
 the accumulation $A_{\Gamma_1}$
 under $\models$ is larger than (under $\models$) 
the smallest element of $A_{\Gamma_2}$.

(ii). Suppose that $\Gamma_1\ll \Gamma$ but $\Gamma_2\sqsubset \Gamma$,
 define $\Gamma_1$ is larger than $\Gamma_2$ under $\vdash$, i.e. 
$\Gamma_1\vdash \Gamma_2$.
\end{defin}

  The heuristic motivation for such a new linear ordering is that the precedence
 among the sequence of 
 blowing ups should be determined by the corresponding 
precedence of the smallest graph under $\models$
 in the set $A_{\Gamma_1}$.

In the new linear ordering $\vdash$, those $\Gamma'$ with $\Gamma'\ll \Gamma$ 
accumulate at the
 larger end of $\bar{I}_{\Gamma}$.

\begin{defin}\label{defin; ggset}
Define $\bar{I}_{\Gamma}^{\gg}\subset \bar{I}_{\Gamma}$ 
to be the set of all elements $\Gamma'$ in 
$\bar{I}_{\Gamma}$ such that the following condition
$(\Gamma, \sum_{e_i\cdot (C-{\bf M}(E)E)<0}e_i)\gg 
(\Gamma', \sum_{e_i'\cdot (C-{\bf M}(E)E)<0}e_i')$ holds.
\end{defin}

 Then definition \ref{defin; newordering}
  implies that an arbitrary element in 
$\bar{I}_{\Gamma}^{\gg}$ is greater than (under the newly defined 
$\vdash$) an arbitrary element in
 $\bar{I}_{\Gamma}-\bar{I}_{\Gamma}^{\gg}$.

\medskip

 The revised 
linear ordering $\vdash$ will play an essential role in the enumeration
 of the algebraic family Seiberg-Witten invariants in the next section.

\medskip 

\section{The Inductive Proof and the Identification with the Modified
 Family Invariants}\label{section; proof}

\bigskip

 The goal of this section is to finish up the proof of the main theorem
 by identifying the localized contribution of the original algebraic family
Seiberg-Witten invariant of $C-{\bf M}(E)E$ 
over $D_{\Gamma}$ with the modified algebraic
 family Seiberg-Witten invariant defined in subsection \ref{subsection; modinv}.

 In subsection \ref{subsection; comb}, we have
 introduced the reduced index sets $\bar{I}_{\Gamma}$ (see definition
 \ref{defin; reduce}), its subset 
$\bar{I}_{\Gamma}^{\gg}$ (see definition \ref{defin; ggset}) 
and the new linear ordering $\vdash$ (see definition \ref{defin; newordering}).
 We may modify the original blowing up
 sequence by blowing ups the (strict transforms) of all the
 ${\cal M}_{C-{\bf M}(E)E}\times_{M_n}\times Y(\Gamma')\subset 
 {\bf P}({\bf V}_{canon})$ for $\Gamma'\in \bar{I}_{\Gamma}$ instead of 
$I_{\Gamma}$,
 starting from the smallest element under $\models$ or 
$\vdash$ and along the reversed
$\models$ or $\vdash$ orderings. 

 In section \ref{subsection; independence}, we prove the
 independence of the localized top Chern class contribution along
 $D_{\Gamma}$ to the detailed history of the blowing ups performed
 ahead of $\Gamma$ under the reversed ordering of $\models$ or 
 $\vdash$.

 In section \ref{subsection; identifyMA}, we present the key argument to identify
 the integral of the cap product of top power of $c_1({\bf H})$ with the
 localized contribution of top Chern class along $D_{\Gamma}$ with the
 modified algebraic family Seiberg-Witten invariants. 

 Then in section \ref{subsection; proofmain} 
we finish the proof of the main theorem in the paper by combining the
discussion in section \ref{section; bundle}, \ref{section; discrepancy}, 
\ref{section; main} and the current section.

 In section \ref{subsection; transv}, we show with the help of G$\ddot{o}$ttsche's
 argument that ${\cal AFSW}_{M_{\delta+1}\times \{t_L\}\mapsto 
M_{\delta}\times \{t_L\}}^{\ast}(1, C-2\sum E_i)$ can be realized as a counting
 of discrete number of nodal curves in a generic $\delta$-dimensional
 linear sub-system of a $5\delta-1$ very ample line bundle $L\mapsto M$.

\medskip

\subsection{The Independence of the Localized Top Chern Class Contribution
 to the Orderings of the Blowing Ups}\label{subsection; independence}

\bigskip

It makes sense to ask the following question:

\medskip

\noindent {\bf Question}: For a fixed $\Gamma\in \Delta(n)-\{\gamma_n\}$, 
is the localized top Chern class contribution over
 $D_{\Gamma}$ independent \footnote{By 'independence', we do not mean
 for arbitrary blowing ups. See below for more details.}
 to the ``history'' of the sequences of blowing ups we 
had performed on $X$ before the one associated with $\Gamma$? 
Different choices of earlier blowing ups leads to mutually birational divisor
 $D_{\Gamma}$. Thus it becomes a non-trivial question to ask. 
More precisely
 suppose that we blow up the scheme $X={\bf P}({\bf V}_{canon})$ along the
 strict transformations of ${\cal M}_{C-{\bf M}(E)E}\times_{M_n}Y(\Gamma')$ 
following the reversed 
ordering $(I_{\Gamma}, \models)$, $(\bar{I}_{\Gamma}, \models)$ or
 $(\bar{I}_{\Gamma}, \vdash)$, do we get ``identical'' localized top Chern
 class contributions upon the resulting exceptional divisor over 
${\cal M}_{C-{\bf M}(E)E}\times_{M_n}Y(\Gamma)$ in the three different
 cases?

 The answer to this question is affirmative as will be
 shown in the proof of the
 following proposition.

 Before we state our result, let us introduce some notations.

 Denote the consecutively blown ups of $X={\bf P}({\bf V}_{canon})$ along
 the strict transformations of ${\cal M}_{C-{\bf M}(E)E}\times_{M_n}Y(\Gamma')$,
 $\Gamma'\in {I}_{\Gamma}$ following the reversed ordering of $\models$ by
 $X_{\Gamma}$. The 
 $\pi_{X_{\Gamma}}:X_{\Gamma}\mapsto M_n\times T(M)$ is the projection map
 to $M_n\times T(M)$. 
 Some repeated applications of the 
residual intersection formula of top Chern classes, proposition
 \ref{prop; rif}, results in the residual obstruction vector bundle 
$\pi_{{X}_{\Gamma}}^{\ast}{\bf W}_{canon}\otimes 
{\bf H}\otimes_{\Gamma'\in {I}_{\Gamma}}{\cal O}(-{D}_{\Gamma'})$.

 Likewise we denote the consecutively blown ups of $X$ along
 the strict transformations of $Z(s_{canon})\times_{M_n}Y(\Gamma')
={\cal M}_{C-{\bf M}(E)E}\times_{M_n}Y(\Gamma')$,
$\Gamma'\in \bar{I}_{\Gamma}$ with the reversed ordering of $\models$ (or $\vdash$)
 by $\bar{X}_{\Gamma}$ with the projection map
$\pi_{\bar{X}_{\Gamma}}:\bar{X}_{\Gamma}\mapsto M_n\times T(M)$ (or
$\hat{X}_{\Gamma}$ with the projection map
$\pi_{\hat{X}_{\Gamma}}:\hat{X}_{\Gamma}\mapsto M_n\times T(M)$).
 We denote the corresponding exceptional divisors by $\bar{D}_{\Gamma'}$ (or 
$\hat{D}_{\Gamma'}$),
 $\Gamma'\in \bar{I}_{\Gamma}$ and the corresponding residual obstruction
 vector bundle is
 $\pi_{\bar{X}_{\Gamma}}^{\ast}{\bf W}_{canon}\otimes 
{\bf H}\otimes_{\Gamma'\in \bar{I}_{\Gamma}}{\cal O}(-\bar{D}_{\Gamma'})$ (or
$\pi_{\hat{X}_{\Gamma}}^{\ast}{\bf W}_{canon}\otimes 
{\bf H}\otimes_{\Gamma'\in \bar{I}_{\Gamma}}{\cal O}(-\hat{D}_{\Gamma'})$),
 respectively.

 By blowing up the strict transform of 
${\cal M}_{C-{\bf M}(E)E}\times_{M_n}Y(\Gamma)$ in $X_{\Gamma}$, $\bar{X}_{\Gamma}$
 or $\hat{X}_{\Gamma}$, 
denote the blown up schemes by 
 $\tilde{X}_{\Gamma}, \tilde{\bar{X}}_{\Gamma}$ and $\tilde{\hat{X}}_{\Gamma}$,
 respectively. We denote the blown up exceptional divisors from the
 strict transformations of ${\cal M}_{C-{\bf M}(E)E}\times_{M_n}Y(\Gamma)
=Z(s_{canon})\times_{M_n}Y(\Gamma)$ in $X_{\Gamma}$, $\bar{X}_{\Gamma}$, 
 $\hat{X}_{\Gamma}$ by
${D}_{\Gamma}$, $\bar{D}_{\Gamma}$, $\hat{D}_{\Gamma}$, respectively. 

By applying lemma
 \ref{lemm; local}, we get the corresponding 
localized contributions of the top Chern classes
 for $\pi_{X_{\Gamma}}^{\ast}{\cal W}_{canon}\otimes {\bf H} \otimes_{\Gamma'\in 
I_{\Gamma}}{\cal O}(-D_{\Gamma'})$, or
 for $\pi_{\bar{X}_{\Gamma}}^{\ast}
{\cal W}_{canon}\otimes {\bf H} \otimes_{\Gamma'\in 
\bar{I}_{\Gamma}}{\cal O}(-\bar{D}_{\Gamma'})$ or for
$\pi_{\hat{X}_{\Gamma}}^{\ast}{\cal W}_{canon}\otimes {\bf H} \otimes_{\Gamma'\in 
\bar{I}_{\Gamma}}{\cal O}(-\hat{D}_{\Gamma'})$ along 
 $D_{\Gamma}$, $\bar{D}_{\Gamma}$ and $\hat{D}_{\Gamma}$, respectively.

 These three spaces $X_{\Gamma}, \bar{X}_{\Gamma}, \hat{X}_{\Gamma}$ are all
 birational and all of them map onto $X$ through the blowing down projection maps.
The following proposition asserts that the images of the
 localized contributions of top Chern classes along $D_{\Gamma}$, 
$\bar{D}_{\Gamma}$ and $\hat{D}_{\Gamma}$ in ${\cal A}_{\cdot}(X)$
 are all equal and provide an explanation. 

\begin{prop}\label{prop; identical} 
The images in ${\cal A}_{dim_{\bf C}X-rank_{\bf C}{\bf W}_{canon}}(X)$
of the localized contribution of the top Chern class for
$\pi_X^{\ast}{\cal W}_{canon}\otimes {\bf H}\otimes_{\Gamma'\in 
I_{\Gamma}}{\cal O}(-D_{\Gamma'})$ over the strict transformation
\footnote{Along the blown up divisor $D_{\Gamma}\subset \tilde{X}_{\Gamma}$ above
 the strict transformation of ${\cal M}_{C-{\bf M}(E)E}\times_{M_n}Y(\Gamma)$
 in $X_{\Gamma}$.} of 
${\cal M}_{C-{\bf M}(E)E}\times_{M_n}Y(\Gamma)$ in $X_{\Gamma}$,
the localized contribution of top Chern class for
$\pi_{\bar{X}_{\Gamma}}^{\ast}
{\cal W}_{canon}\otimes {\bf H}\otimes_{\Gamma'\in 
\bar{I}_{\Gamma}}{\cal O}(-\bar{D}_{\Gamma'})$ over the strict transform of 
${\cal M}_{C-{\bf M}(E)E}\times_{M_n}Y(\Gamma)$ in $\bar{X}_{\Gamma}$,
and the localized contribution of top Chern class for
$\pi_{\hat{X}_{\Gamma}}^{\ast}
{\cal W}_{canon}\otimes {\bf H}\otimes_{\Gamma'\in 
\bar{I}_{\Gamma}}{\cal O}(-\hat{D}_{\Gamma'})$ over the strict transform of 
${\cal M}_{C-{\bf M}(E)E}\times_{M_n}Y(\Gamma)$ in $\hat{X}_{\Gamma}$
 are all equal to each other.
\end{prop}

\noindent Proof: 
 The key issue is to understand why the re-grouping of elements related by
 $\gg$ and the changing of the orderings from $\models$ to $\vdash$ do not
 affect the image cycle class of the above localized contributions of
 top Chern classes.

 By using the residual intersection formula of the top Chern class,
 i.e. proposition \ref{prop; rif}, it implies that the images in 
 ${\cal A}_{\cdot}(X)$ of the 
three localized contributions
 of the top Chern classes along $D_{\Gamma}$, $\bar{D}_{\Gamma}$, and
 $\hat{D}_{\Gamma}$ can be identified with  
the push-forwards into ${\cal A}_{\cdot}(X)$ of

 $$\hskip -.7in 
c_{rank_{\bf C}{\bf W}_{canon}}(
\pi_{\tilde{X}_{\Gamma}}^{\ast}{\bf W}_{canon}\otimes 
{\bf H}\otimes_{\Gamma'\in {I}_{\Gamma}}{\cal O}(-{D}_{\Gamma'}))
-c_{rank_{\bf C}{\bf W}_{canon}}(\pi_{\tilde{X}_{\Gamma}}^{\ast}
{\bf W}_{canon}\otimes 
{\bf H}(-D_{\Gamma})\otimes_{\Gamma'\in {I}_{\Gamma}}{\cal O}(-{D}_{\Gamma'})),$$

 $$\hskip -.7in
c_{rank_{\bf C}{\bf W}_{canon}}(\pi_{\tilde{\bar{X}}_{
\Gamma}}^{\ast}{\bf W}_{canon}\otimes 
{\bf H}\otimes_{\Gamma'\in \bar{I}_{\Gamma}}{\cal O}(-\bar{D}_{\Gamma'}))
-c_{rank_{\bf C}{\bf W}_{canon}}(
\pi_{\tilde{\bar{X}}_{\Gamma}}^{\ast}{\bf W}_{canon}\otimes 
{\bf H}(-\bar{D}_{\Gamma})\otimes_{\Gamma'\in \bar{I}_{\Gamma}}
{\cal O}(-\bar{D}_{\Gamma'})),$$

 and

$$\hskip -.7in
c_{rank_{\bf C}{\bf W}_{canon}}(\pi_{\tilde{\hat{X}}_{\Gamma}}^{\ast}
{\bf W}_{canon}\otimes 
{\bf H}\otimes_{\Gamma'\in \bar{I}_{\Gamma}}{\cal O}
(-\hat{D}_{\Gamma'}))
-c_{rank_{\bf C}{\bf W}_{canon}}(
\pi_{\tilde{\hat{X}}_{\Gamma}}^{\ast}{\bf W}_{canon}\otimes 
{\bf H}(-\hat{D}_{\Gamma})\otimes_{\Gamma'\in \bar{I}_{\Gamma}}
{\cal O}(-\hat{D}_{\Gamma'})),$$

 respectively.

 On the other hand, because $D_{\Gamma}$, $\bar{D}_{\Gamma}$ and
 $\hat{D}_{\Gamma}$ map into $X\times_{M_n}Y(\Gamma)$ under the
 projections, so the images of the localized top Chern classes in
 ${\cal A}_{\cdot}(X)$ factor through the map
 $i_{\Gamma\ast}:
{\cal A}_{\cdot}(X\times_{M_n}Y(\Gamma))\mapsto {\cal A}_{\cdot}(X)$.

 On the other hand, $j:X-X\times_{M_n}Y(\Gamma)\subset X$ is open in $X$ and by
 proposition 1.8. on page 21 of [F] we have the following exact 
sequence,

$$\hskip -.6in
{\cal A}_{\cdot}(X\times_{M_n}Y(\Gamma))\stackrel{i_{\Gamma\ast}}{\longrightarrow} 
{\cal A}_{\cdot}(X)\stackrel{j^{\ast}}{\longrightarrow}
 {\cal A}_{\cdot}(X-X\times_{M_n}Y(\Gamma))\mapsto 0.$$

 For any cycle $\beta$ in $X-X\times_{M_n}Y(\Gamma)$, its Zariski closure 
$\bar{\beta}$ in $X$
 defines a cycle in $X$. This extension of cycles defines the right inverse
 of $j^{\ast}$. Thus for each $\alpha\in {\cal A}_{\cdot}(X)$, there 
exists a unique \footnote{As $i_{\Gamma\ast}$ may not be injective, the
uniqueness of the class in ${\cal A}_{\cdot}(X\times_{M_n}Y(\Gamma))$ 
is not ensured. Nevertheless its image in ${\cal A}_{\cdot}(X)$ is.}
 cycle class $\alpha|_{X\times_{M_n}Y(\Gamma)}\in Im(i_{\Gamma\ast})$ such that
 $\alpha-\alpha|_{X\times_{M_n}Y(\Gamma)}=\overline{j^{\ast}\alpha}$. From now on
 we refer to $\alpha|_{X\times_{M_n}Y(\Gamma)}$ informally as the component of
 $\alpha$ inside the subspace $X\times_{M_n}Y(\Gamma)$.

 Because the ${\cal A}_{\cdot}(X)$ images of the 
three intersection pairings are in $Im(i_{\Gamma\ast})$,
 it suffices to consider the push-forwards of 
the three top Chern classes in the first group, 
$c_{top}(
\pi_{\tilde{X}_{\Gamma}}^{\ast}{\bf W}_{canon}\otimes 
{\bf H}\otimes_{\Gamma'\in {I}_{\Gamma}}{\cal O}(-{D}_{\Gamma'}))$,
$c_{top}(\pi_{\tilde{\bar{X}}_{
\Gamma}}^{\ast}{\bf W}_{canon}\otimes 
{\bf H}\otimes_{\Gamma'\in \bar{I}_{\Gamma}}{\cal O}(-\bar{D}_{\Gamma'}))$, 
and $c_{top}(\pi_{\tilde{\hat{X}}_{\Gamma}}^{\ast}
{\bf W}_{canon}\otimes 
{\bf H}\otimes_{\Gamma'\in \bar{I}_{\Gamma}}{\cal O}
(-\hat{D}_{\Gamma'}))$, and 
identify their components \footnote{Over here we do not intend to claim that
 the push-forwards of these top Chern classes are all 
equal in ${\cal A}_{\cdot}(X)$.
 In fact their $j^{\ast}$-restrictions in 
${\cal A}_{\cdot}(X-X\times_{M_n}Y(\Gamma))$ may be different. The
 object we really care about is the differences of the push-forwards of
 top Chern classes and the cycle classes extended from 
${\cal A}_{\cdot}(X-X\times_{M_n}Y(\Gamma)$ have to cancel out completely.}
 inside $X\times_{M_n}Y(\Gamma)$.
Then consider the push-forwards of
the second group of three top Chern classes,
$c_{top}(\pi_{\tilde{X}_{\Gamma}}^{\ast}
{\bf W}_{canon}\otimes 
{\bf H}(-D_{\Gamma})\otimes_{\Gamma'\in {I}_{\Gamma}}{\cal O}(-{D}_{\Gamma'}))$,
$c_{top}(
\pi_{\tilde{\bar{X}}_{\Gamma}}^{\ast}{\bf W}_{canon}\otimes 
{\bf H}(-\bar{D}_{\Gamma})\otimes_{\Gamma'\in \bar{I}_{\Gamma}}
{\cal O}(-\bar{D}_{\Gamma'}))$, and $c_{top}(
\pi_{\tilde{\hat{X}}_{\Gamma}}^{\ast}{\bf W}_{canon}\otimes 
{\bf H}(-\hat{D}_{\Gamma})\otimes_{\Gamma'\in \bar{I}_{\Gamma}}
{\cal O}(-\hat{D}_{\Gamma'}))$ in ${\cal A}_{\cdot}(X\times_{M_n}Y(\Gamma)$, and
identify their components inside $X\times_{M_n}Y(\Gamma)$.

\bigskip

\noindent $\diamondsuit$ Case I:
 We identify the ${\cal A}_{\cdot}(X\times_{M_n}Y(\Gamma))$-components
of the push-forwards of the 
top Chern classes within the first group. The identification of
 the second group is similar, and will be handled in Case II below.
 Introduce three divisors \footnote{We introduce the new notations $X_a, D_a$
to avoid writing parallel formulae repeatedly!} 
$D_1=(\cup_{\Gamma'\in I_{\Gamma}}D_{\Gamma'})
\subset X_1=X_{\Gamma}$,
$D_2=(\cup_{\Gamma'\in \bar{I}_{\Gamma}}\bar{D}_{\Gamma'})
\subset X_2=\bar{X}_{\Gamma}$
 and $D_3=(\cup_{\Gamma'\in \bar{I}_{\Gamma}}\hat{D}_{\Gamma'})
\subset X_3=\hat{X}_{\Gamma}$.
 The restriction of the vectors bundles
 involved in the first group are pull-back from $X_1, X_2$, and
$X_3$, respectively. We define $\pi_{X_a}:X_a\mapsto M_n\times T(M)$. 
 The yet-to-be-identified classes are 
 $c_{top}(\pi_{X_i}^{\ast}{\bf W}_{canon}\otimes {\bf H}\otimes {\cal O}(-D_a))$
 for $a=1, 2, 3$. 
By applying the residual intersection formula, i.e. proposition \ref{prop; rif}, 
to the $D_1, D_2$ and $D_3$, these top Chern classes can be
 re-written as 

$$c_{top}(\pi_{X_a}^{\ast}{\bf W}_{canon}\otimes {\bf H})\cap [X_a]-
\sum_{1\leq i\leq rank_{\bf C}{\bf W}_{canon}}(-1)^{i-1}
c_{rank_{\bf C}{\bf W}_{canon}-i}(\pi_{X_a}^{\ast}
{\bf W}_{canon}\otimes{\bf H})(D_a)^{i-1}[D_a],$$

 for $a=1, 2, 3$.

 The image of the first term into ${\cal A}_{\cdot}(X)$
 is apparently independent of $a$ and has a unique component in 
 $X\times_{M_n}Y(\Gamma)$. So it suffices to show that the
 components in ${\cal A}_{\cdot}(X\times_{M_n}Y(\Gamma))$ 
 of the push-forwards of the second terms involving $D_a$ are $a$-independent.  

 Firstly recall that when $D_a$ is a divisor, the total Chern class 
$c_{total}({\cal O}(D_a))=1+D_a$ and the total Segre class
 $s_{total}({\cal O}(D_a))=1+\sum_{j\geq 1}(-1)^j D_a^j$ is the total Segre class
 of the normal cone
 $s_{total}(D_a, X_a)$. 
So for $a=1, 2, 3$, the push-forwards of the
 second terms can be re-expressed as the 
push-forward of 

 $$\eta_a=\{c_{total}(\bigl(\pi_{X_a}^{\ast}{\bf W}_{canon}\otimes 
{\bf H}\bigr)|_{D_a})\cap
 s_{total}(D_a, X_a)\}_{dim_{\bf C}X-rank_{\bf C}{\bf W}_{canon}}$$

 into ${\cal A}_{\cdot}(X)$.

 Define $h_a:X_a\mapsto X$ 
to be the blowing down projection maps for $a=1, 2, 3$.

 Now recall the following proposition 4.2.(a) on page 74 of [F].

\begin{prop}\label{prop; birational}(Fulton)
Let $f:Y'\mapsto Y$ be a morphism of pure-dimensional scheme,
 $Z\subset Y$ a closed sub-scheme, $Z'=f^{-1}(Z)$ the inverse image, 
 $g':Z'\mapsto Z$ the induced morphism.

Suppose that $f$ is proper, $Y$ irreducible and $f$ maps each irreducible
component of $Y'$ onto $Y$, then 

$$g_{\ast}(s(Z', Y'))=deg(Y'/Y)\cdot s(Z, Y).$$
\end{prop}

 In our context, the blowing down map $h_a:X_a\mapsto X$ 
is proper, and $deg(X_a/X)=1$ (because they are birational). 
 The sub-scheme $Z=h_a(D_a)\times_{M_n}Y(\Gamma)$. Because
 $h_a$ are composite blowing down maps, $D_a\times_{M_n}Y(\Gamma)
=h_a^{-1}(h_a(D_a)\times_{M_n}Y(\Gamma))$.
 Both of $X_a$ and $X$ are irreducible and $h_a$ maps $X_a$ onto 
$X$. 

 By applying proposition \ref{prop; birational}, for all $1\leq a\leq 3$ 
the push-forward-images of the
 above classes $\eta_a$ in ${\cal A}_{\cdot}(h_a(D_a))$ are equal to
 $$\{c_{total}(\bigl(\pi_X^{\ast}{\bf W}_{canon}\otimes {\bf H}\bigr)|_{h_a(D_a)})
\cap s(h_a(D_a), X)\}_{dim_{\bf C}X-rank_{\bf C}{\bf W}_{canon}}.$$

Then the fact that the $Im(i_{\Gamma\ast})$ components 
of their push-forwards into ${\cal A}_{\cdot}(X)$ are equal follows from the
 following observation,

\begin{lemm}\label{lemm; thesame}
Let $h_a:X_a\mapsto X$ and $D_a$, $1\leq a\leq 3$ be as described above. 
Let the family moduli space of $C-{\bf M}(E)E$, 
${\cal M}_{C-{\bf M}(E)E}\subset X$, denote the sub-scheme defined by
 the canonical section $s_{canon}$ of 
$\pi_X^{\ast}{\bf W}_{canon}\otimes {\bf H}$. Then
for $a=1, 2, 3$, the sub-schemes $h_a(D_a)\cap \bigl(X\times_{M_n}Y(\Gamma)\bigr)=
h_a(D_a)\times_{M_n}Y(\Gamma)\subset X$ all coincide and are all equal
 to the finite union 
$\cup_{\Gamma'\in I_{\Gamma}}{\cal M}_{C-{\bf M}(E)E}\times_{M_n}\bigl(
Y(\Gamma')\cap Y(\Gamma)\bigr)$.
\end{lemm}
 
\noindent Proof: 
 It is easy to see that the change of the linear ordering from $\models$ to
 $\vdash$ in $\bar{I}_{\Gamma}$ does not affect the total locus which is blown up.
 Thus we know that $h_2(D_2)=h_3(D_3)=\cup_{\Gamma'\in \bar{I}_{\Gamma}}
 {\cal M}_{C-{\bf M}(E)E}\times_{M_n}Y(\Gamma')$. So their intersections with
 $X\times_{M_n}Y(\Gamma)$ are equal.
 
 On the other hand, to argue that $h_1(D_1)\times_{M_n}Y(\Gamma)=h_2(D_2)
\times_{M_n}Y(\Gamma)=
\cup_{\Gamma'\in \bar{I}_{\Gamma}}
 {\cal M}_{C-{\bf M}(E)E}\times_{M_n}\bigl(
Y(\Gamma')\cap Y(\Gamma)\bigr)$, it suffices to show that for all
 $\Gamma'\in I_{\Gamma}-\bar{I}_{\Gamma}$, the corresponding sub-scheme
 ${\cal M}_{C-{\bf M}(E)E}\times_{M_n}\bigl(Y(\Gamma')\cap Y(\Gamma)\bigr)$ has
 been included in the union of closed sub-schemes $\cup_{\Gamma''\in 
\bar{I}_{\Gamma}}{\cal M}_{C-{\bf M}(E)E}\times_{M_n}\bigl(Y(\Gamma'')\cap
 Y(\Gamma)\bigr)$ already.

 We may assume that $Y(\Gamma')\cap Y(\Gamma)\not=\emptyset$ or the
 statement is trivial to prove. By lemma \ref{lemm; common}, we know that there are
 three exclusive possibilities (a). $\Gamma'\succ \Gamma$, (b). 
 $\Gamma\succ \Gamma'$, (c). $\exists \Gamma''\in \Delta(n)$ such that 
 both of $\Gamma, \Gamma'\succ \Gamma''$.

Suppose that $S_{\Gamma'}\cap Y(\Gamma)=\emptyset$. We argue that we may 
replace $\Gamma'$ by some $\Gamma''$ with $S_{\Gamma''}\cap 
Y(\Gamma)\not=\emptyset$. We already have $S_{\Gamma'}\cap Y(\Gamma)=\emptyset$
 by our assumption.
 We know that $Y(\Gamma')\cap S_{\Gamma}=\emptyset$, too.
 If not, the hypothesis $Y(\Gamma')\cap S_{\Gamma}\not=\emptyset$ and lemma
 \ref{lemm; suff}
 imply $\Gamma'\succ \Gamma$ and therefore $\Gamma'\models \Gamma$.
 Then such a $\Gamma'$ cannot be in the index set
 $I_{\Gamma}$ at all. As both (a). and
 (b). fail, it falls into the
 situation (c). that $S_{\Gamma}\cap Y(\Gamma')=S_{\Gamma'}\cap 
Y(\Gamma)=\emptyset$. From the proof of lemma \ref{lemm; common} we know that
 for ``all'' $b\in Y(\Gamma)\cap Y(\Gamma')$, there exists a $\Gamma''\in
 \Delta(n)$ such that $b\in S_{\Gamma''}\cap Y(\Gamma)\not=\emptyset$.

\medskip

  From this digestion
 we learn that $\cup_{\Gamma'\in I_{\Gamma}-\bar{I}_{\Gamma}}
{\cal M}_{C-{\bf M}(E)E}\times_{M_n}\bigl(Y(\Gamma')\cap Y(\Gamma)\bigr)$
 can be replaced by the union
 $\cup_{\Gamma'\in I_{\Gamma}; S_{\Gamma'}\cap Y(\Gamma)\not=
\emptyset}
{\cal M}_{C-{\bf M}(E)E}\times_{M_n}\bigl(S_{\Gamma'}\cap Y(\Gamma)\bigr)$.
 
 According to definition \ref{defin; reduce} on the reduced index set 
$\bar{I}_{\Gamma}$, any $\Gamma'\in I_{\Gamma}$ with 
$S_{\Gamma'}\cap Y(\Gamma)\not=\emptyset$
 is thrown away to form $\bar{I}_{\Gamma}$ exactly when there exists another
 $\Gamma''\in \bar{I}_{\Gamma}$ with 
 $(\Gamma'', \sum_{e_i''\cdot (C-{\bf M}(E)E)<0}e_i'')
\gg (\Gamma', \sum_{e_i'\cdot (C-{\bf M}(E)E)<0}e_i')$ 
and $(\Gamma, \sum_{e_i\cdot (C-{\bf M}(E)E)<0}e_i)
\sqsupset (\Gamma'', \sum_{e_i''\cdot (C-{\bf M}(E)E)<0}e_i'')$.

 However the $\gg$ relationship (see definition \ref{defin; partial}) 
between $\Gamma''$ and $\Gamma'$ implies $Y(\Gamma'')\supset Y(\Gamma')$.
 This implies that ${\cal M}_{C-{\bf M}(E)E}\times_{M_n}S_{\Gamma'}
\subset {\cal M}_{C-{\bf M}(E)E}\times_{M_n}Y(\Gamma'')$. So we have the inclusion 

$$\hskip -.4in
\cup_{\Gamma'\in I_{\Gamma}; S_{\Gamma'}\cap Y(\Gamma)\not=
\emptyset}
{\cal M}_{C-{\bf M}(E)E}\times_{M_n}\bigl(S_{\Gamma'}\cap Y(\Gamma)\bigr)
\subset \cup_{\Gamma''\in \bar{I}_{\Gamma}}
{\cal M}_{C-{\bf M}(E)E}\times_{M_n}\bigl(Y(\Gamma'')\cap Y(\Gamma)\bigr).$$

By combining the inclusions we know that $h_1(D_1)\times_{M_n}Y(\Gamma)$
 must be included in
 $h_2(D_2)\times_{M_n}Y(\Gamma)$. 
But the reversed inclusion $h_2(D_2)\times_{M_n}Y(\Gamma)
\subset h_1(D_1)\times_{M_n}Y(\Gamma)$ is apparent.
 So we have $h_1(D_1)\times_{M_n}Y(\Gamma)=h_2(D_2)\times_{M_n}Y(\Gamma)=
h_3(D_3)\times_{M_n}Y(\Gamma)$ as sub-schemes of $X$ 
and the lemma is proved. $\Box$

  As usual let $i_{h_a(D_a)}:h_a(D_a)\hookrightarrow X$ be the inclusion maps.
  The above lemma tells us that the restriction of $h_a(D_a)$ to 
 $X\times_{M_n}Y(\Gamma)$ coincide. On the other hand the normal cones
 ${\bf C}_{h_a(D_a)}X$ can always be written as the unions of irreducible
 normal cones supporting over irreducible components of the 
sub-schemes $h_a(D_a)$. By separating
the irreducible components of $h_a(D_a)$ in $X\times_{M_n}Y(\Gamma)$ and
 the zariski-closures of $h_a(D_a)\cap (X-X\times_{M_n}Y(\Gamma))$, we
 may write each ${\bf C}_{h_a(D_a)}X={\bf C}_{h_a(D_a)\times_{M_n}Y(\Gamma)}X
\cup {\bf C}_a'$, with
 ${\bf C}_a'$ supported over the zariski-closure of 
$h_a(D_a)\cap (X-X\times_{M_n}Y(\Gamma))=h_a(D_a)\times_{M_n}(M_n-Y(\Gamma))$
 in $X$.
 Then it is easy to see that
 the components along $X\times_{M_n}Y(\Gamma)$ of 
$i_{h_a(D_a)\ast}\{s({\bf C}_{h_a(D_a)}X)\}=
i_{h_a(D_a)\ast}\{s({\bf C}_{h_a(D_a)\times_{M_n}Y(\Gamma)}X)\}
+i_{h_a(D_a)\ast}\{s({\bf C}_a')\}$ is exactly
 $i_{h_a(D_a)\ast}\{s({\bf C}_{h_a(D_a)\times_{M_n}Y(\Gamma)}X)\}$ while 
 $i_{h_a(D_a)\ast}\{s({\bf C}_a')\}$ are the extension 
 ${\cal A}_{\cdot}(X-X\times_{M_n}Y(\Gamma))\mapsto {\cal A}_{\cdot}(X)$)
 of the restricted Segre class
 $i_{h_a(D_a)\ast}\{s(h_a(D_a)\times_{M_n}(M_n-Y(\Gamma)),
 X-X\times_{M_n}Y(\Gamma))\}$.

 By the above lemma, we know that the $X\times_{M_n}Y(\Gamma)$-components
 of the total Segre classes 
$i_{h_a(D_a)\ast}s(h_a(D_a), X)$ coincide for all $1\leq a\leq 3$.  By capping with
 $c_{total}(\pi_X^{\ast}{\bf W}_{canon}\otimes {\bf H})$, we conclude that
 the components in $X\times_{M_n}Y(\Gamma)$ of 
$i_{h_a(D_a)\ast}
\{c_{total}(\bigl(\pi_X^{\ast}{\bf W}_{canon}\otimes {\bf H}\bigr))
\cap s(h_a(D_a), X)\}_{dim_{\bf C}X-rank_{\bf C}{\bf W}_{canon}}$ are all the 
same. So Case I is proved.

\medskip

\noindent $\diamondsuit$ Case II: \label{case2}
 The identifications of the components in $X\times_{M_n}Y(\Gamma)$ of the 
second groups of three top Chern classes are rather
 parallel to the previous argument in Case I, with some minute difference.
 We define $X_1'=\tilde{X}_{\Gamma}$, $X_2'=\tilde{\bar{X}}_{\Gamma}$ and
 $X_3'=\tilde{\hat{X}}_{\Gamma}$. Define $h_a':X_a'\mapsto X$ to be the
 projection maps.

 Then the total transformations of
 $D_a\subset X_a$ under the pull-backs of the 
blowing down maps $X_a'\mapsto X_a$ define
 Cartier divisors in $X_a'$ and we skip the pull-back notations and denote them
 by the same symbols \footnote{This is consistent with our earlier
 convention.} $D_a$. Set 
 $D_1'=D_1\cup D_{\Gamma}$, $D_2'=D_2\cup\bar{D}_{\Gamma}$ and
 $D_3'=D_3\cup \hat{D}_{\Gamma}$.

 Then we show that the push-forward of 
$c_{top}(\pi_{X_a}^{\ast}{\bf W}_{canon}\otimes {\cal O}(-D_a'))$ to 
 ${\cal A}_{\cdot}(X)$  
have identical components in ${\cal A}_{\cdot}(X\times_{M_n}Y(\Gamma))$ for
 $a=1, 2, 3$.

 Following the previous convention let $i_{h_a(D_a')}:h_a(D_a')\hookrightarrow X$
 be the inclusions into $X$.

\begin{lemm}\label{lemm; iden2}
 Let $D_1'=D_1\cup D_{\Gamma}$, $D_2'=D_2\cup \bar{D}_{\Gamma}$ and
 $D_3'=D_3\cup \hat{D}_{\Gamma}$ be defined above. Then
for $a=1, 2, 3$, the components in $X\times_{M_n}Y(\Gamma)$ of 
the push-forwarded Segre classes $i_{h_a'(D_a')\ast}s(h_a'(D_a'), X)$ are 
all equal.
\end{lemm}

\noindent Proof: From the argument in Case I, we know that we only need to
 prove that $h_a'(D_a')\times_{M_n}Y(\Gamma)=h_a(D_a')\cap \bigl(
X\times_{M_n}Y(\Gamma)\bigr)$
 are all equal.
We notice that for all three $'a'$ we have
 $h_1'(D_1')=h_1'(D_1)\cup h_1'(D_{\Gamma})$, $h_2'(D_2')=h_2'(D_2)\cup
 h_2'(\bar{D}_{\Gamma})$, and $h_3'(D_3')=h_3'(D_3)\cup h_3(\hat{D}_{\Gamma})$.

 Firstly we notice that $h_a'(D_a)$ are nothing but the $h_a(D_a)$ in Case I.
 On the other hand,
 despite that $D_{\Gamma}$, $\bar{D}_{\Gamma}$, $\hat{D}_{\Gamma}$ are
 different exceptional divisors blown up from the strict transformations of
 $Z(s_{canon})\times_{M_n}Y(\Gamma)$ in three mutually birational spaces 
$X_{\Gamma}, \bar{X}_{\Gamma}, \hat{X}_{\Gamma}$, their images under
 $h_1'$, $h_2'$ and $h_3'$ are identical and their common image is
 $Z(s_{canon})\times_{M_n}Y(\Gamma)=
{\cal M}_{C-{\bf M}(E)E}\times_{M_n}Y(\Gamma)$. So by combining these
 conclusions we have
 $h_a'(D_a')=h_a(D_a)\cup \bigl(Z(s_{canon})\times_{M_n}Y(\Gamma)\bigr)$.

 Thus $h_a'(D_a')\times_{M_n}Y(\Gamma)$ are nothing but
 $h_a(D_a)\times_{M_n}Y(\Gamma)\cup \bigl(Z(s_{canon})\times_{M_n}Y(\Gamma)\bigr)$.
 By lemma \ref{lemm; thesame} the sub-schemes 
$h_a(D_a)\times_{M_n}Y(\Gamma)$ has been shown to 
 be $a-$independent, so we conclude that $h_a'(D_a')\times_{M_n}Y(\Gamma)$ are
 $a-$independent as well. $\Box$

\medskip

  Once we identify the $X\times_{M_n}Y(\Gamma)$ components of
 their Segre classes, the rest of the proof
 is almost identical to Case I. We omit the details. The proof of
 proposition \ref{prop; identical} is finished. 
$\Box$

\medskip

\subsection{The Identification of the Localized Contribution
 with the Modified Algebraic Family Invariant}\label{subsection; identifyMA}

\bigskip
 
 In this subsection, we proceed to identify the integral of
 the intersection pairing of the localized contribution of top Chern class
 with the modified algebraic family Seiberg-Witten invariant defined in
 section \ref{subsection; modinv}.

 The push-forward into ${\cal A}_{\cdot}(X)$ of the 
localized contribution of the top Chern class defines a cycle class of grade
 $dim_{\bf C}X-rank_{\bf C}{\bf W}_{canon}$. In order to get a numerical
 invariant $\in {\bf Z}$, 
we can either pair it with the suitable power of the tautological
 class $c_1({\bf H})^{dim_{\bf C}X-
rank_{\bf C}{\bf W}_{canon}}$ on the projective space bundle 
$X={\bf P}({\bf V}_{canon})$ and push-forward the resulting class into 
${\cal A}_0(pt)\cong {\bf Z}$, 
or we may fix a point $t_L\in T(M)$ and
 pair the push-forward of the localized contribution of top Chern class
 with $c_1({\bf H})^{dim_{\bf C}X-
rank_{\bf C}{\bf W}_{canon}-q}\cap [t_L]$ and then
 \footnote{It depends on whether we counts curves in the non-linear or 
linear systems.} push it forward
 into ${\cal A}_0(pt)$. 
Over here $[t_L]$ represents
 the zero dimensional cycle class of the point $t_L$ and the integer 
$q=dim_{\bf C}T(M)$ denotes the irregularity
 of the algebraic surface.

Now we are ready to identify the yet-to-be-enumerated intersection pairing 
involving the localized contribution of the
top Chern class,
 
$$\hskip -1.3in
\sum_{1\leq i\leq rank_{\bf C}{\bf W}_{canon}}(-1)^{i-1}
c_{rank_{\bf C}{\bf W}_{canon}-i}(\pi_{\tilde{X}}^{\ast}{\bf W}_{canon}\otimes 
{\bf H} \otimes_{\Gamma'\in I_{\Gamma}} {\cal O}(-D_{\Gamma'})|_{D_{\Gamma}})
D_{\Gamma}^{i-1}[D_{\Gamma}]\cap c_1({\bf H})^{dim_{\bf C}X-
rank_{\bf C}{\bf W}_{canon}}$$

with the
 modified mixed algebraic family Seiberg-Witten 
invariant $${\cal AFSW}_{M_{n+1}\times_{M_n}Y(\Gamma)\times
T(M)\mapsto Y(\Gamma)\times T(M)}^{\ast}(c_{total}(\tau_{\Gamma}), 
C-{\bf M}(E)E-\sum_{e_i\cdot (C-{\bf M}(E)E)<0}e_i).$$

 And identify 

$$\hskip -1.3in
\sum_{1\leq i\leq rank_{\bf C}{\bf W}_{canon}}(-1)^{i-1}
c_{rank_{\bf C}{\bf W}_{canon}-i}(\pi_{\tilde{X}}^{\ast}{\bf W}_{canon}\otimes 
{\bf H}\otimes_{\Gamma'\in I_{\Gamma}}
{\cal O}(-D_{\Gamma'})|_{D_{\Gamma}})
D_{\Gamma}^{i-1}[D_{\Gamma}]\cap c_1({\bf H})^{dim_{\bf C}X-
rank_{\bf C}{\bf W}_{canon}-q}\cap [t_L]$$ with the $T(M)$-restricted version
 of modified mixed algebraic family Seiberg-Witten invariant 
${\cal AFSW}_{M_{n+1}\times_{M_n}Y(\Gamma)\times
\{t_L\}\mapsto Y(\Gamma)\times \{t_L\}}^{\ast}(c_{total}(\tau_{\Gamma}), 
C-{\bf M}(E)E-\sum_{e_i\cdot (C-{\bf M}(E)E)<0}e_i)$.

Because the identification of the latter objects is completely identical to the
 identification of the former, if we 
replace $c_1({\bf H})^{dim_{\bf C}X-rank_{\bf C}
{\bf W}_{canon}-q}\cap [t_L]$ by $c_1({\bf H})^{dim_{\bf C}X-rank_{\bf C}
{\bf W}_{canon}}$, we will discuss only the former case in the proof.

 Please consult subsection \ref{subsection; modinv} for the definitions of the
 modified algebraic family Seiberg-Witten invariants and the
 construction of $\tau_{\Gamma}$.

\medskip

 The main tools we will adopt are the machineries developed in subsection
 \ref{subsection; lf} (proposition \ref{prop; id})
 and section \ref{section; discrepancy} (proposition \ref{prop; rif}
, \ref{prop; equivalent}, \ref{prop; cap} and lemma \ref{lemm; local}).

\medskip

\begin{prop}\label{prop; =mo}
 Given an $n$-vertex admissible graph $\Gamma\in \Delta(n)\subset adm(n)$, 
 the integration into ${\cal A}_0(pt)\cong {\bf Z}
$ of $c_1({\bf H})^{dim_{\bf C}X-rank_{\bf C}
{\bf W}_{canon}}$ capping with the push-forward
 of the localized contribution of the top Chern class along the 
blown up divisor $D_{\Gamma}$, 
$$\sum_{1\leq i\leq rank_{\bf C}
{\bf W}_{canon}}(-1)^{i-1}
c_{rank_{\bf C}{\bf W}_{canon}-i}(\pi_{\tilde{X}}^{\ast}{\bf W}_{canon}\otimes 
{\bf H}\otimes_{\Gamma'\in I_{\Gamma}}{\cal O}(-D_{\Gamma'})|_{D_{\Gamma}})
D_{\Gamma}^{i-1}[D_{\Gamma}]$$

 in ${\cal A}_{\cdot}(X)$, 
 is equal to the modified mixed algebraic family invariant
${\cal AFSW}_{M_{n+1}\times_{M_n}Y(\Gamma)\times
T(M)\mapsto Y(\Gamma)\times T(M)}^{\ast}(c_{total}(\tau_{\Gamma}), 
C-{\bf M}(E)E-\sum_{e_i\cdot (C-{\bf M}(E)E)<0}e_i)$ defined in subsection 
\ref{subsection; modinv}.

 Likewise for an arbitrary point $t_L\in T(M)$, 
the integral of the pairing of 
$c_1({\bf H})^{dim_{\bf C}X-rank_{\bf C}{\bf W}_{canon}-q}\cap [t_L]$
 with the image of the above localized contribution of top Chern class along 
$D_{\Gamma}$ into ${\cal A}_{\cdot}(X)$,
 is equal to the modified mixed algebraic family invariant
$${\cal AFSW}_{M_{n+1}\times_{M_n}Y(\Gamma)\times
\{t_L\}\mapsto Y(\Gamma)\times \{t_L\}}^{\ast}(c_{total}(\tau_{\Gamma}), 
C-{\bf M}(E)E-\sum_{e_i\cdot (C-{\bf M}(E)E)<0}e_i).$$
\end{prop}

\medskip

 Because the $t_L$-restricted version is completely parallel to the non-restricted
 version, we only offer a proof for the non-restricted version. Please 
consult remark \ref{rem; restrict} on page \pageref{rem; restrict}
 right after the end of the proof.

\medskip

\noindent Proof of proposition \ref{prop; =mo}:\label{proof=mo}
 The proof of the proposition 
involves an induction on the element $\Gamma\in \Delta(n)$ based on the linear
 ordering $\models$ (see page \pageref{models} for the recursive 
 definition of $\models$).

  Firstly, we provide a simple computation on the dimension formula which 
motivates the
 appearance of $\tau_{\Gamma}$ in the modified family invariant.
  Because $X={\bf P}({\bf V}_{canon})\mapsto M_n\times T(M)$, $dim_{\bf C}X=
rank_{\bf C}{\bf V}_{canon}-1+dim_{\bf C}M_n+q$. Thus, 
$dim_{\bf C}X-rank_{\bf C}{\bf W}_{canon}=dim_{\bf C}M_n+q+
rank_{\bf C}({\bf V}_{canon}-{\bf W}_{canon})-1$. Based on the fact that 
 $(\Phi_{{\bf V}_{canon}{\bf W}_{canon}}, {\bf V}_{canon}, {\bf W}_{canon})$ is the
 canonical algebraic family Kuranishi model of the class $C-{\bf M}(E)E$, we know
 that 

$$rank_{\bf C}({\bf V}_{canon}-{\bf W}_{canon})=1-q+p_g+
{(C-{\bf M}(E)E)^2-c_1({\bf K}_{M_{n+1}/M_n})
\cdot (C-{\bf M}(E)E)\over 2}$$
$$=1-q+p_g+{C^2-c_1({\bf K}_M)\cdot C-\sum_{1\leq i\leq n}(m_i^2+m_i)\over 2},$$ 

by surface Riemann-Roch formula. From this we can infer the relationship between
 the raised power of $c_1({\bf H})$ in the intersection pairing (which is also the
 expected algebraic family dimension of ${\cal M}_{C-{\bf M}(E)E}$)
  and the singular multiplicities $m_i, 1\leq i\leq n$.

On the other hand, by a direct computation 
the expected (family algebraic) dimension of the space 
 ${\cal M}_{C-{\bf M}(E)E-\sum_{e_i\cdot (C-{\bf M}(E)E)<0}e_i}
\times_{M_n}Y(\Gamma)$ is
 given by 

$$
p_g+dim_{\bf C}Y(\Gamma)+{(C-{\bf M}(E)E)^2-c_1({\bf K}_{M_{n+1}/M_n})
\cdot (C-{\bf M}(E)E)
\over 2}-\sum_{e_i\cdot (C-{\bf M}(E)E)<0}(C-{\bf M}(E)E)\cdot e_i$$
$$+{(\sum_{e_i\cdot
(C-{\bf M}(E)E)<0}e_i)^2+c_1({\bf K}_{M_{n+1}/M_n})\cdot (\sum_{e_i\cdot
(C-{\bf M}(E)E)<0}e_i) \over 2}.$$

 By using $dim_{\bf C}Y(\Gamma)=dim_{\bf C}M_n+\sum_{e_i\cdot (C-{\bf M}(E)E)<0}
{e_i^2-c_1({\bf K}_{M_{n+1}/M_n})\cdot e_i\over 2}$, the above expression can be
 simplified to 

$$
=p_g+dim_{\bf C}M_n+{(C-{\bf M}(E)E)^2-c_1({\bf K}_{M_{n+1}/M_n})
\cdot (C-{\bf M}(E)E)
\over 2}$$
$$+\{\sum_{e_i\cdot (C-{\bf M}(E)E)<0}e_i\cdot \bigl(e_i+\sum_{j<i; e_j\cdot 
(C-{\bf M}(E)E)<0}e_j-(C-{\bf M}(E)E)\bigr)\}.$$

  A direct 
 comparison with the formula of $rank_{\bf C}\tau_{\Gamma}$ 
shows that this correction term
 matches up with the rank of $\tau_{\Gamma}\otimes {\bf H}$ found in subsection 
\ref{subsection; modinv} lemma \ref{lemm; free}. This explains morally why we need 
to insert $c_{top}(\tau_{\Gamma}\otimes
{\bf H})=\sum_{l\leq rank_{\bf C}\tau_{\Gamma}}c_l(\tau_{\Gamma})\cap 
c_1({\bf H})^{rank_{\bf C}\tau_{\Gamma}-l}$ in the corresponding modified 
algebraic 
family invariant. The dimension count singles out the role of $\tau_{\Gamma}$-as
 a mean to
 {\bf compensate the discrepancy of the 
expected family dimensions} between ${\cal M}_{C-{\bf M}(E)E}$ 
 and ${\cal M}_{C-{\bf M}(E)E-\sum_{e_i\cdot (C-{\bf M}(E)E)<0}e_i}
\times_{M_n} Y(\Gamma)$.
In the latter half of the proof, we will see why a correct choice of 
$\tau_{\Gamma}$ (not only the rank itself) is essential in our identification.

\bigskip

We start from the
 simplest case when $\Gamma\in \Delta(n)$ is a minimal element of $\Delta(n)$ 
under $\models$. Under this assumption, $I_{\Gamma}=\emptyset$ and 
$\Gamma$ is a minimal element under $\succ$.
The minimality assumption of $\Gamma$ under $\succ$ implies that there can be no
 $\Gamma'$ with $S_{\Gamma'}\cap Y(\Gamma)\not=\emptyset$. In such a case
the space $S_{\Gamma}\subset Y(\Gamma)$ (over which the type $I$ exceptional cone 
${\cal C}_{\Gamma}$ is constant) itself is 
a closed subset of $M_n$, and therefore is equal to $Y(\Gamma)$. The consequence 
 $S_{\Gamma}=Y(\Gamma)$ implies that all the type $I$ exceptional curves 
dual to $e_i$, with $e_i\cdot (C-{\bf M}(E)E)<0$, remain smooth and 
irreducible throughout the whole
 $Y(\Gamma)$. In particular, no curves dual to such $e_i$ can break into more than 
one irreducible component over $Y(\Gamma)$.

\medskip

 By lemma \ref{lemm; local}, the sum \footnote{It is simplified as
 there is no blowing up ahead of the one parametrized by $\Gamma$.}
 $\sum_{1\leq i\leq rank_{\bf C}{\bf W}_{canon}}(-1)^{i-1}
c_{rank_{\bf C}{\bf W}_{canon}-i}(\pi_{\tilde{X}}^{\ast}{\bf W}_{canon}\otimes {\bf H}) 
D_{\Gamma}^{i-1}[D_{\Gamma}]$ is nothing but the localized contribution of the 
 top Chern class defined in section 6 of [Liu5].
 By proposition 11 of [Liu5] and our knowledge that the type $I$ exceptional curves dual 
 to $e_i$, with $e_i\cdot (C-{\bf M}(E)E)<0$
 remain irreducible and smooth throughout $Y(\Gamma)$, these imply that
 the natural bundle map $\pi_X^{\ast}{\bf W}_{canon}^{\circ}\otimes 
{\bf H}|_{X\times_{M_n}Y(\Gamma)\times T(M)}\mapsto 
\pi_X^{\ast}{\bf W}_{canon}\otimes {\bf H}|_{X\times_{M_n}Y(\Gamma)\times T(M)}$ 
 (see proposition 9 of [Liu5] for its construction) 
is injective over the whole $X\times_{M_n}Y(\Gamma)\times T(M)$. In terms of the 
notation of proposition \ref{prop; section} of the current paper or 
proposition 12/corollary 3 of [Liu5],
 the union of cones $\cup_{i>0}{\bf C}_{\rho_i}$ corresponding the kernel of the bundle map 
is empty. Thus the {\bf simplifying assumption} in section 6.1 of [Liu5]
 has been satisfied automatically because the restricted family 
moduli space ${\cal M}_{C-{\bf M}(E)E-\sum_{e_i\cdot (C-{\bf M}(E)E)<0}e_i}
\times_{M_n}Y(\Gamma)=Z(s_{canon}^{\circ})\times_{M_n}Y(\Gamma)$ does not intersect
 with $\cup_{i>0}{\bf C}_{\rho_i}=\emptyset$ at all.
 The argument of theorem 4 of [Liu5] is then 
applicable and we may identify the integration of the top intersection 
pairing of the 
localized top Chern class along $D_{\Gamma}$ 
(over $Y(\Gamma)$) and 
$c_1({\bf H})^{dim_{\bf C}M_n+rank_{\bf C}({\bf V}_{canon}-{\bf W}_{canon})+q-1}$
to \footnote{We 
skip the push-forward operation into ${\cal A}_{\cdot}(X)$ on the localized
 top Chern class by interpreting the cap product with the complementary 
power of $c_1({\bf H})$ as capping this natural Chern class pull-back by 
 $D_{\Gamma}\mapsto X$. To simplify our notations, we will adopt the same 
convention afterward. The reader should be able to recover it from the context.}
 be ${\cal AFSW}_{M_{n+1}\times T(M)\times_{M_n}Y(\Gamma)
\mapsto Y(\Gamma)\times  T(M)}(c_{total}(\tau_{\Gamma}), 
C-{\bf M}(E)E-\sum_{e_i\cdot (C-{\bf M}(E)E)<0}e_i)$, which is
 nothing but the modified invariant ${\cal AFSW}_{M_{n+1}\times T(M)\times_{M_n}
Y(\Gamma)\mapsto Y(\Gamma)\times T(M)}^{\ast}(c_{total}(\tau_{\Gamma}), 
C-{\bf M}(E)E-\sum_{e_i\cdot (C-{\bf M}(E)E)<0}e_i)$ by definition 
\ref{defin; initial}.
 Consult section 6.1-6.4 of [Liu5] for the details of the identification.
\footnote{In the following inductive argument,
a specialization of our argument for the general case also provides a proof for the special case.}

\medskip

  Next we consider the general (and a priori more complicated) situation when
 $\Gamma$ is not minimal under $\models$. 

\medskip

\noindent {\bf Induction Hypothesis}: Assuming that for all the $\Gamma'\in I_{\Gamma}$ (i.e. $\Gamma\models \Gamma'$), the 
 integral of the \label{IH}
following top intersection pairing with localized contribution of top Chern class

$$
\sum_{1\leq i\leq rank_{\bf C}{\bf W}_{canon}}(-1)^{i-1}
c_{rank_{\bf C}{\bf W}_{canon}-i}(\pi_{\tilde{X}}^{\ast}{\bf W}_{canon}\otimes {\bf H}
\otimes_{\Gamma''\in I_{\Gamma'}}
{\cal O}(-D_{\Gamma''})|_{D_{\Gamma'}})$$
$$\cap D_{\Gamma'}^{i-1}[D_{\Gamma'}]\cap c_1({\bf H})^{dim_{\bf C}M_n+
rank_{\bf C}({\bf V}_{canon}-{\bf W}_{canon})+q-1}$$

 have been identified with the modified algebraic family invariant,

$${\cal AFSW}_{M_{n+1}\times_{M_n}Y(\Gamma')\times
T(M)\mapsto Y(\Gamma')\times T(M)}^{\ast}(c_{total}(\tau_{\Gamma'}), 
C-{\bf M}(E)E-\sum_{e_i'\cdot (C-{\bf M}(E)E)<0}e_i').$$

 As usual 
$e_i'$ are the type $I$ exceptional classes over $Y_{\Gamma'}$ and $\tau_{\Gamma'}$
 is the associated tau class defined for $\Gamma'$ by definition \ref{defin; tau} on
 page \pageref{defin; tau}.

\medskip

 By proposition \ref{prop; identical}, one may ``collapse'' the blowing up sequence indexed
 by $I_{\Gamma}$ (following the reversed ordering of $\models$) to the new blowing up
 sequence indexed by the reduced index set $\bar{I}_{\Gamma}$ (see definition \ref{defin; reduce} for
 its definition) following the reversed linear ordering of $\vdash$ without
 changing the answer. As was argued, the
 permutation and the collapsing of blowing up centers not affect 
the result of the localized contribution of top Chern class. Thus 
 the yet-to-be-identified intersection pairing is equal to \footnote{The space 
$\tilde{\hat{X}}_{\Gamma}$ with the exceptional divisor
 $\hat{D}_{\Gamma}$ denotes the blowing up of $\hat{X}_{\Gamma}$.
 The hatted divisor $\hat{D}_{\Gamma'}$ have
 been used in the previous subsection already to denote the exceptional divisors blown up following
 the reversed linear ordering of $\vdash$ inside $\bar{I}_{\Gamma}$.}

$$\hskip -1.3in
\{\sum_{i=1}^{i=rank_{\bf C}{\bf W}_{canon}}(-1)^{i-1}
c_{rank_{\bf C}{\bf W}_{canon}-i}(\pi_{\tilde{\hat{X}}_{\Gamma}}^{\ast}{\bf W}_{canon}\otimes 
 {\bf H}\otimes_{\Gamma'\in \bar{I}_{\Gamma}} 
{\cal O}(-\hat{D}_{\Gamma'}))\cap \hat{D}_{\Gamma}^{i-1}\cap [\hat{D}_{\Gamma}]\cap 
c_1({\bf H})^{dim_{\bf C}M_n+rank_{\bf C}({\bf V}_{canon}-{\bf W}_{canon})+q-1}\}.$$

Among the many different blowing ups indexed by the graphs 
$\Gamma'\in \bar{I}_{\Gamma}$, whenever $\Gamma'\in \bar{I}_{\Gamma}^{\gg}$
 the restricted family moduli spaces
 ${\cal M}_{C-{\bf M}(E)E-\sum_{e_i'\cdot (C-{\bf M}(E)E)<0}e_i'}\times_{M_n}
Y(\Gamma')$ can  be viewed naturally as sub-schemes of 
${\cal M}_{C-{\bf M}(E)E-\sum_{e_i\cdot (C-{\bf M}(E)E)<0}e_i}\times_{M_n}
Y(\Gamma)$ (by adjoining curves in ${\cal M}_{C-{\bf M}(E)E-\sum_{e_i'\cdot (C-{\bf M}(E)E)<0}e_i'}$
 to the exceptional curves dual to 
$\sum_{e_i\cdot (C-{\bf M}(E)E)>0; e_i'\cdot (C-{\bf M}(E)E)<0}e_i'$) and therefore
 sub-scheme of ${\cal M}_{C-{\bf M}(E)E}\times_{M_n}Y(\Gamma)$.

Under the reversed $\vdash$ linear ordering of blowing ups, these blowing ups with
$\Gamma'\in \bar{I}_{\Gamma}^{\gg}$ are performed at the very end of the linear chain
 of blowing ups parametrized
 by $\bar{I}_{\Gamma}$. Thus we may decompose $\bar{I}_{\Gamma}=\bar{I}_{\Gamma}^{\gg}\coprod 
 (\bar{I}_{\Gamma}-\bar{I}_{\Gamma}^{\gg})$ and 
use residual intersection formula to re-write the above
 intersection pairing as 

$$
\{\sum_{i=1}^{i=rank_{\bf C}{\bf W}_{canon}}(-1)^{i-1}
c_{rank_{\bf C}{\bf W}_{canon}-i}(\pi_{\tilde{\hat{X}}_{\Gamma}}^{\ast}{\bf W}_{canon}\otimes 
 {\bf H}\otimes_{\Gamma'\in \bar{I}_{\Gamma}-\bar{I}_{\Gamma}^{\gg}} 
{\cal O}(-\hat{D}_{\Gamma'}))$$
$$\cap \hat{D}_{\Gamma}^{i-1}\cap [\hat{D}_{\Gamma}]\cap 
c_1({\bf H})^{dim_{\bf C}M_n+rank_{\bf C}({\bf V}_{canon}-{\bf W}_{canon})+q-1}\}$$

$$
-\sum_{\Gamma'\in \bar{I}_{\Gamma}^{\gg}}
\{\sum_{i=1}^{i=rank_{\bf C}{\bf W}_{canon}}(-1)^{i-1}
c_{rank_{\bf C}{\bf W}_{canon}-i}(\pi_{\tilde{\hat{X}}_{\Gamma'}}^{\ast}{\bf W}_{canon}\otimes 
 {\bf H}\otimes_{\Gamma_1\in \bar{I}_{\Gamma}; \Gamma'\vdash \Gamma_1} 
{\cal O}(-\hat{D}_{\Gamma_1}))$$
$$\cap \hat{D}_{\Gamma'}^{i-1}\cap [\hat{D}_{\Gamma'}]\cap 
c_1({\bf H})^{dim_{\bf C}M_n+rank_{\bf C}({\bf V}_{canon}-
{\bf W}_{canon})+q-1}\}.$$ 

In the second term of the above sum (i.e. when $\Gamma'\in 
\bar{I}_{\Gamma}^{\gg}$),
 the index set restriction ``$\Gamma_1\in \bar{I}_{\Gamma}; \Gamma'\vdash
 \Gamma_1$'' is the same as the alternative restriction 
$''\Gamma_1\in (\bar{I}_{\Gamma}-\bar{I}_{\Gamma}^{\gg})\coprod
 \{\Gamma_1|\Gamma_1\in \bar{I}_{\Gamma}^{\gg}; \Gamma'\gg \Gamma_1\}''$. 
\footnote{It is because $\vdash$ is identical to $\gg$ on the subset 
$\bar{I}_{\Gamma}^{\gg}$ and by the definition of
 $\vdash$ (see definition \ref{defin; newordering}) 
the elements in $\bar{I}_{\Gamma}^{\gg}$ are larger than all elements in 
$\bar{I}_{\Gamma}-\bar{I}_{\Gamma}^{\gg}$)}

 By the definition/construction
 of $\bar{I}_{\Gamma}$ and $\bar{I}_{\Gamma}^{\gg}$, all elements 
 $\Gamma''\in (\bar{I}_{\Gamma}-\bar{I}_{\Gamma}^{\gg})$ satisfy 
 $(\Gamma, \sum_{e_i\cdot (C-{\bf M}(E)E)<0}e_i)\sqsupset (\Gamma'', 
\sum_{e_i''\cdot (C-{\bf M}(E)E)<0}e_i'')$. Thus 
 $(\Gamma', \sum_{e_i'\cdot (C-{\bf M}(E)E)<0}e_i')\sqsupset (\Gamma'', 
\sum_{e_i''\cdot (C-{\bf M}(E)E)<0}e_i'')$ as well since we know 
 $\Gamma'\in \bar{I}_{\Gamma}^{\gg}$. By the discussion of lemma \ref{lemm; restriction},
 the alternative index set restriction on $\Gamma_1$ can then be replaced by the equivalent one 
``$\Gamma_1\in \bar{I}_{\Gamma'}$''. Then by applying proposition \ref{prop; identical} to 
$I_{\Gamma'}$ and $\bar{I}_{\Gamma'}$, we may ``un-collapse'' to restore
 the reduced index set $\bar{I}_{\Gamma'}$ back to $I_{\Gamma'}$ without affecting the
 result of the intersection pairing. This implies that the second term of the above expression
 can be expressed as (observe that the hat of $\hat{D}_{\Gamma'}$ or 
$\hat{D}_{\Gamma_1}$ has been removed)

$$
-\sum_{\Gamma'\in \bar{I}_{\Gamma}^{\gg}}
\{\sum_{i=1}^{i=rank_{\bf C}{\bf W}_{canon}}(-1)^{i-1}
c_{rank_{\bf C}{\bf W}_{canon}-i}(\pi_{\tilde{X}}^{\ast}{\bf W}_{canon}\otimes 
 {\bf H}\otimes_{\Gamma_1\in I_{\Gamma'}} 
{\cal O}(-D_{\Gamma_1}))$$
$$\cap D_{\Gamma'}^{i-1}\cap [D_{\Gamma'}]\cap 
c_1({\bf H})^{dim_{\bf C}M_n+rank_{\bf C}({\bf V}_{canon}-{\bf W}_{canon})+q-1}\}.$$

\medskip

 Then by the {\bf Inductive Hypothesis} on page \pageref{IH} above,
 the integral of 
each of these terms 
is equal to a modified algebraic family 
invariant attached to $\Gamma'$ and the original 
yet-to-be-identified intersection pairing of the 
localized contribution of the top Chern class with
 $c_1({\bf H})^{dim_{\bf C}M_n+rank_{\bf C}({\bf V}_{canon}-{\bf W}_{canon})+q-1}$ 
is equal to \footnote{This inductive pattern 
is exactly why we had defined the modified invariants
 earlier on page \pageref{defin; inductive} in this way.}

$$
\int_{\tilde{\hat{X}}_{\Gamma}}\{\sum_{i=1}^{i=rank_{\bf C}{\bf W}_{canon}}(-1)^{i-1}
c_{rank_{\bf C}{\bf W}_{canon}-i}(\pi_{\tilde{\hat{X}}_{\Gamma}}^{\ast}{\bf W}_{canon}\otimes 
 {\bf H}\otimes_{\Gamma'\in \bar{I}_{\Gamma}-\bar{I}_{\Gamma}^{\gg}} 
{\cal O}(-\hat{D}_{\Gamma'}))$$
$$\cap \hat{D}_{\Gamma}^{i-1}\cap [\hat{D}_{\Gamma}]\cap 
c_1({\bf H})^{dim_{\bf C}M_n+rank_{\bf C}({\bf V}_{canon}-{\bf W}_{canon})+q-1}\}$$

$$-\sum_{\Gamma'\in \bar{I}_{\Gamma}^{\gg}}
{\cal AFSW}_{M_{n+1}\times_{M_n}Y(\Gamma')\times
T(M)\mapsto Y(\Gamma')\times T(M)}^{\ast}(c_{total}(\tau_{\Gamma'}), 
C-{\bf M}(E)E-\sum_{e_i'\cdot (C-{\bf M}(E)E)<0}e_i').$$

 If we can identify the first term in this sum with 
 $${\cal AFSW}_{M_{n+1}\times_{M_n}Y(\Gamma)\times T(M)\mapsto Y(\Gamma)\times
 T(M)}(c_{total}(\tau_{\Gamma}), C-{\bf M}(E)E-\sum_{e_i\cdot (C-{\bf M}(E)E)<0}e_i),$$ 
 then by definition \ref{defin; inductive}, the total expression is exactly
 what was defined to be the modified algebraic family invariant attached to
 $\Gamma$, 

$${\cal AFSW}_{M_{n+1}\times_{M_n}Y(\Gamma)\times T(M)\mapsto Y(\Gamma)\times
 T(M)}^{\ast}(c_{total}(\tau_{\Gamma}), C-{\bf M}(E)E-\sum_{e_i\cdot 
(C-{\bf M}(E)E)<0}e_i)$$ and then the identification is complete.

 In the rest of the proof we identify the first term of the last sum with the specific 
mixed algebraic family Seiberg-Witten invariant associated with $c_{total}(\tau_{\Gamma})$.

\medskip

 Following the same convention as in [Liu5], we let $k_1<k_2<\cdots<k_p$ be the subscripts
 such that $e_{k_i}\cdot (C-{\bf M}(E)E)<0$ for all $1\leq i\leq p$.

\medskip

\noindent Case I: In this case we deal with the more interesting situation that 
$e_{k_i}^2\geq e_{k_i}\cdot (C-{\bf M}(E)E)$ for all $1\leq i\leq p$.

\medskip

Step One:
 Firstly we make use of the assumption $e_{k_i}^2\geq e_{k_i}\cdot (C-{\bf M}(E)E)$, $1\leq i\leq p$
 and show the following Chern classes identity for ${\bf V}_{quot}$,

\begin{lemm}\label{lemm; chernidentity}
Let $\pi_g:X\times_{M_n}Y(\Gamma)\mapsto Y(\Gamma)\times T(M)$ and $\pi_t:Y(\Gamma)\times T(M)\mapsto 
 Y(\Gamma)$ be the natural projection maps.
Let ${\bf V}_{quot}$ be the quotient bundle of ${\bf W}_{canon}|_{Y(\Gamma)\times T(M)}$ as was
 defined on page \pageref{underlinewcanon}. Let $\tau_{\Gamma}$  be the tau class defined 
in definition \ref{defin; tau} on page \pageref{defin; tau}. 

Suppose that $e_{k_i}^2\geq e_{k_i}\cdot (C-{\bf M}(E)E)$ for all $1\leq i\leq p$, then
there is an identity among the top Chern classes over $X\times_{M_n}Y(\Gamma)$, 

$$c_{top}({\bf H}\otimes \pi_g^{\ast}{\bf V}_{quot})=c_{top}({\bf H}\otimes \pi_g^{\ast}\tau_{\Gamma})
\cap c_{top}(\pi_g^{\ast}\pi_t^{\ast}{\bf N}_{Y(\Gamma)}X).$$ 
\end{lemm}

\noindent Proof: Let ${\bf Q}_{k_i}$ 
and ${\bf E}_C$ be 
the line bundles associated to the invertible sheaves ${\cal Q}_{k_i}$ and ${\cal E}_C$ which
 appeared in definition \ref{defin; tau} (see also proposition \ref{prop; exact}). 

Define the vector bundle ${\bf G}_{\Gamma}=
{\bf H}\otimes 
\pi_g^{\ast}\oplus_{1\leq i\leq p}\pi_t^{\ast}
{\bf N}_{Y(\Gamma_{e_{k_i}})}M_n|_{Y(\Gamma)}\otimes {\bf Q}_{k_i}\otimes {\bf E}_C
\mapsto X\times_{M_n}Y(\Gamma)$ to be a 
rank $dim_{\bf C}M_n-dim_{\bf C}Y(\Gamma)=codim_{\bf C}\Gamma$
 vector bundle over $X\times_{M_n}Y(\Gamma)$
\footnote{Notice that ${\bf G}_{\Gamma}$ is constructed from 
 ${\bf N}_{X\times_{M_n}Y(\Gamma)}X$ twisted by ${\bf Q}_{k_i}$ and ${\bf E}_C$
 on direct direct factors.}.

By prop \ref{prop; id} on page \pageref{prop; id}
 in subsection \ref{subsection; lf}, definition \ref{defin; replace} 
on page \pageref{defin; replace} and the tensor product formula
 of the top Chern classes, we know $c_{top}({\bf H}\otimes \pi_g^{\ast}{\bf V}_{quot})
=c_{top}({\bf H}\otimes \pi_g^{\ast}\tilde{\bf V}_{quot})$.  

 Thus we have

 $$\hskip -1.3in
c_{top}({\bf H}\otimes \pi_g^{\ast}{\bf V}_{quot})=c_{top}(({\bf H}\otimes 
\pi_g^{\ast}\tilde{\bf V}_{quot}-{\bf G}_{\Gamma})\oplus {\bf G}_{\Gamma}).$$

 Recall that by proposition \ref{prop; transversal} in 
section \footnote{Or equivalently by proposition 4.7 on page 426 of 
[Liu1].} \ref{section; strata} the space 
$Y(\Gamma)=\cap_{1\leq i\leq p}Y(\Gamma_{e_{k_i}})$ is a transversal 
intersection of smooth loci $Y(\Gamma_{e_i})$. 
We also know from lemma 9 of [Liu5] that all the
$Y(\Gamma_{e_{k_i}})$, being the family moduli space of $e_{k_i}$,
 are defined by regular global sections of 
 ${\cal R}^0\pi_{\ast}\bigl({\cal O}_{\sum_{j_{k_i}}E_{j_{k_i}}}(E_{k_i})\bigr)$ 
over $M_n$ determined by the morphism of locally free sheaves 
\footnote{The former is invertible, while the latter is the canonical
obstruction bundle of $e_{k_i}$.} 
$${\cal R}^0\pi_{\ast}{\cal O}_{E_{k_i}}\mapsto {\cal R}^0\pi_{\ast}\bigl(
{\cal O}_{\sum_{j_{k_i}}E_{j_{k_i}}}(E_{k_i})\bigr)$$ 
and ${\cal R}^0\pi_{\ast}\bigl({\cal 
O}_{\sum_{j_{k_i}}E_{j_{k_i}}}(E_{k_i})|_{Y(\Gamma_{e_{k_i}})}\bigr)\cong 
{\cal R}^1\pi_{\ast}\bigl({\cal O}_{\tilde{\Xi}_{k_i}}(E_{k_i}-
\sum_{j_{k_i}}E_{j_{k_i}})\bigr)\cong {\cal N}_{Y(\Gamma_{e_{k_i}})}M_n$, the normal sheaf of
 $Y(\Gamma_{e_{k_i}})$ in $M_n$.

  By definition \ref{defin; tau}, lemma \ref{lemm; free}
 and definition \ref{defin; replace}, the class
$\tilde{\bf V}_{quot}-\oplus_{1\leq i\leq p}
{\bf N}_{Y(\Gamma_{e_{k_i}})}M_n|_{Y(\Gamma)}
\otimes {\bf Q}_{k_i}\otimes {\bf E}_C$ is equal to 
$\tau_{\Gamma}$ (expressed here by
 a difference of vector bundles rather than the 
corresponding locally free sheaves), and is
 represented by a vector bundle of rank 
$rank_{\bf C}\tilde{\bf V}_{quot}-codim_{\bf C}Y(\Gamma)$. This and the Whitney sum formula of
 the top Chern class imply that 
$c_{top}({\bf H}\otimes \pi_g^{\ast}{\bf V}_{quot})=c_{top}({\bf H}\otimes \pi_g^{\ast}\tau_{\Gamma})
\cap c_{top}({\bf G}_{\Gamma})$.

Because for
 the type $I$ exceptional class $e_{k_i}=E_{k_i}-\sum_{j_{k_i}}E_{j_{k_i}}$
 the smooth locus $Y(\Gamma_{e_{k_i}})\subset M_n$ has been
 the regular zero locus defined by the canonical global
 section of ${\cal W}_{e_{k_i}}=
{\cal R}^0\pi_{\ast}\bigl({\cal O}_{\sum_{j_{k_i}}E_{j_{k_i}}}(E_{k_i})\bigr)$, 
defined by the canonical algebraic Kuranishi model of $e_{k_i}$ 
(see section 6.2 of [Liu5]). 
 So ${\cal W}_{e_{k_i}}|_{Y(\Gamma_{e_{k_i}})}\cong
 {\cal N}_{Y(\Gamma_{e_{k_i}})}M_n$. By lemma 10 of
\footnote{This lemma implies that the restriction of the top Chern class of 
a vector bundle ${\bf E}$ to the regular zero locus 
$Z(s)$ with codimension $rank_{\bf C}{\bf E}$ 
of its regular section $s$ is equal to the restriction of the top Chern class of 
 ${\bf E}\otimes {\bf Q}$, twisted by a line bundle ${\bf Q}$. Notice that
 it holds only because we are working in ${\cal A}_{\cdot}(Z(S))$ instead of the
 whole space $X$.}
 [Liu5], this implies that the top Chern class 
$c_{top}(\pi_g^{\ast}(\pi_t^{\ast}{\bf N}_{Y(\Gamma_{e_{k_i}})}M_n|_{Y(\Gamma)}\otimes {\bf Q}_{k_i})
\otimes {\bf H}\otimes {\cal E}_C)$ on $X\times_{M_n}Y(\Gamma)$, 
where $Y(\Gamma)=\cap_{1\leq i\leq p}Y(\Gamma_{e_{k_i}})$, is equal to
 the top Chern class of the un-twisted version $c_{top}({\bf N}_{\pi_g^{\ast}(\pi_t^{\ast}
Y(\Gamma_{e_{k_i}})}M_n|_{Y(\Gamma)})$.

 Because this is applicable to all $1\leq i\leq p$, we find that 
$c_{top}({\bf G}_{\Gamma})=c_{top}(\pi_g^{\ast}\pi_t^{\ast}{\bf N}_{Y(\Gamma)}M_n)$.
The proof of this lemma is complete. 
$\Box$

  Secondly we consider the normal cone (into its compactification) 
${\bf C}_{\hat{X}_{\Gamma}\times_{M_n}Y(\Gamma)}\hat{X}_{\Gamma}\subset 
{\bf P}({\bf C}_{\hat{X}_{\Gamma}\times_{M_n}Y(\Gamma)}\hat{X}_{\Gamma}\oplus 1)$ of
 the closed embedding $\hat{X}_{\Gamma}\times_{M_n}Y(\Gamma)\subset \hat{X}_{\Gamma}$.

 Consider \footnote{following chapter 5 of [F].}
the blowing up along $\hat{X}_{\Gamma}\times_{M_n}Y(\Gamma)\times \{0\}\subset 
\hat{X}_{\Gamma}\times
 {\bf C}$. The exceptional divisor of this blowing up is isomorphic to 
${\bf P}({\bf C}_{\hat{X}_{\Gamma}\times_{M_n}Y(\Gamma)}\hat{X}_{\Gamma}\oplus 1)$.

 As it projects onto $X\times_{M_n}Y(\Gamma)\times \{0\}\subset X\times {\bf C}$ under 
 $\hat{X}_{\Gamma}\times {\bf C}\mapsto X\times {\bf C}$, by the universal
 property of the blowing up (proposition 7.14. on page 164 of [Ha]), 
${\bf P}({\bf C}_{\hat{X}_{\Gamma}\times_{M_n}Y(\Gamma)}\hat{X}_{\Gamma}\oplus 1)$ maps onto the
 exceptional divisor \footnote{The embedding is regular, so we use
 ${\bf N}_{X\times_{M_n}Y(\Gamma)}X$ and ${\bf C}_{X\times_{M_n}Y(\Gamma)}X$
interchangeably.}
${\bf P}({\bf N}_{X\times_{M_n}Y(\Gamma)}X\oplus 1)$ 
of the blowing up along $X\times_{M_n}Y(\Gamma)\times 
\{0\}\subset X\times {\bf C}$.
Then it induces a 
surjection of the normal cones
 ${\bf C}_{\hat{X}_{\Gamma}\times_{M_n}Y(\Gamma)}\hat{X}_{\Gamma}
\mapsto {\bf N}_{X\times_{M_n}Y(\Gamma)}X$ 
 and we have the following commutative diagram,

\[
\hskip -.5in
\begin{array}{ccccc}\label{proj}
\hat{X}_{\Gamma}\times_{M_n}Y(\Gamma) & \longrightarrow & 
{\bf C}_{\hat{X}_{\Gamma}\times_{M_n}Y(\Gamma)}\hat{X}_{\Gamma} & \longrightarrow 
 & {\bf P}({\bf C}_{\hat{X}_{\Gamma}\times_{M_n}Y(\Gamma)}\hat{X}_{\Gamma}\oplus 1) \\
 \Big\downarrow\vcenter{%
\rlap{$\scriptstyle{\mathrm{\pi_h}}\, $}}
 & & \Big\downarrow\vcenter{%
\rlap{$\scriptstyle{\mathrm{\pi_{\bf C}}}\, $}} & & \Big\downarrow\vcenter{%
\rlap{$\scriptstyle{\mathrm{\pi_{\bf P}}}\, $}}\\
X\times_{M_n}Y(\Gamma) & \longrightarrow & 
{\bf C}_{X\times_{M_n}Y(\Gamma)}X & \longrightarrow 
 & {\bf P}({\bf C}_{X\times_{M_n}Y(\Gamma)}X\oplus 1)
\end{array}
\]

 As the exceptional divisor $\hat{D}_{\Gamma}$ maps into $Y(\Gamma)$ under 
 $\tilde{\hat{X}}_{\Gamma}\mapsto M_n$, we have 
 $\hat{D}_{\Gamma}\subset \tilde{\hat{X}}_{\Gamma}\times_{M_n}Y(\Gamma)$. 
Our original intersection pairing involving $\hat{D}_{\Gamma}$ can be pushed-forward by 
 $(\pi_r)_{\ast}$ (here $\pi_r: \tilde{\hat{X}}_{\Gamma}\mapsto \hat{X}_{\Gamma}$ is the
 blowing down map) into
 ${\cal A}_0(\hat{X}_{\Gamma}\times_{M_n}Y(\Gamma))$ 
and it defines a zero dimensional cycle class lying inside 
$\hat{X}_{\Gamma}\times_{M_n}Y(\Gamma)$. 

Either by a direct computation on the localized contribution of top Chern class (involving the
 Segre class of some normal cone), or by the technique of the 
deformation to the normal cone (consult
 chapter 5 of [F]) from 
$\hat{X}_{\Gamma}\times_{M_n}Y(\Gamma)\subset \hat{X}_{\Gamma}$ to 
$\hat{X}_{\Gamma}\times_{M_n}Y(\Gamma)\subset {\bf C}_{\hat{X}_{\Gamma}
\times_{M_n}Y(\Gamma)}\hat{X}_{\Gamma}$,
 one may replace the total space $\hat{X}_{\Gamma}$ by a 
linearized object of the same dimension, namely the projectified normal cone 
${\bf P}({\bf C}_{\hat{X}_{\Gamma}\times_{M_n}Y(\Gamma)}\hat{X}_{\Gamma}\oplus 1)$ or
 its affine part \footnote{We prefer the former if we want the space to be complete.} 
${\bf C}_{\hat{X}_{\Gamma}\times_{M_n}Y(\Gamma)}\hat{X}_{\Gamma}$. 

Define $\pi_f:{\bf P}({\bf C}_{\hat{X}_{\Gamma}\times_{M_n}Y(\Gamma)}\hat{X}_{\Gamma}\oplus 1)
\mapsto \hat{X}_{\Gamma}\times_{M_n}Y(\Gamma)$ to be the projection map.
To get a global intersection pairing on 
${\bf P}({\bf C}_{\hat{X}_{\Gamma}\times_{M_n}Y(\Gamma)}\hat{X}_{\Gamma}\oplus 1)$ which is refined
 to our localized 
intersection pairing, we twist the bundle $\pi_f^{\ast}\pi_g^{\ast}{\bf W}_{canon}$ by
 ${\cal O}(\hat{\bf P}_{\infty})$, where $\hat{\bf P}_{\infty}=
{\bf P}({\bf C}_{\hat{X}_{\Gamma}\times_{M_n}Y(\Gamma)}\hat{X}_{\Gamma})$ 
$\subset {\bf P}({\bf C}_{\hat{X}_{\Gamma}\times_{M_n}Y(\Gamma)}\hat{X}_{\Gamma}
\oplus 1)$ is the divisor at 
 infinity. As our intersection pairing is localized at the zero section 
 $\hat{X}_{\Gamma}\times_{M_n}Y(\Gamma)$ (totally disjoint from
 $\hat{\bf P}_{\infty}$), the fact that ${\cal O}(\hat{\bf P}_{\infty})$ is trivialized over
 the affine cone ${\bf C}_{\hat{X}_{\Gamma}\times_{M_n}Y(\Gamma)}\hat{X}_{\Gamma}$
 allows us to remove the
 ${\cal O}(\hat{\bf P}_{\infty})$ tensor product effectively 
in our calculation
\footnote{as far as our intersection cycle does not overlap 
with $\hat{\bf P}_{\infty}$.}.

\medskip

 Consider the pull-back of the ${\bf H}$-twisted short exact sequence 
$0\mapsto \underline{\bf W}_{canon}|\mapsto
 {\bf W}_{canon}|_{Y(\Gamma)\times T(M)} 
\mapsto {\bf V}_{quot}\mapsto 0$ (consult page \pageref{underlinewcanon} in
 section \ref{section; bundle}) by $(\pi_g\pi_h\pi_f)^{\ast}$,
\footnote{The map $\pi_h$ was defined in the above commutative diagram on 
page \pageref{proj}.} 
the short exact sequence exists on the whole space 
${\bf P}({\bf C}_{\hat{X}_{\Gamma}\times_{M_n}Y(\Gamma)}\hat{X}_{\Gamma}\oplus 1)$. 
\footnote{This is the benefit of adopting
 the projectified normal cone 
${\bf P}({\bf C}_{\hat{X}_{\Gamma}\times_{M_n}Y(\Gamma)}\hat{X}_{\Gamma}\oplus 1)$
 than the original space $\hat{X}_{\Gamma}$, a replacement of tubular 
neighborhood in the ${\cal C}^{\infty}$ category.}

 By proposition \ref{prop; cap} we know that the residual intersection formula of 
${\bf W}_{canon}$ and of $\underline{\bf W}_{canon}$ are compatible. We may replace the above 
 top intersection pairing of the localized Chern class by 

$$\hskip -1in
\{c_{total}((\pi_h\pi_f)^{\ast}(\pi_g^{\ast}\underline{\bf W}_{canon}\otimes {\bf H})
\otimes_{\Gamma'\in \bar{I}_{\Gamma}-\bar{I}_{\Gamma}^{\gg}} 
{\cal O}(-\pi_f^{\ast}\hat{D}_{\Gamma'}))\cap s_{total}(\pi_r(\hat{D}_{\Gamma}), 
{\bf C}_{\hat{X}_{\Gamma}\times_{M_n}Y(\Gamma)}\hat{X}_{\Gamma})\}_{dim_{\bf C}X-
rank_{\bf C}\underline{\bf W}_{canon}}$$
$$\cap c_{top}((\pi_h\pi_f)^{\ast}({\bf H}\otimes \pi_g^{\ast}{\bf V}_{quot}))\cap 
c_1((\pi_h\pi_f)^{\ast}{\bf H})^{dim_{\bf C}M_n+rank_{\bf C}({\bf V}_{canon}-{\bf W}_{canon})+q-1}.$$

 By lemma \ref{lemm; chernidentity} we may replace the top Chern class
$c_{top}((\pi_h\pi_f)^{\ast}({\bf H}\otimes \pi_g^{\ast}{\bf V}_{quot}))$ by 
 $c_{top}((\pi_h\pi_f)^{\ast}({\bf H}\otimes 
\pi_g^{\ast}\tau_{\Gamma}))
\cap c_{top}((\pi_t\pi_g\pi_h\pi_f)^{\ast}
{\bf N}_{Y(\Gamma)}M_n)$, which is the same as
 $c_{top}((\pi_h\pi_f)^{\ast}({\bf H}\otimes 
\pi_g^{\ast}\tau_{\Gamma}))
\cap c_{top}((\pi_h\pi_f)^{\ast}{\bf N}_{X\times_{M_n}Y(\Gamma)}X)$.

 On the other hand the projection of normal cones 
$\pi_{\bf C}:{\bf C}_{\hat{X}_{\Gamma}\times_{M_n}Y(\Gamma)}\hat{X}_{\Gamma}
\mapsto {\bf N}_{X\times_{M_n}Y(\Gamma)}X$ induces (by pull-back) 
a tautological section of 
 $(\pi_h\pi_f)^{\ast}{\bf N}_{X\times_{M_n}Y(\Gamma)}X$ over 
${\bf C}_{\hat{X}_{\Gamma}\times_{M_n}Y(\Gamma)}\hat{X}_{\Gamma}$ and
\footnote{It extends to a section of 
$(\pi_h\pi_f)^{\ast}{\bf N}_{X\times_{M_n}Y(\Gamma)}X\otimes {\cal O}(\hat{\bf P}_{\infty})$
 on ${\bf P}({\bf C}_{\hat{X}_{\Gamma}\times_{M_n}Y(\Gamma)}\hat{X}_{\Gamma}\oplus 1)$.}
 its zero locus is exactly 
 $\hat{X}_{\Gamma}\times_{M_n}Y(\Gamma)\subset {\bf C}_{\hat{X}_{\Gamma}\times_{M_n}Y(\Gamma)}
\hat{X}_{\Gamma}$.
 
 Therefore the cycle class $[\hat{X}_{\Gamma}\times_{M_n}Y(\Gamma)]\subset 
{\cal A}_{\cdot}({\bf C}_{\hat{X}_{\Gamma}\times_{M_n}Y(\Gamma)}\hat{X}_{\Gamma})$ represents
 the cap product of the fundamental class with the 
top Chern class $c_{top}((\pi_h\pi_f)^{\ast}{\bf N}_{X\times_{M_n}Y(\Gamma)}X)$.

 Thus we may replace this top Chern class 
$c_{top}((\pi_h\pi_f)^{\ast}{\bf N}_{X\times_{M_n}Y(\Gamma)}X)$ in our pairing 
by the zero section cycle class 
$[\hat{X}_{\Gamma}\times_{M_n}Y(\Gamma)]$ of the compactification of the
 normal cone.
 Consequently, we can restrict the intersection pairing to the zero section 
$\hat{X}_{\Gamma}\times_{M_n}Y(\Gamma)$ of its normal cone in $\hat{X}_{\Gamma}$ and get 

$$
\{c_{total}(\pi_h^{\ast}(\pi_g^{\ast}\underline{\bf W}_{canon}\otimes {\bf H})
\otimes_{\Gamma'\in \bar{I}_{\Gamma}-\bar{I}_{\Gamma}^{\gg}} 
{\cal O}(-\hat{D}_{\Gamma'}))\cap s_{total}(\pi_r(\hat{D}_{\Gamma}), 
{\bf C}_{\hat{X}_{\Gamma}\times_{M_n}Y(\Gamma)}\hat{X}_{\Gamma})\}_{dim_{\bf C}X-
rank_{\bf C}\underline{\bf W}_{canon}}$$
$$\cap c_{top}(\pi_h^{\ast}({\bf H}\otimes \pi_g^{\ast}\tau_{\Gamma}))\cap 
c_1(\pi_h^{\ast}{\bf H})^{dim_{\bf C}M_n+rank_{\bf C}({\bf V}_{canon}-{\bf W}_{canon})+q-1}
\in {\cal A}_0(\hat{X}_{\Gamma}\times_{M_n}Y(\Gamma)).$$

$\diamondsuit$ 
Let us summarize: In step one we have succeeded in restricting the top 
intersection pairing to the subspace 
 $\hat{X}_{\Gamma}\times_{M_n}Y(\Gamma)$ by using some novel property of 
$\tilde{\bf V}_{quot}$. The
 class $c_{top}({\bf H}\otimes \pi_g^{\ast}\tau_{\Gamma})$ has appeared because of
 lemma 
\ref{lemm; chernidentity}.

\bigskip

Step Two: Consider the bundle map 
$\pi_h^{\ast}(\pi_g^{\ast}{\bf W}_{canon}^{\circ}|_{Y(\Gamma)\times T(M)}
\otimes {\bf H})
\mapsto \pi_h^{\ast}(\pi_g^{\ast}\underline{\bf W}_{canon}\otimes {\bf H})$ 
induced by ${\bf H}$-twisted
 version of 
the $\pi_g^{\ast}-$pulled-back 
vector bundle map ${\bf W}_{canon}^{\circ}\mapsto \underline{\bf W}_{canon}$ over $Y(\Gamma)\times
 T(M)$
(see page \pageref{underlinewcanon}) pull-back by $\pi_h^{\ast}$. 
Our goal is to explain why we may use the top Chern class of 
$\pi_h^{\ast}(\pi_g^{\ast}{\bf W}_{canon}^{\circ}|_{Y(\Gamma)\times T(M)}
\otimes {\bf H})$
 to replace the complicated 
bundle $\pi_h^{\ast}(\pi_g^{\ast}\underline{\bf W}_{canon}\otimes {\bf H})\otimes
{\cal O}(-\sum_{\Gamma'\in \bar{I}_{\Gamma}-\bar{I}_{\Gamma}^{\gg}}\hat{D}_{\Gamma'})$
 in the localized contribution of top Chern class.

\medskip 

Observe that for all $\Gamma'\in \bar{I}_{\Gamma}-\bar{I}_{\Gamma}^{\gg}$, 
 $(\Gamma, \sum_{e_i\cdot (C-{\bf M}(E)E)<0}e_i)\sqsupset (\Gamma'', \sum_{e_i''\cdot 
(C-{\bf M}(E)E)<0}e_i'')$. This condition $\sqsupset$ (consult definition \ref{defin; sq}
 for its definition) implies that the sub-scheme ${\cal M}_{C-{\bf M}(E)E-
\sum_{e_i\cdot (C-{\bf M}(E)E)<0}
e_i}\times_{M_n}Y(\Gamma)\cap Y(\Gamma')$ can be embedded into
 ${\cal M}_{C-{\bf M}(E)E-\sum_{e_i'\cdot (C-{\bf M}(E)E)<0}
e_i'}\times_{M_n}Y(\Gamma)\cap Y(\Gamma')$.

Then by proposition \ref{prop; union2} and the remark \ref{rem; gg}
 right after its proof, 
we may decompose ${\cal M}_{C-{\bf M}(E)E}\times_{M_n}Y(\Gamma)$
 into the union of the natural image of
 ${\cal M}_{C-{\bf M}(E)E-\sum_{e_i\cdot (C-{\bf M}(E)E)<0}
e_i}\times_{M_n}Y(\Gamma)$
 and the image of the union $\cup_{\Gamma'\in \bar{I}_{\Gamma}}
{\cal M}_{C-{\bf M}(E)E-\sum_{e_i'\cdot (C-{\bf M}(E)E)<0}
e_i'}\times_{M_n}Y(\Gamma)\cap Y(\Gamma')$. \footnote{We have 
changed the index set of the union of sub-schemes
 from $I_{\Gamma}$ to $\bar{I}_{\Gamma}$ by
 remark \ref{rem; gg} on page \pageref{rem; gg}.}

On the other hand, by using the induced bundle map 
 $\pi_g^{\ast}{\bf W}_{canon}^{\circ}|_{Y(\Gamma)\times T(M)}
\otimes {\bf H}\mapsto \pi_g^{\ast}\underline{\bf W}_{canon}
\otimes {\bf H}$ (consult page \pageref{underlinewcanon} in
 section \ref{section; bundle}) and by using $X\supset Z(s_{canon}^{\circ})={
\cal M}_{C-{\bf M}(E)E-\sum_{e_i\cdot (C-{\bf M}(E)E)<0}e_i}$, we realize that
 the image of
 the union 

$\cup_{\Gamma'\in \bar{I}_{\Gamma}}{\cal M}_{C-{\bf M}(E)E-
\sum_{e_i'\cdot (C-{\bf M}(E)E)<0}
e_i'}\times_{M_n}Y(\Gamma)\cap Y(\Gamma')$ in ${\cal M}_{C-{\bf M}(E)E}\subset X$, the
 excess component, is nothing but
 the projection image 
of the intersection of the section $s_{canon}^{\circ}$ and the kernel cone, 
 $\pi_{\pi_g^{\ast}{\bf W}_{canon}^{\circ}\otimes {\bf H}}
(s_{canon}^{\circ}\cap (\cup_{i>0}{\bf C}_{\rho_i}))$.\footnote{The kernel cone
 means 
the algebraic sub-cone associated to the kernel semi-bundle of the map
 $\pi_g^{\ast}{\bf W}_{canon}^{\circ}\otimes {\bf H}\mapsto 
\pi_g^{\ast}\underline{\bf W}_{canon}
\otimes {\bf H}$.} 

 Therefore the blowing ups of these loci into the union of divisors 
$\cup_{\Gamma'\in \bar{I}_{\Gamma}-\bar{I}_{\Gamma}^{\gg}}
\hat{D}_{\Gamma'}$ has fitted into the framework of proposition \ref{prop; equivalent} under
 the identification $\pi_g^{\ast}{\bf W}_{canon}^{\circ}|_{X\times_{M_n}Y(\Gamma)}\otimes {\bf H}=
{\bf E}\mapsto {\bf F}=\pi_g^{\ast}\underline{\bf W}_{canon}\otimes {\bf H}$
 over $\hat{X}_{\Gamma}\times_{M_n}Y(\Gamma)$. 

\medskip

Then by proposition \ref{prop; equivalent}, the push-forward of the  
localized contribution of the top Chern class into $X\times_{M_n}Y(\Gamma)$ under $\pi_h^{\ast}$ 

$$\hskip -.7in
\pi_{h\ast}\{c_{total}(\pi_h^{\ast}(\pi_g^{\ast}\underline{\bf W}_{canon}\otimes {\bf H})
\otimes_{\Gamma'\in \bar{I}_{\Gamma}-\bar{I}_{\Gamma}^{\gg}} 
{\cal O}(-\hat{D}_{\Gamma'}))\cap s_{total}(\pi_r(\hat{D}_{\Gamma}), 
{\bf C}_{\hat{X}_{\Gamma}\times_{M_n}Y(\Gamma)}\hat{X}_{\Gamma})\}_{dim_{\bf C}X-
rank_{\bf C}\underline{\bf W}_{canon}}$$

 into \footnote{Remember that 
$rank_{\bf C}{\bf W}_{canon}^{\circ}=rank_{\bf C}\underline{\bf W}_{canon}$!} 
${\cal A}_{\cdot}(X\times_{M_n}Y(\Gamma))$
is numerically equivalent ($\stackrel{n}{=}$) to the top Chern class of 
$\pi_g^{\ast}{\bf W}_{canon}^{\circ}|_{Y(\Gamma)\times T(M)}\otimes {\bf H}$,

 By the definition of numerical equivalence 
(see page \pageref{defin; numerical}, definition \ref{defin; numerical}), the push-forward to 
 ${\cal A}_0(pt)$ of their pairings to arbitrary
 complementary dimension cycle classes in $X\times_{M_n}Y(\Gamma)$ are identical. So the original
 top intersection pairing can be replaced by the much simplified version,

$$\hskip -.7in
=c_{top}(\pi_g^{\ast}{\bf W}_{canon}^{\circ}|_{Y(\Gamma)\times T(M)}\otimes {\bf H})\cap
c_{top}(\pi_g^{\ast}\tau_{\Gamma}\otimes {\bf H})\cap
 c_1({\bf H})^{dim_{\bf C}M_n+rank_{\bf C}({\bf V}_{canon}-{\bf W}_{canon})+q-1},$$

 for the purpose of evaluating their push-forward to ${\cal A}_0(pt)$.

\bigskip

$\diamondsuit$
 To summarize, we have succeeded in casting the original intersection pairing to one on
 the smooth space $X\times_{M_n}Y(\Gamma)$ which only involves the top Chern classes of
 $\pi_g^{\ast}{\bf W}_{canon}^{\circ}|_{Y(\Gamma)\times T(M)}\otimes {\bf H}$ and 
 $\pi_g^{\ast}\tau_{\Gamma}\otimes {\bf H}$ and the cycle class $c_1({\bf H})$.

\medskip

Step Three: Finally
 we are ready to identify the last expression in Step Two with the mixed algebraic family invariant.
 Recall the tensor product formula of the top Chern class,
 $$c_{top}(\pi_g^{\ast}\tau_{\Gamma}\otimes {\bf H})=\sum_{0\leq t\leq rank_{\bf C}\tau_{\Gamma}}
  c_t(\pi_g^{\ast}\tau_{\Gamma})\cap c_1({\bf H})^{rank_{\bf C}\tau_{\Gamma}-t}.$$

  If we insert this identity into the final expression in Step Two, we get 

$$\hskip -1.3in
=\sum_{0\leq t\leq rank_{\bf C}\tau_{\Gamma}}
c_{top}(\pi_g^{\ast}{\bf W}_{canon}^{\circ}|_{Y(\Gamma)\times T(M)}\otimes {\bf H})\cap
c_{t}(\pi_g^{\ast}\tau_{\Gamma})\cap 
 c_1({\bf H})^{rank_{\bf C}\tau_{\Gamma}-t+
dim_{\bf C}M_n+rank_{\bf C}({\bf V}_{canon}-{\bf W}_{canon})+q-1}.$$

 Recall that from lemma 6 of [Liu5], 
$((\Phi_{{\bf V}_{canon}^{\circ}{\bf W}_{canon}^{\circ}}, {\bf V}_{canon}^{\circ}, 
{\bf W}_{canon}^{\circ})$ with ${\bf V}_{canon}^{\circ}={\bf V}_{canon}$ is the canonical
 algebraic family Kuranishi model of the class $C-{\bf M}(E)E-
\sum_{e_i\cdot (C-{\bf M}(E)E)<0}e_i$ over the space $M_n\times T(M)$.
 As we have pointed out on page \pageref{proof=mo} 
at the beginning of the current proof that after adding the ``correction term''
 $rank_{\bf C}\tau_{\Gamma}$, 

$$rank_{\bf C}\tau_{\Gamma}+dim_{\bf C}M_n+rank_{\bf C}({\bf V}_{canon}-{\bf W}_{canon})+q-1$$
$$=dim_{\bf C}Y(\Gamma)+rank_{\bf C}({\bf V}_{canon}^{\circ}-{\bf W}_{canon}^{\circ})+q-1,$$

 is nothing but the expected family dimension of the new class $C-{\bf M}(E)E-
\sum_{e_i\cdot (C-{\bf M}(E)E)<0}e_i$ over $Y(\Gamma)\times T(M)$. 

 Therefore the push-forward of top intersection pairing 

$$c_{top}(\pi_g^{\ast}{\bf W}_{canon}^{\circ}|_{Y(\Gamma)\times T(M)}\otimes {\bf H})\cap
c_{t}(\pi_g^{\ast}\tau_{\Gamma})\cap 
 c_1({\bf H})^{rank_{\bf C}\tau_{\Gamma}-t+
dim_{\bf C}M_n+rank_{\bf C}({\bf V}_{canon}-{\bf W}_{canon})+q-1}$$

 in $X\times_{M_n}Y(\Gamma)={\bf P}({\bf V}_{canon})\times_{M_n}Y(\Gamma)=
{\bf P}({\bf V}_{canon}^{\circ})\times_{M_n}Y(\Gamma)$ to ${\cal A}_0(pt)$ is equal to
 the mixed algebraic family Seiberg-Witten invariant 
 $${\cal AFSW}_{M_{n+1}\times_{M_n}Y(\Gamma)\times T(M)
\mapsto Y(\Gamma)\times T(M)}(c_t(\tau_{\Gamma}), 
 C-{\bf M}(E)E-\sum_{e_i\cdot (C-{\bf M}(E)E)<0}e_i).$$

 Thus the total summation over $t$ is

$$\sum_{0\leq t}{\cal AFSW}_{M_{n+1}\times_{M_n}Y(\Gamma)\times T(M)\mapsto
 Y(\Gamma)\times T(M)}(c_t(\tau_{\Gamma}), C-{\bf M}(E)E-\sum_{e_i\cdot 
(C-{\bf M}(E)E)<0}e_i)$$

$$={\cal AFSW}_{M_{n+1}\times_{M_n}Y(\Gamma)\times T(M)\mapsto
 Y(\Gamma)\times T(M)}(c_{total}(\tau_{\Gamma}), C-{\bf M}(E)E-\sum_{e_i\cdot 
(C-{\bf M}(E)E)<0}e_i).$$

We are done with the Case I!

\bigskip

\noindent Case II: If there exists a type $I$ exceptional class $e_{k_i}$
 such that $0>e_{k_i}\cdot (C-{\bf M}(E)E)>e_{k_i}^2$, $\tau_{\Gamma}$ 
 has been defined to be zero in section \ref{subsection; modinv} on page \pageref{tau0}. 
In this case we derive a vanishing result on the top ($=
dim_{\bf C}X-rank_{\bf C}{\bf W}_{canon}$)
 intersection pairing of $c_1({\bf H})$ with the
 localized contribution of top Chern class. It is well known that if the total grading of
an intersection pairing of characteristic classes exceeds the dimension of the space, the intersection
 pairing is equal to zero. Our goal is to show that the cap product of 
the localized contribution of the top Chern class with
 $c_1({\bf H})^{dim_{\bf C}M_n+rank_{\bf C}({\bf V}_{canon}-
{\bf W}_{canon})+q-1}$ vanishes due to dimension count.
 
For notational simplicity, we assume that $e_{k_1}\cdot (C-{\bf M}(E)E)>e_{k_1}^2$.
 That is to say,
 we take $i=1$. Because our argument 
 only makes usage of the dimension count, we do not
 lose any generality in adopting this convention.

\medskip

Step One: Firstly we derive a lemma which will be used later.

\begin{lemm}\label{lemm; cs=}
 Let $\pi_{\bf F}: {\bf F}\mapsto B$ be a finite rank vector bundle over $B$. Let $s_{\bf F}$ denote
 the zero section embedding $s_{\bf F}:B\mapsto {\bf F}$. Let $r\geq rank_{\bf C}{\bf F}$ be 
a positive integer.  
Then for all $\beta\in {\cal A}_{\ast}(B)$ such that $s_{total}({\bf F})\cap \beta$ has no grade 
$<r$ components, the following identity holds, 

$$s_{\bf F\ast}\{\beta\cap s_{total}({\bf F})\}_r=\{\pi_{\bf F}^{\ast}\beta\}_r.$$
\end{lemm}

\noindent Proof of the lemma: 
For all $\alpha\in {\cal A}_r(B)$, where $r$ is a fixed natural number $\geq rank_{\bf C}{\bf F}$,
 we have (see example 3.3.2. on page 67 of [F])

$$s_{\bf F}^{\ast}s_{{\bf F}\ast}\{\alpha\}_r=c_{rank_{\bf C}{\bf F}}({\bf F})\cap \{\alpha\}_r
=\{c_{total}({\bf F})\cap \alpha\}_{r-rank_{\bf C}{\bf F}}.$$

  One can extend this equality trivially 
to all $\alpha\in {\cal A}_{\geq r}(Y(\Gamma_{e_{k_1}}))$ with grading
 $\geq r$ as they do not contribute to both sides of the identity. 

Therefore by the reciprocity property of the total Segre class and the total Chern class and by
 taking $\alpha=s_{total}({\bf F})\cap \beta$, we find  
$s_{\bf F}^{\ast}s_{\bf F\ast}\{s_{total}({\bf F})\cap \beta\}_r
=\{\beta\}_{r-rank_{\bf C}{\bf F}}$ for all 
$\beta$ satisfying the grading assumption in the lemma. 

And therefore
 $s_{\bf F\ast}\{s_{total}({\bf F})\cap \beta\}_r=\{\pi_{\bf F}^{\ast}\beta\}_r$
 because the Gysin homomorphism satisfies $s_{\bf F}^{\ast}=(\pi_{\bf F}^{\ast})^{-1}$ (please consult page
 65, definition 3.3. of [F]). The lemma is proved. $\Box$

Step Two: The yet-to-be-identified intersection pairing 

$$
\{\sum_{i=1}^{i=rank_{\bf C}{\bf W}_{canon}}(-1)^{i-1}
c_{rank_{\bf C}{\bf W}_{canon}-i}(\pi_{\tilde{\hat{X}}_{\Gamma}}^{\ast}{\bf W}_{canon}\otimes 
 {\bf H}\otimes_{\Gamma'\in \bar{I}_{\Gamma}-\bar{I}_{\Gamma}^{\gg}} 
{\cal O}(-\hat{D}_{\Gamma'}))$$
$$\cap \hat{D}_{\Gamma}^{i-1}\cap [\hat{D}_{\Gamma}]\cap 
c_1({\bf H})^{dim_{\bf C}M_n+rank_{\bf C}({\bf V}_{canon}-{\bf W}_{canon})+q-1}\}$$

 can be pushed-forward
 as a zero dimensional cycle class into ${\cal A}_0({\hat{X}}_{\Gamma}\times_{M_n}
Y(\Gamma))$. Because that $Y(\Gamma)\subset Y(\Gamma_{e_{k_1}})$, 
 $\hat{X}_{\Gamma}\times_{M_n}Y(\Gamma)\subset \hat{X}_{\Gamma}\times_{M_n}Y(\Gamma_{e_{k_1}})$.
Similar to what was done in step one of Case I, we may deform to the projectified normal cone 
 and replace
 $\hat{X}_{\Gamma}\times_{M_n}Y(\Gamma)\subset \hat{X}_{\Gamma}$ by the inclusion into the 
zero section of the projectified normal cone of $\hat{X}\times_{M_n}Y(\Gamma_{e_{k_1}})$, 

$$\hat{X}_{\Gamma}\times_{M_n}Y(\Gamma)\subset \hat{X}_{\Gamma}\times_{M_n}
Y(\Gamma_{e_{k_1}})\subset {\bf P}({\bf C}_{\hat{X}_{\Gamma}\times_{M_n}
Y(\Gamma_{e_{k_1}})}\hat{X}_{\Gamma}\oplus 1).$$

Correspondingly, we twist the obstruction vector bundle $\pi_f^{\ast}\pi_g^{\ast}{\bf W}_{canon}$ by 
\footnote{The twisting of this line bundle ${\cal O}(\hat{\bf P}_{\infty})$ does not
 play an essential role in our argument, its presence 
only makes the notations slightly more
complicated. Nevertheless we do not remove it in Case II 
as we do not always keep our
cycle disjoint from $\hat{\bf P}_{\infty}$.} ${\cal O}(\hat{\bf P}_{\infty})$.

Then the derived exact sequence of locally free sheaves
$$\hskip -.3in
{\cal R}^0\pi_{\ast}\bigl({\cal O}_{{\bf M}(E)E+\Xi_{k_1}}\otimes {\cal E}_C\bigr)\mapsto 
{\cal R}^0\pi_{\ast}\bigl({\cal O}_{{\bf M}(E)E}\otimes {\cal E}_C\bigr)\mapsto 
{\cal R}^1\pi_{\ast}\bigl({\cal O}_{\Xi_{k_1}}\otimes {\cal E}_{C-{\bf M}(E)E}\bigr)$$

 induces a short exact sequence analogous to the short exact sequence on
 page \pageref{underlinewcanon}.
One can interpret the construction of this new sequence as a special case of the
 general construction once we ``formally'' consider $e_{k_1}$ to be the unique
 type $I$ exceptional class which pairs negatively with $C-{\bf M}(E)E$. For notational simplicity,
 we still denote the corresponding sequence of vector bundles by the same notation 
\footnote{In our current argument only the ranks of these bundles matter. As far as we
 do not use specific properties of these bundles, the slight abuse of
 notations does not cause trouble.}
as before,  
$0\mapsto \underline{\bf W}_{canon}\mapsto {\bf W}_{canon}\mapsto
 {\bf V}_{quot}\mapsto 0$. 
 This sequence
 breaks ${\bf W}_{canon}|_{Y(\Gamma)\times T(M)}$ into the factors $\underline{\bf W}_{canon}$ and
 ${\bf V}_{quot}$. In the current context, the symbol ${\bf V}_{quot}$
 means the vector bundle associated with the locally free
summand of ${\cal R}^1\pi_{\ast}\bigl({\cal O}_{\Xi_{k_1}}\otimes {\cal E}_{C-{\bf M}(E)E}\bigr)$.
Then by proposition \ref{prop; cap},  
we may identify the above expression with 

$$\hskip -1in
\{c_{total}((\pi_g\pi_f)^{\ast}
(\underline{\bf W}_{canon}\otimes {\cal O}(\hat{\bf P}_{\infty})
\otimes {\bf H})\otimes_{\Gamma'\in \bar{I}_{\Gamma}-\bar{I}_{\Gamma}^{\gg}} 
{\cal O}(-\pi_f^{\ast}\hat{D}_{\Gamma'}))\cap s_{total}(\pi_r(\hat{D}_{\Gamma}), 
{\bf C}_{\hat{X}_{\Gamma}\times_{M_n}Y(\Gamma)}\hat{X}_{\Gamma})\}_{dim_{\bf C}X-
\underline{\bf W}_{canon}}$$
$$\hskip -1in
\cap c_{top}(\pi_f^{\ast}(\pi_g^{\ast}{\bf V}_{quot}\otimes 
{\cal O}(\hat{\bf P}_{\infty})\otimes {\bf H}))
\cap c_1({\bf H})^{dim_{\bf C}M_n+rank_{\bf C}({\bf V}_{canon}-{\bf W}_{canon})+q-1}.$$

 Now we push forward this grade zero cycle class into 
${\cal A}_0({\bf P}({\bf N}_{X\times_{M_n}Y(\Gamma_{e_{k_1}})}X\oplus 1))$ by 
\footnote{Here $\pi_{{\bf P}}:{\bf P}({\bf C}_{\hat{X}_{\Gamma}\times_{M_n}
Y(\Gamma_{e_{k_1}})}\hat{X}\oplus 1)\mapsto {\bf P}({\bf N}_{X\times_{M_n}
Y(\Gamma_{e_{k_1}})}X\oplus 1)$ is the projection map.}
 $\pi_{{\bf P}\ast}$. Take 
${\bf P}_{\infty}={\bf P}({\bf N}_{X\times_{M_n}Y(\Gamma_{e_{k_1}})}X)$ to be the
 divisor at infinity in ${\bf P}({\bf N}_{X\times_{M_n}Y(\Gamma_{e_{k_1}})}X\oplus 1)$.
By proposition \ref{prop; rif} and the observation/argument used 
in proposition \ref{prop; identical}, the 
$\pi_{{\bf P}\ast}$ 
push-forward of the localized contribution of top Chern class $\{c_{total}((\pi_g\pi_f)^{\ast}
(\underline{\bf W}_{canon}\otimes {\cal O}({\bf P}_{\infty})
\otimes {\bf H})\otimes_{\Gamma'\in \bar{I}_{\Gamma}-\bar{I}_{\Gamma}^{\gg}} 
{\cal O}(-\pi_f^{\ast}\hat{D}_{\Gamma'}))\cap s(\pi_r(\hat{D}_{\Gamma}), 
{\bf C}_{\hat{X}_{\Gamma}\times_{M_n}Y(\Gamma)}\hat{X}_{\Gamma})\}_{dim_{\bf C}X-
\underline{\bf W}_{canon}}$
 can be written as the difference of two localized contributions of top Chern class,

$$
\{c_{total}(\pi_g^{\ast}
\underline{\bf W}_{canon}\otimes {\cal O}({\bf P}_{\infty})
\otimes {\bf H})\cap s_{total}(Z(\underline{s}_{canon})\times_{M_n}Y(\Gamma), 
{\bf N}_{X\times_{M_n}Y(\Gamma_{e_{k_1}})}X)$$

$$\hskip -1.1in
-c_{total}(\pi_g^{\ast}\underline{\bf W}_{canon}\otimes {\cal O}({\bf P}_{\infty})
\otimes {\bf H})
\cap s_{total}(\cup_{\Gamma'\in \bar{I}_{\Gamma}-\bar{I}_{\Gamma}^{\gg}}
Z(\underline{s}_{canon})\times_{M_n}(Y(\Gamma')\cap Y(\Gamma)), 
{\bf N}_{X\times_{M_n}Y(\Gamma_{e_{k_1}})}X)\}_{dim_{\bf C}X-
\underline{\bf W}_{canon}},$$

 where the section $\underline{s}_{canon}$
 is the induced 
section of $\pi_g^{\ast}\underline{\bf W}_{canon}\otimes {\bf H}$ by

 $$\pi_X^{\ast}{\bf W}_{canon}^{\circ}|_{X\times_{M_n}Y(\Gamma)}
\mapsto \underline{\bf W}_{canon}$$ and 
 $s_{canon}^{\circ}$ and the 
sub-schemes $Z(\underline{s}_{canon})\times_{M_n}Y(\Gamma)$
 and $\cup_{\Gamma'\in \bar{I}_{\Gamma}-\bar{I}_{\Gamma}^{\gg}}
Z(\underline{s}_{canon})\times_{M_n}(Y(\Gamma')\cap Y(\Gamma))$ are embedded in the zero cross section 
 $X\times_{M_n}Y(\Gamma_{e_{k_1}})\subset {\bf N}_{X\times_{M_n}Y(\Gamma_{e_{k_1}})}X$.

Set $Z_1=Z(\underline{s}_{canon})\times_{M_n}Y(\Gamma)$ 
and $Z_2=\cup_{\Gamma'\in \bar{I}_{\Gamma}-\bar{I}_{\Gamma}^{\gg}}
Z(\underline{s}_{canon})\times_{M_n}(Y(\Gamma')\cap Y(\Gamma))$ and 
 $i_{Z_a}:Z_a\subset {\bf P}({\bf N}_{X\times_{M_n}Y(\Gamma_{e_{k_1}})}X\oplus 1)$ for $a=1, 2$.
 
We will give a uniform argument for both
 $a=1, 2$ that the
 top intersection pairings

$$
i_{Z_a\ast}\{c_{total}(\pi_g^{\ast}\underline{\bf W}_{canon}\otimes
 {\cal O}({\bf P}_{\infty})\otimes 
{\bf H})\cap s_{total}(Z_a, {\bf N}_{X\times_{M_n}
Y(\Gamma_{e_{k_1}})}X)\}_{dim_{\bf C}X-rank_{\bf C}
\underline{\bf W}_{canon}}$$
$$\cap c_{top}(\pi_{{\bf N}_{X\times_{M_n}Y(\Gamma_{e_{k_1}})}X}^{\ast}
\pi_g^{\ast}{\bf V}_{quot}\otimes {\cal O}({\bf P}_{\infty})\otimes {\bf H})\cap 
c_1(\pi_{{\bf N}_{X\times_{M_n}Y(\Gamma_{e_{k_1}})}X}^{\ast}
{\bf H})^{dim_{\bf C}M_n+
rank_{\bf C}({\bf V}_{canon}-{\bf W}_{canon})+q-1}$$

 vanish identically.

\medskip

Step Three: Because $e_{k_1}^2<e_{k_1}\cdot (C-{\bf M}(E)E)$, 
one observes that the expected dimension of
 the family moduli space of $C-{\bf M}(E)E-e_{k_1}$ over $Y(\Gamma_{e_{k_1}})$,

$$
dim_{\bf C}Y(\Gamma_{e_{k_1}})+p_g+{(C-{\bf M}(E)E-e_{k_1})^2-
c_1({\bf K}_{M_{n+1}/M_n})\cdot
 (C-{\bf M}(E)E-e_{k_1})\over 2}$$
$$\hskip -.9in =dim_{\bf C}M_n+p_g
+{e_{k_1}^2-c_1({\bf K}_{M_{n+1}/M_n})\cdot e_{k_1}
\over 2}+{(C-{\bf M}(E)E-e_{k_1})^2-c_1({\bf K}_{M_{n+1}/M_n})\cdot
 (C-{\bf M}(E)E-e_{k_1})\over 2}$$

$$<dim_{\bf C}M_n+p_g+{(C-{\bf M}(E)E)^2-c_1({\bf K}_{M_{n+1}/M_n})\cdot
 (C-{\bf M}(E)E)\over 2},$$

 strictly smaller than 
the expected family dimension of the class $C-{\bf M}(E)E$ over the family $M_n$. 
We use this observation to derive the vanishing result.

 Define $B=X\times_{M_n}Y(\Gamma_{e_{k_1}})$, 
${\bf F}={\bf N}_{X\times_{M_n}Y(\Gamma_{e_{k_1}})}X$,
 $\iota_{\bf F}:{\bf F}\subset {\bf P}({\bf F}\oplus 1)$
 and $j_{Z_a}:Z_a\subset X\times_{M_n}Y(\Gamma_{e_{k_1}})=B$. Then 
 $i_{Z_a}$ can be factorized as $\iota_{\bf F}\circ s_{\bf F}\circ j_{Z_a}$.
 Because both $Z_a\subset X\times_{M_n}Y(\Gamma)\subset X\times_{M_n}Y(\Gamma_{e_{k_1}})$, 
there is a short 
exact sequence of normal cones (see example 4.1.6. of [F] for its definition),

$$0\mapsto j_{Z_a}^{\ast}{\bf N}_{X\times_{M_n}Y(\Gamma_{e_{k_1}})}X\mapsto 
{\bf C}_{Z_a}{\bf N}_{X\times_{M_n}Y(\Gamma_{e_{k_1}})}X\mapsto {\bf C}_{Z_a}X\times_{M_n}
Y(\Gamma_{e_{k_1}})\mapsto 0.$$

By the product property of the total Segre classes for short exact sequences
 of cones, 
the final expression in Step Two can be re-casted into

$$\hskip -1in
\iota_{{\bf F}\ast}s_{{\bf F}\ast}j_{Z_a\ast}\{c_{total}(\pi_g^{\ast}\underline{\bf W}_{canon}
\otimes {\cal O}({\bf P}_{\infty})\otimes 
{\bf H})\cap s_{total}(Z_a, X\times_{M_n}Y(\Gamma_{e_{k_1}}))
\cap s_{total}(j_{Z_a}^{\ast}{\bf F})\}_{dim_{\bf C}X-rank_{\bf C}
\underline{\bf W}_{canon}}$$
$$\cap c_{top}(\pi_{{\bf P}({\bf F}\oplus 1)}^{\ast}
(\pi_g^{\ast}{\bf V}_{quot}\otimes {\bf H}) \otimes {\cal O}({\bf P}_{\infty}))\cap 
c_1(\pi_{{\bf P}({\bf F}\oplus 1)}^{\ast}
{\bf H})^{dim_{\bf C}M_n+rank_{\bf C}({\bf V}_{canon}-{\bf W}_{canon})+q-1}$$

$$\hskip -.8in
=\iota_{{\bf F}\ast}\{s_{{\bf F}\ast}\{j_{Z_a\ast}(c_{total}(
\pi_g^{\ast}\underline{\bf W}_{canon}\otimes 
{\bf H})\cap s_{total}(Z_a, X\times_{M_n}Y(\Gamma_{e_{k_1}})))\cap s_{total}({\bf F})
\}_{dim_{\bf X}-rank_{\bf C}\underline{\bf W}_{canon}}$$
$$\cap c_{top}(\pi_{\bf F}^{\ast}(\pi_g^{\ast}{\bf V}_{quot}\otimes {\bf H})\otimes 
{\cal O}({\bf P}_{\infty}))\cap 
c_1(\pi_{\bf F}^{\ast}{\bf H})^{dim_{\bf C}M_n+rank_{\bf C}({\bf V}_{canon}-{\bf W}_{canon})+q-1}\}.$$

Define $\beta=j_{Z_a\ast}(c_{total}(\pi_g^{\ast}\underline{\bf W}_{canon}\otimes 
{\bf H})\cap s_{total}(Z_a, X\times_{M_n}Y(\Gamma_{e_{k_1}})))\in {\cal A}_{\cdot}(X\times_{M_n}
Y(\Gamma_{e_{k_1}}))$ and set $r=dim_{\bf C}X-rank_{\bf C}\underline{\bf W}_{canon}$.
  Then $\{\beta\cap s_{total}({\bf F})\}_r$ is the push-forward of the localized contribution
 of top Chern class over $Z_a\subset {\bf F}$ 
into ${\cal A}_r(X\times_{M_n}Y(\Gamma_{e_{k_1}}))$. Then by e.g. 
proposition 13 of \footnote{We had used this fact earlier in the proof of 
 proposition \ref{prop; cap}, too.}
 [Liu5], $\{\beta\cap s_{total}({\bf F})\}_s=0$ for all $s<r$. Namely,
 the localized contribution of top Chern class is the leading (lowest grading) term of the 
intersection pairing.

 Then the assumption of lemma \ref{lemm; cs=} of Step One is applicable and we know
 $s_{{\bf F}\ast}\{\beta\cap s_{total}({\bf F})\}_r=\{\pi_{\bf F}^{\ast}\beta\}_r$.
 Because the bundle projection $\pi_{\bf F}:{\bf F}\mapsto B=X\times_{M_n}Y(\Gamma_{e_{k_1}})$
 is flat of relative dimension $rank_{\bf C}{\bf F}$, 
and the inclusion $\iota_{\bf F}$ is a proper morphism, 
 and by theorem 3.2.(c)-(d). on pages 50-51 of [F],
 and the fact that the flat pull-back $\pi_{\bf F}^{\ast}:{\cal A}_{r-rank_{\bf C}{\bf F}}(B)\mapsto
 {\cal A}_r({\bf F})$ lifts the gradings up by $rank_{\bf C}{\bf F}$, 
we may rewrite the above intersection pairing as
\footnote{After we introduce $\pi_{\bf F}^{\ast}$ into our formulae,
 the cycle class is not ``refined'' in $X\times_{M_n}Y(\Gamma_{e_{k_1}})$
 any more!}

$$\hskip -1.5in
\iota_{{\bf F}\ast}\{\pi_{\bf F}^{\ast}\{\beta\}_{
dim_{\bf C}X-rank_{\bf C}\underline{\bf W}_{canon}-rank_{\bf C}{\bf F}}\cap 
c_1(\pi_{\bf F}^{\ast}{\bf H})^{dim_{\bf C}M_n+rank_{\bf C}({\bf V}_{canon}-{\bf W}_{canon})+q-1}\cap 
c_{top}(\iota_{\bf F}^{\ast}\bigl(\pi_{{\bf P}({\bf F}\oplus 1)}^{\ast}(\pi_g^{\ast}{\bf V}_{quot}\otimes
 {\bf H})\otimes 
{\cal O}({\bf P}_{\infty})\bigr))\}$$

$$\hskip -1.3in
=\iota_{{\bf F}\ast}\{\pi_{\bf F}^{\ast}(\{\beta\}_{
dim_{\bf C}X-rank_{\bf C}\underline{\bf W}_{canon}-rank_{\bf C}{\bf F}}\cap 
c_1({\bf H})^{dim_{\bf C}M_n+rank_{\bf C}({\bf V}_{canon}-{\bf W}_{canon})+q-1})\cap 
c_{top}(\iota_{\bf F}^{\ast}
\pi_{{\bf P}({\bf F}\oplus 1)}^{\ast}(\pi_g^{\ast}{\bf V}_{quot}\otimes {\bf H})\otimes 
{\cal O}({\bf P}_{\infty}))\}$$

$$\hskip -1.3in
=\iota_{{\bf F}\ast}\{\pi_{\bf F}^{\ast}(\{\beta\}_{
dim_{\bf C}Y(\Gamma_{e_{k_1}})-rank_{\bf C}\underline{\bf W}_{canon}}\cap 
c_1({\bf H})^{dim_{\bf C}X-rank_{\bf C}{\bf W}_{canon}})
\cap c_{top}(i_{\bf F}^{\ast}
\pi_{{\bf P}({\bf F}\oplus 1)}^{\ast}(\pi_g^{\ast}{\bf V}_{quot}\otimes {\bf H})\otimes
{\cal O}({\bf P}_{\infty}))\}$$

$$\hskip -1.3in
=\iota_{{\bf F}\ast}\{
\pi_{\bf F}^{\ast}\{\beta\cap c_1({\bf H})^{dim_{\bf C}X-rank_{\bf C}{\bf W}_{canon}}\}_{
dim_{\bf C}Y(\Gamma_{e_{k_1}})-rank_{\bf C}\underline{\bf W}_{canon}-
(dim_{\bf C}X-rank_{\bf C}{\bf W}_{canon})}\}\cap 
c_{top}(\pi_{{\bf P}({\bf F}\oplus 1)}^{\ast}(\pi_g^{\ast}{\bf V}_{quot}\otimes {\bf H})\otimes
{\cal O}({\bf P}_{\infty})).$$

We have used $dim_{\bf C}X-rank_{\bf C}{\bf F}=dim_{\bf C}Y(\Gamma_{e_{k_1}})$
 in the above derivation.

 Yet the grading $dim_{\bf C}Y(\Gamma_{e_{k_1}})-rank_{\bf C}
\underline{\bf W}_{canon}-
(dim_{\bf C}X-rank_{\bf C}{\bf W}_{canon})$ of $\{\bullet\}$ 
is exactly the difference between
 the expected family dimension of the class $C-{\bf M}(E)E-e_{k_1}$ over $Y(\Gamma_{e_{k_1}})$
 and \footnote{We have used 
the observation that $rank_{\bf C}{\bf W}_{canon}^{\circ}=
rank_{\bf C}\underline{\bf W}_{canon}$ implicitly.}
the expected family dimension of the class $C-{\bf M}(E)E$ over $M_n$.

As we assume $e_{k_1}^2<e_{k_1}\cdot (C-{\bf M}(E)E)<0$, we have already shown at the beginning of 
Step Three that this grading is negative. Therefore 

$$\pi_{\bf F}^{\ast}\{\beta\cap c_1({\bf H})^{dim_{\bf C}X-rank_{\bf C}{\bf W}_{canon}}\}_{
{\bf C}Y(\Gamma_{e_{k_1}})-rank_{\bf C}\underline{\bf W}_{canon}-
(dim_{\bf C}X-rank_{\bf C}{\bf W}_{canon})}=0$$

 and therefore the whole intersection pairing vanishes. In particular its push-forward into
 ${\cal A}_0(pt)\cong {\bf Z}$ is zero. As this holds for both $Z_1$ and $Z_2$, the
 original intersection pairing (their difference) is also zero.
 We are done with Case II.

 As we have finished the identification with the mixed algebraic family Seiberg-Witten invariants
 in both cases, we have finished the proof of proposition \ref{prop; =mo}.

$\Box$

\medskip

\begin{rem}\label{rem; restrict}
In the proof of the proposition \ref{prop; =mo}, 
we only discuss the non-restricted case. 
 If one specifies a point $t_L\in T(M)$ and would like to count curves in 
 ${\cal M}_{C-{\bf M}(E)E}$ whose images in $M$ are in
 the linear system $|L|$ specified by the point $t_L$,
 there are two viewpoints one can
 adopt. 

(i). By restricting to a point $t_L \in T(M)$, effectively one shrinks 
 $T(M)$ to a point. One can think of this procedure as a formal reduction 
of the irregularity $q\rightarrow 0$ and
 the rest of the deduction is identical to the $q=0$ case, where $T(M)$
 does not play any role here.

\medskip

(ii). Alternatively, one may replace the family moduli space 
${\cal M}_{C-{\bf M}(E)E}$, the total projective space bundle 
$X={\bf P}({\bf V}_{canon})$, etc., by their $t_L-$restricted counterparts,
${\cal M}_{C-{\bf M}(E)E}\times_{T(M)}\{t_L\}$ and $X\times_{T(M)}\{t_L\}$, respectively.
One may insert the cycle class $[t_L]\in {\cal A}_0(T(M))$ into the
 intersection theory product and therefore replace the power
$c_1({\bf H})^{dim_{\bf C}M_n+rank_{\bf C}({\bf V}_{canon}-{\bf W}_{canon})+q-1}$  
by $[t_L]\cap c_1({\bf H})^{dim_{\bf C}M_n+rank_{\bf C}({\bf V}_{canon}-{\bf W}_{canon})-1}$ 
and the rest of the discussion goes through without any change.
 By either angles the reader should be able to make the suitable adjustments in all the formulae
 and finish the proof. We do not repeat the redundant details here
\footnote{Please compare with remark \ref{rem; insert}, located right 
after the definition of the  modified algebraic family invariants.}.
\end{rem}

\medskip

\subsection{The Proof of the Main Theorem}\label{subsection; proofmain}

\bigskip

 We are ready to combine all the results proved in the paper to
 prove the main theorem of the paper,

\begin{theo}\label{theo;  main1}
Let $\delta\in {\bf N}$ denote \footnote{In the proof of the main theorem, 
we switch from $n$ to $\delta$, fitting to G$\ddot{o}$ttsche's convention.}
 the number of nodal singularities.
Let $L$ be a $5\delta-1$ very-ample line bundle on an algebraic surface $M$,
 then the number of $\delta$ nodes nodal singular curves in a generic 
 $\delta$ dimensional linear sub-system of $|L|$ can be expressed as 
 a universal polynomial (independent to $M$) 
of $c_1(L)^2$, $c_1(L)\cdot c_1(M)$, $c_1(M)^2$, 
$c_2(M)$ of degree $\delta$.
\end{theo}

\medskip

 For the invertible sheaf ${\cal L}={\cal E}_C\mapsto M\times T(M)$ 
parametrized by a cohomology class 
$C\in H^{1, 1}(M, {\bf Z})$, we may extend the definition of $k-$very ampleness
 by assuming the surjectivity of the restriction morphism 
${\cal R}^0\pi_{T(M)\ast}\bigl({\cal L}\bigr)\mapsto
 {\cal R}^0\pi_{T(M)\ast}\bigl({\cal L}\otimes 
{\cal O}_{Z\times T(M)}\bigr)$ for all
 length $k+1$ sub-schemes $Z\subset M$.

\begin{rem}\label{rem; main2}
Let $\delta\in {\bf N}$ denote the number of nodal singularities.
Let $C$ be a cohomology class in $H^{1, 1}(M, {\bf Z})$ and let 
 ${\cal L}\mapsto M\times T(M)$ be the invertible sheaf with $c_1(i_M^{\ast}{\cal L})=C$, where
 $i_M:M\times \{0\}\subset M\times T(M)$. Suppose that ${\cal L}$ is 
$5\delta-1-$very ample and one can find generic $\delta$ dimensional non-linear sub-system of
 the projective space bundle 
${\bf P}(\pi_{T(M)\ast}\bigl({\cal L}\bigr)^{\ast}) \cong {\bf P}({\bf V}_{canon})$, then
 one may formulate a corresponding theorem for ${\cal L}$, parallel to theorem \ref{theo; main1}. 
The universal polynomial 
 associated to ${\cal L}$ is the product of the universal polynomial found in 
theorem \ref{theo; main1}
 and\footnote{Refer to remark \ref{rem; poincare}.} ${\cal ASW}(C)$. 
\end{rem}

\noindent Proof of the main theorem: Let $L$ be a line bundle over $M$ with $c_1(L)=C$, then
 $L$ determines a unique point $t_L\in T(M)$ in the connected component of the Picard variety.
 As usual $T(M)$ represents the component of 
 Picard group of $M$ parametrizing the line bundles with first Chern class $C$. Let $m_1=m_2=\cdots=
m_{\delta}=2$ and let ${\cal M}_{C-{\bf M}(E)E}$ denote the algebraic family moduli space of 
 curves dual to $C-2\sum_{i\leq \delta}E_i$ which projects to $M_{\delta}\times
T(M)$. Then ${\cal M}_{C-{\bf M}(E)E}\times_{T(M)}\{t_L\}$ is the sub-moduli 
space of curves whose projection into $M$ lie in $|L|$. Let ${\cal M}_V$ denote the pre-image of $V$, a
 general $\delta$ dimensional linear subsystem of $|L|$ under the
 projection map ${\cal M}_{C-{\bf M}(E)E}\times_{T(M)}\{t_L\}\mapsto |L|$, 
 then ${\cal M}_V$ can be viewed as the
${L^2-{\bf K}_M\cdot L\over 2}-q(M)+p_g(M)-\delta$-fold generic hyperplane intersection
 of $|L|={\bf P}({\bf V}_{canon})\times_{T(M)}\{t_L\}$, intersecting with
 the family moduli space of $C-{\bf M}(E)E$,
${\cal M}_{C-{\bf M}(E)E}$. It can be also viewed as
 the set theoretical intersection resembling the 
following intersection theoretical
 product,

$$c_{top}(\pi_X^{\ast}{\bf W}_{canon}\otimes {\bf H})\cap 
c_1({\bf H})^{dim_{\bf C}M_{\delta}+{L^2-{\bf K}_M\cdot L\over 2}-q(M)+p_g(M)-3\delta}\cap
 [X\times_{T(M)}t_L]$$

where ${\cal M}_{C-{\bf M}(E)E}$ is represented by the top Chern class of the
 canonical algebraic obstruction vector bundle $\pi_X^{\ast}{\bf W}_{canon}\otimes
 {\bf H}$ over $X={\bf P}({\bf V}_{canon})$ and $[X\times_{T(M)}t_L]\in {\cal A}_{dim_{\bf C}X-q}(X)$ 
is the fiber cycle class determined by the point $t_L$.

 This object is nothing but the mixed algebraic family Seiberg-Witten invariant of 
 $C-\sum 2E_i$, with an additional $[t_L]$ inserted, to restrict
 ${\cal L}$ to $L$, i.e. ${\cal AFSW}_{M_{\delta+1}\times \{t_L\}\mapsto 
 M_{\delta}\times \{t_L\}}(1, C-2\sum_{1\leq i\leq \delta}E_i)$.
 
By the discussion presented in subsection \ref{subsection; transv} below, if we choose 
 the linear subsystem $V$ generically, then the set ${\cal M}_V$ 
 can be decomposed into a portion over $S_{\gamma_{\delta}}$, 
${\cal M}_V\times_{M_{\delta}}(M_{\delta}-
\cup_{\Gamma\in \Delta(\delta)-\{\gamma_{\delta}\}}Y_{\Gamma})$ 
and the excess component ${\cal M}_V\times_{M_{\delta}}
(\cup_{\Gamma\in \Delta(\delta)-\{\gamma_{\delta}\}}Y_{\Gamma})$.

 By proposition \ref{prop; tran} the $5\delta-1$-very ampleness condition on $L$ implies that
 ${\cal M}_V\times_{M_{\delta}}(M_{\delta}-
\cup_{\Gamma\in \Delta(\delta)-\{\gamma_{\delta}\}}
Y_{\Gamma})$ has the structure of a finite scheme which maps into
 the generic stratum $Y_{\gamma_{\delta}}$. 
I.e. its image will miss all those $Y_{\Gamma}$ associated with fan-like 
$\Gamma\in adm_2(\delta)$. In particular, both ${\cal M}_V\times_{M_{\delta}}(M_{\delta}-
\cup_{\Gamma\in \Delta(\delta)}Y_{\Gamma})$ and ${\cal M}_V\times_{M_{\delta}}
(\cup_{\Gamma\in \Delta(\delta)}Y_{\Gamma})$ are closed sub-schemes of $X$.

 We emphasize that we do {\bf NOT} use the very ampleness condition
on $L$ to gain any regularity of the sub-scheme 
 ${\cal M}_V\times_{M_{\delta}}
(\cup_{\Gamma\in \Delta(\delta)-\{\gamma_{\delta}\}}Y_{\Gamma})$. Instead, the
 machineries developed earlier in this paper, namely 
residual intersection formula of
 top Chern classes, recursive blowing ups of $X={\bf P}({\bf V}_{canon})$ 
 in subsection \ref{subsection; blowup}
 and proposition \ref{prop; =mo}, remark \ref{rem; restrict}, etc.,
 allows us to identify through an induction argument
 the intersection numbers represented by ${\cal M}_V\times_{M_{\delta}}Y_{\gamma_{\delta}}$, 
i.e. the top intersection pairing of the push-forward of the
 localized top Chern class of $\pi_X^{\ast}{\bf W}_{canon}\otimes {\bf H}$ over
 ${\cal M}_V\times_{M_{\delta}}Y_{\gamma_{\delta}}$ with a complementary power of $c_1({\bf H})$,
 with the modified family invariant ${\cal AFSW}^{\ast}$ of $C-{\bf M}(E)E$ over $M_{\delta}\times \{t_L\}$,
 namely the difference of ${\cal AFSW}_{M_{\delta+1}\times \{t_L\}\mapsto M_{\delta}\times \{t_L\}}
(1, C-{\bf M}(E)E)$ and the sum of a hierarchy of the
modified mixed algebraic
 family Seiberg-Witten invariants of 
$C-{\bf M}(E)E-\sum_{e_i\cdot (C-{\bf M}(E)E)<0}e_i$ above $Y(\Gamma)$, for various
 $\Gamma\in \Delta(\delta)-\{\gamma_{\delta}\}$.
 Thus we may identify the degree of the 
 finite cycle class $[{\cal M}_V\times_{M_{\delta}}(M_{\delta}-
\cup_{\Gamma\in \Delta(\delta)-\{\gamma_{\delta}\}}Y_{\Gamma})]$ in ${\cal A}_0(pt)\cong {\bf Z}$
 with the modified algebraic family Seiberg-Witten invariant,
 ${\cal AFSW}_{M_{\delta+1}\times T(M)\mapsto 
 M_{\delta}\times T(M)}^{\ast}(\{t_L\}, C-2\sum_{1\leq i\leq \delta}E_i)$, defined following 
 remark \ref{rem; tL} and parallel to definition \ref{defin; generic}.

 According to proposition \ref{prop; polynomial}, remark \ref{rem; insert} and remark \ref{rem; restrict},
 this modified family invariant can be expressed
 as ${\cal ASW}(\{t_L\}, C)$ times a universal degree $\delta$ polynomial of 
 $C^2, C\cdot c_1(M), c_1(M)^2$ and  $c_2(M)$.  Because ${\cal ASW}(\{t_L\}, C)=c_1({\bf H})^{
 p_g-q+{C^2-c_1({\bf K}_M)\cdot C\over 2}}[{\bf P}(H^0(M, L))] \cong 1$ and
 $C^2=L\cdot L$, $C\cdot c_1(M)=-L\cdot {\bf K}_M$, 
$c_1(M)^2={\bf K}_M\cdot {\bf K}_M$, the integer can be expressed as a
 universal polynomial of $L^2, L\cdot {\bf K}_M, {\bf K}_M\cdot {\bf K}_M$
 and $c_2(M)$.
 
On the other hand, the symmetric group of $\delta$ elements, 
${\bf S}_{\delta}$, acts naturally
 and freely upon the open stratum $Y_{\gamma_{\delta}}$, whose underlying set is the set of
 all ordered distinct $\delta$-tuples of points in $M$. 
This free action induces free actions upon
the sub-scheme ${\cal M}_V\times_{M_{\delta}}Y_{\gamma_{\delta}}$ and the smooth ambient
 space $X\times_{M_{\delta}}Y_{\gamma_{\delta}}$. We have the following proposition, whose
proof will be postponed after we have finished the proof of our main theorem.

\begin{prop}\label{prop; factorial}
Assuming that ${\cal M}_V\times_{M_{\delta}}Y_{\gamma_{\delta}}\subset V\times M_{\delta}\times \{t_L\}$ 
($\subset X$) is a finite sub-scheme. 
Then the push-forward of the 
zero cycle $[{\cal M}_V\times_{M_{\delta}}Y_{\gamma_{\delta}}]\in 
{\cal A}_0(X)$ into ${\cal A}_0(pt)$
 is equal to $\delta!$ times ``the number of $\delta$-node nodal curves'', $d_{\delta}(L)$,
 defined \footnote{Consult the discussion in next subsection.} by G$\ddot{o}$ttsche [Got].
\end{prop}

\medskip

When $L$ is $3\delta-1$-very ample, the sub-scheme $W\subset V\times M^{\delta}_{2, 0}$ in
 the proof of proposition \ref{prop; Got} is a finite
 scheme for a generic choice of $\delta$ dimensional 
linear-sub-system $V\subset |L|$. Under such an assumption,
 ${\cal M}_V\times_{M_{\delta}}Y_{\gamma_{\delta}}$ is a finite scheme as well. From proposition
 \ref{prop; factorial},
 we know that the degree of $[{\cal M}_V\times_{M_{\delta}}Y_{\gamma_{\delta}}]$ is equal to 
$\delta!\cdot d_{\delta}(L)$. 
As we have assumed that $L$ is $5\delta-1$-very ample, G$\ddot{o}$ttsche (in proposition 
\ref{prop; Got})
 has shown that $d_{\delta}(L)$ actually represents the number of $\delta-$node nodal singular curves
 in a generic $\delta$ dimensional linear-subsystem $V$.

 As we have identified the degree of $[{\cal M}_V\times_{M_{\delta}}Y_{\gamma_{\delta}}]$ by
 two different ways, we find that
 the number of $\delta$-nodes nodal curves in a generic $\delta$ dimensional $V\subset |L|$ 
(counted with multiplicities), $d_{\delta}(L)$, 
 is equal to ${1\over \delta!}{\cal AFSW}_{M_{\delta+1}\times \{t_L\}\mapsto 
 M_{\delta}\times \{t_L\}}^{\ast}(1, C-2\sum_{1\leq i\leq \delta}E_i)$. By proposition 
\ref{prop; polynomial} and its ending remark \ref{rem; insert}
 it is a universal degree $\delta$ polynomial in terms of the four variables
$L^2$, $L\cdot {\bf K}_M$, ${\bf K}_M\cdot {\bf K}_M$ and $c_2(M)$. 
So we have finished the proof of our main theorem. $\Box$

\bigskip

\begin{rem}\label{rem; generalm}
 If we replace the singular multiplicities $2$ by $m_1=m_2=\cdots=m_{\delta}=m>2$, and
 replace the $5\delta-1-$very-ampleness condition on $L\mapsto M$ by an
$({(m+1)(m+2)\over 2}-1)\delta-1-$very-ampleness condition,
our main theorem can be generalized to count curves with $\delta$ 
ordinary multiplicities $m$ singularities
\footnote{in a general $({m(m+1)\over 2}-2)\delta$
 dimensional linear sub-system of $|L|$.}.
And the argument is completely parallel to the above argument.
\end{rem}

 At the end of this subsection, we offer a proof of proposition \ref{prop; factorial} cited above.

\noindent Proof of proposition \ref{prop; factorial}:
 We observe that the canonical algebraic obstruction vector bundle 
$\pi_X^{\ast}{\bf W}_{canon}\otimes 
{\bf H}$ of the class $C-{\bf M}(E)E$ 
restricts to an ${\bf S}_{\delta}$ invariant vector bundle over 
 $X\times_{M_{\delta}\times T(M)}
(Y_{\gamma_{\delta}}\times \{t_L\})$. This is easy to check by using the definition of
 ${\bf W}_{canon}$ (see section 5.1, definition 5.3 of [Liu3]) and the fact that different 
 exceptional ${\bf CP}^1$s are completely 
disjoint and are permuted transitively under an induced ${\bf S}_{\delta}$ action.
 Then $\pi_X^{\ast}{\bf W}_{canon}\otimes 
{\bf H}$ descends to a vector bundle on the free 
quotient $(X\times_{M_{\delta}}Y_{\gamma_{\delta}})/{\bf S}_{\delta}
\times_{T(M)} \{t_L\}$, denoted by 
${\bf W}_{descend}$.

 As usual let $M^{[3\delta]}$ denote the Hilbert scheme of $M$ parametrizing the length
 $3\delta$ sub-schemes of $M$.
 Consider the universal sub-scheme $Z_{3\delta}(M)\subset
M\times M^{[3\delta]}$ and the projection maps,

\[
\begin{array}{ccc}
Z_{3\delta}(M) & \stackrel{q_{3\delta}}{\longrightarrow} &  M^{[3\delta]} \\
 \Big\downarrow\vcenter{%
\rlap{$\scriptstyle{\mathrm{p_{3\delta}}}\,$}} &  &  \\
 M &  &  
\end{array}
\]

 The fibration of universal divisors (curves) $D\mapsto |L|$ of the linear system 
forms a divisor in $M\times |L|$, and is called the universal divisor.
 The line bundle ${\cal O}_{|L|\times M}(D)$ for the universal divisor 
$D\subset |L|\times M$ is equivalent
 to $\pi_{|L|\times M\mapsto M}^{\ast}L\otimes \pi_{|L|\times M\mapsto |L|}^{\ast}{\bf H}$, 
twisted by the hyperplane line bundle ${\bf H}\mapsto |L|$.
 Then we may consider 
$q_{3\delta\ast}p_{3\delta}^{\ast}(L\otimes {\bf H})=
 {\bf H}\otimes q_{3\delta\ast}p_{3\delta}^{\ast}L={\bf H}\otimes L_{3\delta}$,
 of rank $3\delta$ over $|L|\times M^{[3\delta]}$.

The smooth quotient space $Y_{\gamma_{\delta}}/{\bf S}_{\delta}$ parametrizes
 the un-ordered $\delta$-tuples of distinct points in $M$ and is embedded naturally
 onto the top open stratum of $M^{[\delta]}$.  On the other hand,
let $x_1, x_2, \cdots, x_{\delta}$ be distinct $\delta$ points on $M$.
 Then $\coprod_{1\leq i\leq \delta}Spec({\cal O}_{M, x_i}/m_{M, x_i}^2)$ is a 
length $3\delta$ sub-scheme of $M$. This enables us to embed the
 top open stratum of $M^{[\delta]}$ into $M^{[3\delta]}$.
Denote this composite inclusion by
 $\Pi_{\delta}:Y_{\gamma_{\delta}}/{\bf S}_{\delta}\mapsto M^{[3\delta]}$. Then 
 we first notice that $\Pi_{\delta}^{\ast}(L_{3\delta}\otimes {\bf H})=
{\bf W}_{descend}$, i.e. the descend of our canonical algebraic family obstruction bundle 
coincides with the obstruction bundle defined by G$\ddot{o}$ttsche 
\footnote{Consult subsection \ref{subsection; transv}.}
 when they are both restricted to the
 top open strata. Again it is because when the $\delta$ 
 blowing up points are distinct in $M$, the corresponding 
exceptional divisors $E_i\subset M_{\delta+1}\times_{M_{\delta}}Y(\gamma_{\delta})$, $1\leq i\leq
 \delta$, are all disjoint. 

 Then we have the following short exact sequence \footnote{Recall $f_{\delta, 1}:M_{\delta+1}\mapsto M_1=M$
 is the composition of $f_{\delta}$, $f_{\delta-1}$, $\dots f_1$, where $f_i:M_{i+1}\mapsto M_i$
 are the projection maps of the universal spaces, introduced in section \ref{section; strata}.} 

$$\hskip -1.4in
0\mapsto {\cal O}_{M_{\delta+1}\times_{M_{\delta}}Y_{\gamma_{\delta}}\times \{t_L\}}
(-2\sum_{1\leq i\leq \delta}E_i)\otimes f_{\delta, 1}^{\ast}{\cal L}
\mapsto {\cal O}_{M_{\delta+1}\times_{M_{\delta}}Y_{\gamma_{\delta}}\times \{t_L\}}\otimes 
f_{\delta, 1}^{\ast}{\cal L}
\mapsto {\cal O}_{2\sum_{1\leq i\leq \delta}E_i|_{M_{\delta+1}\times_{M_{\delta}}Y_{\gamma_{\delta}}}}
\otimes f_{\delta, 1}^{\ast}{\cal L}
\mapsto 0,$$

which is the fundamental building block of the $t_L$-restricted version of 
 the canonical algebraic
 family Kuranishi model of $C-{\bf M}(E)E$.

The push-forward of ${\cal O}_{M_{\delta+1}\times_{M_{\delta}}Y_{\gamma_{\delta}}}
(-2\sum_{1\leq i\leq \delta}E_i)\subset {\cal O}_{M_{\delta+1}\times_{M_{\delta}}Y_{\gamma_{\delta}}}$ 
to $M\times Y_{\gamma_{\delta}}$ defines an ideal sheaf of a universal sub-scheme.
 It is invariant under a free ${\bf S}_{\delta}$ action and
 we denote its free quotient under ${\bf S}_{\delta}$ by the 
new notation $Z_{\gamma_{\delta}}$.

On the other hand $M_{\delta+1}\mapsto M\times M_{\delta}$ projects to the trivial bundle
 $M\times M_{\delta}$ over $M_{\delta}$.
 So the $M_{\delta+1}\times_{M_{\delta}}Y_{\gamma_{\delta}}\longrightarrow M\times 
Y_{\gamma_{\delta}}$-push-forward of the above short exact sheaf 
sequence results in an ${\bf S}_{\delta}$ invariant short exact 
sequence which descends to a short exact sequence on $M\times (Y_{\gamma_{\delta}}/{\bf S}_{\delta})$, 

$$({\ast}) 0\mapsto {\cal I}_{Z_{\gamma_{\delta}}}\otimes \pi_{\gamma}^{\ast}L
\mapsto {\cal O}_{M\times (Y_{\gamma_{\delta}}/{\bf S}_{\delta})} \otimes \pi_{\gamma}^{\ast}L
\mapsto {\cal O}_{Z_{\gamma_{\delta}}}\otimes p_{\gamma_{\delta}}^{\ast}L\mapsto 0.$$

Over here $\pi_{\gamma}:M\times (Y_{\gamma_{\delta}}/{\bf S}_{\delta})\mapsto M$ and
 $p_{\gamma_{\delta}}:Z_{\gamma_{\delta}}\mapsto M$ are the natural projection maps.

On the other hand \footnote{For the definitions of the maps $q_{3\delta}$, $p_{3\delta}$ and 
the Hilbert scheme $M^{[3\delta]}$, 
 please consult the beginning of subsection \ref{subsection; transv}.}
 $q_{3\delta}|_{Z_{3\delta}}:Z_{3\delta}\mapsto M^{[3\delta]}$ and the pre-image of
the subset $\Pi_{\delta}(Y_{\gamma_{\delta}}/{\bf S}_{\delta})\subset M^{[3\delta]}$
 under $(q_{3\delta}|_{Z_{3\delta}})^{-1}$ (inside the universal sub-scheme 
 $Z_{3\delta}(M)$) splits into a disjoint union of the form $\coprod_{i\leq \delta}Z_i$, where
 each $Z_i$ represents a non-reduced sub-scheme of 
relative length $3$ over the base $\Pi_{\delta}(Y_{\Gamma_{\delta}}/{\bf S}_{\delta})$.
 
 And there is a corresponding short exact sequence,

$$\hskip -.6in ({\ast\ast})
0\mapsto \otimes_{1\leq i\leq \delta}{\cal I}_{Z_i}\otimes p_{3\delta}^{\ast}L
\mapsto {\cal O}_{Z_{3\delta}(M)\cap q_{3\delta}^{-1}(\Pi_{\delta}(Y_{\Gamma_{\delta}}/{\bf S}_{\delta}))}
\otimes p_{3\delta}^{\ast}L \mapsto {\cal O}_{\coprod_{1\leq i\leq \delta}Z_i}
\otimes p_{3\delta}^{\ast}L
\mapsto 0.$$

We claim that the push-forward the 
former short exact sequence $({\ast})$ of $Z_{\gamma_{\delta}}$ along 
$M\times (Y_{\gamma_{\delta}}/{\bf S}_{\delta})\mapsto (Y_{\gamma_{\delta}}/{\bf S}_{\delta})$ is
isomorphic to the $\Pi_{\delta}^{\ast}$ pull-back of the $q_{3\delta\ast}$-push-forward of the 
short exact sequence $({\ast\ast})$ 
on $\coprod_{1\leq i\leq \delta}Z_i$, due to the following commutative diagram of maps,

\[
\begin{array}{ccc}
 Z_{\gamma_{\delta}}  & \stackrel{q_{\gamma_{\delta}}}
{\longrightarrow} &  Y_{\gamma_{\delta}}/{\bf S}_{\delta}\\
  \Big\downarrow\vcenter{%
\rlap{$\scriptstyle{\mathrm{\check{\pi}_{\delta}}}\, $}}
 &  & \Big\downarrow\vcenter{%
\rlap{$\scriptstyle{\mathrm{\Pi_{\delta}}}\, $}}
 \\  Z_{3\delta}(M) & \stackrel{q_{3\delta}}{\longrightarrow} & M^{[3\delta]}  
\end{array}
\]

, with $\check{\pi}_{\delta}:Z_{\gamma_{\delta}}\mapsto Z_{3\delta}(M)$ being the canonical 
inclusion.

 Therefore we can identify the descend bundle ${\bf W}_{descend}\cong q_{\gamma_{\delta}\ast}
p_{\gamma_{\delta}}^{\ast}L\otimes {\bf H}$,
 with $\Pi_{\delta}^{\ast}q_{3\delta\ast}p_{3\delta}^{\ast}
 (L)\otimes {\bf H}$. Moreover, because of the following commutative diagram on the bundle maps,

\[
\begin{array}{ccc}
  {\bf H}^{\ast} & \longrightarrow & \pi_{|L|\times Y_{\gamma_{\delta}}/{\bf S}_{\delta}\mapsto 
Y_{\gamma_{\delta}}/{\bf S}_{\delta}}^{\ast}q_{\gamma_{\delta}\ast}
p_{\gamma_{\delta}}^{\ast}L\\
 \Big\downarrow & & \Big\downarrow\\ 
 {\bf H}^{\ast} & \longrightarrow & \pi_{|L|\times M^{[3\delta]}\mapsto M^{[3\delta]}}^{\ast}
q_{3\delta\ast}p_{3\delta}^{\ast}(L)
\end{array}
\]

 the descend of the canonical section $s_{canon}|_{|L|\times Y_{\gamma_{\delta}}}$ corresponds to
 the $\Pi_{\delta}^{\ast}$-pull-back of a section of $q_{3\delta\ast}p_{3\delta}^{\ast}
 (L)\otimes {\bf H}$ and a ray ${\bf l}$ in the projective space 
$|L|$ represents an algebraic curve in $M$ singular along the sub-scheme
 $\coprod_{i\leq \delta}Spec({\cal O}_{x_i, M}/m_{x_i}^2)$ iff the
 values of the canonical section $s_{canon}$ at each of the points ${\bf l}\times \sigma(x_1\times
 x_2\times \cdots \times x_{\delta})\in {\bf P}({\bf V}_{canon})$, 
$\sigma\in {\bf S}_{\delta}$, vanishes.

 By our assumption in the proposition, ${\cal M}_V\times_{M_{\delta}}Y_{\gamma_{\delta}}$
 has been assumed to be a finite sub-scheme of $V\times M_{\delta}\subset X\times_{T(M)}\{t_L\}$.
 Because ${\cal M}_V\times_{M_{\delta}}
Y_{\gamma_{\delta}}$ can be identified with the zero locus of $s_{canon}$ in $V\times Y_{\gamma_{\delta}}$, 
 then according to Section 14.1 of [F], one may define a localized top Chern class of 
$\pi_X^{\ast}{\bf W}_{canon}\otimes {\bf H}|_{V\times Y_{\gamma_{\delta}}}$ with respect to
 $s_{canon}|_{V\times Y_{\gamma_{\delta}}}$ inside 
${\cal A}_0({\cal M}_V\times_{M_{\delta}}Y_{\gamma_{\delta}})$. Because we have identified the
descend bundle ${\bf W}_{descend}$ of 
$\pi_X^{\ast}{\bf W}_{canon}\otimes {\bf H}|_{X\times_{M_{\delta}}Y_{\gamma_{\delta}}}$ 
with $\Pi_{\delta}^{\ast}q_{3\delta\ast}p_{3\delta}^{\ast}(L)\otimes {\bf H}$, the image of
 ${\cal M}_V\times_{M_{\delta}}Y_{\gamma_{\delta}}$ in 
$(X\times_{M_{\delta}}Y_{\gamma_{\delta}})/{\bf S}_{\delta}$ can be identified with a finite sub-scheme
$\subset |L|\times M^{\delta}_{2, 0}$, denoted as $W$ in the proof of 
proposition \ref{prop; Got}.

Since the localized top Chern class is defined by the local datum, i.e.
 the total Segre class of the normal cone of the zero locus
and the restriction of the total Chern class of the vector bundle to the zero locus, the 
localized top Chern class of ${\bf W}_{descend}$ along 
${\cal M}_V\times_{M_{\delta}}Y_{\gamma_{\delta}}/{\bf S}_{\delta}$ is equal to the
 localized top Chern class of $q_{3\delta\ast}p_{3\delta}^{\ast}(L)\otimes {\bf H}$ over $W$.

Since the quotient map
$X\times_{M_{\delta}}Y_{\gamma_{\delta}}\mapsto X\times_{M_{\delta}}Y_{\gamma_{\delta}}/
{\bf S}_{\delta}$ is an un-ramified covering map of degree $\delta!$, by proposition 14.1. (d).(iii). of
 [F], the degree of the 
 localized top Chern class of $\pi_X^{\ast}{\bf W}_{canon}\otimes {\bf H}$ along
 the zero locus ${\cal M}_V\times_{M_{\delta}}Y_{\gamma_{\delta}}$ is $\delta!$ times the
 localized top Chern class of ${\bf W}_{descend}$ along the quotient of zero locus 
 ${\cal M}_V\times_{M_{\delta}}Y_{\gamma_{\delta}}/{\bf S}_{\delta}$. 

 Because $W\mapsto pt$ factors through $W\subset V\times M^{\delta}_{2, 0}\mapsto pt$, the degrees of the
 localized Chern class of ${\bf W}_{descend}|_{V\times Y_{\gamma_{\delta}}}$ along 
${\cal M}_V\times_{M_{\delta}}Y_{\gamma_{\delta}}/{\bf S}_{\delta}$ and the localized 
Chern class of $q_{3\delta\ast}p_{3\delta}^{\ast}(L)\otimes {\bf H}$ along $W$ are equal and their
 common value is equal to 

$$\hskip -.6in
\int_{|V|\times M^{\delta}_{2, 0}}c_{3\delta}({\bf H}\otimes L_{3\delta})
=\int_{|V|\times M^{\delta}_{2, 0}}c_1({\bf H})^{\delta}\cap 
 c_{2\delta}(L_{3\delta})=\int_{M^{\delta}_{2, 0}}c_{2\delta}(L_{3\delta})=d_{\delta}(L).$$

Thus the degree of $[{\cal M}_V\times_{M_{\delta}}Y_{\gamma_{\delta}}]$ is $\delta!d_{\delta}(L)$.
The proof of proposition \ref{prop; factorial} is finished.
$\Box$

\medskip

\subsection{The Finiteness Result of 
${\cal M}_{C-{\bf M}(E)E}\times_{M_n}Y_{\gamma_{\delta}}$}
\label{subsection; transv}

\bigskip

 In this subsection, we survey the finiteness result based on 
G$\ddot{o}$ttsche's argument in [Got]. 

 Recall the following definition of $k$-very ampleness of a line bundle on $M$, 
due to [BS].

\begin{defin}\label{defin; kvery}
A line bundle $L$ on an algebraic manifold is $k$-very ample if
 for all length $k+1$ sub-scheme $Z\subset M$, the following 
 restriction map $H^0(M, L)\mapsto H^0(M, {\cal O}_Z\otimes L)$ is surjective.
\end{defin}

 The $1-$very ampleness is equivalent to the usual very ample condition, by
 page 120, prop. 7.3. on page 152 and remark 7.8.2. on page 158 of [Ha].

\medskip

\begin{defin}\label{defin; s}
Let $M^{\delta}_2\subset M^{[3\delta]}$ be the closure (with the
 reduced induced structure) of the locally closed subset 
 $M^{\delta}_{2, 0}$ which parametrizes sub-schemes of the form
 $\coprod_{i=1}^{\delta}Spec({\cal O}_{M, x_i}/m_{M,x_i}^2)$, 
where $x_1, x_2, \cdots , x_{\delta}$ are distinct points in $M$. 
\end{defin}

 The symbol $M^{[n]}$ denote the Hilbert scheme of finite sub-schemes of
 length $n$ on $M$, and let $Z_n\subset M\times M^{[n]}$ denote the
 universal family of sub-schemes
 with projection $p_n:Z_n(M)\mapsto M$ and $q_n:Z_n(M)\mapsto 
 M^{[n]}$. Then $L_n=(q_n)_{\ast}(p_n)^{\ast}L$ is locally free of rank $n$
 on $M^{[n]}$.
 It is easy to see that $M^{\delta}_2$ is birational to $M^{[\delta]}$ and
 we set $d_{\delta}(L)=\int_{M^{\delta}_2}c_{2\delta}(L_{3\delta})$.

 Recall the following proposition due to G$\ddot{o}$ttsche. It is a word by word
 duplication of proposition 5.3. of [Got]. We include it here for the convenience
 of the readers.

\begin{prop}\label{prop; Got}({\bf G$\ddot{o}$ttsche})
Assume $L$ is $3\delta-1$ very ample, then a general $\delta$ dimensional
 linear subsystem $V\subset |L|$ contains only finitely many curves
 $C_1, C_2, C_3, \cdots, C_s$ with $\geq \delta$ singularities.
There exists positive integers $n_1, n_2, \cdots, n_s$ such that 
 $\sum_i n_i=d_{\delta}(L)$. If furthermore $L$ is $(5\delta-1)-$very
 ample ( $5$-very ample if $\delta=1$), then the $C_i$ have precisely 
$\delta$ nodes as singularities.
\end{prop}

For completeness, we include its original proof here.

\noindent Proof (due to {\bf G$\ddot{o}$ttsche}): 
Assume that $L$ is $(3\delta-1)$-very
 ample. We apply the Thom-Porteous formula to the restrictions of the
 evaluation map $H^0(M, L)\otimes {\cal O}_{M^{[3\delta]}}\mapsto
 L_{3\delta}$ to $M^{\delta}_2$ and to $M^{\delta}_2-M^{\delta}_{2, 0}$.
 As $L$ is $(3\delta-1)$-very ample, the evaluation map is surjective. 
 Then ([F] ex. 14.3.2) applied to $M^{\delta}_2$ gives that for a 
general $\delta-$dimensional sub-linear system $V\subset |L|$ the class
 $d_n(L)$ is represented by the class of the finite scheme $W$ of $Z\in 
 M^{\delta}_2$ with $Z\subset D$ for $D\in V$. The scheme structure of 
$W$ might be non-reduced. The application of ([F] ex. 14.3.2) to
 $M^{\delta}_2\backslash M^{\delta}_{2, 0}$ and a dimension count give that $W$ 
lies entirely in $M^{\delta}_{2, 0}$. 

 Now assume that $L$ is $(5\delta-1)-$very ample. Let $V\subset |L|$ 
again be general $\delta-$dimensional subsystem of $|L|$. The Porteous
 formula applied to the restriction of $L_{3\delta+3}$ to
 $M^{\delta+1}_2$ and a dimension count shows that there will be no curves
in $V$ with more than $\delta$ singularities.

 Let $M^{\delta}_{3, 0}\subset M^{[5\delta]}$ be the locus of schemes 
 of the form $Z_1 \cup Z_2 \cup Z_3 \cdots \cup Z_s$, where each $Z_i$ is
 of the form $Spec({\cal O}_{M, x_i}/(m^3+xy))$ with $x, y$ local
 parameters at $x_i$ and let $M^{\delta}_3$ be the closure. 
  If a curve $C$ with precisely $\delta$ singularities does not contain
 a sub-scheme corresponding to a point in $M^{\delta}_3\\M^{\delta}_{3, 0}$,
 then it has $\delta$ nodes as only singularities. It is easy to see that
 $M^{\delta}_{3, 0}$ is smooth of dimension $4\delta$. Applying the Porteous
formula to the restriction of $L_{5\delta}$ to 
 $M^{\delta}_3\\M^{\delta}_{3, 0}$ and a dimension count we see that all the
 curves in $V$ with $\delta$ singularities have precisely $\delta$ nodes.
$\Box$

 In the following, we generalize G$\ddot{o}$ttsche's argument to our context.
 Let ${\bf M}(E)E=\sum_{i\leq \delta}2E_i$. Namely, $m_i=2$ for all
 integers $i$, $1\leq i\leq \delta$.
 A line bundle $L\mapsto M$ with $c_1(L)=C$ determines a unique point
 in the connected component of Picard group, 
$T(M)$, denoted as $t_L$. The fiber product 
${\cal M}_{C-{\bf M}(E)E}\times_{T(M)}\{t_L\}$ is the algebraic family 
moduli sub-space of curves in the fibers of the family 
$M_{\delta+1}\mapsto M_{\delta}$ 
projecting onto curves in $|L|$.
 Then there exists a natural map 
${\cal M}_{C-{\bf M}(E)E}\times_{T(M)}\{t_L\}\mapsto |L|$.

\begin{prop}\label{prop; tran}    
Let $L$ be a $3\delta-1$-very ample line bundle over $M$. Let
 ${\bf M}(E)E=\sum_{1\leq i\leq \delta}2E_i$, and let ${\cal M}_V$ 
denote the pre-image of
 a general $\delta$ dimensional linear subsystem $V\subset |L|$ under
 ${\cal M}_{C-{\bf M}(E)E}\times_{T(M)}\{t_L\}\mapsto |L|$. Then the fiber product
 ${\cal M}_V\times_{M_{\delta}}(M_{\delta}-\cup_{\Gamma\in 
\Delta(\delta)-\{\gamma_{\delta}\}}S_{\Gamma})$
 of ${\cal M}_V\mapsto M_{\delta}$ and $M_{\delta}-\cup_{\Gamma\in 
\Delta(\delta)-\{\gamma_{\delta}\}}S_{\Gamma}\subset
 M_{\delta}$ is a finite scheme and its image under the projection map to 
 $M_{\delta}$ lies in the generic stratum $Y_{\gamma_{\delta}}$.
\end{prop}

\noindent Proof: It is not hard to see that the image 
${\cal M}_V\times_{M_n}Y_{\gamma_{\delta}}\mapsto V$ corresponds to all the curves in the 
linear sub system $V$ which has at least $\delta$ distinct singularities.
The space $Y_{\gamma_{\delta}}$ parametrizes all the ordered distinct $\delta$ points
 on $M$, and the image of ${\cal M}_V\times_{M_{\delta}}Y_{\gamma_{\delta}}$ in $Y_{\gamma_{\delta}}$
 are the ordered $\delta$-tuples of singular points of the curves. According to
 proposition \ref{prop; Got} (by G$\ddot{o}$ttsche), there are at most a finite number
 of singular curves in $V$ with exactly $\delta$ distinct singularities.

 This implies ${\cal M}_V\times_{M_{\delta}}Y_{\gamma_{\delta}}$ to be a finite
 scheme. On the other hand, the image of ${\cal M}_V\mapsto M_{\delta}$ 
may have non-trivial intersections
 with the various subsets $S_{\Gamma}, \Gamma\in 
\Delta(\delta)-\{\gamma_{\delta}\}$.

 To prove the proposition, it suffices to show that the image
 ${\cal M}_V\mapsto M_{\delta}$ intersect with $Y_{\Gamma}$ trivially for
 all the chain-like \footnote{Consult
 definition \ref{defin; special} and the comment afterward.}
$\Gamma\in adm_2(\delta)$.
  Then ${\cal M}_V\times_{M_{\delta}}(M_{\delta}-\cup_{\Gamma\in \Delta(\delta)
-\{\gamma_{\delta}\}}S_{\Gamma})$ can be identified with the space 
${\cal M}_V\times_{M_{\delta}}Y_{\gamma_{\delta}}$ and
 it has been shown to be
 a finite scheme which projects into $Y_{\gamma_{\delta}}$ automatically.

 To show that the image of ${\cal M}_V\mapsto M_{\delta}$ avoids all the
 $Y_{\Gamma}$ for chain-like $\Gamma$, 
we fix an arbitrary $\bar{b}\in Y_{\Gamma}$ and show
 that $\bar{b}$ is not included in the image of ${\cal M}_V\mapsto M_{\delta}$
 for generic choices of $V$.

\begin{lemm}\label{lemm; saturate}
Let $\Gamma\in adm_2(\delta)$ be a chain-like admissible graph.
The fiber above the point $\bar{b}\in Y(\Gamma)\subset M_{\delta}$ of 
$M_{\delta+1}\mapsto M_{\delta}$ determines a $\delta-$consecutive 
blowing ups of $M$, denoted by $\tilde{M}$. As usual, let $E_1, E_2, \cdots,
 E_{\delta}$ denote the $\delta$ exceptional divisors in $\tilde{M}$ of the
 blowing down map $f:\tilde{M}\mapsto M$. Then 
$f_{\ast}{\cal O}_{\tilde{M}}(-2\sum_{i\leq \delta}E_i)$ is an ideal sheaf
 $\subset {\cal O}_M$ of a finite sub-scheme of $M$ of length $3\delta$.
\end{lemm}
 
\noindent Proof of Lemma \ref{lemm; saturate}: Firstly we prove that
 $f_{\ast}{\cal O}_{\tilde{M}}(-2\sum_{i\leq \delta}E_i)$ is an ideal sheaf by
 showing that it is a sub-sheaf of ${\cal O}_M$.

 To see this, we notice that ${\cal O}_{\tilde{M}}(-2\sum_{i\leq \delta}E_i)\subset
 {\cal O}_{\tilde{M}}$ and thus 
 $f_{\ast}{\cal O}_{\tilde{M}}(-2\sum_{i\leq \delta}E_i)\subset 
f_{\ast}{\cal O}_{\tilde{M}}$. On the other hand, the exceptional divisors of
 $f:\tilde{M}\mapsto M$ are all rational, this implies that 
$f_{\ast}{\cal O}_{\tilde{M}}={\cal O}_M$. Thus, 
$f_{\ast}{\cal O}_{\tilde{M}}(-2\sum_{i\leq \delta}E_i)$ is a sub-sheaf 
of ${\cal O}_M$. Let $Z$ be
\footnote{The $Z$ has nothing to do with the various $Z$ used in the
 previous sections.} the sub-scheme of $M$ defined by the ideal sheaf
 ${\cal I}_Z=f_{\ast}{\cal O}_{\tilde{M}}(-2\sum_{i\leq \delta}E_i)$.

 Secondly, we prove that the length of $Z$ is bounded by $3\delta$ from above.
 This is achieved by an induction argument on $\delta$.  For $\delta=1$, there is a
 unique exceptional divisor $E_1$. Let $x\in M$ be the blowing up point.
 It is easy to see that ${\cal O}_Z\cong{\cal O}_M/{\cal I}_Z\cong 
{\cal O}_{x}/m_x^2$ and $Z$ is of length $3=3\cdot 1=3\cdot\delta$.

 For $\delta>1$, suppose that for all the smooth algebraic surfaces 
$M$ and all $\delta$-consecutive blowing ups, $\tilde{M}$, of $M$, the ideal sheaf 
$f_{\ast}{\cal O}_{\tilde{M}}(-2\sum_{i\leq \delta}
E_i)$ is known to define a length $\leq 3\delta$ sub-scheme of $M$, we would
like to show that for $\delta+1$, and the $\delta+1$-consecutive blowing ups
$\check{M}$ of $M$, the ideal sheaf $\check{f}_{\ast}{\cal O}_{\check{M}}(-2
\sum_{i\leq \delta+1}E_i)$ defines a sub-scheme of $M$ of length $\leq 3(\delta+1)$.

  We notice that $\check{f}:\check{M}\mapsto M$ can be factored into
 $\check{\bar{f}}:\check{M}\mapsto \bar{M}$ and $\bar{f}:\bar{M}\mapsto M$, where
$\bar{M}$ is a one-point blowing up of $M$, with the exceptional divisor $E_1$, 
and $\check{M}$ can be constructed
 from $\bar{M}$ by $\delta$-consecutive blowing ups.

 By induction hypothesis, $\check{\bar{f}}_{\ast}{\cal O}_{\check{M}}(-2
\sum_{2\leq i\leq \delta+1}E_i)$ defines an ideal sheaf on $\bar{M}$ of
 the sub-scheme, $\bar{Z}\subset \bar{M}$ of length $\leq 3\delta$.

 Then 
$$\check{\bar{f}}_{\ast}{\cal O}_{\check{M}}(-2
\sum_{1\leq i\leq \delta+1}E_i)=\check{\bar{f}}_{\ast}\bigl({\cal O}_{\check{M}}(-2
\sum_{2\leq i\leq \delta+1}E_i)\otimes 
\check{\bar{f}}^{\ast}{\cal O}_{\check{M}}(-2E_1)\bigr)={\cal I}_{\bar{Z}}\otimes
 {\cal O}_{\bar{M}}(-2E_1),$$

 and 

$$\check{f}_{\ast}{\cal O}_{\check{M}}(-2
\sum_{1\leq i\leq \delta+1}E_i)=\bar{f}_{\ast}\bigl(
\check{\bar{f}}_{\ast}{\cal O}_{\check{M}}(-2
\sum_{1\leq i\leq \delta+1}E_i)\bigr)=\bar{f}_{\ast}
\bigl({\cal I}_{\bar{Z}}\otimes 
 {\cal O}_{\bar{M}}(-2E_1)\bigr).$$

 By using the exactness of the sequence,
$0\mapsto \bar{f}_{\ast}{\cal I}_{\bar{Z}}
\mapsto \bar{f}_{\ast}{\cal O}_{\bar{M}}
\mapsto \bar{f}_{\ast}{\cal O}_{\bar{Z}}$,

 and the fact $\bar{f}_{\ast}{\cal O}_{\bar{M}}(-2E_1)$ being an ideal sheaf of $M$
 of co-length $3$, the length of the sub-scheme defined by
 $\bar{f}_{\ast}\bigl({\cal I}_{\bar{Z}}\otimes 
 {\cal O}_{\bar{M}}(-2E_1)\bigr)$ is bounded by $length(\bar{Z})+3=
 3\delta+3=3(\delta+1)$ from above.

 Thirdly, we show that the equality is saturated, namely
 $length(Z)=3\delta$ for all $\delta\in {\bf N}$, 
when $\Gamma\in \Delta(\delta)$ is a chain-like
 admissible graph $\in adm_2(\delta)$.
 We prove this by contradiction. We know that for $\delta=1$ the
 equality always saturates. Suppose that there is some chain-like admissible
graph $\Gamma$ and for some $\bar{b}\in Y_{\Gamma}$, the fiber $\tilde{M}$ of
 $M_{\delta+1}\mapsto M_{\delta}$ above $\bar{b}$, $f:\tilde{M}\mapsto M$,
 such that the ideal sheaf 
$f_{\ast}{\cal O}_{\tilde{M}}(-2\sum_{i\leq \delta}E_i)={\cal I}_Z$
 is of length $<3\delta$. We may assume additionally that the smallest 
natural number 
$\delta_0>1$ satisfying the above condition has been chosen. I.e. for
 all natural numbers $\delta$ smaller than $\delta_0$ and any chain-like
 admissible graphs $\in adm_2(\delta)$, the ideal sheaf always defines length
 $3\delta$ sub-schemes in $M$. 

 Apparently the question is of local nature, so we may assume that
 we are blowing up consecutively at the origin ${\bf 0}\in 
{\bf C}^2=M$. Suppose that the connected graph $\Gamma\in adm_2(\delta_0)$ is
 a linear chain and the $i-$th vertex is the unique 
direct descendent of the $i-1-$th vertex, for 
all $2\leq i\leq \delta_0$. Let $x, y$ be the affine coordinates
 around the origin ${\bf 0}$.

 Consider the blowing down of $\tilde{M}$ along the last exceptional divisor 
$E_{\delta_0}$, $\acute{f}: 
\tilde{M}\mapsto \acute{M}$. Then $\acute{M}$ is a
 $\delta_0-1$ consecutive blowing ups of $M$ at ${\bf 0}$. 
Define ${\cal I}_{Z_0}=
{\cal O}_{\acute{M}}$, 
 ${\cal I}_{Z_1}=\acute{f}_{\ast}{\cal O}_{\tilde{M}}(-E_{\delta_0})$ and
 ${\cal I}_{Z_3}=\acute{f}_{\ast}{\cal O}_{\tilde{M}}(-2E_{\delta_0})$. Let
 ${\cal I}_{Z_2}$ be an ideal sheaf of co-length $2$ 
in-between ${\cal I}_{Z_1}$ and
 ${\cal I}_{Z_3}$, i.e. ${\cal I}_{Z_1}\supset {\cal I}_{Z_2}\supset 
{\cal I}_{Z_3}$. Then for $\breve{f}:\acute{M}\mapsto M$, we have 

$$f_{\ast}{\cal O}_{\tilde{M}}(-2\sum_{i\leq \delta_0-1}E_i-E_{\delta_0})
\cong \breve{f}_{\ast}
\bigl({\cal O}_{\acute{M}}(-2\sum_{i\leq \delta_0-1}E_i)\otimes 
{\cal I}_{Z_1}\bigr)\supset
\breve{f}_{\ast}
\bigl({\cal O}_{\acute{M}}(-2\sum_{i\leq \delta_0-1}E_i)\otimes 
{\cal I}_{Z_2}\bigr)$$

$$\supset 
\breve{f}_{\ast}\bigl({\cal O}_{\acute{M}}(-2\sum_{i\leq \delta_0-1}E_i)
\otimes {\cal I}_{Z_3}\bigr)\cong 
\breve{f}_{\ast}{\cal O}_{\tilde{M}}(-2\sum_{i\leq \delta_0}E_i).$$

 The minimality of $\delta_0$ implies that 
$f_{\ast}{\cal O}_{\tilde{M}}(-2\sum_{i\leq \delta_0-1}E_i)$ is of co-length
 $3(\delta_0-1)$ in ${\cal O}_M$.
Because 
$f_{\ast}{\cal O}_{\tilde{M}}(-2\sum_{i\leq \delta_0-1}E_i)\supset 
f_{\ast}{\cal O}_{\tilde{M}}(-2\sum_{i\leq \delta_0}E_i)$ is of co-length 
$<3$, there must be some $a\in \{0, 1, 2\}$ such that
$\breve{f}_{\ast}\bigl({\cal O}_{\acute{M}}(-2\sum_{i\leq \delta_0-1}E_i)\otimes
 {\cal I}_{Z_a}\bigr)= 
\breve{f}_{\ast}
\bigl({\cal O}_{\acute{M}}(-2\sum_{i\leq \delta_0-1}E_i)\otimes
 {\cal I}_{Z_{a+1}}\bigr)$.

 Define 
$\breve{f}_{\ast}\bigl({\cal O}_{\acute{M}}(-2\sum_{i\leq \delta_0-1}E_i)\otimes
 {\cal I}_{Z_a}\bigr)={\cal I}_Z$ 
for such an $a$. Consider a polynomial $g(x, y)\in {\bf C}[x, y]$
 vanishing along the sub-scheme 
$Z$. 

 The sheaf identification 
$\breve{f}_{\ast}\bigl({\cal O}_{\acute{M}}(-2\sum_{i\leq \delta_0-1}E_i)\otimes
 {\cal I}_{Z_a}\bigr)={\cal I}_Z$ induces an identification
 $\psi_a:\Gamma(M, {\cal I}_Z)\stackrel{\cong}{\longrightarrow}\Gamma(\acute{M}, 
\breve{f}_{\ast}\bigl({\cal O}_{\acute{M}}(-2\sum_{i\leq \delta_0-1}E_i)\otimes
 {\cal I}_{Z_a})$.
 
 Then the image $\psi_a(g)$ of $g$ in $\Gamma(\acute{M}, 
\breve{f}_{\ast}\bigl({\cal O}_{\acute{M}}(-2\sum_{i\leq \delta_0-1}E_i)\otimes
 {\cal I}_{Z_a})$ defines a zero locus in $\acute{M}$. Suppose that $g$ 
 can be chosen such that $\psi_a(g)$ does not vanish identically on
 $E_1, E_2, \cdots, E_{\delta_0-1}$, then the zero locus 
$Z(\psi_a(g))=\{x|\psi_a(g)(x)=0, x\in \acute{M}\}$ in
 $\acute{M}$ is nothing but the strict transform of $Z(g)=\{x|g(x)=0, 
 x\in M\}$ under the $\delta_0-1-$consecutive blowing ups.

 By the choice of $a$
 we know that $\breve{f}_{\ast}
{\cal O}_{\acute{M}}(-2\sum_{i\leq \delta_0-1}E_i)\otimes
 {\cal I}_{Z_{a+1}}={\cal I}_Z$ as well. So the defining equation $\psi_a(g)$
 vanishes 
along the length $a+1$ sub-scheme $Z_{a+1}\supset Z_a$ automatically. 

 We demonstrate the existence of 
 some counter-example violating the above assertion in the following lemma
 \ref{lemm; existence}. After it is achieved,
 then the co-length of $f_{\ast}{\cal O}_{\tilde{M}}(-2\sum_{1\leq i\leq 
\delta_0}E_i)$ has to be exactly $3\delta_0$ and therefore the minimal
 $\delta_0$ violating the saturation condition can never exist. 
 Then the proof of lemma \ref{lemm; saturate} is finished. $\Box$

\medskip

 The following lemma supports the counter-example needed in the proof of
 lemma \ref{lemm; saturate}.

\medskip

\begin{lemm}\label{lemm; existence}
 Let $M={\bf C}^2$ and let $\Gamma$ be a connected 
chain-like admissible graph in 
 $adm_2(\delta_0)$. As above fix a point $\bar{b}\in Y_{\Gamma}$ and therefore
 a $\delta_0$-consecutive blowing up of $M$ at its origin ${\bf 0}$. 
Let $a$ be chosen as above and let 
 $Z_a, Z_{a+1}$ be the length $a$ and $a+1$ sub-schemes of $\acute{M}$ defined
 above. Given any nonzero $g\in H^0(M, {\cal I}_Z)$, let 
$\acute{g}\in H^0(\acute{M}, {\cal O}_{\acute{M}})$ be the
 defining equation of the strict transform of the locus $Z(g)$
 in $\acute{M}$ (well-defined up to a ${\bf C}^{\ast}$ multiplication). Then
 there exists a nonzero 
${\bf g}\in H^0({\bf C}^2, {\cal I}_Z)$ such that $\acute{{\bf g}}\in
 H^0(\acute{M}, {\cal O}_{\acute{M}}(-2\sum_{1\leq i\leq \delta_0-1}E_i)\otimes 
{\cal I}_{Z_a})$ but $\acute{{\bf g}}\not\in
 H^0(\acute{M}, {\cal O}_{\acute{M}}(-2\sum_{1\leq i\leq \delta_0-1}E_i)\otimes 
{\cal I}_{Z_{a+1}})$.
\end{lemm}

\noindent Proof of lemma \ref{lemm; existence}:
By the embedded resolutions of singular curves in algebraic surfaces,
 (for example consult 8B, page 160-166 of [Mum]), the rational double
 point in a singular curve in ${\bf C}^2$ 
defined by the equation $x^2=y^{2(\delta_0-1)}$ can be
 resolved into smooth points by consecutively blowing up $\delta_0-1$ times,
 each upon the unique singular point of the intermediate strict transforms. 

 Algebraically blowing up a point
 corresponds to replacing the coordinates $(x, y)$ by 
$(x', y')(={x\over y}, y)$ (or
 $(x', y')=(x, {y\over x})$)
 in the defining equations.  And under the first set of
 change of variables the equation becomes

 $$x^2-y^{2(\delta_0-1)}=({x\over y}\cdot y)^2-y^{2(\delta_0-1)}=
y'^2(({x\over y})^2-
y'^{2(\delta_0-2)})=y'^2(x'^2-y'^{2(\delta_0-2)}).$$
 
 Firstly, the strict transformation of the zero locus 
in the one-point blowing up of $M$ at the origin, 
defined by $x'^2=y'^{2(\delta_0-2)}$, has a
 rational double point ($A_{2\delta_0-5}$ singularity) at $x'=y'=0$ and
 it intersects with the exceptional divisor (defined locally by $y'=0$ here) with
 a singular 
multiplicity $\mu=2$. By induction one realizes that the original singular 
curve gets resolved
 into a smooth curve after $\delta_0-1$ consecutively blowing ups
 and the resolved smooth curve 
intersects with the last exceptional ${\bf CP}^1$ (dual to $E_{\delta_0-1}$)
 at two distinct points.
 This can be seen by observing that the $\delta_0=2$ case corresponds to
 nothing but the ordinary 
double (nodal singular) point. The blowup sequence determined by $\bar{b}\in 
 Y_{\Gamma}$ 
determines a sequence of $\delta_0-1$
 points in a sequence of the first $\delta_0-1$
exceptional ${\bf P}^1$, each representing an exceptional divisor in the 
intermediate blowing
 ups. Apparently the blowing-up centers in resolving
 $x^2=y^{2(\delta_0-1)}$ to a smooth curve
 may not be identical to the first $\delta_0-1$
 blowing-up centers determined by
 $\bar{b}\in Y_{\Gamma}$. On the other hand, the change of variables on page
 161, in the subsection (8.6) of [Mum] allows us to move the locations of
 the intermediate singularities that are blown up.  We proceed as the following.

 Firstly notice that it requires at least $2\delta_0-1$ affine coordinate
 systems to cover the $\delta_0-1$ distinct exceptional ${\bf P}^1$ of 
 $\acute{M}\mapsto M$ and the
 punctured neighborhood of $M$ around ${\bf 0}$. Let $({\bf x}, {\bf y})$ 
 be the affine coordinate of $M$ at ${\bf 0}$ and let 
 $(x_{2i-1}, y_{2i-1})$, $(x_{2i}, y_{2i})$, $1\leq i\leq \delta_0-1$ be the
 dual affine coordinates on the neighborhoods of the 
$i-$th exceptional ${\bf P}^1$. For a fixed $i$,
 they satisfy the following transition rules $x_{2i-1}y_{2i-1}=x_{2i}$, 
 ${1\over x_{2i-1}}=y_{2i}$. The locus $y_{2i-1}$ or $x_{2i}=0$ corresponds to the
 $i-$th ${\bf P}^1$.
To determine the transitions of coordinates among
 different $i$, it suffices to work out the transition for the adjacent pairs.

The blowing up sequence determined by $\bar{b}$ determines $\delta_0-1$
points ${\bf v}_1, {\bf v}_2, \cdots, {\bf v}_{\delta_0-1}$ 
in the $\delta_0-1$ exceptional ${\bf P}^1$ of $\acute{M}\mapsto M$, respectively.

Firstly focus on the exceptional ${\bf P}^1$ representing the $\delta_0-1$-th 
divisor $E_{\delta_0-1}$.
 Because $Z_1={\bf v}_{\delta_0-1}$ 
is a point in this ${\bf P}^1$, either it is covered by the
 coordinate system $(x_{2\delta_0-3}, y_{2\delta_0-3})$ with $y_{2\delta_0-3}=0$,
 or it is at the origin of $(x_{2\delta_0-2}, y_{2\delta_0-2})$ coordinate 
 system.

 By choosing either $(u, v)=(x_{2\delta_0-3}, y_{2\delta_0-3})$ 
or $(u, v)=(y_{2\delta_0-2}, x_{2\delta_0-2})$, we assume that
 $(u, v)$ is a coordinate system around the $\delta_0-1-$th exceptional
 ${\bf P}^1$ containing th point ${\bf v}_{\delta_0-1}$ 
such that $v=0$ defines the exceptional curve locally and 
$u$ is a local uniformizer of the exceptional curve.
 We choose the constants $\alpha, \beta\in {\bf C}$, and
 $A\in {\bf C}$ 
in the quadratic polynomial 
$g_{\delta_0-1}(u, v)=(u-\alpha)(u-\beta)+Av$
according to the value of $a\in \{0, 1, 2\}$.

 Notice that when $\alpha\not=\beta$, 
the equation $v^2g_{\delta_0-1}({u\over v}, v)=u^2-(\alpha+\beta)uv+
\alpha\beta v^2+Av^3=0$
 represents a curve with a rational double point at the origin of $(u, v)$.
 We know that $u=\alpha$, $u=\beta$ are the affine coordinates of
 the two intersection points of the
 resolved smooth curve with $E_{\delta_0-1}$, locally defined by $v=0$.
 
\noindent Case $0$: Suppose that $a=0$, then set $A=0$ and choose generic 
$\alpha$ and $\beta$ to  move the two intersection points of the resolved smooth 
curve with the
 exceptional curve $E_{\delta_0-1}$ away from the blown up
point $Z_1={\bf v}_{\delta_0-1}$ of$\acute{f}:\tilde{M}\mapsto \acute{M}$.

\noindent Case $1$: Suppose that $a=1$. Firstly choose $\beta=\beta_0$
 such that $u=\beta_0$ is the affine coordinate of blown up point $Z_1$ 
in $\acute{M}$. 
 Thus the smooth curve resolved from the nodal curve locally defined by 
$v^2g_{\delta_0-1}({u\over v}, v)=
u^2-(\alpha+\beta)uv+\alpha\beta v^2+Av^3=0$ 
vanishes along the length one sub-scheme $Z_1$.
 Notice that in terms of the local uniformizers $(u-\beta_0), v$ at
 $(\beta_0, 0)$ the
 first jets of $g_{\delta_0-1}(u, v)$ are determined by $\alpha-\beta_0$ and $A$. 
 We choose a generic $\alpha$, $\alpha\not=\beta_0$ and $A$ 
 such that the resolved smooth curve
 does not vanish along the sub-scheme $Z_2$. This is
 possible because the length $2$ sub-scheme $Z_2\subset \acute{M}$ determines
 a tangent direction of $\acute{M}$ at $Z_1$ and the generic choices
 of $\alpha$ and $A$ can prevent the locus $g_{\delta_0-1}(u, v)=0$
 from being tangent to this tangent direction specified by $Z_2$
 at $Z_1$.

\medskip

 \noindent Case $2$: Suppose that $a=2$. As before we choose $\beta=\beta_0$
 such that $u=\beta_0$ is the affine coordinate of the blown up point
 $Z_1\subset {\bf P}^1$. The rest of the discussion depends on $Z_2$ explicitly.
  If the length two sub-scheme $Z_2$ represents the tangent direction to
 $E_{\delta_0-1}$ at $Z_1$, then we take $\alpha=\beta_0$ and $A\not=0$.
   If $Z_2$ determines a tangent direction of $\acute{M}$ at $Z_1$ other than
 the tangent direction of $E_{\delta_0-1}$ at $Z_1$, then we choose
 a pair of non-identically zero $\alpha\not=\beta_0$ and $A$ such that 
 the conic determined by the equation $g_{\delta_0-1}(u, v)=0$ is tangent to this
 given direction specified by $Z_2$. Because the first jets of $g_{\delta_0-1}$ 
 at $(u, v)=(\beta_0, 0)$ are not identically zero, it is apparent that
 the polynomial $g_{\delta_0-1}$ does not vanish along $Z_3$.

\begin{defin}\label{defin; lead}
A polynomial $f(x,y)\in {\bf C}[x, y]$ 
is said to be leaded by the variable $x$ of degree 
two if it only contains monomials
 $x^sy^t$ with $0\leq s\leq 2$. It is said to be
 leaded by the variable $y$ of degree two if it only contains monomials
 of the type $x^sy^t$ with $0\leq t\leq 2$.
\end{defin}

 It is easy to observe that if a polynomial is leaded by $x$ (or by $y$)
 of degree two, then $y^2f({x\over y}, y)$ (or $x^2f(x, {y\over x})$),
 $f(x+a\cdot y, y)$ (or $f(x, y+bx)$) 
are still leaded by $x$ (or by $y$) as well.

 We employ the following inductive procedure with decreasing $i$, $1\leq i\leq 
 \delta_0-2$, 
to determine the transition
 maps between different coordinate charts and $g_i$.  
Starting from $g_{\delta_0-1}=g_{i+1}$ with $i=\delta_0-2$.

\noindent Case I:
 If $g_{i+1}$ is leaded by $x_{2i+1}$ or $y_{2i+2}$, then consider the
 following transition rule.

\noindent (i). Suppose that the point 
${\bf v}_i$ is in the open subset of the $i-th$
 ${\bf P}^1$ covered by 
the coordinate system $(x_{2i-1}, 0)$ with an affine coordinate 
 $(\alpha_{2i-1}, 0)$, then set $(x_{2i-1}-\alpha_{2i-1})=y_{2i+1}x_{2i+1}$, 
$y_{2i-1}=y_{2i+1}$; 
 $(x_{2i-1}-\alpha_{2i-1})=x_{2i+2}$, $y_{2i-1}=x_{2i+2}y_{2i+2}$
 for the coordinate transitions. 

\noindent (ii). Suppose that the point 
${\bf v}_i$ is not in the open subset of ${\bf P}^1$ covered
 by the coordinate system $(x_{2i-1}, 0)$, then it must be covered by
 the affine coordinate system $(0, y_{2i})$ with an affine coordinate 
$y_{2i}=0$.
 We set $y_{2i}=x_{2i+2}y_{2i+2}$, $x_{2i}=x_{2i+2}$; 
 $y_{2i}=y_{2i+1}$, $x_{2i}=x_{2i+1}y_{2i+1}$ for the coordinate transitions. 

 It is apparent that our choices of transition maps are consistent with 
the transitions of dual coordinates 
$(x_{2i+1}, y_{2i+1})\leftrightarrow (x_{2i+2}, y_{2i+2})$ defined earlier.

 Define $g_i$ as $(x_{2i-1}-\alpha_i)^2g_{i+1}({x_{2i-1}-\alpha_i\over y_{2i-1}},
 y_{2i-1})$ (in alternative (i))
 or $g_i=y_{2i}^2g_{i+1}({x_{2i}\over y_{2i}}, y_{2i})$ (in alternative (ii))
 if $g_{i+1}$ is leaded by $x_{2i+1}$ of degree two.

 Define $g_i$ as $(x_{2i-1}-\alpha_i)^2g_{i+1}(x_{2i-1}-\alpha_i, {y_{2i-1}\over 
 x_{2i-1}-\alpha_i})$ (in the alternative (i)) or 
$g_i=(x_{2i})^2g_{i+1}(x_{2i}, {y_{2i}\over x_{2i}})$ (in the alternative (ii))
 if $g_i$ is leaded by $y_{2i+2}$ of degree two.

\medskip

\noindent Case II:

 If $g_{i+1}$ is leaded by 
$y_{2i+1}$ or $x_{2i+2}$, we consider the following alternative
 scheme instead.

\noindent (i)'. If ${\bf v}_i$ is in the open subset of the $i-$th 
${\bf P}^1$ covered by 
the coordinate system $(0, y_{2i})$ with an affine coordinate 
 $(0, \beta_{2i})$, then set $(y_{2i}-\beta_{2i})=x_{2i+1}y_{2i+1}$, 
$x_{2i}=x_{2i+1}$;
 $(y_{2i}-\beta_{2i})=y_{2i+2}$, $x_{2i}=x_{2i+2}y_{2i+2}$
 for the coordinate transitions. 

\noindent (ii)'. If the point 
${\bf v}_i$ is not in the open subset of ${\bf P}^1$ covered
 by the coordinate system $(0, y_{2i})$, it must be covered by
 the affine coordinate system $(x_{2i-1}, 0)$ with an affine coordinate 
$x_{2i-1}=0$.
 Then we set $y_{2i}=x_{2i+1}y_{2i+1}$, $x_{2i}=x_{2i+1}$; 
 $y_{2i}=y_{2i+2}$, $x_{2i}=x_{2i+2}y_{2i+2}$ for the coordinate transitions. 

\medskip

Define $g_i$ to be $(y_{2i}-\beta_{2i})^2
g_{i+1}({x_{2i}\over y_{2i}-\beta_{2i}}, y_{2i}-\beta_{2i})$  
(in the alternative (i)') 
or $y_{2i}^2g_{i+1}({x_{2i}\over y_{2i}}, y_{2i})$ (in alternative (ii)') if 
 $g_{i+1}$ is leaded by $x_{2i+2}$.

Define $g_i$ to be $x_{2i}^2g_{i+1}(x_{2i}, {y_{2i}-\beta_{2i}\over  x_{2i}})$ 
(in the alternative (i)') or 
$x_{2i}^2g_{i+1}(x_{2i}, {y_{2i}\over x_{2i}})$ (in alternative (ii)')
 if $g_{i+1}$ is leaded by
 $y_{2i+1}$.

 It is easy to check that the two-variable polynomial $g_i$ is still 
leaded by one of its variables of degree two.

\medskip

  If $i\geq 1$, decrease $i$ by one, $i\mapsto i-1$, and repeat the above process
 until $i=0$.

Finally when $i=0$, jump out of the defining loops and 
define\footnote{Because we blow up ${\bf 0}\in M$.}
 $({\bf x}, {\bf y})=(x_1y_1, y_1)$ and 
$({\bf x}, {\bf y})=(x_2, x_2y_2)$. Define 
$g_0({\bf x}, {\bf y})={\bf y}^2g_1({{\bf x}\over {\bf y}}, {\bf y})$ if 
 $g_1$ is leaded by $x_1$ or $x_2$ of degree two. 
Define $g_0({\bf x}, {\bf y})={\bf x}^2g_1({\bf x}, {{\bf y}\over {\bf x}})$
 if $g_1$ is leaded by $y_1$ or $y_2$ of degree two.

 The union of the 
zero loci in $\acute{M}$ defined by $g_i=0$, $0\leq i\leq \delta_0-1$ on
 the $\delta_0$ different coordinate charts form 
 an algebraic curve intersecting $E_i$, $1\leq i\leq \delta_0-1$, with
 multiplicity two. By our construction of $g_{\delta_0-1}$ above, 
it vanishes along $Z_a$ but not along
 $Z_{a+1}$. Because 
 this curve is a divisor in $\acute{M}$, it is defined by a global section
$\acute{g}_0\in H^0(\acute{M}, {\cal O}_{\acute{M}}(-2\sum_{1\leq i\leq 
\delta_0-1}E_i)
\otimes {\cal I}_{Z_a})$. By our construction we know that
 $\acute{g}_0\not\in H^0(\acute{M}, 
{\cal O}_{\acute{M}}(-2\sum_{1\leq i\leq \delta_0-1}E_i)
\otimes {\cal I}_{Z_{a+1}})$.

  The explicit form of $\acute{g}_0$ on all $2\delta_0-1$ coordinate charts
   can be determined by the $g_i, 0\leq i\leq \delta_0-1$ (on the
 $\delta_0$ charts) and the transition maps
 among dual charts covering $E_i$. 
On the other hand, $\acute{f}|_{\acute{M}-\cup_{1\leq i\leq 
\delta_0-1}E_i}:\acute{M}-\cup_{1\leq i\leq 
\delta_0-1}E_i\mapsto M-{\bf 0}$ is an isomorphism under the blowing down map.
 Under this
 identification $({\bf x}, {\bf y})=\acute{f}^{\ast}(x, y)$, $g_0(x, y)$ defines a polynomial 
$\in H^0(M, {\cal O}_M)\cong {\bf C}[x, y]$. Because $\acute{g}_0$ is
 a global section of ${\cal O}_{\acute{M}}(-2\sum_{1\leq i\leq \delta_0-1}E_i)
\otimes {\cal I}_{Z_a}$, $g_0\in H^0(M, 
\acute{f}_{\ast}({\cal O}_{\acute{M}}(-2\sum_{1\leq i\leq \delta_0-1}E_i)
\otimes {\cal I}_{Z_a}))=H^0(M, {\cal I}_Z)$.

 The pair $(\acute{g}_0, g_0)$ satisfy the requirement in lemma 
\ref{lemm; existence}. So the proof of 
lemma \ref{lemm; existence} is finished. $\Box$

 From lemma \ref{lemm; saturate}, we know that for all chain-like 
 $\Gamma_0\in adm_2(\delta)$, the scheme $Z$ defined by the ideal sheaf 
${\cal I}_Z=
f_{\ast}{\cal O}_{\tilde{M}}(-2\sum_{i\leq \delta}E_i)$ is of length $3\delta$.

 This implies that for any chain-like $\Gamma_0\in adm_2(\delta)$, 
a point $\bar{b}\in Y_{\Gamma_0}\subset M_{\delta_0}$ 
is in the image of ${\cal M}_V\times_{M_{\delta}}
(M_{\delta}-\cup_{\Gamma\in \Delta(\delta)}S_{\Gamma})$ if there exists 
 a point $c\in{\cal M}_V$ above $\bar{b}$ such that its
 image inside the $\delta$ dimensional linear sub-system 
 $V$ under ${\cal M}_V\mapsto V$ lies in the kernel
 $H^0(M, {\cal I}_Z\otimes L)$ of $H^0(M, L)\mapsto H^0(M, {\cal O}_Z\otimes L)$.
 I.e. the corresponding curve represented by the point $c$ vanishes along
 the length $3\delta$ sub-scheme $Z$.

 By the $3\delta-1$-very ampleness condition on $L$ and the argument of 
G$\ddot{o}$ttsche's proposition \ref{prop; Got}, for generic choices of
 $V\subset |L|$ there can be no such curve.  So the proof of 
proposition \ref{prop; tran} is complete. $\Box$

\medskip

\section{Appendix:
The Relationship with the Gromov-Witten Invariant}\label{section; compare}

\bigskip

 In the previous section, we have given
 an algebraic proof that the ``number of
 $\delta$-nodes nodal curves'' in a general $\delta$ dimensional 
linear-subsystem of $|L|$ can be expressed as a universal polynomial
 of $L\cdot L, L\cdot c_1({\bf K}_M), c_1({\bf K}_M)^2, c_2(M)$. This
 ``number of nodal curves'' is understood in the sense of G$\ddot{o}$ttsche (see
 proposition \ref{prop; Got}) and our proof involves the various 
constructions in
 the algebraic family Seiberg-Witten theory. The reader working on
 Gromov-Witten invariant may desire to understand the relationship between the
 ``family Seiberg-Witten invariant count'' and the Gromov-Witten invariant 
count. As sometimes it may lead to some misunderstanding of the result,
so we offer some clarification here.

  Firstly, for the difference between the usual ``algebraic'' Seiberg-Witten 
invariant (over $B=pt$) and the topological version of Seiberg-Witten invariant (which
 is equivalent to the ``right genus''
 Gromov-Witten invariant of an algebraic 
surface by [T1], [T2], [T3] and [IP]), please consult sub-section 4.3.1 of [Liu3].
 So we will focus on the difference of our ``number of nodal curves'' and 
the usual ``wrong genera'' Gromov-Witten invariant.

 Given a holomorphic line bundle $L$ on $M$, the adjunction formula,
 $C^2+c_1({\bf K}_M)\cdot C=2g(C)-2$ with $C=c_1(L)$, predicts a preferred genus
 of curves in $|L|$. An identical adjunction formula holds in the
 pseudo-holomorphic category as well.
Taubes had used the pseudo-holomorphic 
curve counting of the preferred 
genus in developing his ``SW=Gr'' theorem, [T1], [T2], [T3] etc.

 On the other hand, in the standard $GW$ invariant the genus of the
 source curve is not pre-determined by the class $C=c_1(L)\in H^2(M, {\bf Z})$.
 In fact, for all $g\in {\bf N}$, it makes sense to define the genus $g$
 Gromov-Witten invariant $GW_g(C)$ 
which enumerates the virtual number of (pseudo-)holomorphic
 maps representing $C\in H^2(M, {\bf Z})$ from source curves with genus $g$.

 The fundamental observation which relates the ``nodal curve counting'' with
the number $GW_g(C)$ is that a genus $g<g(C)$ (pseudo)-holomorphic curve tends to
 develop $g(C)-C$ nodal singularities. It is because a (pseudo)-holomorphic
 map from a genus $g$, $g<g(C)$ curve, $\Sigma_g$ mapping into $M$ cannot be
 embedded (or it violates the adjunction formula) and tends to develop
 isolated singularities in its image (if it is not multiple-covered 
or bubbling off multiple coverings of two spheres). 
 The nodal curve singularities are
 preferred because of dimension reason. Suppose that 
the image of $\Sigma_g$ has
 developed singularities at $x_1, x_2, \cdots, x_k\in M$ with
 singular multiplicities $m_1, m_2, \cdots, m_k$ for some $k\in {\bf N}$.

 Then the adjunction formula for singular curves, exercise 1.3 and 
corollary 3.7 of chapter V of [Ha], implies that

 $$2g-2+\sum_{i\leq k}m_i(m_i-1)=C^2+c_1({\bf K}_M)\cdot C.$$

 On the other hand for $g\geq 2$,
 the expected dimension of the Gromov-Witten moduli space
 is equal to 

$$\hskip -1.3in
\int_{\Sigma_g}c_1(M)-2(g-1)+dim_{\bf C}{\cal M}_g
=-c_1({\bf K}_M)\cdot C+
{C^2+c_1({\bf K}_M)\cdot C\over 2}-\sum_{i\leq k}{m_i(m_i-1)\over 2}=
{C^2-c_1({\bf K}_M)\cdot C\over 2}-\sum_{i\leq k}{m_i(m_i-1)\over 2},$$

 where
 ${C^2-c_1({\bf K}_M)\cdot C\over 2}$ is both (i). the expected dimension of
 Gromov-Taubes moduli space (see Taubes [T3]) and
(ii). The $C$ dependent term of the
 surface Riemann-Roch formula and closely related to 
the dimension of the (non-)linear
 system associated to a given $C\in H^2(M, {\bf Z})$. On the other hand,
 the expected dimension of algebraic curves carrying $k$ singularities with
 multiplicities $m_1, m_2, \cdots, m_k$ is at most 
 ${C^2-c_1({\bf K}_M)\cdot C\over 2}-
\sum_{1\leq i\leq \delta}({m_i^2+m_i\over 2}-2)$. Because 
 $m(m-1)\leq m(m+1)-4$ for $m\geq 2$ and the equality saturates only when $m=2$,
 the curves with singular 
multiplicities $>2$ are of lower dimensions in the moduli space
 of genus $g$ curves. A closer look at the type of the curve singularity shows
 that any double point other than nodal ($A_1$) singularity drops the complex
 dimension of the deformation space of curves by at least two.
 Therefore a generic genus $g$ immersed curve dual to $C$ develop 
 $\delta=g(C)-g$ nodal singularities.

  When one works in the $C^{\infty}$ category and perturbs the almost complex
 structures of the algebraic surface $M$ to a generic one, one expects the
 pseudo-holomorphic maps to develop nodal singularities. Thus both family
 Seiberg-Witten invariant of $C-2\sum_{1\leq i\leq \delta}E_i$ and 
$Gr_{g(C)-\delta}(C)$ are objects enumerating $\delta$-nodes curve dual to $C$.

 The fundamental question we have to clarify and answer is,

\noindent {\bf Question}: Do ${1\over \delta!}{\cal AFSW}_{M_{\delta+1}\times
 T(M)\mapsto M_{\delta}\times T(M)}^{\ast}(1, 
C-\sum_{1\leq i\leq \delta}2E_i)$ and the Ruan-Tian version of
 $GW_{g(C)-\delta}(C)$ ``always'' 
enumerate the ``number of nodal curves'' in a totally identical way?

Certainly there are many important cases that they do enumerate the same
 numbers, e.g. when $C$ is a primitive cohomology class of a $K3$ or $T^4$, etc.
  
But the general answer of this question is ``No''. In the following
we offer an explanation of the causes. 

\medskip

\noindent 1. The Gromov-Witten invariant enumerates the 
``number of (pseudo)-holomorphic maps'' instead of immersed curves (viewed as
 divisors in $M$). When the cohomology class $C$ is primitive, i.e. it is not
 a multiple of any other element in $H^2(M, {\bf Z})$, each holomorphic map
 determines uniquely a nodal curve and vice versa. But when $C$ becomes
 non-primitive, sometimes there can be multiple coverings of holomorphic maps
 such that the image (without counting multiplicity) is dual to ${1\over m}C$ 
, for some $m\in {\bf N}$. These multiple coverings of holomorphic
 maps contribute to $GW_{g(C)-\delta}(C)$ as well.  But they do not 
correspond to immersed nodal curve dual to $C$ and are mostly ignored by the scheme
 of family invariant. In general, the object $GW_g(C)$ is ${\bf Q}$ valued,
 reflecting the orbifold structure of the space $\overline{{\cal M}_{g, n}}$.
On the other hand, either 
${1\over \delta!}{\cal AFSW}_{M_{\delta+1}\times T(M)\mapsto
 M_{\delta}\times T(M)}^{\ast}(1, C-2\sum_{i\leq
 \delta} E_i)$ or
${1\over \delta!}{\cal AFSW}_{M_{\delta+1}\times \{t_L\}\mapsto
 M_{\delta}\times \{t_L\}}^{\ast}(1, C-2\sum_{i\leq
 \delta} E_i)$ is always integral valued.

\medskip

\noindent 2. In the algebraic category, it is highly non-trivial to
make the appropriated moduli space of curves defined as the zero locus 
of a transversal algebraic section.  This applies to both family Seiberg-Witten
 invariant and Gromov-Witten invariant. When the moduli space is 
 not transversal, one interprets the invariants as some types of virtual
 number counts. Without any transversality result, the correspondence between
 algebraic family Seiberg-Witten invariant and Gromov-Witten invariant is
 not transparent at all.
 If one decides to work instead in the symplectic (pseudo-holomorphic) 
category, it is usually easier to get the transversality result of the appropriated
 moduli spaces by choosing generic almost complex structures tamed by a
 simplectic structure on $M$. On the other hand, at this moment it is not
 clear how to define the ``number of nodal singular curves'' of a general
 symplectic four-manifold without selecting special classes $C$ or
 almost complex structures $J$. In the situation when one can make sure the
 cut down moduli space consists of a finite number pseudo-holomorphic
 nodal singular curves, one has to identify the algebraic family Seiberg
Witten invariant with its topological cousin ``Family Seiberg-Witten invariant''
 and employ the technique of Taubes' ``SW=Gr'' to compare the
 solutions to family Seiberg-Witten equations and the smooth curve
 resolved from the nodal singular curves in $M$.

  Unluckily the gluing machineries of Taubes in his seminal long papers [T1], [T2],
falls out of the algebraic category. Thus one may desire a purely algebraic
 method to determine the family algebraic Seiberg-Witten invariants or relate
 them with the Gromov-Witten invariant.

\medskip

\noindent 3. In G$\ddot{o}$ttsche's definition, the ``number of nodal curves'' is defined
 for $L$ to be $5\delta-1-$very ample.  Under this assumption, there is a
 well defined integer, gotten by counting the discrete number 
(with multiplicities) of nodal curves.
On the other hand, when $L$ fails to be $5\delta-1-$very ample, generally speaking
 we do not expect ${1\over \delta}{\cal AFSW}_{M_{\delta+1}\times \{t_L\}\mapsto
 M_{\delta}\times \{t_L\}}^{\ast}(1, C-2\sum_{i\leq
 \delta} E_i)$ to calculate the number of nodal curves. In general, 
 we have to subtract all the contributions from type $II$ exceptional curves [Liu6] 
 as well. The result is usually manifold dependent and involves some
 generalization of the technique used in this paper. 
 
 One exception is the case that $M=K3$ or $T^4$ when all the contributions
 from type $II$ exceptional curves are known to vanish due to the fact
 $SW(2C)=0$ for any $C\in H^2(M, {\bf Z})$, $M=K3$ or $M=T^4$. 
 The details about the contribution of the type $II$ curve multiple-coverings
 will be developed in a separated article and we do not elaborate it here.
In these cases, the numbers ${1\over \delta!}
{\cal AFSW}_{M_{\delta+1}\times T(M)\mapsto
 M_{\delta}\times T(M)}^{\ast}(1, C-2\sum_{i\leq \delta} E_i)$
  are actually equal to ``the number of nodal singular curves'' dual to $C$,
 understood as a virtual intersection number. When $C$ is primitive, this 
number coincides with the usual Gromov-Witten invariant ([BL1], [BL2]).
 When $C$ is not primitive, the number of immersed nodal singular curves 
 dual to $C$ differ from the usual Gromov-Witten invariant count as the
 former does not count the multiple covering maps of curves from the
 fractional multiples of $C$. This fact can be seen from the discrepancy between
 the Yau-Zaslow formula and Gathman's calculation [Gat] of $GW_{g}(C)$ of 
 a $2$-multiple of a primitive class in the $K3$ lattice.

 In Gathman's calculation, he considers an algebraic K3 surface which
 is the double-covering of ${\bf P}^2$ ramified along a generic sextic
 curve. He considers the class $C$, $C^2=5$, to be twice of the pull-back of
 the hyperplane class from ${\bf P}^2$ and the answer he got for $GW_0(C)$ was 
 $N_5+{1\over 8}N_2$, where $\sum_{\delta\geq 0}N_{\delta}q^{\delta}=
\prod_{i>0}{1\over (1-q^i)^{24}}$ is the generating function of the
 Yau-Zaslow formula. In his calculation, he had used a degenerated complex
 structure to enumerate curves and the $176256$ rational 
curves from his (i)., (ii)., and (iii)(b). can be thought to be the 
 degenerations from nodal pseudo-holomorphic rational curves of generic
 ${\bf S}^2$ families of almost complex structures. On the other
 hand, the $324$ rational curves from his (iii).(a). are honest double coverings
 of primitive rational curves and can not be degenerated from immersed
  pseudo-holomorphic nodal curves.

 This example indicates clearly that Yau-Zaslow conjecture is not
 about the prediction of Gromov-Witten invariant at all. The prediction
 of Yau-Zaslow conjecture coincides with the Gromov-Witten calculation
 only for the primitive classes when the multiple coverings addressed in
 1. do not occur. Therefore it is dangerous to mix up ``the number of
 {\bf nodal} rational curves'' with ``the number of holomorphic maps
 to rational curves'' and identify them conceptually.

Moreover when the class $C$ is not primitive, in Gathman's calculation 
the contributions from the multiple coverings
 does show up in the correction term and it is desirable to find
 out the relationship between them explicitly.

 In the symplectic category 
if one can prevent the discrete number of pseudo-holomorphic maps
 with nodal curve images to converge
 to some multiple covering of curves in the most general content (this is
not achieved at this moment), 
 then a multiple covering formula (presumably determined by
 certain intersection numbers on ${\cal M}_{g, n}$) should allow us to
 relate the number of nodal singular curves dual to $C$ with the numbers 
$GW_g(C)$ and the formula should be of combinatorial nature.

\medskip

\noindent 4. When the geometric genus $p_g$ is greater than zero. The usual
 $GW_g(C)$ counts are mostly zero except for a finite number of
 classes (direct related to the so-called Seiberg-Witten basic classes for
 K$\ddot{a}$hler surfaces). Yet ${\cal AFSW}^{\ast}$ still picks up
 non-trivial contributions. This can be seen by the lower $\delta$ formula
 calculated by Vainsencher/Kleiman\&Piene[][]. The reason behind the
 discrepancy is that $GW_g(C)$ (like Taubes' version of Gromov-Taubes invariant)
 are symplectic invariants. The algebraic surfaces are Seiberg-Witten
 simple type that classes with positive moduli space dimension
 ${C^2-c_1({\bf K}_M)\cdot C\over 2}$ have vanishing invariants. Algebraically
 it is reflected in the fact that there is a $p_g$ difference between the
  dimension of the projective space 
$p_g+{C^2-C\cdot c_1({\bf K}_M)\over 2}$ (assuming $q(M)=0$) and the
 expected Gromov-Taubes dimension of the curve ${C^2-C\cdot 
c_1({\bf K}_M)\over 2}$.

 For a very ample line bundle $L$ such that ${\bf K}_M\otimes L$ is ample,
 this gives a trivial rank $p_g$ complex obstruction sub-bundle above 
 ${\cal M}_L=|L|$, which causes the usual Seiberg-Witten invariant of $c_1(L)=C$ 
 and the sub-sequent family invariant of $C-\sum_{i\leq n}m_i E_i$ to vanish.

 This is why the hyperwinding families of $K3$ or $T^4$ had been used to
 count the nodal curves in [Liu1].
 The algebraic Seiberg-Witten invariants and the family algebraic Seiberg-Witten
 invariants are defined such that for such $L$, they remove the trivial
 rank $p_g$ obstruction bundle from $|L|$ and shift the expected dimension
 of the moduli space up by $p_g$. This causes the algebraic family
 Seiberg-Witten invariants to be ``enumeration invariants'' but not the usual
 symplectic invariants under deformation. The way that we realize ${\cal ASW}$ or
 ${\cal AFSW}^{\ast}$, etc. to be invariants is through a different route: 
 After calculating ${\cal ASW}$, ${\cal AFSW}$ or the modified invariants 
 ${\cal AFSW}^{\ast}$,
 they can be expressed in terms of some datum which depend on only the homotopic
 type of the algebraic surface and $C$.

\medskip

\noindent 5. Another significant difference between the algebraic family
Seiberg-Witten invariant of immersed curves 
and the Gromov-Witten invariant of maps is that (assuming $p_g=0$ for 
simplicity) when the appropriated
 moduli spaces are transversal and the dimensions of the
 compactifying boundary components drop, the former object enumerates all the
irreducible as well as reducible nodal singular curves dual to $C$ while
 $GW_g(C)$ encodes only the irreducible nodal curves dual to $C$. The rough 
reason is that according to Taubes, $SW(2C-c_1({\bf K}_M))=Gr(C)$ enumerates
  the connected as well as disconnected smooth curves dual to $C$. Based
 on this philosophy that the enumeration of the family invariants of a given family
 should include
connected as well as dis-connected smooth curves within a family, 
the modified family invariants of the universal families also enumerate 
 reducible nodal curves where two curves intersect at normal crossing 
singularities (and the normal crossing singularities get resolved after
 the blowing ups). On the other hand, reducible nodal curves with multiple
 irreducible components can not be the pseudo-holomorphic 
image of irreducible Riemann surface. Thus, reducible nodal curves can only
 be viewed as the images of semi-stable maps while the source curve is parametrized
 by a point in the boundary point 
$\overline{{\cal M}_{g, n}}-{\cal M}_{g, n}$ for some pair of $(g, n)$, 
and is of lower expected dimension.

 This symptom is purely of combinatorial nature and should
 be cured by re-developing the generating series.

{}
\end{document}